Noè A. Caruso
Alessandro Michelangeli

# Inverse linear problems on Hilbert space and their Krylov solvability

# Preface

This monograph is centred at the intersection of three mathematical topics, that are theoretical in nature, yet with motivations and relevance deep rooted in applications: the linear inverse problems on abstract, in general infinite-dimensional Hilbert space; the notion of Krylov subspace associated to an inverse problem, i.e., the cyclic subspace built upon the datum of the inverse problem by repeated application of the linear operator; the possibility to solve the inverse problem by means of Krylov subspace methods, namely projection methods where the finite-dimensional truncation is made with respect to the Krylov subspace and the approximants converge to an exact solution to the inverse problem.

For subject classification purposes, the issue of the *Krylov solvability of abstract inverse linear problems* lies at the edge between functional analysis, measure theory, operator theory, and spectral theory on the one hand, and numerical analysis and approximation theory on the other hand.

The underlying motivation is the ubiquitous occurrence of linear phenomena that produce an *output g* from an *input f* according to a *linear law A*, so that from the exact or approximate measurement of $g$ one tries to recover exact or approximate information on $f$. A vast generality of such phenomena may be encompassed into a unified abstract framework, where the datum $g$ and the unknown $f$ belong to an abstract, possibly infinite-dimensional Hilbert space $\mathcal{H}$, and $A$ is a linear, bounded or unbounded, operator acting on $\mathcal{H}$. One then speaks of an *abstract linear inverse problem* of the form $Af = g$.

As in countless other instances in mathematics, such an abstraction is aimed at recognising and analysing general features and patterns that are in common with the very many different concrete realisations of linear inverse problems in applications, where there is a huge spectrum of methodologies and techniques, tailored on this or that class of inverse problems.

In fact, scientific computing demands a suitable discretisation of a given inverse problem, be it finite-dimensional but with huge size, or even infinite-dimensional. Greatly widespread and popular methods for the numerical solution of a linear inverse problem consist of solving truncated versions of the given problem, obtained by projecting the operator $A$ and the datum $g$ onto a finite-dimensional subspace of $\mathcal{H}$: under suitable conditions, the sequence of approximated solutions obtained by increasing the size of the projection subspace do converge to an actual solution $f$ to $Af = g$.

An ample variety of such *projection methods* is available nowadays for disparate applications – think of the very common Petrov-Galerkin approximation schemes, the finite element methods, and the like. Among the most popular, sophisticated, and efficient, are the *Krylov subspace methods*, where broadly speaking the approximants to the solution $f$ are sought for within the finite linear combinations of the vectors $g, Ag, A^2g, A^3g, \ldots$ The linear span $\mathcal{K}(A,g)$ of such vectors is called the *Krylov subspace* associated to the inverse problem $Af = g$ and the possibility of approximating $f$ with vectors from the Krylov subspace, or in other words the occurrence that $f$ is a limit point of $\mathcal{K}(A,g)$ and hence belongs to the closure in $\mathcal{H}$ of $\mathcal{K}(A,g)$, is referred to as the *Krylov solvability* of the inverse problem.

When we recently started investigating the general question of Krylov solvability of abstract inverse problems on Hilbert space, we soon faced two amazing circumstances, that we tried at our best to implement also in the present monograph.

To begin with, it is hard to underestimate the *conceptual relevance* of such a theoretical question in applications. Here the matter isn't, of course, to improve the rate of convergence of a particular Krylov-based algorithm, or to make a Krylov subspace method best tailored for a specific application: the present abstract setting surely does not allow such ad hoc analyses. Yet, to have access to general theoretical tools, formulated in abstract, versatile terms, which allow to decide a priori what categories of inverse problems are indeed Krylov solvable, and to understand the structural mechanisms that make an inverse problem Krylov solvable or non-solvable, provides a most valuable information prior to attacking the problem with the actual computational weaponry. If the solution $f$ to $Af = g$ is not a Krylov solution, the norm distance between $f$ and its finite-dimensional Krylov approximants, however accurate their computation, cannot be made smaller than a finite threshold – in other words, the algorithm does not converge to a solution. The lack of Krylov solvability is a possibility that may indeed occur, depending of the features of the operator $A$ and of the datum $g$ of the inverse problem, and is a sort of drawback of this type of projection methods where the advantage is instead the "easy" construction of the approximation space, simply spanned by the "nice" vectors $g, Ag, A^2g, A^3g, \ldots$ (think, in contrast, of how accurate and subtle must the choice be of the approximating functions in finite element methods).

There was another fascinating circumstance at the beginning of our investigation: apart from specific studies of the convergence of certain Krylov-based methods for infinite-dimensional inverse problems, the general study of Krylov solvability at the abstract level we are referring to appeared to be an essentially uncharted territory. Previous related studies can be tracked down in the literature from subject areas that are often considerably far apart and mutually independent, thus denoting the different perspective on several questions connected with Krylov subspaces, existence of Krylov solutions, and the like. This is probably due also to the dominance of other mainstream lines, in particular the interesting subject of ill-posed and ill-conditioned inverse problems.

It was therefore natural to address in rapid sequence a first corpus of fundamental questions – an initial "taxonomy" – through paradigmatic examples and counter-examples, of Krylov-solvable and non-solvable inverse problems; the identification of intrinsic mechanisms for a solution $f$ to be a Krylov solution; the possible uniqueness of a Krylov solution; the identification of special classes of operators inducing Krylov-solvable inverse problems; the rigorous formulation and the analysis of the counterpart questions from a bounded to an unbounded operator $A$ on infinite-dimensional Hilbert space; the behaviour of a Krylov subspace under suitable small perturbations; the phenomena related to Krylov solvability or lack-of under perturbations of the inverse problems; the possibility and genericity of weaker forms of convergence.

This resulted in a series of works that we published over the last couple of years, which produced a first stable amount of knowledge at a sufficiently mature level, together with adequately neat and sharply formulated open questions, that then became natural to conceive the present monograph. It contains a polished and harmonised version of all our recent findings, and the background of the related literature, now displayed within a rational unifying discussion, and also in a self-contained and rather compact form.

The compactness and the self-contained character of the presentation is in fact an additional virtue we particularly cared of in light of the rather diverse range of mathematical tools employed here. There is indeed a core apparatus from abstract functional analysis, operator theory on Hilbert space, spectral theory, with also techniques from measure theory, point-set topology, theory of orthogonal polynomials. Such a mixture of tools and approaches is surely due to the very nature of certain questions that naturally redirect to different subject areas. Thus, for instance, the study of the possible density in $\mathcal{H}$ of the Krylov subspace $\mathcal{K}(A, g)$ (an occurrence that makes the Krylov solvability of $Af = g$ trivial) concerns the classical question whether the vector $g$ is cyclic for the operator $A$. Or also: the characterisation of the Krylov subspace of a self-adjoint inverse problem is naturally dealt with in the language of the spectral theorem and the functional calculus. At the same time the miscellaneous presence of different techniques is presumably the signature of the fact that for further investigation new ad hoc methods need be developed.

At any rate, whereas all the various background materials are given for granted, throughout this monograph we made scrupulous reference to all the standard sources, so as to immediately orient the reader towards the several subject areas involved in our discussion.

We kept a style of presentation at the graduate level and in fact we could present a selection of our results in recent graduate seminars and international workshops. The whole monograph is implicitly designed so as to serve as a single-subject thematic graduate course on the one hand, and as a reference guide for experts of neighbouring fields, including operator theorists, applied and numerical analysts, with the hope that the corpus of established results and the many questions we posed and left open may indeed spark interest and generate a structured research agenda in the forthcoming future.

We should like to acknowledge the very useful inputs received at various stages by expert colleagues, particularly by M. Erceg (Zagreb), L. Grubišić (Zagreb), L. Heltai (Trieste), M. Ligabò (Bari), A. S. Nemirovskiy (Atlanta), P. Novati (Trieste). The results contained in this monograph have been developed at an initial stage at the International School for Advanced Studies (SISSA) Trieste, and later refined in the course of our scientific activities, respectively, at the Gran Sasso Science Institute L'Aquila (N.C.), and at the Institute of Applied Mathematics, and the Hausdorff Centre for Mathematics Bonn (A.M.). For part of this research we also benefited from the support of the Italian National Institute for Advanced Mathematics INdAM and the Alexander von Humboldt Foundation.

L'Aquila (Italy) and Bonn (Germany),                                           *Noè A. Caruso*
April 2021                                                                     *Alessandro Michelangeli*

# Contents

# Acronyms

The following list contains the most frequently used symbols used throughout this book along with their brief description.

| | |
|---|---|
| $\mathbb{N}, \mathbb{Z}, \mathbb{R}, \mathbb{C}$ | Natural, integer, real, complex numbers |
| $\mathbb{N}_0$ | Non-negative integers |
| $\mathcal{H}$ | Abstract Hilbert space |
| $\langle \cdot, \cdot \rangle, \langle \cdot, \cdot \rangle_{\mathcal{H}}$ | Hilbert space scalar product, antilinear in the first argument |
| $\|\cdot\|, \|\cdot\|_{\mathcal{H}}$ | Hilbert space norm |
| $\perp$ | Orthogonality in Hilbert space |
| $S^\perp$ | Orthogonal complement to the linear span of the subset $S \subset \mathcal{H}$ |
| $\overline{V}$ | Hilbert norm closure of the subspace $V \subset \mathcal{H}$ |
| $\|\cdot\|_{\mathrm{op}}$ | Operator norm |
| $\mathcal{B}(\mathcal{H})$ | $C^*$-algebra of everywhere defined, bounded linear maps on $\mathcal{H}$ |
| $A^*$ | Adjoint of a linear operator $A$ |
| $\mathbb{1}$ | Identity operator on $\mathcal{H}$ |
| $\mathbb{O}$ | Zero operator on $\mathcal{H}$ |
| $\rho(A)$ | Resolvent set of a linear operator $A$ |
| $\sigma(A)$ | Spectrum of a linear operator $A$ |
| $\mathcal{D}(A)$ | Domain of a linear operator $A$ |
| $\mathrm{ran}\, A$ | Range of a linear operator $A$ |
| $\ker A$ | Kernel of a linear operator $A$ |
| $\|\cdot\|_A$ | Graph norm of a linear operator $A$ |
| $C^\infty(A)$ | Subspace of $A$-smooth vectors |
| $\mathcal{D}^a(A), \mathcal{D}^{qa}(A)$ | Sets of $A$-analytic, respectively, $A$-quasi-analytic vectors |
| $E^{(A)}$ | Spectral projection measure of the self-adjoint operator $A$ |
| $|\varphi\rangle\langle\psi|$ | Rank-one orthogonal projection $\xi \mapsto \langle\psi, \xi\rangle\varphi$ |
| $\mathcal{K}(A, g)$ | Krylov subspace relative to the operator $A$ and vector $g$ |
| $\mathcal{K}_N(A, g)$ | $N$-th order Krylov subspace relative to $A$ and $g$ |
| $\mathcal{K}^{(\Xi)}(A, g)$ | Rational Krylov subspace relative to $A, g$ for the sequence $\Xi$ |
| $\mathcal{I}(A, g)$ | Krylov intersection relative to $A$ and $g$ |
| $\widehat{d}_w(U, V)$ | Weak gap metric on weakly closed subsets of the unit ball of $\mathcal{H}$ |

# Chapter 1
# Introduction and motivation

## 1.1 Abstract inverse linear problems on Hilbert space

The primary focus of this book is inverse linear problems and the possibility of expressing their solutions in terms of convenient approximants produced by iterative algorithms.

As the "inverse problems" jargon appears in an ample spectrum of mathematical contexts, let us first of all clarify the meaning it has in the present framework.

The underlying motivation is the ubiquitous occurrence of linear phenomena that produce an *output g* from an *input f* according to a *linear law A*, so that from the exact or approximate measurement of $g$ one tries to recover exact or approximate information on $f$.

As an elementary example, knowing (measuring) the accelerations $a_1, a_2, a_3$ of three classical point particles constrained on a straight line, of masses respectively $m_1, m_2, m_3$, and interacting through two-body forces $F_{jk}$ between the $j$-th and the $k$-th particle, hence with $F_{jk} = -F_{kj}$, one can determine such forces $F_{12}, F_{13}, F_{23}$ by exploiting Newton's law

$$\begin{aligned} F_{12} + F_{13} &= m_1 a_1 \\ -F_{12} + F_{23} &= m_2 a_2 \\ -F_{13} - F_{23} &= m_3 a_3, \end{aligned}$$

i.e.,

$$\begin{pmatrix} m_1^{-1} & m_1^{-1} & 0 \\ -m_2^{-1} & 0 & m_2^{-1} \\ 0 & -m_3^{-1} & -m_3^{-1} \end{pmatrix} \begin{pmatrix} F_{12} \\ F_{13} \\ F_{23} \end{pmatrix} = \begin{pmatrix} a_1 \\ a_2 \\ a_3 \end{pmatrix}.$$

This is the problem of finding the unknown vector $f \equiv (F_{12}, F_{13}, F_{23})$ from the knowledge of the datum $g \equiv (a_1, a_2, a_3)$ through the law $Af = g$, where $A$ is the above matrix with mass coefficients. In fact, here $\det A = -2m_1 m_2 m_3 \neq 0$, hence in this case the solution $f$ exists and is unique for given datum $g$.

Another example on the same conceptual footing, but beyond the above finite-dimensional setting, is the determination of the electric potential $V(x)$ at each point $x\in\mathbb{R}^3$ generated by a known (measured) charge density $\rho(x)$: this amounts to solving Poisson's equation

$$\Delta V \;=\; -\frac{\rho}{\varepsilon_0}$$

(where $\Delta \equiv \partial^2/\partial x_1^2 + \partial^2/\partial x_2^2 + \partial^2/\partial x_3^2$ is the Laplace differential operator on $\mathbb{R}^3$, $x \equiv (x_1, x_2, x_3) \in \mathbb{R}^3$, and $\varepsilon_0$ is the vacuum permittivity), which can be again interpreted as the the determination of the unknown $f \equiv V$ (this time a function, not a finite-dimensional vector) from the datum $g \equiv \rho$ through the law $Af = g$, where $A = -\varepsilon_0 \Delta$. As long as $\rho$ is sufficiently regular and fast decaying at infinity, it is well known that the solution $V$ exists and is unique, given by

$$V(x) \;=\; \frac{1}{4\pi\varepsilon_0} \int_{\mathbb{R}^2} \frac{\rho(x')}{|x-x'|}\, \mathrm{d}x\,.$$

Observe that in this second example, essentially as a consequence of the infinite-dimensionality, the problem $Af = g$ must be supplemented with the information of which functional space the datum $\rho$ belongs to and hence over which space the unknown $V$ must be searched for.

These and many other examples can be cast into a unified abstract framework where the possible inputs and outputs constitute an abstract linear vector space $\mathcal{H}$ on which an action is performed by a linear operator $A$. We shall then refer to the problem

$$Af \;=\; g \qquad (1.1)$$

in the **unknown** $f \in \mathcal{H}$ as the **inverse linear problem on the space** $\mathcal{H}$ associated with the law (the **linear operator**) $A$ and the **datum** $g \in \mathcal{H}$. We shall often drop the adjective 'linear', as throughout our discussion we will not deal with non-linear maps.

Borrowing familiar nomenclature, we shall refer to the inverse problem (1.1) as **solvable** when a solution $f$ exists, namely when $g \in \mathrm{ran} A$, the range of the operator $A$ in $\mathcal{H}$, and **well-defined** when in addition the existing solution is unique, namely when $A$ is injective. In this case it is customary to refer to $f$ as *the* **exact solution** to the problem (otherwise one speaks of *an* exact solution). By linearity, should the inverse problem (1.1) admit a multiplicity of solutions, they are all of the form $f + f_0$ with $f_0 \in \ker A$, the kernel (null space) of $A$. If $A$ is injective, then $A$ is invertible on its range and one writes $A^{-1}$ for its inverse, so that the unique solution $f$ to $Af = g$ for given $g \in \mathrm{ran} A$ is $f = A^{-1}g$.

Inverse linear problems in the above sense arise in a vast multitude of applications. Even when one has the a priori information that the problem is well-defined, hence when the main theoretical issues of existence and uniqueness of solution are answered in the affirmative, an explicit, closed solution formula for $f$ may well be lacking. This frequent occurrence motivates, in all the realm of numerical analysis and scientific computing, the determination of $f$ through subsequent, finer and finer approximants that are explicitly produced by suitable approximation algorithms.

For the above abstract setting of inverse problem on the vector space $\mathcal{H}$ to include the notion of a sequence $(f_n)_{n\in\mathbb{N}}$ of approximants such that $f_n \to f$, or $f_n \approx f$, as $n \to \infty$, the space $\mathcal{H}$ must be equipped with additional structure, a topology in the first place.

In fact, a vast variety of phenomena are encompassed by such a mathematical abstraction when $\mathcal{H}$ is taken to be a (possibly infinite-dimensional) inner product and complete vector space, namely a *Hilbert space*, and *A* is a *closed linear operator* acting on $\mathcal{H}$. We give for granted here all the basic pre-requisites on Hilbert spaces and linear (including closed) operators acting on them, referring to standard references such as [90, 67, 88, 116, 11, 29, 30, 97] for details. Besides, in essentially all the discussion of this monograph the choice of the field $\mathbb{K} = \mathbb{R}$ or $\mathbb{C}$ the Hilbert space $\mathcal{H}$ is built upon does not make any difference: for convenience, we shall carry on the notation with $\mathbb{K} = \mathbb{C}$, thus with a *complex* Hilbert space $\mathcal{H}$, making the convention that its inner product $\langle \cdot, \cdot \rangle$ is linear in the second entry and anti-linear in the first, and with associated norm $\|\cdot\|$. (When multiple norms will be used in our discusson, we shall write $\|\cdot\|_{\mathcal{H}}$ to avoid confusion.) Let us also remark that the (somewhat minimal) requirement of operator closedness of *A* is aimed at having a non-trivial notion of the spectrum of *A*, hence to allow for the possible use of spectral methods in solving (1.1).

For such $\mathcal{H}$ we shall refer to '$Af = g$' above as an **abstract inverse linear problem on the Hilbert space** $\mathcal{H}$.

We should like to point out that various levels of abstraction are implemented here: (a) the dimensionality of $\mathcal{H}$, finite or infinite, (b) the boundedness or unboundedness of *A*, (c) the spectral properties of *A* (from a purely discrete spectrum to richer structures with continuous components, separated or not from zero, and so on).

When $\dim \mathcal{H} < \infty$ the inverse problem involves finite matrices, albeit of possibly huge size, and is typically under a very accurate control in all its aspects (algebraic, analytic, numerical, including in applications the control of the rate of convergence of approximants, etc.). For this reason, in this monograph we shall keep the finite-dimensional problem in the background as a comparison playground and source of generalisations, and we shall rather focus on the infinite-dimensional setting, namely infinite-dimensional inverse problems on Hilbert space.

For a large part of our discussion, that is natural to expose first, *A* is going to be an everywhere defined and bounded operator on $\mathcal{H}$, meaning that the domain $\mathcal{D}(A)$, a linear subspace of $\mathcal{H}$, is actually the whole $\mathcal{H}$, and that *A* has finite operator norm $\|A\|_{\mathrm{op}} < +\infty$, recalling that in general

$$\|A\|_{\mathrm{op}} := \sup_{\substack{h \in \mathcal{D}(A) \\ h \neq 0}} \frac{\|Ah\|}{\|h\|}.$$

We refer to the above-mentioned references for the general theory of bounded operators on Hilbert space – in particular let us recall that, by linearity, continuity (in the operator norm) and boundedness are equivalent, and we will use both terms interchangeably. When the inverse problem (1.1) on $\mathcal{H}$ is induced by an operator *A* that has everywhere defined and bounded inverse, not only is the problem (solvable and) well-defined, for $\operatorname{ran} A = \mathcal{H}$ and $\ker A = \{0\}$, but also it is said to be **well-posed**, equivalently, the unique solution $f = A^{-1}g$ depends continuously (in the norm topology of $\mathcal{H}$) on the datum $g$.

Obviously when $\dim \mathcal{H} < \infty$ the finite matrix $A$ of the inverse problem (1.1) is bounded. When instead $\dim \mathcal{H} = \infty$, and $A$ is a *bounded* operator on $\mathcal{H}$ with infinite rank ($\dim \operatorname{ran} A = \infty$), infinite dimensionality brings in new features as far as the approximability of the solution(s) $f$ to (1.1) is concerned – more precisely, as we shall discuss in due time, the actual possibility of approximating $f$ through certain explicit approximants $(f_n)_{n \in \mathbb{N}}$, the topology of the convergence $f_n \to f$, the rate of convergence, etc.

Next, we shall also be interested in the class of inverse problems when $A$ is a (closed) *unbounded* operator on $\mathcal{H}$. Minimal operational requirements on an unbounded $A$ force its domain $\mathcal{D}(A)$ to be a proper, possibly dense linear subspace of $\mathcal{H}$: this introduces domain issues that are new as compared to the bounded case (where $\mathcal{D}(A) = \mathcal{H}$). For instance, when beside the approximation $f_n \to f$, and hence the smallness in some sense of the 'error' vector $f_n - f$, one investigates also the approximation $A f_n \to g$, hence the smallness of the 'residual' vector $A f_n - g$, one needs to ensure first that the approximants $f_n \in \mathcal{D}(A)$.

## 1.2 General truncation and approximation scheme

We should also like to clarify under which perspective we shall discuss in this book the issue of the *approximation* of an exact solution $f$ to the inverse problem $Af = g$ on Hilbert space $\mathcal{H}$. This is indeed a huge field in theoretical numerical analysis, approximation theory, scientific computing, and related disciplines, whereas our discussion concerns a very circumscribed scenario.

It is beneficial to recall first the general idea underlying truncation and approximation schemes. We stress "general", as we are in an abstract framework on Hilbert space, which has the virtue of encompassing a whole variety of concrete inverse problems, but has at the same time the limit of not including the specific peculiarities of this or that model, for which there are often ad hoc, tailored methods to find convenient approximate solutions.

The essential point is that in the infinite-dimensional setting scientific computing demands a suitable discretisation of the problem (1.1) to *truncated* finite-dimensional Hilbert spaces, in order to run the numerics with actual matrices. Thus, one *general* and natural approach to a genuinely infinite-dimensional inverse problem is to reduce it to a finite-dimensional auxiliary problem, under conditions guaranteeing that the solution to the latter is an approximation of the solution to the former.

We referred to genuine infinite-dimensionality to mean that not only $\dim \mathcal{H} = \infty$, but also that $A$ is *not* reduced to $A = A_1 \oplus A_2$ by an orthogonal direct sum decomposition (see, e.g., [97, Sect. 1.4]) $\mathcal{H} = \mathcal{H}_1 \oplus \mathcal{H}_2$ of $\mathcal{H}$ into two Hilbert subspaces with $\dim \mathcal{H}_1 < \infty$, $\dim \mathcal{H}_2 = \infty$, and $A_2 = \mathbb{O}$, for otherwise the inverse problem would be in practice induced by the finite matrix $A_1$.

Typically the finite-dimensional truncation is performed through what is customarily referred to as **projection method**, the essence of which is the following (we refer to Appendix A for details and further discussion). Two convenient, a priori known orthonormal systems $(u_n)_{n \in \mathbb{N}}$ and $(v_n)_{n \in \mathbb{N}}$ of $\mathcal{H}$ are taken, and for each (large enough) integer $N$ the corresponding finite-dimensional **solution space** $\mathrm{span}\{u_1, \ldots, u_N\}$ and **trial space** $\mathrm{span}\{v_1, \ldots, v_N\}$ are constructed. Denoting by $P_N$ and $Q_N$ the orthogonal projection maps onto, respectively, the solution and the trial space, one replaces the original problem $Af = g$ with its $N$-dimensional truncation

$$Q_N A P_N \widehat{f^{(N)}} = Q_N g \qquad (1.2)$$

in the unknown $\widehat{f^{(N)}}$ in the solution space. In (1.2) the operator $Q_N A P_N$ on $\mathcal{H}$, called the **compression** of $A$ between the above solution and trial space, acts non-trivially as a $N \times N$ matrix on $\mathbb{C}^N$. Thus, one searches for $\widehat{f^{(N)}} \in \mathrm{span}\{u_1, \ldots, u_N\}$, solving

$$Q_N(A\widehat{f^{(N)}} - g) = 0, \qquad (1.3)$$

an inverse problem on $\mathbb{C}^N$. For large enough $N$, the resulting $\widehat{f^{(N)}}$'s are expected to produce reasonable approximations for the 'exact' solution(s) $f$ to $Af = g$.

The same idea applies of course to a finite-dimensional inverse problem ($\dim \mathcal{H} < \infty$), by reducing it to truncated problems whose smaller size, say, the dimension $N$ of the truncation (solution and trial) space, is gradually increased: in this case, instead of controlling a limit $N \to \infty$, one performs a finite number of steps until when $N$ reaches the dimension of $\mathcal{H}$.

Surely the most popular class of truncation schemes are the celebrated **Petrov-Galerkin projection methods** [69, Chapter 4], [96, Chapter 5], [9, Chapter 9] (and one simply speaks of a *Galerkin method* when $u_n = v_n$ for all $n$). They are all based on the above general idea and differ among each other in the choice of the truncation spaces, the procedure for the computation of the approximants $\widehat{f^{(N)}}$, often iteratively in $N$, as well as in the scope of their applicability, in particular the assumptions that the problem must satisfy in order for the projection method to display good convergence features, and of course they also differ in their effectiveness at the computational level (rate of convergence, computational costs, implementation time, etc.). It is also in view of such differences that in the following we shall use various (similar) notations for the approximants, beside the symbol $\widehat{f^{(N)}}$ used above.

In this monograph, which is more operator-theoretic and functional-analytic in nature, we keep this beautiful and prolific subject in the background, with no pretension of surveying it systematically, given the already existing numerous excellent references such as [69, 83, 112, 56, 35, 24, 57, 73, 96, 37, 9, 106, 87]. In view of the materials developed in the following chapters, in retrospect we shall reflect in Appendix A on a broader perspective on projection methods, in particular on the possibility that certain operational assumptions that are crucial in Petrov-Galerkin methods may not be applicable, and that consequently the projection scheme fails to converge, or at least converges in weaker senses than the traditional ones. For instance, one prototypical occurrence that is going to be discussed at length in this book is the possibility that the orthonormal systems $(u_n)_{n\in\mathbb{N}}$ and $(v_n)_{n\in\mathbb{N}}$ chosen for the truncation are not *complete* in $\mathcal{H}$, in which case a solution $f$ to $Af = g$ that has non-trivial component on the orthogonal complement of the truncation space cannot have approximants from such space with arbitrarily good precision in the Hilbert space norm.

## 1.3 Krylov subspace and Krylov solvability

We are finally in the condition to outline the main subject of this book. We are concerned with a special and quite popular class of truncation and projection methods for the approximation of solutions to (infinite-dimensional) inverse linear problems of the form $Af = g$, and precisely those schemes where a solution $f$ is searched over a subspace of distinguished approximants in $\mathcal{H}$ explicitly constructed from $A$ and $g$ as finite linear combinations of $g, Ag, A^2g, A^3g, \ldots$.

Given a linear operator $A$ acting on the Hilbert space $\mathcal{H}$ and given a vector $g \in \mathcal{H}$, the subspace

$$\mathcal{K}(A,g) := \mathrm{span}\{g, Ag, A^2g, \ldots\} \subset \mathcal{H} \tag{1.4}$$

is called the **Krylov subspace** associated with $A$ and $g$. Here by the **(linear) span** of a collection of vectors one means the subspace of all their finite linear combinations. In other mathematical contexts $\mathcal{K}(A,g)$ is also called the **cyclic subspace** for $A$ relative to the vector $g$.

This definition is completely innocent when $\mathcal{H} \cong \mathbb{C}^N$ and $A$ is a $N \times N$ matrix, or also when $\dim \mathcal{H} = \infty$ and $A$ an everywhere defined and bounded operator on $\mathcal{H}$, because in these cases the generic vector $A^k g$ makes sense for arbitrary $k \in \mathbb{N}_0$. If instead $A$ is unbounded with a domain $\mathcal{D}(A)$ that is a proper linear subspace of $\mathcal{H}$, then $A^k g$ makes only sense if $g$ belongs to the domain of $A^k$, hence the space (1.4) is only meaningfully defined provided that $g$ is taken in the intersection of the domains of all positive integer powers of $A$. Let us defer these technical aspects to the forthcoming discussion of the unbounded case, and for the sake of this introductory presentation let us carry on here the notion of Krylov subspace with no additional specifications.

Next to (1.4) it is also convenient to introduce the **$N$-th order Krylov subspace**

$$\mathcal{K}_N(A,g) := \mathrm{span}\{g, Ag, \ldots, A^{N-1}g\}. \tag{1.5}$$

For each $N \in \mathbb{N}$, (1.5) above defines a finite-dimensional subspace of $\mathcal{H}$ with dimension $\dim \mathcal{K}_N(A,g) \leqslant N$.

Now, projection methods for an inverse problem $Af = g$ where the truncation (in the sense sketched in the previous Section) is performed with respect to solu-tion and trial spaces that are taken from the Krylov subspace $\mathcal{K}(A,g)$ are called **Krylov subspace methods**. These are schemes where, in the above notation, the $N$-th order approximant $\widehat{f^{(N)}}$ is searched for over the subspace $\mathcal{K}_N(A,g)$. To say that the algorithm converges to a solution $f$, namely $\|\widehat{f^{(N)}} - f\| \to 0$ as $N \to \infty$, means therefore that $f$ admits approximants in the form of finite linear combinations of vectors $g, Ag, A^2g, \ldots$

When this happens, namely when the inverse problem (1.1) admits a solution $f$ that is arbitrarily close in norm to a finite linear combination of $g, Ag, A^2g, \ldots$, one says that the problem is **Krylov solvable** and that $f$ is a **Krylov solution**. In other words, a Krylov solution $f$ is such that $f \in \overline{\mathcal{K}(A,g)}$, the closure in $\mathcal{H}$ of the Krylov subspace.

Krylov subspace methods, at least in their finite-dimensional realisation (i.e., when $\dim \mathcal{H} < \infty$) constitute a class of extremely popular and efficient numerical schemes, even counted among the 'Top 10 Algorithms' of the 20th century [31, 26]. In finite dimension this is by now a classical and deeply understood framework (see, e.g., the above-mentioned monographs [96, 73] or also [95]). They are naturally exported to infinite dimension, although the latter scenario is less systematically studied and is better understood through special sub-classes of interest [66, 27, 64, 81, 117, 61, 22, 20].

One of the convenient features of an approximation scheme based on Krylov subspace, although at this abstract level this cannot be visualised concretely, is precisely the "easy" construction of the spanning vectors $g, Ag, A^2g, A^3g, \ldots$ for the subspace of approximants. One should regard them as "good" vectors, determined by the sole knowledge of $g$ and $A$ – the comparison, for instance, with '*finite element methods*' extensively discussed in the references listed in the previous Section, should provide an immediate flavour of how different and more laborious the construction of the trial space may actually be.

From the theoretical viewpoint, which is in fact our primary perspective throughout, the relevance of Krylov subspace methods is connected with profound aspects of the theory of cyclic subspaces: we shall elaborate on that repeatedly over the next chapters.

*This book analyses the general question of Krylov solvability of an inverse linear problem.*

Sometimes this is a trivial issue: if $A$ and $g$ are such that the subspace $\mathcal{K}(A,g)$ is *dense* in $\mathcal{H}$, namely $\overline{\mathcal{K}(A,g)} = \mathcal{H}$, it is obvious that any solution $f$ to $Af = g$ is a Krylov solution, hence admits approximants in the form of finite linear combinations of vectors of the type $A^k$, $k \in \mathbb{N}_0$. In this case $g$ is said to be a **cyclic vector** for $A$ in $\mathcal{H}$.

In many other circumstances the "nice" choice of the "good" vectors $g, Ag, A^2g, A^3g, \ldots$ for the subspace of approximants for the considered inverse problem comes with the price that $\mathcal{K}(A,g)$ has a non-trivial orthogonal complement $\mathcal{K}(A,g)^\perp$ in $\mathcal{H}$, equivalently, that $g$ is non-cyclic for $A$. In may then happen that the solution(s) $f$ have non-trivial component in $\mathcal{K}(A,g)^\perp$, i.e., $f = f_\mathcal{K} + f_{\mathcal{K}^\perp}$ with $f_\mathcal{K} \in \overline{\mathcal{K}(A,g)}$, $f_{\mathcal{K}^\perp} \in \mathcal{K}(A,g)^\perp$, and $f_{\mathcal{K}^\perp} \neq 0$. In this case, the best one can obtain with Krylov approximants $\widehat{f^{(N)}}$ of $f$ is

$$\|\widehat{f^{(N)}} - f\|^2 \;=\; \|\widehat{f^{(N)}} - f_{\mathcal{K}}\|^2 + \|f_{\mathcal{K}^\perp}\|^2 \;\geqslant\; \|f_{\mathcal{K}^\perp}\|^2 \;>\; 0\,,$$

that is, the displacement between $\widehat{f^{(N)}}$ and $f$ cannot be made arbitrarily small in norm. This is the scenario of *lack of Krylov solvability*.

It goes without saying that having tools to decide *a priori* whether an inverse problem of the form (1.1) is or is not Krylov solvable has conceptual and practical relevance in applications, because it informs on whether a "brute force" (and possibly costly!) attempt to compute $f$ through suitable Krylov approximants is bound to fail.

To our knowledge, and in fact to our surprise when we recently started investigating these topics, the study of Krylov solvability of inverse problems, at the abstract level of formalisation described so far, and in infinite dimension, is indeed a new subject, although several related questions attracted recent and past consideration in other mathematical contexts, often quite far apart from each other. This monograph is aimed at presenting a first comprehensive picture of rigorous results and open problems – in fact more the latter than the former, so that in many respects the results that we demonstrate here only scratch the surface of a much deeper field that deserves being explored in the future.

For a given abstract inverse problem (1.1) on Hilbert space, hence assigned $g$ and $A$ on $\mathcal{H}$, natural questions are, in view of what discussed above:

- the characterisation of the Krylov subspace $\mathcal{K}(A,g)$;
- the study of its density in $\mathcal{H}$ (i.e., cyclicity or non-cyclicity of $g$ in $\mathcal{H}$);
- the characterisation of the Krylov solvability, or lack of thereof;
- how the answers to the above questions depend on the type of operator $A$;
- the role of boundedness vs unboundedness of $A$;
- the Krylov solvability of special classes of Krylov subspace methods;
- the effects on Krylov solvability from perturbations of $A$ or $g$;

and these are of course just the first, most natural ones.

In this book we present non-trivial partial answers to all the questions above and we highlight many related and important difficulties. In fact, even the identification of examples and counterexamples, on which we shall insist a lot, is often arduous, because for concrete $A$ and $g$ it is normally hard to compute precisely the subspace $\mathcal{K}(A,g)$, which in general requires a combination of algebraic, analytic, and spectral methods.

## 1.4 Structure of the book

The materials of the present monograph consist of a comprehensive update, a reworked version, and a coherent unification of several results that we have published over the last three years, some of which in collaboration with P. Novati, supplemented with additional discussion and reference to related literature.

In Chapter 2, mainly based on [22] and the references therein, the general set-up of Krylov subspace and the main results and examples concerning Krylov solvability of bounded linear inverse problems is discussed.

In Chapter 3, that stems from our recent work [20] and from the precursors surveyed therein, we export such discussion to a first, paradigmatic case of unbounded inverse problems: the case where the Krylov-based scheme is the unbounded version of the very popular conjugate gradient projection method. This has interest per se, given the important history of investigations on the infinite-dimensional, bounded, conjugate gradient method, for which we established a non-trivial generalisation, and it is also preparatory for the following analysis on the general unbounded scenario.

Indeed, in Chapter 4, mainly based on [19] and the references therein, we develop the counterpart of the study of Krylov solvability of an ample generality of unbounded inverse problems, with a special focus on the self-adjoint problem.

Then in Chapter 5, originating from our recent analysis [21], we present a first study on the effect that controlled perturbations of a given inverse problem have on its Krylov solvability, and we also outline the most relevant questions that deserve being investigated in this uncharted direction, together with a first set of meaningful examples and counterexamples.

As side material, in Appendix A, mainly based on our [23], a supplementary reflection on projection methods is discussed (of which Krylov subspaces methods are only one representative) from the point of view of general mechanisms of finite-dimensional truncation and of convergence of the approximants in weaker senses than the usual Hilbert space norm.

# Chapter 2
# Krylov solvability of bounded linear inverse problems

We start in this chapter the mathematical investigation of (genuinely) infinite-dimensional inverse linear problems in Hilbert space with respect to the associated Krylov subspace, in the scenario where the linear operator is (everywhere defined and) *bounded*.

After recalling the standard setting in infinite dimension (Section 2.1), we address the general question of the Krylov solvability, first of all through several paradigmatic examples and counter-examples.

Most importantly, we demonstrate necessary and sufficient conditions, for certain relevant classes of bounded operators, in order for the solution to be a Krylov solution (Section 2.3).

In this respect, we identify a somewhat 'intrinsic' notion associated to the operator $A$ and the datum $g$, a subspace that we call the '*Krylov intersection*', that turns out to qualify the operator-theoretic mechanism for the Krylov solvability of the problem.

For one of the most investigated study cases, namely self-adjoint bounded inverse linear problems, this mechanism takes a more explicit form, that we shall refer to as the '*Krylov reducibility*' of the operator $A$.

Last, in Section 2.6 we control the main features discussed theoretically through a series of numerical tests on inverse problems in infinite-dimensional Hilbert space, suitably truncated and analysed by increasing the size of the truncation.

## 2.1 Krylov subspace of a Hilbert space

Let us recall for convenience the main definition from Section 1.3. Throughout this chapter, $A : \mathcal{H} \to \mathcal{H}$ is an everywhere-defined bounded linear operator on a given Hilbert space $\mathcal{H}$ (on the complex field, for concreteness). For short, $A \in \mathcal{B}(\mathcal{H})$, the $C^*$-algebra of everywhere-defined bounded linear maps on $\mathcal{H}$.

To any such given $A$ and any given $g \in \mathcal{H}$ one associates the subspaces

$$\mathcal{K}_N(A,g) := \mathrm{span}\{g, Ag, \ldots, A^{N-1}g\}, \qquad N \in \mathbb{N}, \tag{2.1}$$

and

$$\mathcal{K}(A,g) := \mathrm{span}\{A^k g \mid k \in \mathbb{N}_0\}, \tag{2.2}$$

called respectively the ***N*-th order Krylov subspace** and the **Krylov subspace** (or **cyclic subspace**) relative to (or associated with) $A$ and $g$. Obviously $\mathcal{K}_N(A,g) \subset \mathcal{K}(A,g)$, with equality if $A^N g$ depends linearly on $g, Ag, \ldots, A^{N-1}g$. Such a definition is insensitive of $\dim \mathcal{H}$: as commented already, it is the infinite-dimensional setting that we are primarily focussed on. Clearly, if $\mathcal{H} \cong \mathbb{C}^d$ for some $d \in \mathbb{N}$, then $1 \leqslant \dim \mathcal{K}_N(A,g) \leqslant N$, and there always exists $N_0 \leqslant d$ such that all $N$-th order spaces $\mathcal{K}_N(A,g)$ are the same whenever $N \geqslant N_0$, and therefore $\mathcal{K}(A,g) = \mathcal{K}_{N_0}(A,g)$.

Observe that assuming $\|A\|_{\mathrm{op}} < +\infty$ and $\mathcal{D}(A) = \mathcal{H}$ makes each vector $A^k g$ well-defined.

It is also straightforward to see that, when $\dim \mathcal{K}(A,g) = \infty$, the subspace $\mathcal{K}(A,g)$ is neither open nor closed in $\mathcal{H}$, and its closure can either be a proper closed subspace of $\mathcal{H}$, or even the whole $\mathcal{H}$ itself.

In purely operator-theoretical contexts, where $\mathcal{K}(A,g)$ is rather referred to as the cyclic subspace relative to $A$ and $g$, the spanning vectors $g, Ag, A^2 g, \ldots$ are said to form the **orbit** of $g$ under $A$ and the possible density of $\mathcal{K}(A,g)$ in $\mathcal{H}$ is called the **cyclicity** of $g$, in which case $g$ is called a **cyclic vector** for $A$. An operator $A \in \mathcal{B}(\mathcal{H})$ that admits cyclic vectors in $\mathcal{H}$ is called **cyclic operator**.

For completeness of information, let us recall a few well-known facts about cyclic vectors and cyclic operators [55, Chapter 18].

(I) In non-separable Hilbert spaces there are no cyclic vectors.
(II) The set of (bounded) cyclic operators on a Hilbert space $\mathcal{H}$ is dense in $\mathcal{B}(\mathcal{H})$, with respect to the $\|\ \|_{\mathrm{op}}$-norm, if $\dim \mathcal{H} < \infty$. Instead, it is not dense in $\mathcal{B}(\mathcal{H})$ if $\dim \mathcal{H} = \infty$.
(III) The set of cyclic operators on a separable Hilbert space $\mathcal{H}$ is not closed in $\mathcal{B}(\mathcal{H})$. It is open in $\mathcal{B}(\mathcal{H})$ if $\dim \mathcal{H} < \infty$, it is not open if $\dim \mathcal{H} = \infty$.
(IV) If $\dim \mathcal{H} = \infty$ and $\mathcal{H}$ is separable, then the set of non-cyclic operators on $\mathcal{H}$ is dense in $\mathcal{B}(\mathcal{H})$ [39] (whereas, instead, the set of cyclic operators is not).
(V) It is *not known* whether there exists a bounded operator on a separable Hilbert space $\mathcal{H}$ such that *every* non-zero vector in $\mathcal{H}$ is cyclic.
(VI) The set of cyclic vectors for a bounded operator $A$ on a Hilbert space $\mathcal{H}$ is either empty or a dense subset of $\mathcal{H}$ (the Geher theorem [49]): in particular, if $g$ is a cyclic vector for $A$, then all vectors

$$g_\alpha^{(n)} := (1-\alpha A)^n g, \qquad |\alpha| \in (0, \|A\|_{\mathrm{op}}^{-1}), \qquad n \in \mathbb{N} \tag{2.3}$$

are cyclic too and for each fixed $\alpha$ the $g_\alpha^{(n)}$'s span the whole $\mathcal{H}$.
(VII) A bounded operator $A$ on the separable Hilbert space $\mathcal{H}$ is cyclic if and only if there is an orthonormal basis $(e_n)_n$ of $\mathcal{H}$ with respect to which the matrix elements $a_{ij} := \langle e_i, A e_j \rangle$ are such that $a_{ij} = 0$ for $i > j+1$ and $a_{ij} \neq 0$ for $i = j+1$ (thus, $A$ is an upper Hessenberg infinite-dimensional matrix).

Let us survey a few instructive examples that will be recurrent in this book. Further examples are presented in Section 2.3.1.

*Example 2.1* Let $\mathcal{H} = \ell^2(\mathbb{N})$, the space of square summable sequences $a \equiv (a_k)_{k \in \mathbb{N}}$ in $\mathbb{C}$ with inner product $\langle a, b \rangle_{\ell^2} = \sum_{k \in \mathbb{N}} \overline{a_k} b_k$. One customarily denotes by $(e_n)_{n \in \mathbb{N}}$ the canonical orthonormal basis of $\ell^2(\mathbb{N})$, namely with $e_n$ being the sequence consisting of all zero entries but the $n$-th one, which amounts to 1. The **right shift** operator $R$ on $\ell^2(\mathbb{N})$ is the map acting on the canonical basis as $R e_n = e_{n+1}$ $\forall n \in \mathbb{N}$ and then extended by linearity and density. In other words,

$$R = \sum_{n=1}^{\infty} |e_{n+1}\rangle\langle e_n| \tag{2.4}$$

as a strongly convergent series in the operator sense. $R$ is an isometry (i.e., it is norm-preserving) with closed range $\operatorname{ran} R = \{e_1\}^\perp$ and norm $\|R\|_{\mathrm{op}} = 1$. $R$ is non-compact, injective, and invertible on its range, with bounded inverse

$$R^{-1} : \operatorname{ran} R \to \mathcal{H}, \qquad R^{-1} = \sum_{n=1}^{\infty} |e_n\rangle\langle e_{n+1}| \tag{2.5}$$

(as a strongly convergent series in the operator sense). The adjoint of $R$ is the **left shift** operator $L = R^*$ on $\ell^2(\mathbb{N})$: $L$ acts as $R^{-1}$ on $\operatorname{ran} R$, in particular $L e_{n+1} = e_n$ $\forall n \in \mathbb{N}$, and $L e_1 = 0$. Thus,

$$LR = \mathbb{1}, \qquad RL = \mathbb{1} - |e_1\rangle\langle e_1|, \qquad \ker L = \operatorname{span}\{e_1\}. \tag{2.6}$$

$R$ and $L$ have the same spectrum $\sigma(R) = \sigma(L) = \{z \in \mathbb{C} \mid |z| \leqslant 1\}$, but $R$ has no eigenvalue, whereas the eigenvalue of $L$ form the open unit ball $\{z \in \mathbb{C} \mid |z| < 1\}$ [88, Section VI.3]. Now, for the right shift operator $R$ on $\ell^2(\mathbb{N})$ and the vector $g = e_{m+1}$ (one of the canonical basis vectors),

$$\overline{\mathcal{K}(R, e_{m+1})} = \operatorname{span}\{e_1, \ldots, e_m\}^\perp, \tag{2.7}$$

which is a proper subspace of $\ell^2(\mathbb{N})$ if $m \geqslant 1$, and instead is the whole $\ell^2(\mathbb{N})$ if $g = e_1$. The vector $e_1$ is therefore a cyclic vector for $R$ in $\mathcal{H}$. In comparison, $\mathcal{K}(L, e_m) = \operatorname{span}\{e_1, \ldots, e_m\} \cong \mathbb{C}^m$. Let us also point out for later reference that *R admits a dense of non-cyclic vectors*, for instance, the subspace $\operatorname{span}\{f_1, e_2, e_3, e_4, \ldots\}$ with $f_1 := e_1 - 2e_2$: density in $\mathcal{H}$ is obvious, as well as the non-cyclicity of $e_n$ with $n \geqslant 2$; concerning $f_1$, each $R^k f_1$, $k \in \mathbb{N}$, is orthogonal to the vector $w := (1, \frac{1}{2}, \frac{1}{4}, \frac{1}{8}, \ldots) \in \ell^2(\mathbb{N})$, meaning $w \in \mathcal{K}(R, f_1)^\perp$, whence the non-cyclicity of $f_1$. In fact, *R also admits a dense set of cyclic vectors*: indeed, $e_1$ is a cyclic vector, and a whole dense of cyclic vectors then exists on account of the above-mentioned Geher theorem.

*Example 2.2* Let $\mathcal{H} = \ell^2(\mathbb{Z})$, defined in complete analogy to $\ell^2(\mathbb{N})$ of Example 2.1, the orthonormal canonical basis being denoted again as $(e_n)_{n \in \mathbb{Z}}$. The **right shift** $R$ and the **left shift** $L$ on $\ell^2(\mathbb{Z})$ are the operators now expressed by the (strongly convergent, in operator sense) series

$$R = \sum_{n\in\mathbb{Z}} |e_{n+1}\rangle\langle e_n|, \qquad L = \sum_{n\in\mathbb{Z}} |e_n\rangle\langle e_{n+1}|. \tag{2.8}$$

One has $R^* = R^{-1} = L$ and $L^* = L^{-1} = R$, thus both $R$ and $L$ are unitary. For each canonical basis vector $e_m$,

$$\overline{\mathcal{K}(R, e_m)} = \overline{\text{span}\{e_k \,|\, k \in \mathbb{Z}, \, k \geqslant m\}}. \tag{2.9}$$

This shows that $R$ (as well as $L$) admits a *dense of non-cyclic vectors*, the canonical basis. It is also true that $R$ (as well as $L$) admits a *dense of cyclic vectors* [60, 98].

*Example 2.3* Examples 2.1-2.2 involve non-compact operators on $\ell^2$, which in fact have a compact counterpart with the same features as far as cyclicity and Krylov subspaces are concerned. Let us consider for concreteness $\mathcal{H} = \ell^2(\mathbb{Z})$ and let $\sigma \equiv (\sigma_n)_{n\in\mathbb{N}_0}$ be a bounded sequence in $\mathbb{R}^+$ with $0 < \sigma_{n+1} < \sigma_n \; \forall n \in \mathbb{N}_0$ and $\lim_{n\to\infty} \sigma_n = 0$. The **compact (weighted) right shift** on $\ell^2(\mathbb{Z})$ with weight $\sigma$ is the operator

$$R_\sigma = \sum_{n\in\mathbb{Z}} \sigma_{|n|} |e_{n+1}\rangle\langle e_n|, \tag{2.10}$$

the above series now converging in operator norm. In fact, (2.10) gives the singular value decomposition of $R_\sigma$ (see (A.11) from Section A.3). In particular, $R_\sigma e_n = \sigma_{|n|} e_{n+1} \; \forall n \in \mathbb{Z}$. $R_\sigma$ is injective and compact, has dense range, norm $\|R_\sigma\|_{\text{op}} = \sigma_0$, and adjoint

$$R_\sigma^* = \sum_{n\in\mathbb{Z}} \sigma_{|n|} |e_n\rangle\langle e_{n+1}|, \tag{2.11}$$

the **compact (weighted) left shift** on $\ell^2(\mathbb{Z})$ with weight $\sigma$. Formula (2.9) is valid also with $R_\sigma$ instead of $R$ and indeed $R_\sigma$ admits both a dense of cyclic vectors [60, 98] and a dense of non-cyclic vectors (the canonical basis, for example).

*Example 2.4* For a given bounded sequence $a \equiv (a_n)_{n\in\mathbb{N}}$ in $\mathbb{C}$, the operator

$$M^{(a)} = \sum_{n=1}^{\infty} a_n |e_n\rangle\langle e_n| \tag{2.12}$$

(as a strongly convergent series in operator sense) is called the **multiplication operator** by $a$ on $\mathcal{H} = \ell^2(\mathbb{N})$. $M^{(a)}$ has norm $\|M^{(a)}\|_{\text{op}} = \sup_n |a_n|$ and spectrum $\sigma(M^{(a)})$ given by the closure in $\mathbb{C}$ of the set $\{a_1, a_2, a_3 \dots\}$. Its adjoint is the multiplication by $\bar{a}$, and therefore $M^{(a)}$ is normal. $M^{(a)}$ is self-adjoint whenever $a$ is real and it is compact if $\lim_{n\to\infty} a_n = 0$. $M^{(a)}$ admits a dense of non-cyclic vectors (the canonical basis, consisting of the eigenvectors of $M^{(a)}$).

*Example 2.5* The **Volterra operator** on the Hilbert space $\mathcal{H} = L^2[0,1]$ is the map $V \in \mathcal{B}(\mathcal{H})$ such that

$$(Vf)(x) = \int_0^x f(y)\,\mathrm{d}y, \qquad x \in [0,1]. \tag{2.13}$$

$V$ is compact and injective with spectrum $\sigma(V) = \{0\}$ (thus, the spectral point 0 is not an eigenvalue) and norm $\|V\|_{\mathrm{op}} = \frac{2}{\pi}$. Its adjoint $V^*$ acts as

$$(V^* f)(x) = \int_x^1 f(y)\,dy, \qquad x \in [0,1], \tag{2.14}$$

therefore $V + V^*$ is the rank-one orthogonal projection

$$V + V^* = |\mathbf{1}\rangle\langle\mathbf{1}| \tag{2.15}$$

onto the function $\mathbf{1}(x) = 1$. The singular value decomposition of $V$ is

$$V = \sum_{n=0}^{\infty} \sigma_n |\psi_n\rangle\langle\varphi_n|, \qquad \begin{aligned} \sigma_n &= \tfrac{2}{(2n+1)\pi} \\ \varphi_n(x) &= \sqrt{2}\cos\tfrac{(2n+1)\pi}{2}x \\ \psi_n(x) &= \sqrt{2}\sin\tfrac{(2n+1)\pi}{2}x, \end{aligned} \tag{2.16}$$

where both $(\varphi_n)_{n\in\mathbb{N}_0}$ and $(\psi_n)_{n\in\mathbb{N}_0}$ are orthonormal bases of $L^2[0,1]$. Thus, $\operatorname{ran} V$ is dense, but strictly contained in $\mathcal{H}$: for example, $\mathbf{1} \notin \operatorname{ran} V$. In fact, $V$ is invertible on its range, but does not have (everywhere defined) bounded inverse; yet $V - z\mathbb{1}$ does, for any $z \in \mathbb{C}\setminus\{0\}$ (recall that $\sigma(V) = \{0\}$), and

$$(z\mathbb{1} - V)^{-1}\psi = z^{-1}\psi + z^{-2}\int_0^x e^{\frac{x-y}{z}}\psi(y)\,dy \qquad \forall \psi \in \mathcal{H},\ z \in \mathbb{C}\setminus\{0\}. \tag{2.17}$$

The explicit action of the powers of $V$ is

$$(V^n f)(x) = \frac{1}{(n-1)!}\int_0^x (x-y)^{n-1} f(y)\,dy, \qquad n \in \mathbb{N}. \tag{2.18}$$

Now, for $g := \mathbf{1}$ (the constant function with value 1), the functions $Vg, V^2 g, V^3 g, \dots$ are (multiples of) the polynomials $x, x^2, x^3, \dots$, therefore $\mathcal{K}(V,g)$ is the space of polynomials on $[0,1]$, which is dense in $L^2[0,1]$. In fact, more generally,

$$\overline{\mathcal{K}(V, x^k)} = L^2[0,1] \qquad \forall k \in \mathbb{N}_0. \tag{2.19}$$

Indeed, $\mathcal{K}(V, x^k)$ is spanned by the monomials $x^k, x^{k+1}, x^{k+2}\dots$, therefore if $h \in \mathcal{K}(V, x^k)^\perp$, then $0 = \int_0^1 \overline{h(x)} x^k p(x)\,dx$ for any polynomial $p$; the $L^2$-density of polynomials on $[0,1]$ implies necessarily that $x^k h = 0$, whence also $h = 0$. This proves that $\mathcal{K}(V,g)^\perp = \{0\}$, whence (2.19).

While the above examples are meant to provide a first intuition on various concrete Krylov subspaces, the reader should be warned that in general the explicit identification of $\mathcal{K}(A,g)$ for given $A$ and $g$ is a hard question!

## 2.2 Krylov reducibility and Krylov intersection

For given $A \in \mathcal{B}(\mathcal{H})$ and $g \in \mathcal{H}$, the space $\mathcal{H}$ undergoes the orthogonal decomposition

$$\mathcal{H} \;=\; \overline{\mathcal{K}(A,g)} \oplus \mathcal{K}(A,g)^\perp \qquad (2.20)$$

that we shall often refer to as the **Krylov decomposition** of $\mathcal{H}$ relative to $A$ and $g$. The corresponding Krylov space is invariant under $A$ and its orthogonal complement is invariant under $A^*$, that is,

$$A\,\mathcal{K}(A,g) \subset \mathcal{K}(A,g)\,, \qquad A^*\mathcal{K}(A,g)^\perp \subset \mathcal{K}(A,g)^\perp. \qquad (2.21)$$

The first statement is obvious and the second follows from $\langle A^*w,z\rangle = \langle w,Az\rangle = 0$ $\forall z \in \mathcal{K}(A,g)$, where $w$ is a generic vector in $\mathcal{K}(A,g)^\perp$. Owing to the evident relations

$$\begin{aligned} A\mathcal{V} \subset A\overline{\mathcal{V}} \subset \overline{A\mathcal{V}} &= \overline{A\overline{\mathcal{V}}}\,, \\ A\overline{\mathcal{V}} = \overline{A\mathcal{V}} \quad &\text{if } A^{-1} \in \mathcal{B}(\mathcal{H})\,, \end{aligned} \qquad (2.22)$$

valid for any subspace $\mathcal{V} \subset \mathcal{H}$ and any $A \in \mathcal{B}(\mathcal{H})$, (2.21) implies also

$$A\,\overline{\mathcal{K}(A,g)} \subset \overline{\mathcal{K}(A,g)}\,. \qquad (2.23)$$

A relevant occurrence for our purposes is when the operator $A$ is *reduced* by the Krylov decomposition (2.20), meaning [97, Section 1.4] that both $\overline{\mathcal{K}(A,g)}$ and $\mathcal{K}(A,g)^\perp$ are invariant under $A$. We shall refer to this occurrence as the **Krylov reducibility** of $A$ with given $g$, or also $\mathcal{K}(A,g)$**-reducibility**, to avoid ambiguity.

It follows by this definition that if $A$ is $\mathcal{K}(A,g)$-reducible, so is $A^*$, and vice versa, as one sees from the following elementary Lemma, whose proof is omitted.

**Lemma 2.1** *If $A$ is a bounded operator on a Hilbert space $\mathcal{H}$ and $\mathcal{V} \subset \mathcal{H}$ is a closed subspace, then properties (i) and (ii) below are equivalent:*

(i) $A\mathcal{V} \subset \mathcal{V}$ and $A\mathcal{V}^\perp \subset \mathcal{V}^\perp$;
(ii) $A^*\mathcal{V} \subset \mathcal{V}$ and $A^*\mathcal{V}^\perp \subset \mathcal{V}^\perp$.

*Example 2.6*

(i) For generic $A \in \mathcal{B}(\mathcal{H})$ and $g \in \mathcal{H}$, $A$ may fail to be $\mathcal{K}(A,g)$-reducible. All bounded self-adjoint operators are surely Krylov reducible, owing to (2.21).
(ii) Yet, Krylov reducibility is not a feature of self-adjoint operators only. To see this, let $A, B \in \mathcal{B}(\mathcal{H})$ and $\widetilde{A} := A \oplus B : \mathcal{H} \oplus \mathcal{H} \to \mathcal{H} \oplus \mathcal{H}$. If $g \in \mathcal{H}$ is a cyclic vector for $A$ in $\mathcal{H}$ and $\widetilde{g} := g \oplus 0$, then $\overline{\mathcal{K}(\widetilde{A},\widetilde{g})} = \mathcal{H} \oplus \{0\}$, implying that $\mathcal{K}(\widetilde{A},\widetilde{g})^\perp = \{0\} \oplus \mathcal{H}$. Therefore, $\widetilde{A}$ is Krylov reduced. On the other hand, $\widetilde{A}$ is self-adjoint (on $\mathcal{H} \oplus \mathcal{H}$) if and only if so are both $A$ and $B$ on $\mathcal{H}$.

For *normal operators* we have the following equivalent characterisation of Krylov reducibility.

**Proposition 2.1** *Let A be a bounded normal operator on a Hilbert space $\mathcal{H}$ and let $g \in \mathcal{H}$. Then A is reduced with respect to the Krylov decomposition* (2.20) *if and only if $A^*g \in \overline{\mathcal{K}(A,g)}$.*

*Proof* If $A$ is Krylov reducible, then $\overline{\mathcal{K}(A,g)}$ is invariant under $A^*$ (Lemma 2.1), hence in particular $A^*g \in \overline{\mathcal{K}(A,g)}$. Conversely, if $A^*g \in \overline{\mathcal{K}(A,g)}$, then

$$\overline{\mathcal{K}(A,A^*g)} \;=\; \overline{\mathrm{span}\{A^k A^* g \,|\, k \in \mathbb{N}_0\}} \;\subset\; \overline{\mathcal{K}(A,g)}\,,$$

and moreover, since $A$ is normal, $A^*\mathcal{K}(A,g) = \mathcal{K}(A,A^*g)$; therefore (using (2.22)),

$$A^* \overline{\mathcal{K}(A,g)} \;\subset\; \overline{\mathcal{K}(A,A^*g)} \;\subset\; \overline{\mathcal{K}(A,g)}\,.$$

The latter property, together with $A^*\mathcal{K}(A,g)^\perp \subset \mathcal{K}(A,g)^\perp$ (from (2.21) above) imply that $A^*$ is reduced with respect to the Krylov decomposition (2.20), and so is $A$ itself, owing to Lemma 2.1. $\square$

For $A \in \mathcal{B}(\mathcal{H})$ and $g \in \mathcal{H}$, an obvious consequence of $A$ being Krylov reducible is that $\overline{\mathcal{K}(A,g)} \cap (A\mathcal{K}(A,g)^\perp) = \{0\}$. For its relevance in the following, we shall call the intersection

$$\mathcal{I}(A,g) \;:=\; \overline{\mathcal{K}(A,g)} \cap (A\mathcal{K}(A,g)^\perp) \tag{2.24}$$

the **Krylov intersection** in $\mathcal{H}$ for the given $A$ and $g$.

*Example 2.7* The Krylov intersection may be trivial also in the absence of Krylov reducibility. This is already clear for finite-dimensional matrices: for example, taking (with respect to the Hilbert space $\mathbb{C}^2$)

$$A_\theta = \begin{pmatrix} 1 & \cos\theta \\ 0 & \sin\theta \end{pmatrix} \qquad \theta \in (0, \tfrac{\pi}{2}], \qquad g = \begin{pmatrix} 1 \\ 0 \end{pmatrix},$$

one sees that $A_\theta$ is Krylov reducible only when $\theta = \tfrac{\pi}{2}$, whereas the Krylov intersection (2.24) is trivial for any $\theta \in (0, \tfrac{\pi}{2})$.

## 2.3 Krylov solutions for a bounded linear inverse problem

### 2.3.1 Krylov solvability. Examples.

Let us get back to the abstract bounded linear inverse problem of the type (1.1): given $A \in \mathcal{B}(\mathcal{H})$ and $g \in \mathrm{ran}\,A$, to find solution(s) $f \in \mathcal{H}$ to

$$Af \;=\; g\,. \tag{2.25}$$

The general question we are studying here is when a solution $f$ to (2.25) admits arbitrarily close (in the norm of $\mathcal{H}$) approximants expressed by finite linear combi-nations of the spanning vectors $A^k g$'s, $k \in \mathbb{N}_0$, or equivalently when $f$ belongs to the closure $\overline{\mathcal{K}(A,g)}$ of the Krylov subspace relative to $A$ and $g$.

Let us recall from Section 1.3 that a solution $f$ satisfying the above property is referred to as a **Krylov solution**, and that when a Krylov solution exists the inverse problem (2.25) is said to be **Krylov solvable**.

*Example 2.8* The self-adjoint multiplication operator $M : L^2[1,2] \to L^2[1,2]$, $f \mapsto xf$ is (everywhere defined and) bounded and invertible on the Hilbert space $\mathcal{H} = L^2[1,2]$, with an everywhere defined bounded inverse. The unique solution to $Mf = \mathbf{1}$ is the function $f(x) = \frac{1}{x}$. Moreover, $\mathcal{K}(M,\mathbf{1})$ is the space of polynomials on $[1,2]$, hence it is *dense* in $L^2[1,2]$. Therefore, $f$ is a Krylov solution.

*Example 2.9* The multiplication operator $M_z : L^2(\Omega) \to L^2(\Omega)$, $f \mapsto zf$, acting on functions on the open disk $\Omega := \{z \in \mathbb{C} \,|\, |z-2| < 1\}$, is a normal bounded bijection on $L^2(\Omega)$, and the solution to $M_z f = g$ for given $g \in L^2(\Omega)$ is the function $f(z) = z^{-1}g(z)$. Moreover, $\mathcal{K}(M_z,g) = \{pg \,|\, p \text{ a polynomial in } z \text{ on } \Omega\}$. One can see that $f \in \overline{\mathcal{K}(M_z,g)}$ and hence the problem $M_z f = g$ is Krylov-solvable. Indeed, $\Omega \ni z \mapsto z^{-1}$ is holomorphic and hence is realised by a uniformly convergent power series (e.g., the Taylor expansion of $z^{-1}$ about $z=2$). If $(p_n)_n$ is such a sequence of polynomial approximants, then

$$\begin{aligned}\|f - p_n g\|_{L^2(\Omega)} &= \|(z^{-1} - p_n)g\|_{L^2(\Omega)} \\ &\leqslant \|z^{-1} - p_n\|_{L^\infty(\Omega)} \|g\|_{L^2(\Omega)} \xrightarrow{n \to \infty} 0.\end{aligned}$$

*Example 2.10* The left shift operator $L$ on $\ell^2(\mathbb{N}_0)$ (Example 2.1) is bounded and with $\ker L = \mathrm{span}\{e_0\}$, $\mathrm{ran}\, L = \ell^2(\mathbb{N}_0)$. A solution $f$ to $Lf = g$ with $g := \sum_{n \in \mathbb{N}_0} \frac{1}{n!} e_n$ is $f = \sum_{n \in \mathbb{N}_0} \frac{1}{n!} e_{n+1}$. Moreover, as we shall now see, $\mathcal{K}(L,g)$ is *dense* in $\ell^2(\mathbb{N}_0)$: therefore, $f$ is a Krylov solution and in fact *all* other solutions, namely the vectors $f + c e_0$ for $c \in \mathbb{C}$, are Krylov solutions too. For the density of $\mathcal{K}(L,g)$: the vector $e_0$ belongs to $\overline{\mathcal{K}(L,g)}$ because

$$\begin{aligned}\|k! L^k g - e_0\|_{\ell^2}^2 &= \|(1, \tfrac{1}{k+1}, \tfrac{1}{(k+2)(k+1)}, \cdots) - (1,0,0,\ldots)\|_{\ell^2}^2 \\ &= \sum_{n=1}^\infty \left(\frac{k!}{(n+k)!}\right)^2 \xrightarrow{k \to \infty} 0.\end{aligned}$$

As a consequence, $(0, \frac{1}{k!}, \frac{1}{(k+1)!}, \frac{1}{(k+2)!}, \cdots) = L^{k-1} g - \frac{1}{(k-1)!} e_0 \in \mathcal{K}(L,g)$, therefore the vector $e_1$ too belongs to $\overline{\mathcal{K}(L,g)}$, because

$$\left\| k! \left(L^{k-1} g - \frac{1}{(k-1)!} e_0\right) - e_1 \right\|_{\ell^2}^2 = \sum_{n=1}^\infty \left(\frac{k!}{(n+k)!}\right)^2 \xrightarrow{k \to \infty} 0.$$

Repeating inductively the above two-step argument proves that any $e_n \in \overline{\mathcal{K}(L,g)}$, whence the cyclicity of $g$.

*Example 2.11* For the right shift operator $R$ on $\ell^2(\mathbb{N})$ (Example 2.1), which is bounded, injective, and with non-dense range, the unique solution $f$ to $Rf = e_2$ is $f = e_1$. However, $f$ is not a Krylov solution, for $\mathcal{K}(R, e_2) = \overline{\text{span}\{e_2, e_3, \dots\}}$. The problem $Rf = e_2$ is not Krylov solvable.

*Example 2.12* The compact, $\sigma$-weighted right-shift operator $R_\sigma$ on $\ell^2(\mathbb{Z})$ (Example 2.3) is injective and with dense range, and the unique solution $f$ to $R_\sigma f = \sigma_1 e_2$ is $f = e_1$. However, $f$ is not a Krylov solution, for $\mathcal{K}(R_\sigma, e_2) = \overline{\text{span}\{e_2, e_3, \dots\}}$. The problem $R_\sigma f = \sigma_1 e_2$ is not Krylov solvable.

*Example 2.13* Let $A$ be a bounded injective operator on a Hilbert space $\mathcal{H}$ with cyclic vector $g \in \operatorname{ran} A$ and let $\varphi_0 \in \mathcal{H} \setminus \{0\}$. Let $f \in \mathcal{H}$ be the solution to $Af = g$. The operator $\widetilde{A} := A \oplus |\varphi_0\rangle\langle\varphi_0|$ on $\widetilde{\mathcal{H}} := \mathcal{H} \oplus \mathcal{H}$ is everywhere defined and bounded. One solution to $\widetilde{A}\widetilde{f} = \widetilde{g} := g \oplus 0$ is $\widetilde{f} = f \oplus 0$ and $\widetilde{f} \in \mathcal{H} \oplus \{0\} = \overline{\mathcal{K}(\widetilde{A}, \widetilde{g})}$. Another solution is $\widetilde{f}_\xi = f \oplus \xi$, where $\xi \in \mathcal{H} \setminus \{0\}$ and $\xi \perp \varphi_0$. Clearly, $\widetilde{f}_\xi \notin \overline{\mathcal{K}(\widetilde{A}, \widetilde{g})}$.

*Example 2.14* If $V$ is the Volterra operator on $L^2[0,1]$ (Example 2.5) and $g(x) = \frac{1}{2}x^2$, then $f(x) = x$ is the unique solution to $Vf = g$. As

$$\mathcal{K}(V, g) \;=\; \{x^2 p(x) \,|\, p \text{ is a polynomial on } [0,1]\},$$

$f \notin \mathcal{K}(V, g)$, because $f(x) = x^2 \cdot \frac{1}{x}$ and $\frac{1}{x} \notin L^2[0,1]$. Yet, $f \in \overline{\mathcal{K}(V, g)} = L^2[0,1]$, as seen in (2.19) above.

### 2.3.2 General conditions for Krylov solvability: case of injectivity

Even stringent assumptions on $A$ such as the simultaneous occurrence of compactness, injectivity, and density of the range do *not* ensure, in general, that the solution $f$ to $Af = g$, for given $g \in \operatorname{ran} A$, is a Krylov solution (Example 2.12).

Let us consider first the case when $A$ is injective (i.e., when the problem (2.25) is well-defined, according to the nomenclature of Section 1.1.

A *necessary* condition for the solution to a well-defined bounded linear inverse problem to be a Krylov solution, which becomes necessary and sufficient if the linear map is a bounded bijection, is the following. (Recall that for $A \in \mathcal{B}(\mathcal{H})$ these three properties are *equivalent*: $A$ is a bijection; $A$ is invertible with everywhere defined, bounded inverse on $\mathcal{H}$; the spectral point 0 belongs to the resolvent set of $A$.)

**Proposition 2.2** *Given a Hilbert space $\mathcal{H}$, let $A \in \mathcal{B}(\mathcal{H})$ be injective, $g \in \operatorname{ran} A$, and $f \in \mathcal{H}$ be the unique solution to $Af = g$.*

(i) *If $f \in \overline{\mathcal{K}(A, g)}$, then $A\overline{\mathcal{K}(A, g)}$ is dense in $\overline{\mathcal{K}(A, g)}$.*
(ii) *Assume further that $A$ is invertible with everywhere defined, bounded inverse on $\mathcal{H}$. Then $f \in \overline{\mathcal{K}(A, g)}$ if and only if $A\overline{\mathcal{K}(A, g)}$ is dense in $\overline{\mathcal{K}(A, g)}$.*

***Proof*** One has $A\overline{\mathcal{K}(A,g)} \supset A\mathcal{K}(A,g) = \text{span}\{A^k g \,|\, k \in \mathbb{N}\}$, owing to the definition of Krylov space and to (2.22). If $f \in \overline{\mathcal{K}(A,g)}$, then $A\overline{\mathcal{K}(A,g)} \ni Af = g$, in which case $A\overline{\mathcal{K}(A,g)} \supset \text{span}\{A^k g \,|\, k \in \mathbb{N}_0\}$; the latter inclusion, by means of (2.22) and (2.23), implies $\overline{\mathcal{K}(A,g)} \supset \overline{A\overline{\mathcal{K}(A,g)}} \supset \overline{\mathcal{K}(A,g)}$, whence $\overline{A\overline{\mathcal{K}(A,g)}} = \overline{\mathcal{K}(A,g)}$. This proves part (i) and the 'only if' implication in part (ii). For the converse, let us now assume that $A^{-1} \in \mathcal{B}(\mathcal{H})$ and that $A\overline{\mathcal{K}(A,g)}$ is dense in $\overline{\mathcal{K}(A,g)}$. Let $(Av_n)_{n \in \mathbb{N}}$ be a sequence in $A\overline{\mathcal{K}(A,g)}$ of approximants of $g \in \overline{\mathcal{K}(A,g)}$ for some $v_n$'s in $\overline{\mathcal{K}(A,g)}$. Since $A^{-1}$ is bounded on $\mathcal{H}$, then $(v_n)_{n \in \mathbb{N}}$ is a Cauchy sequence, hence $v_n \to v$ as $n \to \infty$ for some $v \in \overline{\mathcal{K}(A,g)}$. By continuity, $Af = g = \lim_{n \to \infty} Av_n = Av$, and by injectivity $f = v \in \overline{\mathcal{K}(A,g)}$, which proves also the 'if' implication of part (ii). □

On the other hand, a *sufficient* condition for the Krylov solvability of a well-defined bounded linear inverse problem is the triviality of the Krylov intersection; it actually becomes a necessary and sufficient if $A$ is a bounded bijection.

**Proposition 2.3** *Given a Hilbert space $\mathcal{H}$, let $A \in \mathcal{B}(\mathcal{H})$ be injective, $g \in \text{ran} A$, and $f \in \mathcal{H}$ be the unique solution to $Af = g$.*

*(i) If the Krylov intersection is trivial, namely $\mathcal{I}(A,g) = \overline{\mathcal{K}(A,g)} \cap (A\mathcal{K}(A,g)^\perp) = \{0\}$, then $f \in \overline{\mathcal{K}(A,g)}$.*

*(ii) Assume further that $A$ is invertible with everywhere defined, bounded inverse on $\mathcal{H}$. Then $f \in \overline{\mathcal{K}(A,g)}$ if and only if $\mathcal{I}(A,g) = \overline{\mathcal{K}(A,g)} \cap (A\mathcal{K}(A,g)^\perp) = \{0\}$.*

***Proof*** (i) Let $P_{\mathcal{K}} : \mathcal{H} \to \mathcal{H}$ be the orthogonal projection onto $\overline{\mathcal{K}(A,g)}$. On the one hand, $A(\mathbb{1} - P_{\mathcal{K}})f \in \overline{\mathcal{K}(A,g)}$, because $AP_{\mathcal{K}}f \in \overline{\mathcal{K}(A,g)}$ (from (2.23) above) and $Af = g \in \overline{\mathcal{K}(A,g)}$. On the other hand, obviously, $A(\mathbb{1} - P_{\mathcal{K}})f \in A\mathcal{K}(A,g)^\perp$. As the Krylov intersection is trivial, necessarily $A(\mathbb{1} - P_{\mathcal{K}})f = 0$, and by injectivity $f = P_{\mathcal{K}}f \in \overline{\mathcal{K}(A,g)}$.

(ii) If $f \in \overline{\mathcal{K}(A,g)}$ and $A^{-1} \in \mathcal{B}(\mathcal{H})$, then $A\overline{\mathcal{K}(A,g)}$ is dense in $\overline{\mathcal{K}(A,g)}$ (Proposition 2.2(ii)). Let now $z \in \overline{\mathcal{K}(A,g)} \cap (A\mathcal{K}(A,g)^\perp)$, say, $z = Aw$ for some $w \in \mathcal{K}(A,g)^\perp$, and based on the density proved above let $(Ax_n)_{n \in \mathbb{N}}$ be a sequence in $A\overline{\mathcal{K}(A,g)}$ of approximants of $z$ for some $x_n$'s in $\overline{\mathcal{K}(A,g)}$. From $Ax_n \to z = Aw$ and $\|A^{-1}\|_{\text{op}} < +\infty$ one has $x_n \to w$ as $n \to \infty$. Since $x_n \perp w$, then

$$0 = \lim_{n \to \infty} \|x_n - w\|_{\mathcal{H}}^2 = \lim_{n \to \infty} \left(\|x_n\|_{\mathcal{H}}^2 + \|w\|_{\mathcal{H}}^2\right) = 2\|w\|_{\mathcal{H}}^2,$$

whence necessarily $w = 0$ and $z = 0$. This proves the 'only if' implication of (ii). □

**Corollary 2.1** *Well-defined self-adjoint inverse problems are always Krylov solvable, i.e., if $A$ is bounded, injective, and self-adjoint, and $g \in \text{ran} A$, then $Af = g$ implies $f \in \overline{\mathcal{K}(A,g)}$.*

***Proof*** Being self-adjoint, $A$ is $\mathcal{K}(A,g)$-reducible (Example 2.6(i)), whence $\overline{\mathcal{K}(A,g)} \cap (A\mathcal{K}(A,g)^\perp) = \{0\}$. Then one applies Proposition 2.3. □

Propositions 2.2(ii) and 2.3(ii) provide equivalent conditions to the Krylov solvability of linear inverse problems on Hilbert space when the linear maps are bounded bijections. (As such, these results do not apply to compact operators on infinite-dimensional Hilbert spaces.)

In particular, Proposition 2.3(ii) shows that *for such linear inverse problems the Krylov solvability is tantamount as the triviality of the Krylov intersection*, which was the actual reason to introduce the subspace (2.24).

### 2.3.3 Krylov reducibility and Krylov solvability

Concerning the relation between the Krylov reducibility and the Krylov solvability, we know that when $A$ is injective the former implies the latter (Proposition 2.3(i)).

Moreover, there are classes of operators in $\mathcal{B}(\mathcal{H})$ for which the two notions coincide, as the following remark shows.

*Remark 2.1* For *unitary operators*, the Krylov solvability of the associated inverse problem is *equivalent* to the Krylov reducibility. The fact that the latter implies the former is the general property of Proposition 2.3. Conversely, when $U : \mathcal{H} \to \mathcal{H}$ is unitary and $f = U^*g$ is the solution to $Uf = g$ for some $g \in \mathcal{H}$, then the assumption $f \in \overline{\mathcal{K}(U,g)}$ implies $U^*g \in \overline{\mathcal{K}(U,g)}$, which by Proposition 2.1 is the same as the fact that $U$ is reduced with respect to $\mathcal{H} = \overline{\mathcal{K}(U,g)} \oplus \mathcal{K}(U,g)^\perp$.

There are also Krylov solvable problems that are not Krylov reduced, even among well-defined inverse linear problem, namely when $A$ is bounded and injective and $g \in \operatorname{ran} A$. This is the case of Example 2.7 when $\theta \neq \frac{\pi}{2}$.

Even in the relevant class of (bounded, injective) *normal* operators (the operator $A_\theta$ of Example 2.7 is *not* normal), Krylov solvability does not necessarily imply Krylov reducibility. Let us discuss a counterexample that elucidates that.

We need first to recall the following fact from complex functional analysis.

**Proposition 2.4** *Let $\mathcal{U} \subset \mathbb{C}$ be an open subset of the complex plane, and denote by $H(\mathcal{U})$ the space of holomorphic functions on $\mathcal{U}$. Then the space $H(\mathcal{U}) \cap L^2(\mathcal{U})$ is closed in $L^2(\mathcal{U})$.*

Proposition 2.4 is a consequence of the well-known property that for a sequence of holomorphic functions on a compact set convergence in $L^2$ implies uniform convergence, combined with the other standard property [91, Theorem 10.28] that uniform convergence of such a sequence on all compact subsets of $\mathcal{U}$ makes the limit function holomorphic on the whole $\mathcal{U}$. Details are worked out, for instance, in [18, Theorem 4.4.10].

*Example 2.15* Consider the multiplication operator from Example 2.9, namely $M_z : L^2(\Omega) \to L^2(\Omega)$, $M_z f = z f$, with $\Omega = \{z \in \mathbb{C} \,|\, |z - 2| < 1\}$, and let $g \in L^2(\Omega)$ be such that $\varepsilon \leqslant |g(z)| \leqslant \varepsilon' \ \forall z \in \Omega$, for given $\varepsilon, \varepsilon' > 0$. Then:

(i) $M_z$ is bounded, injective, and normal;
(ii) the inverse linear problem $M_z f = g$ is Krylov-solvable: $f \in \overline{K(M_z, g)}$;

(iii) however, $M_z$ is not Krylov-reduced.

Parts (i) and (ii) were discussed in Example 2.9. It was also observed therein that $\mathcal{K}(M_z,g) = \{pg \,|\, p \in \mathbb{P}_\Omega[z]\}$, where $\mathbb{P}_\Omega[z]$ denotes the polynomials in $z \in \Omega$ with complex coefficients. Let us show that

$$\overline{\mathcal{K}(M_z,g)} \;=\; \left\{ \phi g \,\Big|\, \phi \in \overline{\mathbb{P}_\Omega[z]}^{\|\,\|_2} \right\}. \tag{*}$$

Indeed, if $w \in \overline{\mathcal{K}(M_z,g)}$, then $w \xleftarrow{\|\,\|_2} p_n g$ for a sequence $(p_n)_{n\in\mathbb{N}}$ in $\mathbb{P}_\Omega[z]$, and since

$$\|p_n - p_m\|_{L^2(\Omega)} \;\leqslant\; \varepsilon^{-1} \|p_n g - p_m g\|_{L^2(\Omega)}$$

then $(p_n)_{n\in\mathbb{N}}$ is a Cauchy sequence in $L^2(\Omega)$ with $p_n \xrightarrow{\|\,\|_2} \phi$ for some $\phi \in \overline{\mathbb{P}_\Omega[z]}^{\|\,\|_2}$; whence $w = \phi g$. Conversely, if $w = \phi g$ for $\phi \in \overline{\mathbb{P}_\Omega[z]}^{\|\,\|_2}$, then $\phi \xleftarrow{\|\,\|_2} p_n$ for a sequence $(p_n)_{n\in\mathbb{N}}$ of approximants in $\mathbb{P}_\Omega[z]$ and

$$\|w - p_n g\|_{L^2(\Omega)} \;=\; \|\phi g - p_n g\|_{L^2(\Omega)} \;\leqslant\; \varepsilon' \|\phi - p_n\|_{L^2(\Omega)} \xrightarrow{n\to\infty} 0$$

shows that $w \in \overline{\mathcal{K}(M_z,g)}$. The identity (*) is therefore established. Now, if for contradiction $M_z$ was reduced with respect to the decomposition $L^2(\Omega) = \overline{\mathcal{K}(M_z,g)} \oplus \mathcal{K}(M_z,g)^\perp$, then $\bar{z}g = M_{\bar{z}}g = M_z^* g \in \overline{\mathcal{K}(M_z,g)}$ (Proposition 2.1), and the identity (*) would imply that the function $\Omega \ni z \mapsto \bar{z}$ belongs to $\overline{\mathbb{P}_\Omega[z]}^{\|\,\|_2}$; however, the latter space, owing to Proposition 2.4, is formed by holomorphic functions, and the function $z \mapsto \bar{z}$ clearly is not. Part (iii) is thus proved.

### 2.3.4 More on Krylov solutions in the lack of injectivity

To conclude this Section let us consider more generally inverse problems of the form (2.25) that are solvable ($g \in \mathrm{ran}\,A$) but with possibly multiple solutions (i.e., $A$ is possibly non-injective).

The triviality of the Krylov intersection (and in particular Krylov reducibility) still guarantees the existence of Krylov solutions: indeed, Proposition 2.3 has this immediate counterpart valid also in the lack of injectivity.

**Proposition 2.5** *Given a Hilbert space $\mathcal{H}$, let $A \in \mathcal{B}(\mathcal{H})$, $g \in \mathrm{ran}\,A$, and $f \in \mathcal{H}$ with $Af = g$. If $\mathcal{I}(A,g) = \overline{\mathcal{K}(A,g)} \cap (A\mathcal{K}(A,g)^\perp) = \{0\}$, then there exists a Krylov solution to the problem $Af = g$. For example, if $f_\circ \in \mathcal{H}$ satisfies $Af_\circ = g$ and $P_{\mathcal{K}}$ is the orthogonal projection onto $\overline{\mathcal{K}(A,g)}$, then $f := P_{\mathcal{K}} f_\circ$ is a Krylov solution.*

*Proof* One has $A(\mathbb{1} - P_{\mathcal{K}})f_\circ = 0$, owing to the very same argument as in the above proof of Proposition 2.3. Thus, $AP_{\mathcal{K}} f_\circ = g$, that is, $f := P_{\mathcal{K}} f_\circ$ is a Krylov solution. □

Generic bounded inverse problems may or may not admit a Krylov solution (Example 2.8-2.14), and when they do there may exist further non-Krylov solutions (Example 2.13), and even cases when the whole infinity of solutions consist of Krylov solutions (Example 2.10). For a fairly general class of such problems, however, the Krylov solution, when it exists, is *unique*.

**Proposition 2.6** *Given a Hilbert space $\mathcal{H}$, let $A \in \mathcal{B}(\mathcal{H})$ be a normal operator and $g \in \operatorname{ran} A$. There exists at most one $f \in \mathcal{H}$ such that $Af = g$ and $f \in \overline{\mathcal{K}(A,g)}$. More generally, the same conclusion holds if $A$ is bounded with $\ker A \subset \ker A^*$.*

*Proof* If $f_1, f_2 \in \overline{\mathcal{K}(A,g)}$ and $Af_1 = g = Af_2$, then $f_1 - f_2 \in \ker A \cap \overline{\mathcal{K}(A,g)}$. By normality, $\ker A = \ker A^*$, and moreover obviously $\overline{\mathcal{K}(A,g)} \subset \overline{\operatorname{ran} A}$. Therefore, $f_1 - f_2 \in \ker A^* \cap \overline{\operatorname{ran} A}$. But $\ker A^* \cap \overline{\operatorname{ran} A} = \{0\}$, whence $f_1 = f_2$. The second statement is then obvious. □

The above result is similar to comments made in [41, 14, 47] about Krylov solutions to singular systems in finite dimensions.

The general tool of Proposition 2.5 requires the identification of the Krylov intersection $\mathcal{I}(A,g) = \overline{\mathcal{K}(A,g)} \cap (A\mathcal{K}(A,g)^\perp)$, a task that similarly to the identification of the Krylov subspace $\mathcal{K}(A,g)$ itself may be practically hard. Certain classes of operators, such as the self-adjoint operators (Section 2.4), or the $\mathcal{K}$-class operators (Section 2.5), and possibly others, display features that make $\mathcal{I}(A,g)$ easily controllable. Otherwise, deciding whether a bounded inverse problem admits Krylov solution(s) seems to require a case-by-case ad hoc approach.

*Example 2.16* For given $k \in L^2[0,1]$, it is standard that the operator

$$A : L^2[0,1] \to L^2[0,1]$$
$$(Au)(x) := \int_0^1 k(x-y)u(y)\,\mathrm{d}y \qquad (2.26)$$

is a Hilbert-Schmidt normal operator [88, Section VI.6], with integral kernel $k$ and norm $\|A\|_{\mathrm{op}} \leqslant \|k\|_{L^2}$. Actually,

$$(A^*u)(x) = \int_0^1 \overline{k(y-x)}u(y)\,\mathrm{d}y, \qquad (2.27)$$

hence $A$ is self-adjoint if and only if $k(x) = \overline{k(-x)}$ for almost every $x \in [0,1]$. With the concrete choice

$$k(x) := \frac{e \sin \pi x}{(1+e)\pi e^x} - \frac{1}{1+\pi^2} \qquad (2.28)$$

$A$ is not self-adjoint, and the Krylov solvability of the inverse problem $Af = g$ for given $g \in \operatorname{ran} A$ must be investigated through an ad hoc analysis. Let us introduce the orthonormal basis $\{\varphi_n \,|\, n \in \mathbb{Z}\}$ of $L^2[0,1]$, where

$$\varphi_n(x) = e^{2\pi i n x}. \qquad (2.29)$$

Then

$$k = \sum_{n \in \mathbb{Z}} c_n \varphi_n, \qquad c_n := \langle \varphi_n, k \rangle_{L^2}, \qquad (2.30)$$

and a straightforward explicit computation yields

$$c_n = \begin{cases} 0 & \text{if } n = 0 \\ \dfrac{1}{1 + 4\mathrm{i}n\pi + (1 - 4n^2)\pi^2} & \text{if } n \in \mathbb{Z} \setminus \{0\}. \end{cases} \qquad (2.31)$$

As a consequence,

$$\begin{aligned}(Au)(x) &= \int_0^1 k(x-y)u(y)\,\mathrm{d}y = \sum_{n \in \mathbb{Z}} c_n \int_0^1 \varphi_n(x-y)u(y)\,\mathrm{d}y \\ &= \sum_{n \in \mathbb{Z}} c_n\, \varphi_n(x) \int_0^1 \overline{\varphi_n(y)}\, u(y)\,\mathrm{d}y, \end{aligned} \qquad (2.32)$$

that is,

$$A = \sum_{n \in \mathbb{Z}} c_n |\varphi_n\rangle\langle\varphi_n|. \qquad (2.33)$$

In practice, up to the Hilbert space isomorphism $\varphi_n \mapsto e_n$ that maps the orthonormal basis $\{\varphi_n \,|\, n \in \mathbb{Z}\}$ of $L^2[0,1]$ into the canonical basis of $\ell^2(\mathbb{Z})$, $A$ is unitarily equivalent to the multiplication operator $M^{(c)}$ considered in Example 2.4, with $c \equiv (c_n)_{n \in \mathbb{Z}}$. We can thus draw a number of conclusions.

- $A$ is not injective: $\ker A = \mathrm{span}\{\varphi_0\}$.
- If $g \in \mathrm{ran}\,A$ (that is, if $g$ is not a constant function), and $J \subset \mathbb{Z} \setminus \{0\}$ is the subset of non-zero integers $n$ such that $g_n := \langle \varphi_n, g \rangle_{L^2} \neq 0$, then

$$g = \sum_{n \in J} g_n\, \varphi_n$$

and the inverse linear problem $Af = g$ admits an infinity of solutions of the form $f = \alpha \varphi_0 + f_{\mathcal{K}}$ for arbitrary $\alpha \in \mathbb{C}$, where

$$f_{\mathcal{K}} := \sum_{n \in J} \frac{g_n}{c_n}\, \varphi_n.$$

- The vectors $g, Ag, A^2 g, \ldots$ have non-zero components only of order $n \in J$; this, together with the fact that the $c_n$'s are all distinct, implies that

$$\mathcal{K}(A, g) \subset \mathrm{span}\{\varphi_n \,|\, n \in J\}.$$

- The functions $f = \alpha \varphi_0 + f_{\mathcal{K}}$ with $\alpha \in \mathbb{C} \setminus \{0\}$ are non-Krylov solutions to the problem $Af = g$, whereas only $f_{\mathcal{K}}$ can be a Krylov solution, consistently with Proposition 2.6.

## 2.4 Krylov solvability and self-adjointness

Propositions 2.5 and 2.6 yield a noticeable behaviour for the class of self-adjoint inverse problems.

**Proposition 2.7** *If $A \in \mathcal{B}(\mathcal{H})$ is self-adjoint, then the inverse problem $Af = g$ with $g \in \mathrm{ran}\, A$ admits a unique Krylov solution.*

*Proof* $A$ is Krylov reducible (Example (2.6)(i)), hence the induced inverse problem admits a Krylov solution (Proposition 2.5). Such a solution is then necessarily unique (Proposition 2.6). □

It is worth noticing that the self-adjoint case has always deserved a special status in this context, theoretically and in applications: the convergence of Krylov techniques for self-adjoint operators are the object of an ample literature – see, e.g., [66, 27, 64, 117, 81, 61, 84].

When in particular $A \in \mathcal{B}(\mathcal{H})$ is self-adjoint and *positive definite*, a subtle alternative analysis by Nemirovskiy and Polyak [81] (for a more recent discussion of which we refer to [35, Section 7.2] and [56, Section 3.2]) proved that the associated inverse problem $Af = g$ with $g \in \mathrm{ran}\, A$ is actually Krylov solvable. In particular it was proved that the sequence of Krylov approximations from the conjugate gradient algorithm converges strongly to the exact solution. We shall discuss that approach in Chapter 3, in connection to its non-trivial generalisation to unbounded operators.

*Example 2.17* The test problems

<div align="center">

blur    deriv2   foxgood  gravity  heat
i_laplace  parallax  phillips   shaw    ursell

</div>

of Hansen's REGULARIZATION TOOLS Matlab package [58] correspond to integral operators $A_K$ on some $L^2[a,b]$ whose integral kernels $K(x,y)$ are square-integrable and have the property $K(x,y) = \overline{K(y,x)}$, namely they are Hilbert-Schmidt and self-adjoint operators. Owing to Proposition 2.7, all such inverse problems admit a unique Krylov solution (in fact they are Krylov solvable, as long as $g \in \mathrm{ran}\, A_K$ and $A_K$ is injective).

*Example 2.18* The `PRdiffusion` two-dimensional test problem of Gazzola-Hansen-Nagy's IR Tools [48] consists of reconstructing, from the heat diffusion problem

$$\begin{cases} \dfrac{\partial u}{\partial t} = \Delta_N u \\ u(0) = u_0 \end{cases}$$

with unknown $u \equiv u(x,y;t)$ in the Hilbert space $L^1([0,1] \times [0,1], \mathrm{d}x\mathrm{d}y)$ and with Neumann Laplacian $\Delta_N$, the initial datum $u_0$ starting from the knowledge of the function $u(t)$ at time $t > 0$. By standard functional-analytic arguments one has $u_0 = e^{-t\Delta_N} u(t)$, that is, for given $t$ the inverse problem $u(t) \mapsto u_0$ is self-adjoint and hence (Proposition 2.7) it admits a unique Krylov solution.

A useful abstract tool for the study of Krylov solvability of inverse problems induced by bounded self-adjoint operators is provided by the following isomorphism. In Chapter 4 (Theorem 4.3) we will lift this result to the *unbounded* self-adjoint case.

**Proposition 2.8** *Consider a Hilbert space $\mathcal{H}$, a self-adjoint operator $A \in \mathcal{B}(\mathcal{H})$, and a vector $g \in \mathcal{H}$. Denote by $E^{(A)}$ the spectral measure for A, and by $\mu_g^{(A)}$ the scalar spectral measure associated with A and g, i.e., $\mu_g^{(A)}(\Omega) = \langle g, E^{(A)}(\Omega)g \rangle$ for every Borel subset $\Omega \subset \mathbb{R}$. There is the Hilbert space isomorphism*

$$L^2(\sigma(A), \mathrm{d}\mu_g^{(A)}) \overset{\cong}{\longrightarrow} \overline{\mathcal{K}(A,g)} \tag{2.34}$$
$$\varphi \longmapsto \varphi(A)g$$

*where $\varphi(A)$ is the operator constructed with the bounded functional calculus of A.*

A standard reference for the spectral measure $E^{(A)}$ is [97, Section 5.2] and for the functional calculus of $A$ is [97, Section 5.3].

***Proof (Proof of Proposition 2.8)*** Since $A$ is bounded and self-adjoint, the spectrum $\sigma(A)$ is a compact subset of $\mathbb{R}$ (in fact, $\sigma(A) \subset [-\|A\|_{\mathrm{op}}, \|A\|_{\mathrm{op}}]$). Therefore, the space $\mathbb{P}_{\sigma(A)}[\lambda]$ of complex coefficient polynomials in the variable $\lambda \in \sigma(A)$ is dense in the space $C(\sigma(A), \mathbb{C})$ of complex-valued continuous functions on $\sigma(A)$ (the Stone-Weierstrass theorem, [40, Section 4.7], [88, Theorem IV.10]). In turn, it is well known that $C(\sigma(A), \mathbb{C})$ is dense in $L^2(\sigma(A), \mathrm{d}\mu_g^{(A)})$ [91, Theorem 3.14], and since on the compact $\sigma(A)$ the uniform norm is stronger then the $L^2$-norm, then $\mathbb{P}_{\sigma(A)}[\lambda]$ is dense in $L^2(\sigma(A), \mathrm{d}\mu_g^{(A)})$. Now, for $p \in \mathbb{P}_{\sigma(A)}[\lambda]$, the map $p \mapsto p(A)g$ defines an isometric linear bijection $\mathbb{P}_{\sigma(A)}[\lambda] \to \mathcal{K}(A,g)$: linearity and bijectivity are obvious, and norm preservation follows from

$$\|p(A)g\|_{\mathcal{H}}^2 = \int_{\sigma(A)} |p(\lambda)|^2 \, \mathrm{d}\mu_g^{(A)}(\lambda) = \|p\|_{L^2(\sigma(A), \mathrm{d}\mu_g^{(A)})}^2.$$

Exploiting the above-mentioned $L^2$-density of polynomials, by standard linear bounded extension [88, Theorem I.7] such a map admits a unique bounded extension $L^2(\sigma(A), \mathrm{d}\mu_g^{(A)}) \to \overline{\mathcal{K}(A,g)}$, $\varphi \mapsto \varphi(A)g$. It remains to prove that the extended map is also surjective. Let then $v \in \overline{\mathcal{K}(A,g)}$; there is a sequence $(p_n)_{n \in \mathbb{N}}$ in $\mathbb{P}_{\sigma(A)}[\lambda]$ such that $\mathcal{K}(A,g) \ni p_n(A)g \to v$ in $\mathcal{H}$ as $n \to \infty$. Owing to the above norm preservation, $(p_n)_{n \in \mathbb{N}}$ is a Cauchy sequence in $L^2(\sigma(A), \mathrm{d}\mu_g^{(A)})$ and hence $p_n \to \varphi$, as a $L^2$-limit, for some $\varphi \in L^2(\sigma(A), \mathrm{d}\mu_g^{(A)})$, whence $p_n(A)g \to \varphi(A)g$ in $\mathcal{H}$. Then necessarily $\varphi(A)g = v$. $\square$

As a first application of Proposition 2.8 to the issue of Krylov solvability, let us re-prove Corollary 2.1, now within a completely different conceptual scheme.

***Proof (Alternative proof of Corollary 2.1)*** Since $A$ is injective, it is invertible on $\mathrm{ran}\,A$ and the solution $f$ to $Af = g$, for given $g \in \mathrm{ran}\,A$, is the vector $f = A^{-1}g$. This means that the function $\varphi(\lambda) = \lambda^{-1}$ belongs to $L^2(\sigma(A), \mathrm{d}\mu_g^{(A)})$, indeed

$$\|\varphi\|^2_{L^2(\sigma(A),\mathrm{d}\mu_g^{(A)})} \;=\; \int_{\sigma(A)} |\lambda|^{-2} \mathrm{d}\mu_g^{(A)}) \;=\; \|A^{-1}g\|^2_{\mathcal{H}} \;=\; \|f\|^2_{\mathcal{H}} \;<\; +\infty.$$

As such, one concludes from Proposition 2.8 that $f = \varphi(A)g \in \overline{\mathcal{K}(A,g)}$. □

## 2.5 Special classes of Krylov solvable problems

In the current lack of a satisfactory characterisation of Krylov solvable inverse problems on infinite-dimensional Hilbert space, it is of interest to identify special subclasses of them.

We already examined simple explicit cases in Examples 2.8, 2.9, 2.10, 2.14.

We also concluded (Proposition 2.7) that a whole class of paramount relevance, the bounded self-adjoint operators, induce inverse problems that are Krylov solvable, and Examples 2.17 and 2.18 report on several applications.

Outside the self-adjoint class, there are various types of Krylov solvable problems, in particular whenever the datum $g$ is cyclic for the operator $A$, as is the case for suitable choices of $g$ with shift operators (Examples 2.1-2.3), multiplication operators (Example 2.4), Volterra operators (Example 2.5), and the like.

Furthermore, operators with the structure of $A_\theta$ from Example 2.7 are non-normal, yet give rise to non-trivial Krylov intersection, thereby inducing Krylov solvable problems (Propositions 2.3 and 2.5). On the other hand, in general Krylov solvability is not automatic for compact injective operators (Example 2.12).

All these are particular cases of a rather diverse scenario that deserves being studied systematically.

In this spirit, let us present one further class of well posed inverse linear problems that are Krylov solvable (Corollary 2.2 below).

Given a Hilbert space $\mathcal{H}$, we say that an operator $A \in \mathcal{B}(\mathcal{H})$ is of (or: belongs to the) $\mathcal{K}$-**class** when there exists an open subset $\mathcal{W} \subset \mathbb{C}$ such that

- $\sigma(A) \subset \mathcal{W}$,
- $\overline{\mathcal{W}}$ is compact with $0 \notin \overline{\mathcal{W}}$,
- and $\mathbb{C} \setminus \overline{\mathcal{W}}$ is connected in $\mathbb{C}$.

In particular, if $A$ is of $\mathcal{K}$-class, then $0 \notin \sigma(A)$.

*Example 2.19* The $\mathcal{K}$-class condition is not satisfied by the right shift on $\ell^2$ (Examples 2.1-2.2), the compact multiplication operator on $\ell^2$ (Example 2.4), the Volterra operator (Example 2.5). The multiplication of Example 2.9 is a $\mathcal{K}$-class opeator. The identity operator $\mathbb{1}$ is obviously of $\mathcal{K}$-class, but in general unitary or normal operators are not. The multiplication operator from Example 2.8 is of $\mathcal{K}$-class.

$\mathcal{K}$-class operators have a polynomial approximation of their inverse, which eventually yields Krylov solvability of the associated inverse problem.

**Proposition 2.9** *Let $A$ be an operator of $\mathscr{K}$-class on a Hilbert space $\mathcal{H}$. Then there exists a polynomial sequence $(p_n)_{n\in\mathbb{N}}$ over $\mathbb{C}$ such that $\|p_n(A) - A^{-1}\|_{\mathrm{op}} \to 0$ as $n \to \infty$.*

*Proof* Let $\mathcal{U} \subset \mathbb{C}$ be an open set such that $0 \notin \mathcal{U}$ and $\overline{\mathcal{W}} \subset \mathcal{U}$, where $\mathcal{W}$ is an open set fulfilling the definition of $\mathscr{K}$-class for the given $A$. The function $z \mapsto z^{-1}$ is holomorphic on $\mathcal{U}$ and hence (see, e.g., [91, Theorem 13.7]) there exists a polynomial sequence $(p_n)_{n\in\mathbb{Z}}$ on $\mathcal{W}$ such that

$$\|z^{-1} - p_n(z)\|_{L^\infty(\overline{\mathcal{W}})} \xrightarrow{n\to\infty} 0.$$

On the other hand, there exists a closed curve $\Gamma \subset \mathcal{W} \setminus \sigma(A)$ such that (see, e.g., [91, Theorem 13.5])

$$z^{-1} = \frac{1}{2\pi\mathrm{i}} \int_\Gamma \frac{\mathrm{d}\zeta}{\zeta(\zeta - z)}, \qquad p_n(z) = \frac{1}{2\pi\mathrm{i}} \int_\Gamma \frac{p_n(\zeta)}{(\zeta - z)} \mathrm{d}\zeta,$$

whence also (see, e.g. [90, Chapter XI, Section 151])

$$A^{-1} = \frac{1}{2\pi\mathrm{i}} \int_\Gamma \zeta^{-1}(\zeta\mathbb{1} - A)^{-1} \mathrm{d}\zeta, \qquad p_n(A) = \frac{1}{2\pi\mathrm{i}} \int_\Gamma p_n(\zeta)(\zeta\mathbb{1} - A)^{-1} \mathrm{d}\zeta.$$

Thus,

$$\begin{aligned}
\|A^{-1} - p_n(A)\|_{\mathrm{op}} &= \left\|\frac{1}{2\pi\mathrm{i}} \int_\Gamma (\zeta^{-1} - p_n(\zeta))(\zeta\mathbb{1} - A)^{-1} \mathrm{d}\zeta\right\|_{\mathrm{op}} \\
&\leq \|z^{-1} - p_n(z)\|_{L^\infty(\overline{\mathcal{W}})} \left\|\frac{1}{2\pi\mathrm{i}} \int_\Gamma (\zeta\mathbb{1} - A)^{-1} \mathrm{d}\zeta\right\|_{\mathrm{op}} \\
&= \|z^{-1} - p_n(z)\|_{L^\infty(\overline{\mathcal{W}})}
\end{aligned}$$

(indeed, $(2\pi\mathrm{i})^{-1} \int_\Gamma (\zeta\mathbb{1} - A)^{-1} \mathrm{d}\zeta = \mathbb{1}$), and the conclusion follows. $\square$

**Corollary 2.2** *Let $A$ be an operator of $\mathscr{K}$-class on a Hilbert space $\mathcal{H}$. Then the inverse problem $Af = g$ for given $g \in \mathcal{H}$ is Krylov solvable, i.e., the unique solution $f$ belongs to $\overline{\mathcal{K}(A,g)}$.*

*Proof* As $\|p_n(A) - A^{-1}\|_{\mathrm{op}} \xrightarrow{n\to\infty} 0$ (Proposition 2.9), then also $\|p_n(A)g - f\|_{\mathcal{H}} = \|p_n(A)g - A^{-1}g\|_{\mathcal{H}} \xrightarrow{n\to\infty} 0$, and obviously $p_n(A)g \in \mathcal{K}(A,g)$. $\square$

## 2.6 Some illustrative numerical tests

Let us illustrate certain main features discussed theoretically so far through a series of numerical tests with inverse problems in infinite-dimensional Hilbert space, suitably truncated using the GMRES algorithm, and analysed by increasing the size of the truncation (i.e. the number of iterations of GMRES).

We focus on the behaviour of the truncated problems under these circumstances:

I) the solution to the original problem is or is not a Krylov solution;
II) the linear operator is or is not injective (well-defined vs ill-defined problem).

Let us consider four inverse problems. As a baseline case, where the solution is known a priori to be a Krylov solution, we take the compact, injective, self-adjoint multiplication operator on $\ell^2(\mathbb{N})$ (Example 2.4)

$$M^{(\sigma)} := \sum_{n=1}^{\infty} \sigma_n |e_n\rangle\langle e_n|, \qquad \sigma_n := (5n)^{-1}. \tag{2.35}$$

For comparison, we consider a non-injective version

$$M^{(\sigma')} := \sum_{n=1}^{\infty} \sigma'_n |e_n\rangle\langle e_n|, \qquad \sigma'_n := \begin{cases} 0 & \text{if } n \in \{3,6,9\} \\ \sigma_n & \text{otherwise}, \end{cases} \tag{2.36}$$

as well as the weighted right shift (Example 2.3)

$$R_\sigma := \sum_{n=1}^{\infty} \sigma_n |e_{n+1}\rangle\langle e_n| \tag{2.37}$$

with the same weights as in (2.35). The inverse problems $M^{(\sigma)} f = g$, $M^{(\sigma')} f = g$, and $R_\sigma f = g$ are tested with datum $g$ generated by the a priori chosen solution

$$f := \sum_{n \in \mathbb{N}} f_n e_n, \qquad f_n = \begin{cases} n^{-1} & \text{if } n \leqslant 250 \\ 0 & \text{otherwise}. \end{cases} \tag{2.38}$$

Let us observe that

$$\|f\|_{\ell^2} = \sqrt{\tfrac{\pi^2}{6} - \Psi^{(1)}(251)} \simeq 1.28099, \tag{2.39}$$

where $\Psi^{(k)}$ is the polygamma function of order $k$ [1, Sec. 6.4].

Fourth and last, we consider the inverse problem $Vf = g$ where $V$ is the Volterra operator in $L^2[0,1]$ (Example 2.5) and $g(x) := \tfrac{1}{2}x^2$. The problem has unique solution

$$f(x) = x, \qquad \|f\|_{L^2[0,1]} = \frac{1}{\sqrt{3}} \simeq 0.5774. \tag{2.40}$$

Depending on the context, we shall denote by $\mathcal{H}$ and by $A$, respectively, the Hilbert space ($\ell^2(\mathbb{N})$ or $L^2[0,1]$) and the operator ($M^{(\sigma)}$, $M^{(\sigma')}$, $R_\sigma$, or $V$) under consideration.

The inverse problems in $\mathcal{H}$ associated with $M^{(\sigma)}$ and $M^{(\sigma')}$ are Krylov solvable (Proposition 2.7), and so too is the inverse problem associated with $V$, with $\mathcal{K}(V,g)$ dense in $L^2[0,1]$ (Example 2.14). Instead, the problem associated with $R_\sigma$ is not Krylov solvable, for $\mathcal{K}(R_\sigma,g)^\perp$ always contains the first canonical vector $e_1$.

For each operator $A$, let us consider spanning vectors $g, Ag, A^2 g, \ldots$ of $\mathcal{K}(A,g)$ up to order $N_{\max} = 500$ if $A = M^{(\sigma)}, M^{(\sigma')}, R_\sigma$, and up to order $N_{\max} = 175$ if $A = V$. Such values represent our practical choice of 'infinite' dimension for $\mathcal{K}(A,g)$.

Analogously, when $A = M^{(\sigma)}, M^{(\sigma')}, R_\sigma$ let us allocate for each of the considered vectors $f, g, Ag, A^2 g, \ldots$ an amount of 2500 entries with respect to the canonical basis of $\ell^2(\mathbb{N})$, such a value representing our practical choice of 'infinite' dimension for $\mathcal{H}$. Observe, in particular, that by repeated application of $R_\sigma$ up to 500 times, the vectors $R_\sigma^k g$ have non-trivial entries up to order 251+500=751 (by construction the last non-zero entries of $f$ and of $g$ are the components, respectively, $e_{250}$ and $e_{251}$), and by repeated application of $M^{(\sigma)}$ and $M^{(\sigma')}$ the vectors $(M^{(\sigma)})^k g$ and $(M^{(\sigma')})^k g$ have the last non-zero entry in the component $e_{250}$: all such limits stay well below our 'infinity' threshold of 2500 for $\mathcal{H}$.

From each collection $\{g, Ag, \ldots, A^{N-1} g\}$ we then obtain an orthonormal basis of the $N$-dimensional truncation of $\mathcal{K}(A,g)$, $N \leqslant N_{\max}$, and truncate the 'infinite-dimensional' inverse problem $Af = g$ to a $N$-dimensional one, that we solved by means of the GMRES algorithm, in the same spirit of the general discussion of Appendix A.2.

Denoting by $\widehat{f^{(N)}} \in \mathcal{H}$ the vector of the solution from the GMRES algorithm at the $N$-th iterate, we analysed two natural indicators of the convergence 'as $N \to \infty$', the **infinite-dimensional error** $\mathscr{E}_N$ and the **infinite-dimensional residual** $\mathfrak{R}_N$ as defined in Appendix A.2.3, that is, respectively,

$$\mathscr{E}_N := f - \widehat{f^{(N)}}, \qquad \mathfrak{R}_N := g - A\widehat{f^{(N)}}. \tag{2.41}$$

The behaviour of $\|\mathscr{E}_N\|_\mathcal{H}$, $\|\mathfrak{R}_N\|_\mathcal{H}$, and $\|\widehat{f^{(N)}}\|_\mathcal{H}$ at the $N$-th step of the algorithm is illustrated in Figure 2.1 as a function of $N$. The numerical evidence is the following.

- The error norm of the baseline case and the Volterra case tend to vanish with $N$, and so does the residual norm, consistently with the obvious property $\|\mathfrak{R}_N\|_\mathcal{H} \leqslant \|A\|_{\mathrm{op}} \|\mathscr{E}_N\|_\mathcal{H}$. Moreover, $\|\widehat{f^{(N)}}\|_\mathcal{H}$ stays uniformly bounded and attains asymptotically the theoretical value prescribed by (2.39) or (2.40).
- Instead, the error norm of the forward shift remains of order one indicating a lack of *norm-convergence*, regardless of truncation size. Analogous lack of convergence is displayed in the norm of the finite-dimensional residual. Again, $\|\widehat{f^{(N)}}\|_\mathcal{H}$ remains uniformly bounded, but attains an asymptotic value that is strictly smaller than the theoretical value (2.39).

The asymptotics $\|f - \widehat{f^{(N)}}\|_{\ell^2} \to 1.0$ and $\|g - R_\sigma \widehat{f^{(N)}}\|_{\ell^2} \to 0.2$ found numerically for the problem $R_\sigma f = g$ can be understood as follows. Since $\widehat{f^{(N)}} \in \mathcal{K}(R_\sigma, g)$ and since the latter subspace only contains vectors with zero component along $e_1$, the error vector $\mathscr{E}_N = f - \widehat{f^{(N)}}$ tends to approach asymptotically the vector $e_1$ that gives the first component of $f = (1, \frac{1}{2}, \frac{1}{3}, \ldots)$, and this explains $\|\mathscr{E}_N\|_{\ell^2} \to 1$.

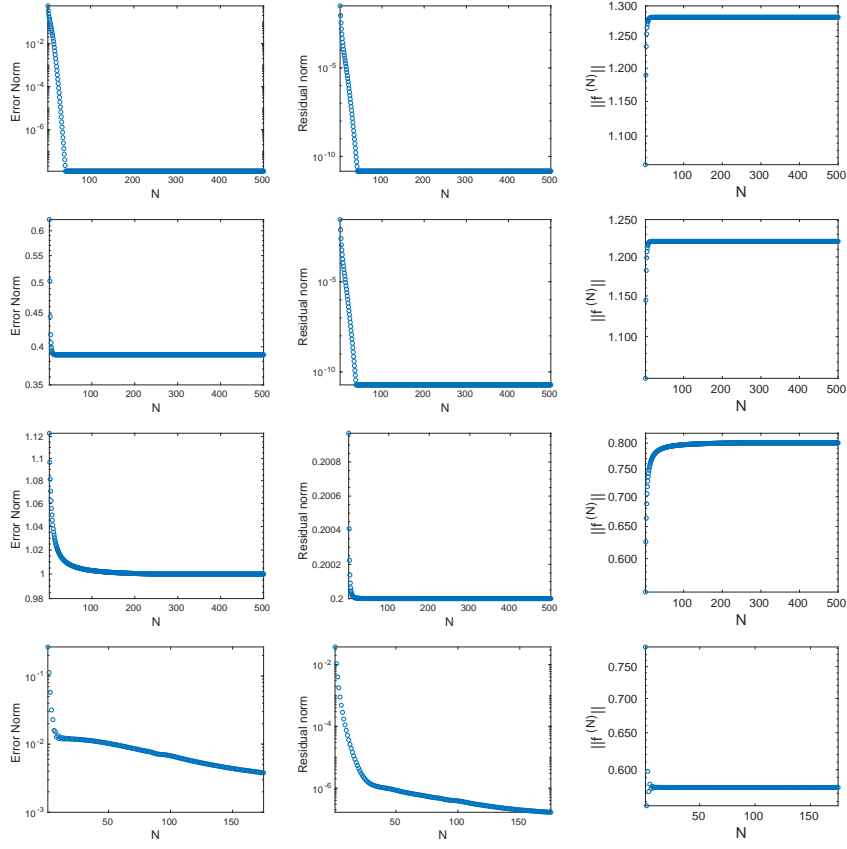

**Fig. 2.1** Error norm, residual norm, and approximate solution magnitude as a function of iterations for the cases of the injective multiplication operator $M^{(\sigma)}$ (top row), the non-injective multiplication operator $M^{(\sigma')}$ (second row from top), the weighted right shift $R_\sigma$ (third row from top), and the Volterra operator $V$ (last row from top).

Analogously, since by construction $g = (0, \frac{1}{5}, \frac{1}{20}, \frac{1}{45}, \dots)$, and since the asymptotics on $\mathscr{E}_N$ implies that each component of $\widehat{f^{(N)}}$ *but the first one* converges to the corresponding component of $f$, then $\widehat{f^{(N)}} \approx (0, \frac{1}{2}, \frac{1}{3}, \dots)$ for large $N$, whence also $R_\sigma \widehat{f^{(N)}} \approx (0, 0, \frac{1}{20}, \frac{1}{45}, \dots)$. Thus $g$ and $R_\sigma \widehat{f^{(N)}}$ tend to differ by only the vector $\frac{1}{5} e_2$, which explains $\|\mathfrak{R}_N\|_{\ell^2} \to \frac{1}{5}$.

In fact, the lack of norm vanishing of error and residual for the problem $R_\sigma f = g$ is far from meaning that the approximants $\widehat{f^{(N)}}$ carry no information about the exact solution $f$: in complete analogy to what we discussed in a more general context in Section A.3 and A.4, $\widehat{f^{(N)}}$ reproduces $f$ *component-wise* for all components but the first.

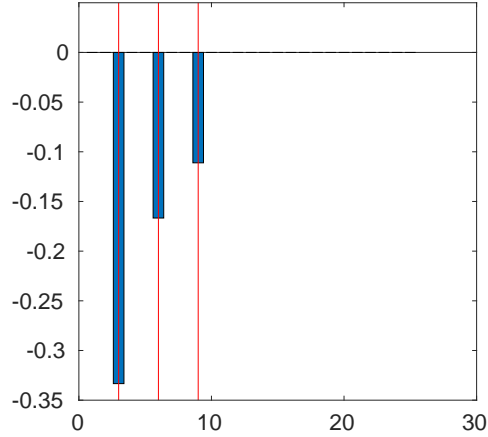

**Fig. 2.2** Support of the error vector (blue bars) for the non-injective problem $M^{(\sigma')}f = g$ at final iteration $N = 500$. The red lines mark the entry positions of the components of the kernel space of $M^{(\sigma')}$.

To summarise the above findings, the Krylov solvable infinite-dimensional problems $M^{(\sigma)}f = g$ and $Vf = g$ display good (i.e., norm-)convergence of error and residual, which is sharper for the multiplication operator $M^{(\sigma)}$ and quite slower for the Volterra operator $V$, indicating that the choice of the Krylov bases is not equally effective for the two problems. This is in contrast with the non-Krylov solvable problem $R_\sigma f = g$, where convergence in norm fails. The uniformity in the size of the solutions produced by the GMRES algorithm appears not to be affected by the presence or the lack of Krylov solvability.

Last, concerning the behaviour of the truncated problems in the absence of injectivity, thus the finite-dimensional truncations of $M^{(\sigma')}f = g$ (which in any case admits a *unique* Krylov solution, on account of Proposition 2.7), the numerical evidence is the following.

- As opposite to the baseline problem $M^{(\sigma)}f = g$, $\|\mathscr{E}_N\|_{\ell^2} = \|f - \widehat{f^{(N)}}\|_{\ell^2}$ does not vanish with the truncation size and remains instead uniformly bounded. Instead, $\|M^{(\sigma')}s\widehat{f^{(N)}} - g\|_{\ell^2(\mathbb{N})}$ displays the same vanishing behaviour as for $M^{(\sigma)}$ (Fig. 2.1).
- $\|\widehat{f^{(N)}}\|_{\ell^2(\mathbb{N}}$ remains uniformly bounded (Fig. 2.1).

The reason as to the observed lack of convergence of the error is unmasked in Figure 2.2. There one can see that the only components in the error vector that are non-zero are the components corresponding to kernel vector entries. This shows that the Krylov algorithm has indeed found a solution to the problem, modulo the kernel components in $f$.

# Chapter 3
# An analysis of conjugate-gradient based methods with unbounded operators

## 3.1 Unbounded positive self-adjoint inverse problems and conjugate gradient approach

We are concerned in this Chapter with the rigorous proof of the convergence, in various meaningful senses, of a particular and well-known iterative algorithm for solving inverse linear problems, the celebrated conjugate gradient method, in the generalised setting of *unbounded* operators on Hilbert space.

The framework of unbounded inverse problems, and the issue of their Krylov solvability, is the actual topic of Chapter 4. Yet, it is advantageous to anticipate the discussion of a particular type of unbounded inverse problem and a specific projection method to solve it. Indeed, on the one hand the special result obtained here is going to provide the basis for the study of Krylov solvability of a larger class of inverse problems (Section 4.2). Moreover, the present particular case requires sophisticated methods and tools that are independent from the general analysis of Chapter 4, and make this separate discussion autonomous and self-consistent.

The object of investigation is the usual inverse linear problem of the form (1.1), except that unlike the framework analysed in Chapter 2, the linear operator $A$ acting on the abstract Hilbert space $\mathcal{H}$ is now self-adjoint and positive, possibly unbounded. In the latter (unbounded self-adjoint) case, the domain $\mathcal{D}(A)$ of $A$ is necessarily a proper dense subspace of $\mathcal{H}$. All our considerations and results will apply also to everywhere defined, bounded, positive, self-adjoint $A$'s, but it is of course the unbounded case that we will have primarily in mind. In all cases, the positivity (equivalently referred to as non-negativity) assumption on $A$ reads $\langle \psi, A\psi \rangle \geqslant 0$ for all $\psi \in \mathcal{D}(A)$.

We shall give for granted an amount of basic properties of (unbounded) self-adjoint operators, for the general theory of which we refer to [88, Chapt. VI-VII] and [97, Chapt. 1-5].

The general setting for unbounded inverse problem is deferred to Section 4.1. The specific set-up for the present case is presented in Section 3.2. In particular, as already implicitly argued in Section 1.3, since now the inverse problem

$$Af = g, \qquad A = A^* \geqslant \mathbb{O} \tag{3.1}$$

involves a (possibly) proper subspace $\mathcal{D}(A) \subset \mathcal{H}$, the very notion (2.2) of Krylov subspace $\mathcal{K}(A,g)$ makes only sense when the datum $g$ belongs to $\mathcal{D}(A^k)$ for every $k \in \mathbb{N}$, an occurrence that we shall call in a moment *A-smoothness* (see (3.6) below, as well as Section 4.1). The datum $g$ will therefore be a vector in ran$A$ that is also *A*-smooth. The extra requirement of *A*-smoothness is trivially satisfied in the bounded-*A* case.

For self-adjoint and positive inverse problems (3.1) in the finite-dimensional case (thus, when $\mathcal{H} \cong \mathbb{C}^d$ for some $d \in \mathbb{N}$ and $A$ is a $d \times d$ positive semi-definite hermitian matrix), a very popular numerical algorithm is available, the **method of conjugate gradients (CG)**. It was first proposed in 1952 by Hestenes and Stiefel [62] and since then, together with its related derivatives (e.g., conjugate gradient method on the normal equations (CGNE), least-square QR method (LSQR), etc.), it has been widely studied in the finite-dimensional setting (we refer to the classical monographs [96, 103, 73]) and also, though to a lesser extent, in the infinite-dimensional Hilbert space setting with bounded, positive, self-adjoint operators.

The conjugate gradient algorithm is an iterative method, essentially of the type of Krylov-based projection methods, that produces approximate solutions $f^{[N]}$ to (3.1) in the form of iterates, namely $f^{[N+1]}$ is constructed in terms of $f^{[N]}$. We shall actually use an equivalent, variational characterisation of $f^{[N]}$, as explained in Section 3.2. For sufficiently large $N$ when $\dim \mathcal{H} < \infty$, or more generally in the limit $N \to \infty$ when $\dim \mathcal{H} = \infty$, such approximants under suitable condition display indeed the norm convergence $f^{[N]} \to f$.

Even though the conjugate gradient convergence theory has been markedly less explored when $\dim \mathcal{H} = \infty$, there are lines of investigation in which yet important works have been produced over the last five decades, both in the scenario where $A$ is bounded with everywhere-defined bounded inverse [27, 28, 61], or at least with bounded inverse on its range [64], and in the scenario where $A$ is bounded with possible unbounded inverse on its range [64, 81, 82, 74, 12].

A fundamental analysis in this respect is that by Nemirovskiy and Polyak [81, 82] from the mid 1980's. In fact, this Chapter is essentially a *non-trivial* generalisation and extension of the Nemirovskiy-Polyak analysis. In [81] the full convergence estimates of the **error** $\|f^{[N]} - f\|$ and **residual** $\|Af^{[N]} - g\|$ (we refer to Appendix A.2.3 for such nomenclature) were proved for *bounded A*, and in the follow-up work [82] the results were shown to be *optimal* in the sense that for the entire class of bounded, ill-posed problems, one can do no better than the estimates provided.

The boundedness of $A$ was crucial in a two-fold way. First, it forced the blow-up of a suitable sum ($\delta_N$, in their notation – see (3.45) below) of the reciprocals of the $N$ zeroes of a polynomial that represents the minimisation (3.3): since no such zero can exceed $\|A\|_{\mathrm{op}}$, the reciprocals cannot vanish and their sum necessarily diverges. As a consequence, the error and the residual, which in turn can be controlled by an inverse power of $\delta_N$, are then shown to vanish as $N \to \infty$, thus establishing convergence. Second, boundedness of $A$ was also determinant to quantify the convergence, as the latter was boiled down to a min-max procedure for polynomials on the *finite* spectral interval containing $\sigma(A)$, then on such an interval (suitable modifications of) the Chebyshev polynomials are recognised to optimise the rate of convergence, and explicit properties of (the extrema of) Chebyshev polynomials finally provide a quantitative version of the vanishing of error and residual.

In contrast, the scenario where $A$ is *unbounded* has been only recently considered from special perspectives, in particular in view of existence [86] (the GMRES algorithm applied to a particular non self-adjoint differential operator), or convergence when $A$ is regularised and made invertible with everywhere-defined bounded inverse [51]. It is in view of the virtually unexplored general convergence theory of the unbounded conjugate gradient method that we recently produced the study [20] this Chapter is based upon.

Our approach bypasses the Nemirovskiy-Polyak restriction of the finiteness of $\|A\|_{\mathrm{op}}$ as far as the convergence alone is concerned. As for the quantitative rate, the min-max strategy of [81] by no means can be adapted to polynomials over the whole $[0,+\infty)$ and in fact a careful analysis of the structure of the proof of [81, 82], as we shall comment in due time, seems to indicate that if $A$ is unbounded with unbounded inverse on its range, then the convergence rate can be arbitrarily small, and possibly not guaranteed in some cases.

As mentioned already, Section 3.2 contains the rigorous set-up and the statement of the convergence result. It is important to stress that for such convergence to hold one needs to further select the datum $g$ within the class of $A$-smooth vectors: the compromise of reasonable generality adopted here is the sub-class of (quasi-)analytic vectors. Technical and preparatory materials of algebraic and measure-theoretic nature are then developed in Section 3.3. The CG-convergence is proved in Section 3.4, and further supplemented with an amount of additional remarks clarifying the actual novelties and differences of the present scheme as compared to the Nemirovskiy-Polyak (bounded case) scenario. Last, in Section 3.5 a selection of numerical tests are presented confirming the main convergence features and corroborating our intuition on certain relevant differences with respect to the bounded case.

## 3.2 Set-up and main results

Here and in the following $A$ is a positive self-adjoint operator on a Hilbert space $\mathcal{H}$, including the possibility that $A$ be unbounded (in which case the domain $\mathcal{D}(A)$ is a proper dense subspace of $\mathcal{H}$) and with a non-trivial kernel.

First, in order to describe the algorithm explicitly, let us introduce the **solution manifold**

$$\mathcal{S} := \{f \in \mathcal{D}(A) \,|\, Af = g\} \qquad (3.2)$$

relative to the problem (3.1). By assumption $\mathcal{S}$ is a convex, non-empty set in $\mathcal{H}$ which is also closed, owing to the fact that $A$, being self-adjoint, is in particular a closed operator. In fact, owing to the linearity of $A$, $\mathcal{S}$ is an affine space with the structure $\mathcal{S} = \{f_0\} + \ker A$, where $f_0$ is a solution to (3.1): we keep the 'solution manifold' terminology for consistency with the nomenclature in the literature. As a consequence of the convexity, non-emptiness, and closedness of $\mathcal{S}$, the **projection map** $P_\mathcal{S} : \mathcal{H} \to \mathcal{S}$ is unambiguously defined and produces, for generic $x \in \mathcal{H}$, the closest-to-$x$ point in $\mathcal{S}$. Observe that $P_\mathcal{S}$ is not a linear map.

Then one needs to ensure that the **conjugate gradient iterates** are *well-defined*. One chooses a datum $g \in \operatorname{ran} A$ and an **initial guess** $f^{[0]} \in \mathcal{H}$, and for some $\theta \geqslant 0$ defines the $\theta$-iterates

$$f^{[N]} := \arg\min_{h \in \{f^{[0]}\} + \mathcal{K}_N(A, \mathfrak{R}_0)} \|A^{\theta/2}(h - P_\mathcal{S} h)\|, \qquad N \in \mathbb{N} \qquad (3.3)$$

with

$$\mathfrak{R}_N := A f^{[N]} - g, \qquad N \in \mathbb{N}_0, \qquad (3.4)$$

$$\mathcal{K}_N(A, \mathfrak{R}_0) := \operatorname{span}\{\mathfrak{R}_0, A\mathfrak{R}_0, \ldots A^{N-1}\mathfrak{R}_0\}, \qquad N \in \mathbb{N}, \qquad (3.5)$$

where $\mathfrak{R}_N$ is the **residual** at the $N$-th step.

In order to apply an arbitrary positive power of $A$ to $Af^{[0]} - g$, both $g$ and $f^{[0]}$ are required to be ***A*-smooth vectors** [89, Sect. X.6], meaning that they belong to the space

$$C^\infty(A) := \bigcap_{N \in \mathbb{N}} \mathcal{D}(A^N). \qquad (3.6)$$

In the applications where $A$ is a differential operator, $A$-smoothness is a regularity requirement.

In turn, $A$-smoothness of $g$ and $f^{[0]}$ implies $\mathcal{K}_N(A, \mathfrak{R}_0) \subset C^\infty(A)$, and obviously $P_\mathcal{S} h \in \mathcal{S} \subset C^\infty(A)$, whereas by interpolation $C^\infty(A) \subset \mathcal{D}(A^{\theta/2})$ for any $\theta \geqslant 0$. This guarantees that in the minimisation (3.3) one is allowed to apply $A^{\theta/2}$ to any vector $h - P_\mathcal{S} h$.

Thus, under the assumptions

$$g \in \operatorname{ran} A \cap C^\infty(A), \qquad f^{[0]} \in C^\infty(A) \qquad (3.7)$$

the corresponding **$\theta$-iterates** $f^{[N]}$ are unambiguously defined by (3.3)-(3.5) above for *any* $\theta \geqslant 0$. If $A$ is bounded, (3.7) simply reduces to $g \in \operatorname{ran} A$. In fact, (3.7) are *minimal assumptions*, inescapable if one wants to give meaning to conjugate gradient iterates in the unbounded case. It will be also commented in due time (Remark 3.6) why the $f^{[N]}$'s exist and are unique, therefore ensuring that they are indeed well-defined as claimed.

Such iterates have three notable properties, whose proof is deferred to Section 3.3.

**Proposition 3.1** *The $\theta$-iterates $f^{[N]}$ defined for a given $\theta \geqslant 0$ by means of* (3.3)-(3.5) *under the assumption* (3.7) *satisfy*

$$f^{[N]} - P_{\mathcal{S}} f^{[N]} \in (\ker A)^\perp \qquad \forall N \in \mathbb{N}_0, \qquad (3.8)$$

$$P_{\mathcal{S}} f^{[N]} = P_{\mathcal{S}} f^{[0]} \qquad \forall N \in \mathbb{N}, \qquad (3.9)$$

$$f^{[N]} - P_{\mathcal{S}} f^{[N]} = p_N(A)(f^{[0]} - P_{\mathcal{S}} f^{[0]}) \qquad \forall N \in \mathbb{N}, \qquad (3.10)$$

where $p_N(\lambda)$ is for each N a polynomial of degree at most N and such that $p_N(0) = 1$.

As, by (3.9), all such $f^{[N]}$'s have the same projection onto the solution manifold $\mathcal{S}$, the approach of $f^{[N]}$ to $\mathcal{S}$ consists explicitly of a convergence $f^{[N]} \to P_{\mathcal{S}} f^{[0]}$. Let us now specify in which sense this convergence is to be monitored.

The underlying idea, as is clear in the typical applications where A is a differential operator on a $L^2$-space, is that $\|f^{[N]} - P_{\mathcal{S}} f^{[0]}\|_{(A)} \to 0$ in some A-dependent Sobolev norm. For this to make sense, clearly one needs enough A-regularity on $f^{[N]} - P_{\mathcal{S}} f^{[0]}$, which eventually is guaranteed by the regularity initially assumed on $f^{[0]}$. Thus, the general indicator of convergence has the form $\|A^{\sigma/2}(f^{[N]} - P_{\mathcal{S}} f^{[0]})\|$, but extra care is needed if one wants to control the convergence in the abstract analogue of a low-regularity, negative-Sobolev norm, which would amount to formally consider $\sigma < 0$, for in general A can have a kernel and hence is only invertible on its range.

Based on these considerations, and inspired by the analogous discussion in [81] for bounded A, let us introduce the class $\mathscr{C}_{A,g}(\theta)$ of vectors of $\mathcal{H}$ defined for generic $\theta \in \mathbb{R}$ as

$$\mathscr{C}_{A,g}(\theta) := \begin{cases} \{x \in \mathcal{H} \,|\, x - P_{\mathcal{S}} x \in \mathcal{D}(A^{\frac{\theta}{2}})\}, & \theta \geqslant 0, \\ \{x \in \mathcal{H} \,|\, x - P_{\mathcal{S}} x \in \mathrm{ran}\,(A^{-\frac{\theta}{2}})\}, & \theta < 0. \end{cases} \qquad (3.11)$$

(The dependence of $\mathscr{C}_{A,g}(\theta)$ on g is implicit through the solution manifold $\mathcal{S}$.) Distinguishing the two cases in (3.11) is needed whenever A has a non-trivial kernel. If instead A is injective, and so too is therefore $A^{-\frac{\theta}{2}}$ for $\theta < 0$, then $A^{-\frac{\theta}{2}}$ is a bijection between the two dense subspaces $\mathcal{D}(A^{-\frac{\theta}{2}}) = \mathrm{ran}\,(A^{\frac{\theta}{2}})$ and $\mathrm{ran}\,(A^{-\frac{\theta}{2}}) = \mathcal{D}(A^{\frac{\theta}{2}})$ of $\mathcal{H}$.

Related to the class $\mathscr{C}_{A,g}(\theta)$ one has two further useful notions. One, for fixed $\theta \in \mathbb{R}$ and $x \in \mathscr{C}_{A,g}(\theta)$, is the vector

$$u_\theta(x) := \begin{cases} A^{\frac{\theta}{2}}(x - P_{\mathcal{S}} x), & \theta \geqslant 0, \\ \text{the minimal norm solution } u \text{ to } A^{-\frac{\theta}{2}} u = x - P_{\mathcal{S}} x, & \theta < 0. \end{cases} \qquad (3.12)$$

The other is the functional $\rho_\theta$ defined on the vectors $x \in \mathscr{C}_{A,g}(\theta)$ as

$$\rho_\theta(x) := \|u_\theta(x)\|^2. \qquad (3.13)$$

Thus,

$$\rho_\theta(x) = \begin{cases} \|A^{\frac{\theta}{2}}(x - P_{\mathcal{S}} x)\|^2, & \theta \geqslant 0, \\ \left\|\left(A^{-\frac{\theta}{2}}\big|_{\mathrm{ran}\,(A^{-\frac{\theta}{2}})}\right)^{-1}(x - P_{\mathcal{S}} x)\right\|^2, & \theta < 0, \end{cases} \qquad (3.14)$$

with an innocent abuse of notation in (3.14) when $\theta < 0$, as *the operator inverse is to be understood for a (self-adjoint, and positive) operator on the Hilbert subspace* $\overline{\operatorname{ran} A}$.

It is worth remarking that in the special case when $A$ is *bounded*, the following interesting properties hold, whose proof is deferred to Section 3.3, which do not have a counterpart in the unbounded case except for the obvious identity $\mathscr{C}_{A,g}(0) = \mathcal{H}$.

**Lemma 3.1** *If $A$ (besides being self-adjoint and non-negative) is bounded, and if $g \in \operatorname{ran} A$, then:*

(i) $\mathscr{C}_{A,g}(\theta) = \mathcal{H}$ whenever $\theta \geqslant 0$;
(ii) $\mathscr{C}_{A,g}(\theta) \subset \mathscr{C}_{A,g}(\theta')$ for $\theta \leqslant \theta'$;
(iii) for $\theta \leqslant \theta'$ and $x \in \mathscr{C}_{A,g}(\theta)$ one has $u_{\theta'}(x) = A^{(\theta'-\theta)/2} u_\theta(x)$, whence also $\rho_{\theta'}(x) \leqslant \|A\|_{\operatorname{op}}^{\theta'-\theta} \rho_\theta(x)$.

Back to the general case where $A$ is possibly *unbounded*, the goal is to evaluate certain $\rho_\sigma$-functionals along the sequence of the $f^{[N]}$'s. This may require an extra assumption on the initial guess $f^{[0]}$, as the following Lemma shows.

**Lemma 3.2** *Consider the $\theta$-iterates $f^{[N]}$ defined for a given $\theta \geqslant 0$ by means of (3.3)-(3.5) under the assumption (3.7). Then:*

(i) $f^{[N]} \in \mathscr{C}_{A,g}(\sigma) \; \forall \sigma \geqslant 0$;
(ii) $f^{[N]} \in \mathscr{C}_{A,g}(\sigma)$ for any $\sigma < 0$ such that, additionally, $f^{[0]} \in \mathscr{C}_{A,g}(\sigma)$, in which case
$$u_\sigma(f^{[N]}) = p_N(A) u_\sigma(f^{[0]}), \tag{3.15}$$

*where $p_N(\lambda)$ is precisely the polynomial mentioned in Proposition 3.1.*

Lemma 3.2 is a direct consequence of (3.10) in Proposition 3.1 above: for completeness its proof is included in Section 3.3.

It then makes sense to control the convergence $f^{[N]} \to \mathcal{P}_S f^{[0]}$ in the $\rho_\sigma$-sense, for $\sigma$ positive or negative, with suitable assumptions on $f^{[0]}$. Explicitly,

$$\rho_\sigma(f^{[N]}) = \|u_\sigma(f^{[N]})\|^2 \\ = \begin{cases} \|A^{\frac{\sigma}{2}}(f^{[N]} - \mathcal{P}_S f^{[0]})\|^2, & \sigma \geqslant 0, \\ \left\|\left(A^{-\frac{\sigma}{2}}\big|_{\operatorname{ran}\left(A^{-\frac{\sigma}{2}}\right)}\right)^{-1}(f^{[N]} - \mathcal{P}_S f^{[0]})\right\|^2, & \sigma < 0, \end{cases} \tag{3.16}$$

having used (3.9). The most typical and meaningful choices in the applications are

$$\begin{aligned} \rho_0(f^{[N]}) &= \|f^{[N]} - \mathcal{P}_S f^{[0]}\|^2, \\ \rho_1(f^{[N]}) &= \langle f^{[N]} - \mathcal{P}_S f^{[0]}, A(f^{[N]} - \mathcal{P}_S f^{[0]})\rangle, \\ \rho_2(f^{[N]}) &= \|A(f^{[N]} - \mathcal{P}_S f^{[0]})\|^2, \end{aligned} \tag{3.17}$$

that is, respectively, the norm of the **error**, the **'energy' (semi-)norm**, and the norm of the **residual**.

The preparation made so far for our forthcoming main result (Theorem 3.1 below) does not account yet for the necessity of one further, restrictive assumption on the datum $g$ and the initial guess $f^{[0]}$ of the algorithm, a restriction needed once again to deal with the possible unboundedness of the operator $A$ (the bounded case being controllable for arbitrary $g \in \mathrm{ran} A$ and $f^{[0]} \in \mathcal{H}$). The actual need of a special, inevitable choice of $g$ and $f^{[0]}$ will be fully clear in the course of the proof; for the time being, let us outline here a short heuristic reasoning.

From the expression of the indicator of convergence $\rho_\sigma(f^{[N]})$, for concreteness the case $\sigma \geqslant 0$ in (3.16), and from the iterates properties (3.9)-(3.10) announced in Proposition 3.1, it is easy to realise, as we shall argue in the next Section, that the actual quantity to control along the limit $N \to \infty$ is an integral of the form

$$\int_{[0,+\infty)} \left| \lambda^{\frac{\sigma}{2}} p_N(\lambda) \right|^2 \mathrm{d}\langle f^{[0]} - P_{\mathcal{S}} f^{[0]}, E^A(\lambda)(f^{[0]} - P_{\mathcal{S}} f^{[0]}) \rangle\,,$$

for suitable polynomials $p_N$ determined by the minimisation (3.3), where the measure is the scalar spectral measure associated to the self-adjoint operator $A$, and hence (by positivity of $A$) the integration runs over $[0, +\infty)$. For bounded $A$'s the integration is actually restricted within the spectrum of $A$, hence within a compact interval of the non-negative half line, and this naturally provides a kind of *uniformity* in $N$ that is crucial in controlling the vanishing of the above integral in the limit. When instead $A$ (and hence the integration domain) is unbounded, some other source of uniformity in $N$ must be implemented, which eventually is to be some kind of uniformity of the measure $\lambda^{2N} \mathrm{d}\langle f^{[0]} - P_{\mathcal{S}} f^{[0]}, E^A(\lambda)(f^{[0]} - P_{\mathcal{S}} f^{[0]}) \rangle$, in other words, a suitable *control of the growth in N of the norm* $\|A^N(f^{[0]} - P_{\mathcal{S}} f^{[0]})\|$. In turn, this requires a control in $N$ of $\|A^N f^{[0]}\|$ and $\|A^N g\|$.

It is with the above heuristics in mind that one appeals to the following classes of vectors [97, Definition 7.1]: the **analytic vectors** for $A$ are the elements of the *subspace*

$$\mathcal{D}^a(A) := \left\{ g \in C^\infty(A) = \bigcap_{n \in \mathbb{N}} \mathcal{D}(A^n) \,\middle|\, \begin{array}{l} \|A^n g\| \leqslant C_g^n n! \\ \text{for any } n \in \mathbb{N} \\ \text{and some } C_g > 0 \end{array} \right\}, \qquad (3.18)$$

and the **quasi-analytic vectors** for $A$ are the elements of the *set*

$$\mathcal{D}^{qa}(A) := \left\{ g \in C^\infty(A) \,\middle|\, \sum_{n \in \mathbb{N}} \|A^n g\|^{-\frac{1}{n}} = +\infty \right\}. \qquad (3.19)$$

Clearly $\mathcal{D}^a(A) \subset \mathcal{D}^{qa}(A) \subset C^\infty(A)$, and the self-adjointness of $A$ ensures that the subspace of its analytic vectors is *dense* in $\mathcal{H}$ (this is the Nelson theorem: see, e.g., [97, Theorem 7.16]). Obviously when $A$ is bounded the whole $\mathcal{H}$ is made of analytic vectors for $A$.

All this finally allows the formulation of the CG-convergence result.

**Theorem 3.1** *Let A be a non-negative self-adjoint operator on the Hilbert space $\mathcal{H}$. Let*

$$g \in \mathcal{D}^a(A) \cap \operatorname{ran} A. \tag{3.20}$$

*Consider the conjugate gradient algorithm associated with A and g where the initial guess vector $f^{[0]}$ satisfies*

$$f^{[0]} \in \mathcal{D}^a(A) \cap \mathscr{C}_{A,g}(\sigma^*), \qquad \sigma^* = \min\{\sigma, 0\} \tag{3.21}$$

*for a given $\sigma \in \mathbb{R}$, and where the iterates $f^{[N]}$, $N \in \mathbb{N}$, are constructed via (3.3) with parameter $\theta = \xi \geqslant 0$ under the condition $\sigma \leqslant \xi$. Then*

$$\lim_{N \to \infty} \rho_\sigma(f^{[N]}) = 0. \tag{3.22}$$

In fact, as a corollary of the proof discussed in Sections 3.3-3.4, one also has:

**Corollary 3.1** *The same conclusion of Theorem 3.1 follows also when the assumptions on g and $f^{[0]}$ are replaced by*

$$\begin{aligned} g &\in C^\infty(A) \cap \operatorname{ran} A, \\ f^{[0]} &\in C^\infty(A) \cap \mathscr{C}_{A,g}(\sigma^*), \\ f^{[0]} - P_{\mathcal{S}} f^{[0]} &\in \mathcal{D}^{qa}(A), \end{aligned} \tag{3.23}$$

*or by*

$$\begin{aligned} g &\in \mathcal{D}^{qa}(A) \cap \operatorname{ran} A, \\ f^{[0]} &= 0. \end{aligned} \tag{3.24}$$

In other words, Theorem 3.1 states that the convergence holds at a given '*A*-regularity level' $\sigma$ for $\xi$-iterates built with *equal or higher* '*A*-regularity level' $\xi \geqslant \sigma$, and with an initial guess $f^{[0]}$ that is *A*-analytic if $\sigma \geqslant 0$, and additionally belongs to the class $\mathscr{C}_{A,g}(\sigma)$ if $\sigma < 0$.

In particular, with no extra assumption on $f^{[0]}$ but its *A*-analyticity, the $\xi$-iterates with $\xi \geqslant 0$ automatically converge in the sense of the error ($\sigma = 0$, see (3.17) above), the $\xi$-iterates with $\xi \geqslant 1$ automatically converge in the sense of the error and of the energy norm ($\sigma = 1$), the $\xi$-iterates with $\xi \geqslant 2$ automatically converge in the sense of the error, energy norm, and residual ($\sigma = 2$).

*Remark 3.1* If, for a finite $N$, $\rho_\sigma(f^{[N]}) = 0$, then the very iterate $f^{[N]}$ *is* a solution to the linear problem $Af = g$, and one says that the algorithm 'has come to convergence' in a finite number ($N$) of steps. Indeed, $\rho_\sigma(f^{[N]}) = 0$ is the same as $A^{\frac{\sigma}{2}}(f^{[N]} - P_{\mathcal{S}} f^{[0]}) = 0$ if $\sigma \geqslant 0$, i.e., $f^{[N]} - P_{\mathcal{S}} f^{[0]} \in \ker A^{\frac{\sigma}{2}} = \ker A$; this, combined with $f^{[N]} - P_{\mathcal{S}} f^{[0]} \in (\ker A)^\perp$ (see (3.8)-(3.9) above), implies that $f^{[N]} = P_{\mathcal{S}} f^{[0]} \in \mathcal{S}$. On the other hand, $\rho_\sigma(f^{[N]}) = 0$ is the same as $u_\sigma(f^{[N]}) = 0$ with $A^{-\frac{\sigma}{2}} u_\sigma(f^{[N]}) = f^{[N]} - P_{\mathcal{S}} f^{[0]}$ if $\sigma < 0$, whence again $f^{[N]} = P_{\mathcal{S}} f^{[0]} \in \mathcal{S}$.

*Remark 3.2*

(i) In the special scenario where $A$ is (everywhere-defined and) bounded, $A$-analyticity is automatically guaranteed, so one only needs to assume that $g \in \operatorname{ran} A$ and $f^{[0]} \in \mathscr{C}_{A,g}(\sigma^*)$ for some $\sigma \in \mathbb{R}$ ($\sigma^* = \min\{\sigma, 0\}$) in order for the convergence of the $\xi$-iterates ($\xi \geqslant \sigma$) to hold in the sense $\rho_\sigma(f^{[N]}) \to 0$. Then, owing to Lemma 3.1, one automatically has also $\rho_{\sigma'}(f^{[N]}) \to 0$ for any $\sigma' \geqslant \sigma$. This is precisely the form of the convergence result originally established by Nemirovskiy and Polyak [81].

(ii) Thus, in the bounded-case scenario, if $\sigma$ is the minimum level of convergence chosen, then not only are the $\xi$-iterates with $\xi \geqslant \sigma$ proved to $\rho_\sigma$-converge, but in addition the *same* $\xi$-iterates also $\rho_{\sigma'}$-converge at any other level $\sigma' \geqslant \sigma$, with no upper bound on $\sigma'$. In particular, it is shown in [81] that

$$\rho_{\sigma'}(f^{[N]}) \leqslant C(\|A\|_{\mathrm{op}}, \xi - \sigma)(2N+1)^{-2(\sigma'-\sigma)} \rho_\sigma(f^{[0]}), \quad \sigma < \sigma' \leqslant \xi, \tag{3.25}$$

for some constant $C(\|A\|_{\mathrm{op}}, \xi - \sigma) > 0$, thus providing an explicit *rate of convergence* of the $\xi$-iterates in a generic $\rho_{\sigma'}$-sense such that $\sigma' \in (\sigma, \xi]$.

(iii) In the general unbounded-case scenario, instead, the $\rho_\sigma$-convergence guaranteed by Theorem 3.1 is not exportable to $\rho_{\sigma'}$-convergence with $\sigma' > \sigma$.

*Remark 3.3* When, in the unbounded case, $A$ has an everywhere-defined bounded inverse, one has $\mathscr{C}_{A,g}(\sigma) = \mathcal{H}$ for any $\sigma \leqslant 0$. Therefore, Theorem 3.1 guarantees the $\rho_\sigma$-convergence of the $\xi$-iterates for any $\sigma \leqslant 0$, provided that $g$ and $f^{[0]}$ are $A$-analytic. Such 'weaker' convergence can be still informative in many contexts. For instance, choosing

$$\begin{aligned} \mathcal{H} &= L^2(\mathbb{R}^d), \\ A &= -\Delta + \mathbb{1} \quad \text{with} \quad \mathcal{D}(A) = H^2(\mathbb{R}^d) \qquad (\operatorname{ran} A = \mathcal{H}), \\ g, f^{[0]} &\in C^\infty(\mathbb{R}^d), \end{aligned}$$

we see that the $\theta$-iterates defined by (3.3) with the above data converge to the unique solution $f$ to the inverse problem $-\Delta f + f = g$ in any negative Sobolev space $H^\sigma(\mathbb{R}^d)$, $\sigma < 0$; in particular, $f^{[N]}(x) \to f(x)$ point-wise almost everywhere.

*Remark 3.4*

(i) Assumptions (3.20)-(3.21) of Theorem 3.1, as well as assumptions (3.23) of Corollary 3.1, are needed to cover the case of our primary interest, the *unboundedness* of $A$.
(ii) Such restrictions still allow the admissible $g$ and $f^{[0]}$ to run over a *dense* of $\mathcal{H}$.
(iii) Assumptions (3.23) are slightly less restrictive than (3.20)-(3.21). Indeed, from (3.20)-(3.21), since $g$ is analytic, so is $P_S f^{[0]}$ (a fact that we shall prove in Lemma 3.4), and by linearity $f^{[0]} - P_S f^{[0]}$ is analytic too, whence (3.23).
(iv) Albeit more general, assumptions (3.23) have the apparent drawback of being formulated in terms of a vector, $P_S f^{[0]}$, that is unknown prior to actually solving the inverse problem. In fact, (3.23) was singled out because, as is going to emerge in the forthcoming discussion, it is precisely the quasi-analyticity of $f^{[0]} - P_S f^{[0]}$ (together with the inevitable operational assumption $g \in C^\infty(A) \cap \mathrm{ran}\, A$) that makes our proof work. In fact, quasi-analyticity of $f^{[0]} - P_S f^{[0]}$ provides a control on the $N$-growth of $\|A^N (f^{[0]} - P_S f^{[0]})\|$, the quantity we heuristically discussed prior to stating Theorem 3.1.

(v) Assumptions (3.24) are in fact a special case of (3.23), as will be clear from Lemma 3.4(i). They refer to the frequent occurrence, in conjugate gradient methods, where the initial guess $f^{[0]}$ is just the zero function.

*Remark 3.5* As a follow-up of Remark 3.4(iii): although the quasi-analyticity of $f^{[0]} - P_S f^{[0]}$ is only a *sufficient* condition, some possibly weaker assumption of that sort, namely some kind of control of the growth in $N$ of $\|A^N(f^{[0]} - P_S f^{[0]})\|$, is surely needed for the conjugate gradient convergence (3.22). That in the regime $g, f^{[0]} \in C^\infty(A)$ the vanishing $\rho_\sigma(f^{[N]}) \to 0$ is *not* guaranteed, is going to be explained in Proposition 3.5, when the technical details of Theorem 3.1 will be clear.

## 3.3 Algebraic and measure-theoretic background properties

The discussion of this Section concerns an amount of technical properties that are needed for the proof of the main Theorem 3.1.

For convenience, let us set for each $N \in \mathbb{N}$

$$\begin{aligned}
\mathbb{P}([0, +\infty)) &:= \{\text{real-valued polynomials } p(\lambda), \lambda \in [0, +\infty)\} \\
\mathbb{P}_N &:= \{ p \in \mathbb{P}([0, +\infty)) \mid \deg p \leqslant N \} \\
\mathbb{P}_N^{(1)} &:= \{ p \in \mathbb{P}_N \mid p(0) = 1 \}.
\end{aligned} \tag{3.26}$$

Let us start with the proof of those properties stated in Section 3.2. The proof of Proposition 3.1 requires the following elementary fact.

**Lemma 3.3** *Let $z \in \mathcal{H}$. For a point $y \in S$ these conditions are equivalent:*

*(i) $y = P_S z$,*
*(ii) $z - y \in (\ker A)^\perp$.*

*Proof* By linearity of $A$, $S = \{y\} + \ker A$. If $z - y \in (\ker A)^\perp$, then for any $x \in \ker A$, and hence for a generic point $y + x \in S$, one has

$$\|z - (y+x)\|^2 = \|z - y\|^2 + \|x\|^2 \geqslant \|z - y\|^2,$$

therefore $y$ is necessarily the closest to $z$ among all points in $S$, i.e., $y = P_S z$. This proves that (ii) $\Rightarrow$ (i). Conversely, if $y = P_S z$, and if by contradiction $z - y$ does not belong to $(\ker A)^\perp$, then $\langle x_0, z - y \rangle > 0$ for some $x_0 \in \ker A$. In this case, let us consider the polynomial

$$p(t) := \|z - y - tx_0\|^2 = \|x_0\|^2 t^2 - 2\langle x_0, z - y \rangle t + \|z - y\|^2.$$

Clearly, $t = 0$ is not a point of minimum for $p(t)$, as for $t > 0$ and small enough one has $p(t) \leqslant p(0)$. This shows that there are points $y + tx_0 \in \mathcal{S}$ for which $\|z - (y + tx_0)\| \leqslant \|z - y\|$, thus contradicting the assumption that $y$ is the closest to $z$ among all points in $\mathcal{S}$. Then necessarily $z - y \in (\ker A)^\perp$, which proves that (i) $\Rightarrow$ (ii). □

*Proof (Proof of Proposition 3.1)* In the minimisation (3.3)

$$h - f^{[0]} = q_{N-1}(A)(Af^{[0]} - g) = q_{N-1}(A)A(f^{[0]} - P_{\mathcal{S}}f^{[0]})$$

for some polynomial $q_{N-1} \in \mathbb{P}_{N-1}$, whence also

$$h - P_{\mathcal{S}}f^{[0]} = q_{N-1}(A)A(f^{[0]} - P_{\mathcal{S}}f^{[0]}) + (f^{[0]} - P_{\mathcal{S}}f^{[0]}).$$

This implies, upon setting $p_N(\lambda) := \lambda q_{N-1}(\lambda) + 1$, that

$$f^{[N]} - P_{\mathcal{S}}f^{[0]} = p_N(A)(f^{[0]} - P_{\mathcal{S}}f^{[0]}) \qquad \forall N \in \mathbb{N}, \tag{*}$$

where $p_N \in \mathbb{P}_N^{(1)}$.

Moreover, $f^{[N]} - P_{\mathcal{S}}f^{[N]} \in (\ker A)^\perp$, as a consequence of Lemma 3.3 applied to the choice $z = f^{[N]}$ and $y = P_{\mathcal{S}}f^{[N]}$. With an analogous argument, also $f^{[0]} - P_{\mathcal{S}}f^{[0]} \in (\ker A)^\perp$. Thus, (3.8) is proved.

Owing to (3.7) and (3.8), $f^{[0]} - P_{\mathcal{S}}f^{[0]} \in (\ker A)^\perp \cap C^\infty(A)$. Now, $(\ker A)^\perp \cap C^\infty(A)$ is invariant under the action of polynomials of $A$, and therefore owing to (*) we deduce that $f^{[N]} - P_{\mathcal{S}}f^{[0]} \in (\ker A)^\perp$.

Next, let us split

$$P_{\mathcal{S}}f^{[N]} - P_{\mathcal{S}}f^{[0]} = (f^{[N]} - P_{\mathcal{S}}f^{[0]}) - (f^{[N]} - P_{\mathcal{S}}f^{[N]}).$$

Obviously, $P_{\mathcal{S}}f^{[N]} - P_{\mathcal{S}}f^{[0]} \in \ker A$. But in the right-hand side, as just shown, both $f^{[N]} - P_{\mathcal{S}}f^{[0]} \in (\ker A)^\perp$ and $f^{[N]} - P_{\mathcal{S}}f^{[N]} \in (\ker A)^\perp$. So $P_{\mathcal{S}}f^{[N]} - P_{\mathcal{S}}f^{[0]} \in (\ker A)^\perp$. The conclusion is necessarily $P_{\mathcal{S}}f^{[N]} - P_{\mathcal{S}}f^{[0]} = 0$.

This establishes (3.9), by means of which formula (*) above takes also the form of (3.10). □

Let us now prove Lemmas 3.1 and 3.2.

*Proof (Proof of Lemma 3.1)* Part (i) is evident from the fact that $\mathcal{D}(A^{\frac{\theta}{2}}) = \mathcal{H}$ for any $\theta \geqslant 0$, as $A$ is (everywhere-defined and) bounded and non-negative.

Part (ii) is therefore obvious if $\theta' \geqslant 0$. If, instead, $\theta \leqslant \theta' < 0$, then $\operatorname{ran}(A^{-\frac{\theta}{2}}) \subset \operatorname{ran}(A^{-\frac{\theta'}{2}})$, owing again to the boundedness and non-negativity of $A$, so part (ii) is actually valid in general.

If $0 \leqslant \theta \leqslant \theta'$, then

$$u_{\theta'}(x) = A^{\theta'/2}(x - P_{\mathcal{S}}x) = A^{(\theta'-\theta)/2}A^{\theta/2}(x - P_{\mathcal{S}}x) = A^{(\theta'-\theta)/2}u_\theta(x).$$

If instead $\theta < 0 \leqslant \theta'$, then $u_{\theta'}(x) = A^{\theta'/2}(x - P_{\mathcal{S}}x)$ and $A^{-\theta/2}u_\theta(x) = x - P_{\mathcal{S}}x$, whence

$$A^{(\theta'-\theta)/2}u_\theta(x) \;=\; A^{\theta'/2}(x-P_{\mathcal{S}}x) \;=\; u_{\theta'}(x)\,.$$

Last, if $\theta \leqslant \theta' < 0$, then $A^{-\xi/2}u_\xi(x) = x - P_{\mathcal{S}}x$ for both $\xi = \theta$ and $\xi = \theta'$, therefore from

$$x - P_{\mathcal{S}}x = A^{-\theta/2}u_\theta(x) = A^{-\theta'/2}A^{(\theta'-\theta)/2}u_\theta(x) \quad \text{and} \quad A^{-\theta'/2}u_{\theta'}(x) = x - P_{\mathcal{S}}x$$

one deduces that $u_{\theta'}(x) = A^{(\theta'-\theta)/2}u_\theta(x)$. In all possible cases the claimed identity is therefore proved. The inequality $\rho_{\theta'}(x) \leqslant \|A\|_{\mathrm{op}}^{\theta'-\theta}\rho_\theta(x)$ then follows at once from (3.13). This completes the proof of part (iii). $\qquad\square$

*Proof (Proof of Lemma 3.2)* Owing to (3.10) and to the $A$-smoothness of $g$ and $f^{[0]}$, $f^{[N]} - P_{\mathcal{S}}f^{[N]} \in C^\infty(A)$, which by interpolation means in particular that $f^{[N]} - P_{\mathcal{S}}f^{[N]} \in \mathcal{D}(A^{\frac{\sigma}{2}})\ \forall \sigma \geqslant 0$. This proves part (i) of the Lemma.

Assume now that $f^{[0]} \in \mathscr{C}_{A,g}(\sigma)$ for some $\sigma < 0$. In this case (3.10) reads

$$f^{[N]} - P_{\mathcal{S}}f^{[N]} \;=\; p_N(A)(f^{[0]} - P_{\mathcal{S}}f^{[0]}) \;=\; p_N(A)A^{-\frac{\sigma}{2}}u_\sigma(f^{[0]}),$$

thanks to the definition (3.12) of $u_\sigma(f^{[0]})$. Therefore $f^{[N]} - P_{\mathcal{S}}f^{[N]} \in \mathrm{ran}\,(A^{-\frac{\sigma}{2}})$ and, again by (3.12), $u_\sigma(f^{[N]}) = p_N(A)u_\sigma(f^{[0]})$. This proves part (ii). $\qquad\square$

Next, an amount of important results that are measure-theoretic in nature will be established. To this aim, with customary notation [97, Chapt. 4-5], let $E^{(A)}$ denote the projection-valued measure associated with the self-adjoint operator $A$, and let $\mathrm{d}\langle x, E^{(A)}(\lambda)x\rangle$ denote the corresponding scalar measure associated with a vector $x \in \mathcal{H}$. Such measures are supported on $\sigma(A) \subset [0, +\infty)$.

A special role is going to be played by the measure

$$\mathrm{d}\mu_\sigma(\lambda) \;:=\; \mathrm{d}\langle u_\sigma(f^{[0]}), E^{(A)}(\lambda)u_\sigma(f^{[0]})\rangle \tag{3.27}$$

defined under the assumption that $f^{[0]} \in \mathscr{C}_{A,g}(\sigma)$ for a given $\sigma \in \mathbb{R}$. Clearly, by definition, $\mu_\sigma$ is a *finite* measure with

$$\int_{[0,+\infty)} \mathrm{d}\mu_\sigma(\lambda) \;=\; \|u_\sigma(f^{[0]})\|^2\,. \tag{3.28}$$

Two relevant properties of $\mu_\sigma$ are the following.

**Proposition 3.2** *For the given self-adjoint and non-negative operator $A$ on $\mathcal{H}$, and for given $g \in C^\infty(A)$, $\sigma \in \mathbb{R}$, $f^{[0]} \in C^\infty(A) \cap \mathscr{C}_{A,g}(\sigma)$, consider the measure $\mu_\sigma$ defined by (3.27). Then:*

*(i) one has*

$$\mathrm{d}\mu_\sigma(\lambda) \;=\; \lambda^\sigma \mathrm{d}\langle f^{[0]} - P_{\mathcal{S}}f^{[0]}, E^{(A)}(\lambda)(f^{[0]} - P_{\mathcal{S}}f^{[0]})\rangle\,; \tag{3.29}$$

*(ii) the spectral value $\lambda = 0$ is not an atom for $\mu_\sigma$, i.e.,*

$$\mu_\sigma(\{0\}) \;=\; 0\,. \tag{3.30}$$

**Proof** The identity (3.29) when $\sigma \geqslant 0$ follows immediately from the definition (3.27) of $d\mu_\sigma$ and from the definition (3.12) of $u_\sigma(f^{[0]}) = A^{\frac{\sigma}{2}}(f^{[0]} - P_{\mathcal{S}} f^{[0]})$, owing to the property

$$d\langle A^\alpha \psi, E^{(A)}(\lambda) A^\alpha \psi \rangle = \lambda^{2\alpha} d\langle \psi, E^{(A)}(\lambda) \psi \rangle, \qquad \alpha \geqslant 0, \qquad \psi \in \mathcal{D}(A^\alpha).$$

If instead $\sigma < 0$, let us consider the auxiliary measures

$$d\widetilde{\mu}_\sigma(\lambda) := \lambda^{-\sigma} d\mu_\sigma(\lambda), \qquad d\widehat{\mu}_\sigma(\lambda) := d\langle f^{[0]} - P_{\mathcal{S}} f^{[0]}, E^{(A)}(\lambda)(f^{[0]} - P_{\mathcal{S}} f^{[0]}) \rangle.$$

On an arbitrary Borel subset $\Omega \subset [0, +\infty)$ one then has

$$\begin{aligned}\widetilde{\mu}_\sigma(\Omega) &= \int_\Omega \lambda^{-\sigma} d\mu_\sigma(\lambda) = \|E^{(A)}(\Omega) A^{-\frac{\sigma}{2}} u_\sigma(f^{[0]})\|^2 \\ &= \|E^{(A)}(\Omega)(f^{[0]} - P_{\mathcal{S}} f^{[0]})\|^2 = \int_\Omega d\widehat{\mu}_\sigma(\lambda) = \widehat{\mu}_\sigma(\Omega),\end{aligned}$$

having used the definition (3.12) in the form $A^{-\frac{\sigma}{2}} u_\sigma(f^{[0]}) = f^{[0]} - P_{\mathcal{S}} f^{[0]}$. This shows that $d\widetilde{\mu}_\sigma(\lambda) = d\widehat{\mu}_\sigma(\lambda)$, whence again (3.29). Part (i) is proved.

Concerning part (ii), let us recall from (3.8) that $f^{[0]} - P_{\mathcal{S}} f^{[0]} \in (\ker A)^\perp$. Therefore, $\widehat{\mu}_\sigma(\{0\}) = 0$. Thus, (3.29) implies that also $\mu_\sigma(\{0\}) = 0$. $\square$

In turn, Proposition 3.2 allows one to discuss one further set of technical ingredients for the proof of Theorem 3.1. They concern the polynomial $p_N$, in the expression (3.10) of the $\xi$-iterates $f^{[N]}$, that corresponds to the actual minimisation (3.3).

**Proposition 3.3** *For the given self-adjoint and non-negative operator $A$ on $\mathcal{H}$, and for given $g \in C^\infty(A)$, $\sigma \in \mathbb{R}$, $f^{[0]} \in C^\infty(A) \cap \mathscr{C}_{A,g}(\sigma)$, and $\xi \geqslant 0$ let $f^{[N]}$ be the N-th $\xi$-iterate defined by (3.3) with initial guess $f^{[0]}$ and parameter $\theta = \xi$, and let*

$$s_N := \arg\min_{p_N \in \mathbb{P}_N^{(1)}} \int_{[0,+\infty)} \lambda^\xi p_N^2(\lambda) \, d\langle f^{[0]} - P_{\mathcal{S}} f^{[0]}, E^{(A)}(\lambda)(f^{[0]} - P_{\mathcal{S}} f^{[0]}) \rangle \qquad (3.31)$$

*for each $N \in \mathbb{N}$. Then the following properties hold.*

(i) *One has*
$$f^{[N]} - P_{\mathcal{S}} f^{[N]} = s_N(A)(f^{[0]} - P_{\mathcal{S}} f^{[0]}) \qquad \forall N \in \mathbb{N}. \qquad (3.32)$$

(ii) *The family $(s_N)_{N \in \mathbb{N}}$ is a set of orthogonal polynomials on $[0, +\infty)$ with respect to the measure*

$$\begin{aligned}d\nu_\xi(\lambda) &:= \lambda^{\xi - \sigma + 1} d\mu_\sigma(\lambda) \\ &= \lambda^{\xi + 1} d\langle f^{[0]} - P_{\mathcal{S}} f^{[0]}, E^{(A)}(\lambda)(f^{[0]} - P_{\mathcal{S}} f^{[0]}) \rangle\end{aligned} \qquad (3.33)$$

*and satisfying*

$$\deg s_N = N, \qquad s_N(0) = 1 \qquad \forall N \in \mathbb{N} \qquad (3.34)$$

*(under the further tacit assumption that the $s_N$'s are all non-vanishing with respect to the measure $\mu_\sigma$).*

*(iii)* One has
$$\rho_\sigma(f^{[N]}) \;=\; \int_{[0,+\infty)} s_N^2(\lambda)\,\mathrm{d}\mu_\sigma(\lambda) \qquad \forall N \in \mathbb{N}. \tag{3.35}$$

***Proof*** Denote temporarily by $\widetilde{s}_N \in \mathbb{P}_N^{(1)}$ the polynomial that qualifies the iterate $f^{[N]}$ in (3.10) by means of the minimisation (3.3) with $\theta = \xi$. Then

$$\begin{aligned}
\min_{h \in \{f^{[0]}\} + \mathcal{K}_N(A,\mathfrak{R}_0)} \|A^{\xi/2}(h - P_\mathcal{S} h)\|^2 &= \|A^{\xi/2}(f^{[N]} - P_\mathcal{S} f^{[N]})\|^2 \\
&= \|A^{\xi/2}\widetilde{s}_N(A)(f^{[0]} - P_\mathcal{S} f^{[0]})\|^2 \\
&= \int_{[0,+\infty)} \lambda^\xi\, \widetilde{s}_N^2(\lambda)\,\mathrm{d}\langle f^{[0]} - P_\mathcal{S} f^{[0]}, E^{(A)}(\lambda)(f^{[0]} - P_\mathcal{S} f^{[0]})\rangle.
\end{aligned}$$

Comparing the above identity with (3.31) we see that $\widetilde{s}_N$ must be precisely the polynomial $s_N$. Therefore, (3.10) takes the form (3.32). This proves part (i).

By means of (3.29) we may re-write (3.31) as
$$s_N \;=\; \arg\min_{p_N \in \mathbb{P}_N^{(1)}} \int_{[0,+\infty)} \lambda^{\xi-\sigma}\, p_N^2(\lambda)\,\mathrm{d}\mu_\sigma(\lambda).$$

The latter minimising property of $s_N$ implies
$$\begin{aligned}
0 &= \frac{\mathrm{d}}{\mathrm{d}\varepsilon}\bigg|_{\varepsilon=0} \int_{[0,+\infty)} \lambda^{\xi-\sigma} \left(s_N(\lambda) + \varepsilon\lambda\, q_{N-1}(\lambda)\right)^2 \mathrm{d}\mu_\sigma(\lambda) \\
&= 2\int_{[0,+\infty)} \lambda^{\xi-\sigma+1}\, s_N(\lambda)\, q_{N-1}(\lambda)\,\mathrm{d}\mu_\sigma(\lambda)
\end{aligned}$$

for any $q_{N-1} \in \mathbb{P}_{N-1}$ (indeed, $s_N + \varepsilon\lambda\, q_{N-1} \in \mathbb{P}_N^{(1)}$). Equivalently, owing to (3.33),
$$\int_{[0,+\infty)} s_N(\lambda)\, q_{N-1}(\lambda)\,\mathrm{d}\nu_\xi(\lambda) \;=\; 0 \qquad \forall q_{N-1} \in \mathbb{P}_{N-1}.$$

Such a condition is valid for each $N \in \mathbb{N}$ and, as well known [107, 25, 68], this amounts to saying that $(s_N)_{N \in \mathbb{N}}$ is a set of orthogonal polynomials on $[0, +\infty)$ with respect to the measure $\mathrm{d}\nu_\xi$. Part (ii) is thus proved.

If $\sigma \geqslant 0$, then (3.9), (3.16), (3.29), and (3.32) yield
$$\begin{aligned}
\rho_\sigma(f^{[N]}) &= \left\|A^{\frac{\sigma}{2}}(f^{[N]} - P_\mathcal{S} f^{[N]})\right\|^2 = \left\|A^{\frac{\sigma}{2}} s_N(A)(f^{[0]} - P_\mathcal{S} f^{[0]})\right\|^2 \\
&= \int_{[0,+\infty)} s_N^2(\lambda)\,\mathrm{d}\mu_\sigma(\lambda).
\end{aligned}$$

If instead $\sigma < 0$, then owing to (3.32) the identity (3.15) reads
$$u_\sigma(f^{[N]}) \;=\; s_N(A)\, u_\sigma(f^{[0]}).$$

The latter identity, together with (3.16) and (3.27), yield

$$\rho_\sigma(f^{[N]}) = \|u_\sigma(f^{[N]})\|^2 = \|s_N(A)u_\sigma(f^{[0]})\|^2 = \int_{[0,+\infty)} s_N^2(\lambda)\,d\mu_\sigma(\lambda).$$

In either case (3.35) is established. This proves part (iii). □

*Remark 3.6* At this point, one can comment on the *existence* and *uniqueness* of the polynomials $s_N$. As the measure $\lambda^\xi d\langle f^{[0]} - P_S f^{[0]}, E^{(A)}(f^{[0]} - P_S f^{[0]})\rangle$ is positive, there exists a sequence $(Q_N)_{N\in\mathbb{N}}$ of polynomials that are orthogonal with respect to this measure [25, Theorem I.3.3]. Moreover [25, Corollary to Theorem I.2.2], any sequence $(\widetilde{Q}_N)_{N\in\mathbb{N}}$ of polynomials that are orthogonal with respect to the same measure must satisfy $Q_N = c_N \widetilde{Q}_N \ \forall N \in \mathbb{N}$ for some non-zero constants $c_N$. This and the additional requirement that $Q_N, \widetilde{Q}_N \in \mathbb{P}_N^{(1)}$, namely that the $N$-th polynomial attains the value one at $\lambda = 0$, yields $Q_N = \widetilde{Q}_N \ \forall N \in \mathbb{N}$.

*Remark 3.7* The measure $\nu_\xi$ too is finite, with

$$\int_{[0,+\infty)} d\nu_\xi = \|A^{\frac{\xi+1}{2}}(f^{[0]} - P_S f^{[0]})\|^2, \qquad (3.36)$$

as is evident from (3.33). In fact, one could define $\nu_\xi$ for arbitrary $\xi \geqslant -1$: we keep the restriction to $\xi \geqslant 0$ because $\xi$ here is the parameter $\theta = \xi$ required in the definition (3.3) of the $\xi$-iterates, and as such must therefore be non-negative.

*Remark 3.8* There is an implicit dependence on $\xi$ in each $s_N$, as is clear from (3.31), analogously to the fact that the iterates $f^{[N]}$'s depend on the choice of the parameter $\xi$. Such a dependence is simply omitted from the notation $s_N$.

Proposition 3.3(iii) thus shows that the control of the convergence of the $f^{[N]}$'s in the $\rho_\sigma$-sense is boiled down to monitoring a precise spectral integral, namely the right-hand side of (3.35). For an efficient estimate of the latter, suitable properties of the polynomials $s_N$ and of the measure $\nu_\xi$ are needed, which are going to be discussed in the remaining part of this Section.

Here is the main result in this context.

**Proposition 3.4** *Consider the set $(s_N)_{N\in\mathbb{N}}$ of orthogonal polynomials on $[0,+\infty)$ with respect to the measure $\nu_\xi$, as defined in (3.31) and (3.33) under the assumptions of Proposition 3.3.*

(i) *For each $N \in \mathbb{N}$, either $s_N(\lambda) = 0$ $\nu_\xi$-almost everywhere, or $s_N$ has exactly $N$ simple zeroes, all located in $(0,+\infty)$.*

*Assume now the $s_N$'s are all non-vanishing with respect to the $\nu_\xi$-measure, and denote by $\lambda_k^{(N)}$ the k-th zero of $s_N$, ordering the zeros as*

$$0 < \lambda_1^{(N)} < \lambda_2^{(N)} < \cdots < \lambda_N^{(N)}. \qquad (3.37)$$

*(ii) (Separation.)* One has

$$\lambda_k^{(N+1)} < \lambda_k^{(N)} < \lambda_{k+1}^{(N+1)} \qquad \forall k \in \{1,2,\ldots,N-1\}, \qquad (3.38)$$

*that is, the zeroes of $s_N$ and $s_{N+1}$ mutually separate each other.*

*(iii) (Monotonicity.) For each integer $k \geqslant 1$,*

$$\begin{aligned}(\lambda_k^{(N)})_{N=k}^{\infty} &\quad \text{is a decreasing sequence,} \\ (\lambda_{N-k+1}^{(N)})_{N=k}^{\infty} &\quad \text{is an increasing sequence.}\end{aligned} \qquad (3.39)$$

*In particular, the limits*

$$\lambda_1 := \lim_{N \to \infty} \lambda_1^{(N)}, \qquad \lambda_\infty := \lim_{N \to \infty} \lambda_N^{(N)} \qquad (3.40)$$

*exist in $[0,+\infty) \cup \{+\infty\}$.*

*(iv) (Orthogonality.) One has*

$$\int_{[0,\lambda_1^{(N)})} s_N^2(\lambda) \frac{\lambda_1^{(N)}}{\lambda_1^{(N)} - \lambda} \, d\nu_\xi(\lambda) = \int_{[\lambda_1^{(N)},+\infty)} s_N^2(\lambda) \frac{\lambda_1^{(N)}}{\lambda - \lambda_1^{(N)}} \, d\nu_\xi(\lambda) \quad (3.41)$$

*for any $N \in \mathbb{N}$.*

*Finally, assume in addition to the assumptions of Proposition 3.3 also assumptions (3.20)-(3.21) of Theorem 3.1, or assumptions (3.23) of Corollary 3.1. In other words, assume in addition that $f^{[0]}, g \in \mathcal{D}^a(A)$, or also that $f^{[0]} - P_S f^{[0]} \in \mathcal{D}^{qa}(A)$.*

*(v) (Representation.) The measure $\nu_\xi$ is only supported on the so-called 'true interval of orthogonality' $[\lambda_1, \lambda_\infty]$. Here and in the following, the symbol $[\lambda_1, \lambda_\infty]$ is understood as the closure of $(\lambda_1, \lambda_\infty)$.*

Observe that, for the first time, in Proposition 3.4(v) the assumption of *A*-analyticity of $g$ and $f^{[0]}$ kicks in, replacing the mere *A*-smoothness. This is the condition prescribed in the final Theorem 3.1. So, prior to presenting the proof of Proposition 3.4, let us highlight in what form we shall exploit the extra condition of *A*-analyticity of $g$ and $f^{[0]}$.

**Lemma 3.4** *Let A be a linear operator on a Hilbert space $\mathcal{H}$.*

*(i) Assume that $g \in \mathcal{D}^{qa}(A) \cap \mathrm{ran}\, A$. Then any $f \in \mathcal{D}(A)$ such that $Af = g$ satisfies $f \in \mathcal{D}^{qa}(A)$.*
*(ii) Assume that $g \in \mathcal{D}^a(A) \cap \mathrm{ran}\, A$. Then any $f \in \mathcal{D}(A)$ such that $Af = g$ satisfies $f \in \mathcal{D}^a(A)$.*
*(iii) Assume that $g \in \mathcal{D}^a(A) \cap \mathrm{ran}\, A$ and $f^{[0]} \in \mathcal{D}^a(A)$. Then $f^{[0]} - P_S f^{[0]} \in \mathcal{D}^a(A)$.*

***Proof*** (i) As $Af = g$, then

$$\sum_{n=1}^\infty \|A^n f\|^{-\frac{1}{n}} = \|g\|^{-1} + \sum_{n=1}^\infty \|A^n g\|^{-\frac{1}{n+1}}.$$

The latter series, by a standard ratio test (d'Alembert's criterion), is asymptotic to $\sum_{n=1}^{\infty} \|A^n g\|^{-\frac{1}{n}}$ and hence diverges because $g$ is quasi-analytic (see (3.19) above). Then also $\sum_{n=1}^{\infty} \|A^n f\|^{-\frac{1}{n}} = +\infty$, whence the quasi-analyticity of $f$.

(ii) As $Af = g$, then $A^{n-1}g = A^n f$ for any integer $n \geqslant 1$. By definition of $A$-analyticity of $g$ (see (3.18) above), there is $C_g > 0$ such that

$$\|A^n f\| = \|A^{n-1} g\| \leqslant C_g^{n-1}(n-1)! \leqslant D_f^n n!, \qquad n \geqslant 2,$$

having set $D_f := \max\{1, C_g\}$. The latter inequality is due to $D_f \geqslant C_g$, whence $D_f^{n-1} \geqslant C_g^{n-1}$, and to $D_f \geqslant 1$, whence $D_f^n \geqslant D_f^{n-1}$. As $\|Af\| = \|g\|$ (the $n=1$ case), then setting $C_f := \max\{\|g\|, D_f\}$ finally yields

$$\|A^n f\| \leqslant C_f^n n!, \qquad n \geqslant 1,$$

which in view of (3.18) expresses the $A$-analyticity of $f$.

(iii) On account of part (ii), any solution $f$ belongs to $\mathcal{D}^a(A)$. Then in particular $P_g f^{[0]} \in \mathcal{D}^a(A)$, and since $\mathcal{D}^a(A)$ is a linear subspace, the conclusion follows by linearity. $\square$

**Lemma 3.5** *Given $A = A^* \geqslant \mathbb{O}$,*

$$A^\beta \mathcal{D}^{qa}(A) \subset \mathcal{D}^{qa}(A) \qquad \forall \beta \geqslant 0. \tag{3.42}$$

***Proof*** We intend to apply this simple property (see, e.g., [97, Lemma 7.17]):

> if $S$ and $T$ are two densely defined operators with common domain $\mathcal{D}$ and such that $T\mathcal{D} \subset \mathcal{D}$, $S\mathcal{D} \subset \mathcal{D}$, and $TS = ST$ on $\mathcal{D}$, then $S\mathcal{D}^{qa}(T) \subset \mathcal{D}^{qa}(T)$. (*)

In the present case let us take

$$S := A^\beta \big|_{C^\infty(A)}, \qquad T := A \big|_{C^\infty(A)}, \qquad \mathcal{D} := C^\infty(A).$$

With this choice, obviously, $C^\infty(T) = C^\infty(A)$, whence also, owing to the definition (3.19), $\mathcal{D}^{qa}(T) = \mathcal{D}^{qa}(A)$. So, provided that all assumptions of (*) are matched, the conclusion $S\mathcal{D}^{qa}(T) \subset \mathcal{D}^{qa}(T)$ amounts precisely to (3.42). Concerning the assumptions of (*), it is clear that both $T$ and $S$ are symmetric and densely defined, with common domain $\mathcal{D}$. The invariance properties $T\mathcal{D} \subset \mathcal{D}$ and $S\mathcal{D} \subset \mathcal{D}$ are tantamount as $A^\tau C^\infty(A) \subset C^\infty(A)$, respectively with $\tau = 1$ and $\tau = \beta$, and in either case they follow from the fact that for every $h \in C^\infty(A)$ and any $k \in \mathbb{N}$, the vector $A^\tau h$ satisfies

$$\begin{aligned}
\|A^k A^\tau h\|^2 &= \int_{[0,+\infty)} \lambda^{2(k+\tau)} \, d\mu_h^{(A)}(\lambda) \\
&\leqslant \int_{[0,1)} d\mu_h^{(A)}(\lambda) + \int_{[1,+\infty)} \lambda^{2(k+\tau)} \, d\mu_h^{(A)}(\lambda) \\
&\leqslant \|h\|^2 + \|A^{k+\lceil \tau \rceil} h\|^2 < +\infty,
\end{aligned}$$

where $\lceil\tau\rceil$ is the smallest integer greater than $\tau$. Last, the commutativity of $S$ and $T$ on $\mathcal{D}$ is obviously tantamount as $AA^\beta h = A^{1+\beta}h = A^\beta Ah$ for $h \in \mathcal{D}$. All assumptions of (*) are verified, and the Lemma is proved. $\square$

On a related note, for completeness and later use, let us also recall this simple property.

**Lemma 3.6** *For any operator $A$ on a Hilbert space $\mathcal{H}$, $A\mathcal{D}^a(A) \subset \mathcal{D}^a(A)$*

***Proof*** Let $f \in \mathcal{D}^a(A)$ and $g := Af$. Then

$$\|A^n g\| = \|A^{n+1}f\| \leqslant C_f^{n+1}(n+1)!$$

for some $C_f > 0$ and for all $n \in \mathbb{N}_0$. Set $C_g := 2(\max\{C_f, 1\})^2$ and take $n \in \mathbb{N}$. Then $C_f^{1+\frac{1}{n}}(1+n)^{\frac{1}{n}} \leqslant (\max\{C_f, 1\})^2 \cdot 2 = C_g$, whence

$$\|A^n g\| \leqslant C_f^{n+1}(n+1)! \leqslant C_g^n n! \qquad \forall n \in \mathbb{N},$$

which shows that $g \in \mathcal{D}^a(A)$. $\square$

***Proof (Proof of Proposition 3.4)*** Part (i) is standard from the theory of orthogonal polynomials (see, e.g., [107, Theorem 3.3.1] or [25, Theorem I.5.2]), owing to the fact that the map

$$\mathbb{P}([0,+\infty)) \ni p \longmapsto \int_{[0,+\infty)} p(\lambda)\,d\nu_\xi(\lambda)$$

is a positive-definite functional on $\mathbb{P}([0,+\infty))$.

Part (ii) is another standard fact in the theory of orthogonal polynomials (see, e.g., [107, Theorem 3.3.2] or [25, Theorem I.5.3]). Part (iii), in turn, is an immediate corollary of part (ii).

Part (iv) follows from the identity

$$\int_{[0,+\infty)} s_N(\lambda)q_{N-1}(\lambda)\,d\nu_\xi(\lambda) = 0 \qquad \forall q_{N-1} \in \mathbb{P}_{N-1}$$

(already considered in the proof of Proposition 3.3, as a consequence of the orthogonality of the $s_N$'s), when the explicit choice

$$q_{N-1}(\lambda) := \frac{\lambda_1^{(N)} s_N(\lambda)}{\lambda_1^{(N)} - \lambda}$$

is made.

For Part (v) let us first recall [25, Definition I.5.2] that the true interval of orthogonality $[\lambda_1, \lambda_\infty]$ is the smallest closed interval containing all the zeroes $\lambda_k^{(N)}$, and moreover [25, Theorem II.3.1] there exists a measure $\eta$ on $[0,+\infty)$ supported only on $[\lambda_1, \lambda_\infty]$ such that the $s_N$'s remain orthogonal with respect to $\eta$ too and

$$\mu_k := \int_{[0,+\infty)} \lambda^k \, d\nu_\xi(\lambda) = \int_{[\lambda_1,\lambda_\infty]} \lambda^k \, d\eta(\lambda), \qquad \forall k \in \mathbb{N}_0.$$

Such $\eta$-measure is actually a Stieltjes measure associated with a bounded, non-decreasing function $\psi$ obtained as point-wise limit of a sub-sequence of $(\psi_N)_{N\in\mathbb{N}}$, where

$$\psi_N(\lambda) := \begin{cases} 0, & \lambda < \lambda_1^{(N)}, \\ A_1^{(N)} + \cdots + A_p^{(N)}, & \lambda \in [\lambda_p^{(N)}, \lambda_{p+1}^{(N)}) \quad \text{for} \quad p \in \{1,\ldots,n-1\}, \\ \mu_0, & \lambda \geq \lambda_N^{(N)} \end{cases}$$

and $A_1^{(n)},\ldots,A_N^{(n)}$ are positive numbers determined by the Gauss quadrature formula

$$\mu_k = \sum_{p=1}^N A_p^{(N)} (\lambda_p^{(N)})^k, \qquad \forall k \in \{0,1,\ldots,2N-1\}.$$

We want to show that $\nu_\xi = \eta$, i.e., that the Hamburger moment problem that guarantees that $(s_N)_{N\in\mathbb{N}}$ is an orthogonal system on $[0,+\infty)$ is *uniquely* solved with the measure $\nu_\xi$. To this aim, let us re-write the *even* moments of $\nu_\xi$ as

$$\mu_{2k} = \int_{[0,+\infty)} \lambda^{2k} \lambda^{\xi+1} \, d\langle f^{[0]} - P_{\mathcal{S}} f^{[0]}, E^{(A)}(\lambda)(f^{[0]} - P_{\mathcal{S}} f^{[0]})\rangle$$
$$= \|A^k A^{\frac{\xi+1}{2}} (f^{[0]} - P_{\mathcal{S}} f^{[0]})\|^2 = \|A^k \phi\|^2,$$

having set $\phi := A^{\frac{\xi+1}{2}}(f^{[0]} - P_{\mathcal{S}} f^{[0]})$. The extra assumptions made for this part ensure that $f^{[0]} - P_{\mathcal{S}} f^{[0]} \in \mathcal{D}^a(A)$, on account of Lemma 3.4, or directly that $f^{[0]} - P_{\mathcal{S}} f^{[0]} \in \mathcal{D}^{qa}(A)$. As a consequence, owing to Lemma 3.5, $\phi \in \mathcal{D}^{qa}(A)$. The quasi-analyticity of $\phi$ then implies (see (3.19) above)

$$\sum_{k=1}^\infty \mu_{2k}^{-\frac{1}{2k}} = \sum_{k=1}^\infty \|A^k \phi\|^{-\frac{1}{k}} = +\infty.$$

Now, the divergence of the above series $\sum_{k=1}^\infty \mu_{2k}^{-\frac{1}{2k}}$ is a well-known sufficient condition (Carleman's criterion, see, e.g., [99, Theorem I.10]) for the uniqueness of the Hamburger moment problem's solution. This shows that $\nu_\xi = \eta$, thus proving that $\nu_\xi$ is supported only on $[\lambda_1, \lambda_\infty]$. □

*Remark 3.9* Analogously to what was already observed in Remark 3.8, there is an implicit dependence on $\xi$ of all the zeroes $\lambda_k^{(N)}$. For a more compact notation, such a dependence is omitted.

In view of Proposition 3.4(i), when the $s_N$'s are not identically zero one can explicitly represent

$$s_N(\lambda) = \prod_{k=1}^N \left(1 - \frac{\lambda}{\lambda_k^{(N)}}\right). \tag{3.43}$$

The integral (3.41) is going to play a central role in the main proof, so the next technical result we need is the following efficient estimate of such a quantity.

**Lemma 3.7** *Consider the set* $(s_N)_{N\in\mathbb{N}}$ *of orthogonal polynomials on* $[0,+\infty)$ *with respect to the measure* $v_\xi$, *as defined in* (3.31) *and* (3.33) *under the assumptions of Proposition 3.3 and with the further restriction* $\xi - \sigma + 1 \geqslant 0$. *Assume that the* $s_N$'s *are non-zero polynomials with respect to the measure* $v_\xi$. *Then, for any* $N \in \mathbb{N}$,

$$\int_{(\lambda_1,\lambda_1^{(N)})} s_N^2(\lambda) \frac{\lambda_1^{(N)}}{\lambda_1^{(N)} - \lambda} \, dv_\xi(\lambda) \leqslant \mu_\sigma((\lambda_1, \lambda_1^{(N)})) \left( \frac{\xi - \sigma + 1}{\delta_N} \right)^{\xi - \sigma + 1}, \quad (3.44)$$

*where*

$$\delta_N := \frac{1}{\lambda_1^{(N)}} + 2 \sum_{k=2}^N \frac{1}{\lambda_k^{(N)}}. \quad (3.45)$$

*Remark 3.10* Estimate (3.44) provides a $(\xi, \sigma)$-dependent bound on a quantity that is $\xi$-dependent only. This is only possible for a constrained range of $\sigma$, namely $\sigma \leqslant \xi + 1$.

*Proof (Proof of Lemma 3.7)* For each $N \in \mathbb{N}$, the function

$$[0, \lambda_1^{(N)}) \ni \lambda \longmapsto a_N(\lambda) := \frac{\lambda_1^{(N)} \lambda^{\xi - \sigma + 1} s_N^2(\lambda)}{\lambda_1^{(N)} - \lambda}$$

$$= \lambda^{\xi - \sigma + 1} \left( 1 - \frac{\lambda}{\lambda_1^{(N)}} \right) \prod_{k=2}^N \left( 1 - \frac{\lambda}{\lambda_k^{(N)}} \right)^2$$

(where the representation (3.43) for $s_N$ was used) is non-negative, smooth, and such that $a_N(0) = \lim_{\lambda \to \lambda_1^{(N)}} a_N(\lambda) = 0$. Let $\lambda_N^* \in (0, \lambda_1^{(N)})$ be the point of maximum for $a_N$. Then $a_N'(\lambda_N^*) = 0$, which after straightforward computations yields

$$\xi - \sigma + 1 \geqslant \lambda_N^* \left( \frac{1}{\lambda_1^{(N)}} + 2 \sum_{k=2}^N \frac{1}{\lambda_k^{(N)}} \right) = \lambda_N^* \delta_N,$$

whence also

$$\lambda_N^* \leqslant \frac{\xi - \sigma + 1}{\delta_N}.$$

Moreover, $0 \leqslant 1 - \lambda/\lambda_k^{(N)} \leqslant 1$ for $\lambda \in [0, \lambda_1^{(N)})$ and for all $k \in \{1, \ldots, N\}$, as $\lambda_1^{(N)}$ is the smallest zero of $s_N$. Therefore,

$$a_N(\lambda) \leqslant a_N(\lambda_N^*) \leqslant (\lambda_N^*)^{\xi - \sigma + 1} \leqslant \left( \frac{\xi - \sigma + 1}{\delta_N} \right)^{\xi - \sigma + 1}, \qquad \lambda \in [0, \lambda_1^{(N)}).$$

One then concludes

$$\int_{(\lambda_1,\lambda_1^{(N)})} s_N^2(\lambda)\frac{\lambda_1^{(N)}}{\lambda_1^{(N)}-\lambda}\,\mathrm{d}\nu_\xi(\lambda) = \int_{(\lambda_1,\lambda_1^{(N)})} a_N(\lambda)\,\mathrm{d}\mu_\sigma(\lambda)$$

$$\leqslant \mu_\sigma((\lambda_1,\lambda_1^{(N)}))\left(\frac{\xi-\sigma+1}{\delta_N}\right)^{\xi-\sigma+1},$$

which completes the proof. □

## 3.4 Proof of CG-convergence and additional observations

Let us present in this Section the proof of the CG-convergence, Theorem 3.1 and Corollary 3.1, based on the intermediate results established in the previous Section.

Owing to Proposition 3.3, one has to control the behaviour for large $N$ of the quantity

$$\rho_\sigma(f^{[N]}) = \int_{[0,+\infty)} s_N^2(\lambda)\,\mathrm{d}\mu_\sigma(\lambda).$$

Obviously, in the following it is assumed that none of the polynomials $s_N$ vanish with respect to the measure $\nu_\xi$ previously introduced in (3.33), otherwise for some $N$ one would have $\rho_\sigma(f^{[N]}) = 0$ and therefore $f^{[N]} = P_S f^{[0]} \in S$ (see Remark 3.1, or also (3.32)), meaning that the conjugate gradient algorithm has come to convergence in a finite number of steps. The conclusion of Theorem 3.1 would then be trivially true.

Let us first observe, from the relation (3.33) between the measures $\mu_\sigma$ and $\nu_\xi$ and from the fact that the latter is supported on the true interval of orthogonality $[\lambda_1,\lambda_\infty]$ (Proposition 3.4(v)), that the measure $\mu_\sigma$ too is supported on such an interval. Thus, in practice,

$$\rho_\sigma(f^{[N]}) = \int_{[\lambda_1,\lambda_\infty]} s_N^2(\lambda)\,\mathrm{d}\mu_\sigma(\lambda). \tag{3.46}$$

(Let us recall that $[\lambda_1,\lambda_\infty]$ is a shorthand for the closure of $(\lambda_1,\lambda_\infty)$, even when $\lambda_\infty = +\infty$.)

It is convenient to split

$$\int_{[\lambda_1,\lambda_\infty]} s_N^2(\lambda)\,\mathrm{d}\mu_\sigma(\lambda) = \mu_\sigma(\{\lambda_1\})s_N^2(\lambda_1) + \int_{(\lambda_1,\lambda_1^{(N)})} s_N^2(\lambda)\,\mathrm{d}\mu_\sigma(\lambda)$$

$$+ \int_{[\lambda_1^{(N)},+\infty)} s_N^2(\lambda)\,\mathrm{d}\mu_\sigma(\lambda) \tag{3.47}$$

$$\leqslant \mu_\sigma(\{\lambda_1\})s_N^2(\lambda_1) + \mu_\sigma((\lambda_1,\lambda_1^{(N)})) + \int_{[\lambda_1^{(N)},+\infty)} s_N^2(\lambda)\,\mathrm{d}\mu_\sigma(\lambda).$$

Here the bound $s_N^2(\lambda) \leqslant 1$, $\lambda \in [0,\lambda_1^{(N)})$, was used, which is obvious from (3.43).

Next, let us show that

$$\int_{[\lambda_1^{(N)},+\infty)} s_N^2(\lambda)\,\mathrm{d}\mu_\sigma(\lambda) \leqslant \frac{1}{(\lambda_1^{(N)})^{\xi-\sigma+1}} \int_{[0,\lambda_1^{(N)})} s_N^2(\lambda)\,\frac{\lambda_1^{(N)}}{\lambda_1^{(N)}-\lambda}\,\mathrm{d}\nu_\xi(\lambda). \quad (3.48)$$

In fact, (3.48) is a consequence of the properties of $s_N$ discussed in Section 3.3. To see that, consider the inequality

$$\begin{aligned}
1 \leqslant \Big(\frac{\lambda}{\lambda_1^{(N)}}\Big)^{\xi-\sigma} &= \frac{1}{(\lambda_1^{(N)})^{\xi-\sigma+1}}\cdot\frac{\lambda_1^{(N)}}{\lambda}\cdot\lambda^{\xi-\sigma+1}\\
&\leqslant \frac{1}{(\lambda_1^{(N)})^{\xi-\sigma+1}}\cdot\frac{\lambda_1^{(N)}}{\lambda-\lambda_1^{(N)}}\cdot\lambda^{\xi-\sigma+1} \qquad (\lambda\geqslant\lambda_1^{(N)}),
\end{aligned} \quad (3.49)$$

which is valid owing to the constraint $\xi-\sigma\geqslant 0$. Then,

$$\begin{aligned}
\int_{[\lambda_1^{(N)},+\infty)} s_N^2(\lambda)\,\mathrm{d}\mu_\sigma(\lambda) &\leqslant \frac{1}{(\lambda_1^{(N)})^{\xi-\sigma+1}} \int_{[\lambda_1^{(N)},+\infty)} s_N^2(\lambda)\,\frac{\lambda_1^{(N)}}{\lambda-\lambda_1^{(N)}}\,\mathrm{d}\nu_\xi(\lambda)\\
&= \frac{1}{(\lambda_1^{(N)})^{\xi-\sigma+1}} \int_{[0,\lambda_1^{(N)})} s_N^2(\lambda)\,\frac{\lambda_1^{(N)}}{\lambda_1^{(N)}-\lambda}\,\mathrm{d}\nu_\xi(\lambda),
\end{aligned}$$

having used (3.33) and (3.49) in the first step, and the orthogonality property (3.41) in the second. Estimate (3.48) is thus proved.

In turn, from (3.48) one gets

$$\begin{aligned}
\int_{[\lambda_1^{(N)},+\infty)} s_N^2(\lambda)\,\mathrm{d}\mu_\sigma(\lambda) &\leqslant \frac{\lambda_1^{(N)} s_N^2(\lambda_1)}{\lambda_1^{(N)}-\lambda_1}\,\frac{\nu_\xi(\{\lambda_1\})}{(\lambda_1^{(N)})^{\xi-\sigma+1}}\\
&\quad + \frac{1}{(\lambda_1^{(N)})^{\xi-\sigma+1}} \int_{(\lambda_1,\lambda_1^{(N)})} s_N^2(\lambda)\,\frac{\lambda_1^{(N)}}{\lambda_1^{(N)}-\lambda}\,\mathrm{d}\nu_\xi(\lambda)\\
&\leqslant \frac{\lambda_1^{(N)} s_N^2(\lambda_1)}{\lambda_1^{(N)}-\lambda_1}\,\mu_\sigma(\{\lambda_1\}) + \Big(\frac{\xi-\sigma+1}{\lambda_1^{(N)}\delta_N}\Big)^{\xi-\sigma+1}\mu_\sigma((\lambda_1,\lambda_1^{(N)}))\\
&\leqslant \frac{\lambda_1^{(N)} s_N^2(\lambda_1)}{\lambda_1^{(N)}-\lambda_1}\,\mu_\sigma(\{\lambda_1\}) + (\xi-\sigma+1)^{\xi-\sigma+1}\mu_\sigma((\lambda_1,\lambda_1^{(N)})),
\end{aligned} \quad (3.50)$$

where formula (3.33) was used to pass from $\nu_\xi$ to $\mu_\sigma$ and Lemma 3.7 was applied, and in the final inequality the property $\lambda_1^{(N)}\delta_N\geqslant 1$ (following from (3.45)) was used.

Thus, (3.46), (3.47) and (3.50) yield

$$\begin{aligned}
\rho_\sigma(f^{[N]}) &\leqslant \Big(s_N^2(\lambda_1) + \frac{\lambda_1^{(N)} s_N^2(\lambda_1)}{\lambda_1^{(N)}-\lambda_1}\Big)\mu_\sigma(\{\lambda_1\})\\
&\quad + \Big(1+(\xi-\sigma+1)^{\xi-\sigma+1}\Big)\mu_\sigma((\lambda_1,\lambda_1^{(N)})),
\end{aligned}$$

whence also, using the factorisation (3.43) for $s_N$,

$$\rho_\sigma(f^{[N]}) \leqslant 2\left(1 - \frac{\lambda_1}{\lambda_1^{(N)}}\right) \prod_{k=2}^{N}\left(1 - \frac{\lambda_1}{\lambda_k^{(N)}}\right)^2 \mu_\sigma(\{\lambda_1\}) \\ + \left(1 + (\xi - \sigma + 1)^{\xi - \sigma + 1}\right) \mu_\sigma((\lambda_1, \lambda_1^{(N)})). \tag{3.51}$$

In the right-hand side of (3.51) one has $\mu_\sigma((\lambda_1, \lambda_1^{(N)})) \xrightarrow{N \to \infty} 0$. Moreover, depending on the value of $\lambda_1$, the quantity

$$\left(1 - \frac{\lambda_1}{\lambda_1^{(N)}}\right) \prod_{k=2}^{N}\left(1 - \frac{\lambda_1}{\lambda_k^{(N)}}\right)^2 \mu_\sigma(\{\lambda_1\})$$

either attains at every $N$ the value $\mu_\sigma(\{0\})$, if $\lambda_1 = 0$, and hence vanishes, owing to (3.30) from Proposition 3.2, or in general is bounded by

$$\frac{\lambda_1^{(N)} - \lambda_1}{\lambda_1^{(N)}} \|u_\sigma(f^{[0]})\|^2,$$

owing to (3.28) and to the ordering $0 < \lambda_1^{(N)} < \lambda_2^{(N)} < \cdots < \lambda_N^{(N)}$ and $\lambda_1 < \lambda_1^{(N)}$, and hence when $\lambda_1 > 0$ it vanishes in the limit $N \to \infty$.

In either case one concludes from (3.51) that $\rho_\sigma(f^{[N]}) \xrightarrow{N \to \infty} 0$, thus completing the proof of Theorem 3.1, and also of Corollary 3.1, as our crucial Proposition 3.4 is proved under the assumptions of either of them.

In the second part of this Section, we intend to highlight a number of important observations.

*Remark 3.11* The assumption that none of the polynomials $s_N$ vanish with respect to the measure $v_\xi$ in the proof of Theorem 3.1 immediately excludes the possibility that $\lambda_1^{(N)} = \lambda_1$ for any $N \in \mathbb{N}$ by considerations in Proposition 3.4(ii). Clearly $\lambda_1^{(N)} \neq 0$ for any $N \in \mathbb{N}$ too owing to Proposition 3.4(i).

*Remark 3.12* In retrospect, the assumption $\xi \geqslant \sigma$ was necessary to establish the bound (3.48) – more precisely, the inequality (3.49). In the step (3.51) (which is an application of Lemma 3.7), only the less restrictive assumption $\xi \geqslant \sigma - 1$ was needed.

*Remark 3.13* Where exactly the true interval of orthogonality lies within $[0, +\infty)$ depends on the behaviour of the zeroes of the $s_N$'s. In particular, in terms of the quantity $\delta_N$ defined in (3.45) we distinguish two alternative scenarios:

CASE I: $\delta_N \to +\infty$ as $N \to \infty$;
CASE II: $\delta_N$ remains uniformly bounded, strictly above 0, in $N$.

If the operator $A$ is bounded, then we are surely in Case I: indeed the orthogonal polynomials $s_N$ are defined on $\sigma(A) \subset [0, \|A\|_{\mathrm{op}}]$, and their zeroes cannot exceed $\|A\|_{\mathrm{op}}$: this forces $\delta_N$ to blow up with $N$. Moreover, $\lambda_\infty = \lim_{N \to \infty} \lambda_N^{(N)} < +\infty$.

If instead $A$ is unbounded, the $\lambda_k^{(N)}$'s fall in $[0,+\infty)$ and depending on their rate of possible accumulation at infinity $\delta_N$ may still diverge as $N \to \infty$ or stay bounded.

Clearly in Case II one has $\lambda_1 > 0$ and $\lambda_\infty = +\infty$, for otherwise the condition $\lambda_1 = \lim_{N\to\infty} \lambda_1^{(N)} = 0$ or $\lambda_\infty = \lim_{N\to\infty} \lambda_N^{(N)} < +\infty$ would necessarily imply $\delta_N \to +\infty$. Thus, in Case II the true interval of orthogonality is $[\lambda_1, +\infty)$ and it is separated from zero.

*Remark 3.14* Estimate (3.51) in the proof and the reasoning thereafter show that the vanishing rate of $\rho_\sigma(f^{[N]})$ is actually controlled by the vanishing rate of the quantity $\mu_\sigma((\lambda_1, \lambda_1^{(N)}))$ if $\lambda_1 = 0$, or more generally of both quantities $(\lambda_1^{(N)} - \lambda_1)$ and $\mu_\sigma((\lambda_1, \lambda_1^{(N)}))$ if $\lambda_1 > 0$. It is however unclear how to possibly quantify, in the above senses, the pace of $\lambda_1^{(N)} \to \lambda_1$. Let us recall (see Remark 3.2 and (3.25) in particular) that the Nemirovskiy-Polyak analysis [81] for the bounded-$A$ case provides an explicit vanishing rate for $\rho_{\sigma'}(f^{[N]})$ for any $\sigma' \in (\sigma, \xi]$, based on a polynomial min-max argument that relies crucially on the *finiteness* of the interval where the orthogonal polynomials $s_N$ are supported on (i.e., it relies on the boundedness of $\sigma(A)$). Therefore, there is certainly no room for applying the same argument to the present setting. In fact, it is reasonable to expect that for generic (unbounded) $A$ the quantity $\rho_\sigma(f^{[N]})$ vanishes with arbitrarily slow pace depending on the choice of the initial guess $f^{[0]}$. A strong indication in this sense comes from the numerical tests discussed in Section 3.5.

*Remark 3.15* It is worth pointing out that removing from the hypotheses of Theorem 3.1 (respectively Corollary 3.1) the $A$-analyticity of $g$ and $f^{[0]}$ (respectively, the quasi-analyticity assumptions (3.23) or (3.24)) and replacing it with just the minimal assumption of $A$-smoothness, one could have only come to the (still non-trivial, yet not-informative) conclusion that $\rho_\sigma(f^{[N]}) \leqslant \kappa$ uniformly in $N$ for some $\kappa > 0$. This is seen as follows. For sure, even if the moment problem for $\nu_\xi$ is indeterminate, the measure $\nu_\xi$ has some support within $[\lambda_1, \lambda_\infty]$ (see, e.g., [25, Theorem II.3.2]), and so does $\mu_\sigma$. However, in the lack of the information that $\mu_\sigma$ is *only* supported in $[\lambda_1, \lambda_\infty]$, in the above proof one should additionally estimate, besides the vanishing quantity (3.47), the extra term

$$\int_{[0,\lambda_1)} s_N^2(\lambda)\,d\mu_\sigma(\lambda)\,.$$

On account of the inequalities $\lambda_1 \leqslant \lambda_1^{(N)}$ (Proposition 3.4(iii)) and $s_N(\lambda) \leqslant 1\ \forall \lambda \in [0, \lambda_1)$ (representation (3.43)), the above integral is controlled by $\int_{[0,\lambda_1)} d\mu_\sigma(\lambda)$, and is therefore bounded uniformly in $N$.

*Remark 3.16 (**Comparison with the proof of [81] valid for bounded** $A$)*

Theorem 3.1 and its proof generalise the Nemirovskiy-Polyak analysis [81] with a crucial technical novelty that is necessary when $A$ is unbounded, and in fact it also yields a subtle improvement of the old argument for the bounded case. More precisely, in [81] one does *not* make use of the very useful property that $\mu_\sigma$ is only supported on $[\lambda_1,\lambda_\infty]$, which is the outcome of the somewhat laborious path that led to Proposition 3.4(v) here. The sole measure-theoretic information used in [81] is that $\lambda = 0$ is not an atom for $\mu_\sigma$. Then in [81], instead of naturally splitting the integration as in (3.47) above, one separates the small and the large spectral values at a threshold $\gamma_N = \min\{\lambda_1^{(N)}, \delta_N^{-1/2}\}$. Clearly $\gamma_N \to 0$, because $\delta_N \to +\infty$ since $A$ is bounded (see Remark 3.13 above), and through a somewhat lengthy analysis of the integration for $\lambda < \gamma_N$ and $\lambda \geqslant \gamma_N$ one reduces both integrations to one over $[0,\gamma_N)$. Then one finally pulls out the upper bound $\mu_\sigma([0,\gamma_N))$, which vanishes as $N \to \infty$ precisely because $\mu_\sigma$ is atom-less at $\lambda = 0$. In the unbounded case such a scheme cannot work: $\delta_N$ does not necessarily diverge and only the information that $\mu_\sigma$ is supported at the right of, and possibly at, $\lambda_1$ makes the final estimate meaningful. Furthermore, in retrospect, by splitting the integration as in (3.47) and not in the old manner of [81], our proof shortens the overall argument and applies both to the bounded and to the unbounded case, with no need to introduce the $\gamma_N$ cut-off.

*Remark 3.17 (**Continuation: comparison with subsequent surveys of [81]**)*

The analysis of conjugate gradients in the bounded case is nicely revisited by Hanke in the monograph [56], both by presenting a version of the same Nemirovskiy-Polyak $\gamma_N$-argument [81], and by relying on a dominated convergence argument for a choice of a sequence of polynomials that vanishes point-wise over $(0,1]$ or on a Banach-Steinhaus uniform boundedness argument. For obvious reasons, none of such schemes are exportable to the present unbounded setting by merely updating the assumption of $g$ and $f^{[0]}$ so as to be $A$-smooth: in order to deal with a measure supported on $[\lambda_1,\lambda_\infty]$, an interval that is possibly infinite and separated from zero, the additional measure-theoretic analysis of Proposition 3.4 is needed. (One should not be misled when in [56] certain spectral integrals appear to run over the whole positive half-line: it is clear from the discussion therein that the zeroes of the considered orthogonal polynomials only fall within a *bounded* interval.) Of course, as commented already, in the present generalised scheme one pays the price that any quantitative bound on the rate of convergence is lost.

We conclude this Section with one further important fact that was alluded to in Remark 3.5.

**Proposition 3.5** *Let $\xi \geqslant 0$, and with respect to the Hilbert space $\mathcal{H} = L^2(\mathbb{R},\mathrm{d}x)$ let $A$ be the self-adjoint multiplication by $x^2$, and let $g := x^2 f$ with*

$$f(x) := \frac{\sqrt{2}\,\mathbf{1}_{\mathbb{R}^+}(x)}{x^{\frac{3}{2}+\xi}(2\pi)^{\frac{1}{4}}}\, e^{-\frac{1}{4}(\log x^2)^2}. \qquad (3.52)$$

*Set further $f^{[0]} \equiv 0$. Then:*

 (i) *$A$ is non-negative and $g \in C^\infty(A) \cap \mathrm{ran}\,(A)$;*
 (ii) *Neither $f$ nor $g$ are quasi-analytic for $A$;*
(iii) *$P_\mathcal{S} f^{[0]} = f$;*

*(iv) the measure $\nu_\xi$ defined in (3.33) is a log-normal distribution, i.e.,*

$$d\nu_\xi(\lambda) \;=\; \frac{\mathbf{1}_{\mathbb{R}^+}(\lambda)}{\lambda\sqrt{2\pi}}\, e^{-\frac{1}{2}(\log\lambda)^2}\,d\lambda\,, \qquad \lambda \geqslant 0.$$

*(v) The Hamburger moment problem for $\nu_\xi$ is indeterminate.*

Proposition 3.5 shows that in general, when $g$ and $f^{[0]}$ are only assumed to be $A$-smooth, and the vector $f^{[0]} - P_S f^{[0]}$ is not necessarily quasi-analytic for $A$, the measure $\nu_\xi$ may fail to be supported *entirely* in the true interval of orthogonality $[\lambda_1, \lambda_\infty]$. (Moreover, let us recall – see, e.g., [25, Exercise II.5.7] – that in the lack of unique solution to the moment problem discussed in the proof of Proposition 3.4, at least one representative measure has a part of its support *outside* $[\lambda_1, \lambda_\infty]$.) As a consequence, the quantity $\rho_\sigma(f^{[N]})$, while staying uniformly bounded (Remark 3.15) is not guaranteed to vanish as $N \to \infty$.

***Proof (Proof of Proposition 3.5)*** The facts that $A \geqslant \mathbb{O}$ and $g \in \operatorname{ran}(A)$ (provided that $f \in \mathcal{D}(A)$) are obvious. Let us prove that $f \in C^\infty(A)$ (whence $f \in \mathcal{D}(A)$ and $g \in C^\infty(A)$). For $n \in \mathbb{N}_0$, and with the change of variable $y := \log x^2$, one computes

$$\begin{aligned}
\|A^n f\|_{L^2}^2 \;&=\; \|x^{2n} f\|_{L^2} \;=\; \frac{2}{\sqrt{2\pi}}\int_0^{+\infty} x^{4n-3-2\xi}\, e^{-\frac{1}{2}(\log x^2)^2}\,dx \\
&=\; \frac{1}{\sqrt{2\pi}}\int_{\mathbb{R}} e^{(2n-1-\xi)y}\, e^{-\frac{1}{2}y^2}\,dy \;=\; e^{\frac{1}{2}(2n-1-\xi)^2} \;<\; +\infty.
\end{aligned}$$

Thus, $f \in C^\infty(A)$ and part (i) is proved.

From the latter computation one also finds

$$\sum_{n\in\mathbb{N}} \|A^n f\|_{L^2}^{-\frac{1}{n}} \;=\; \sum_{n\in\mathbb{N}} e^{-\frac{1}{4n}(2n-1-\xi)^2} \;<\; +\infty,$$

whence $f \notin \mathcal{D}^{qa}(A)$ (on account of definition (3.19)), and also $g \notin \mathcal{D}^{qa}(A)$. Thus, (ii) is proved.

Part (iii) follows from $Af = g$ and from the injectivity of $A$.

Concerning part (iv), one observes first of all that the spectral measure of $A$ is only supported on $\sigma(A) = [0, +\infty)$, and moreover (see, e.g., [97, Example 5.3]), the spectral projections $E^{(A)}(\Omega)$, for given Borel subset $\Omega \subset [0, +\infty)$, are nothing but the multiplication operators by the characteristic functions $\mathbf{1}_{a^{-1}(\Omega)}$, where $x \mapsto a(x) := x^2$. Thus,

$$\begin{aligned}
\langle f, E^{(A)}(\Omega) f\rangle \;&=\; \int_{a^{-1}(\Omega)} |f(x)|^2\,dx \;=\; \int_{\{\lambda\in\mathbb{R}\,|\,\lambda^2\in\Omega\}} |f(x)|^2\,dx \\
&=\; \int_{\{\sqrt{\lambda}\,|\,\lambda\in\Omega\}} |f(x)|^2\,dx + \int_{\{-\sqrt{\lambda}\,|\,\lambda\in\Omega\}} |f(x)|^2\,dx.
\end{aligned}$$

Therefore, with $\Omega = [0, \lambda]$ and $E^{(A)}(\lambda) \equiv E^{(A)}([0, \lambda])$ for $\lambda \geqslant 0$,

$$\begin{aligned}
\frac{d\langle f, E^{(A)}(\lambda)f\rangle}{d\lambda} &= \frac{d}{d\lambda}\int_{-\sqrt{\lambda}}^{\sqrt{\lambda}}|f(\lambda)|^2 d\lambda = \frac{1}{2\sqrt{\lambda}}|f(\sqrt{\lambda})|^2 + \frac{1}{2\sqrt{\lambda}}|f(-\sqrt{\lambda})|^2 \\
&= \frac{1}{2\sqrt{\lambda}}|f(\sqrt{\lambda})|^2 = \frac{\mathbf{1}_{\mathbb{R}^+}(\lambda)}{\lambda^{2+\xi}\sqrt{2\pi}}e^{-\frac{1}{2}(\log\lambda)^2}.
\end{aligned}$$

From this, from part (iii), and from (3.33),

$$d\nu_\xi(\lambda) = \lambda^{\xi+1} d\langle f, E^A(\lambda)f\rangle = \frac{\mathbf{1}_{\mathbb{R}^+}(\lambda)}{\lambda\sqrt{2\pi}}e^{-\frac{1}{2}(\log\lambda)^2} d\lambda.$$

Part (iv) is proved.

Last, concerning part (v), it is well known (see, e.g., [106, Exercise 8.7]) that the space of polynomials on $[0,+\infty)$ is *not dense* in $L^2([0,+\infty),d\nu_\xi)$ when the measure $\nu_\xi$, as in the present case, is a log-normal distribution. As an immediate consequence the solution for the Hamburger moment problem for $\nu_\xi$ is not unique: for any $\nu_\xi$-integrable function $\varphi$ that is $\nu_\xi$-orthogonal to the subspace of polynomials, the distinct measures $\nu_\xi$ and $(1+\varphi)\nu_\xi$ have obviously the same moments (see, e.g., [2, Corollary 2.3.3] for the same conclusion in abstract terms). □

## 3.5 Unbounded CG-convergence tested numerically

In this Section a selection of numerical tests is discussed whose three-fold purpose is to confirm the main features of our convergence result, to corroborate our intuition on certain relevant differences with respect to the bounded case, and to explore the behaviour of the unbounded conjugate gradient algorithm beyond the regime covered by Theorem 3.1.

With respect to the Hilbert space $\mathcal{H} = L^2(\mathbb{R})$ one considers:

$$\begin{aligned}
\text{test-1a:} \quad & A := -\frac{d^2}{dx^2} + \mathbb{1}, \quad \mathcal{D}(A) := H^2(\mathbb{R}), \quad f(x) := e^{-x^2}, \\
\text{test-1b:} \quad & A := -\frac{d^2}{dx^2}, \quad \mathcal{D}(A) := H^2(\mathbb{R}), \quad f(x) := e^{-x^2}, \\
\text{test-2a:} \quad & A := -\frac{d^2}{dx^2} + \mathbb{1}, \quad \mathcal{D}(A) := H^2(\mathbb{R}), \quad f(x) := (1+x^2)^{-1}, \\
\text{test-2b:} \quad & A := -\frac{d^2}{dx^2}, \quad \mathcal{D}(A) := H^2(\mathbb{R}), \quad f(x) := (1+x^2)^{-1},
\end{aligned} \quad (3.53)$$

where $H^2$ denotes the usual Sobolev space of second order. In either case $A$ is an unbounded, injective, non-negative, self-adjoint operator on $\mathcal{H}$; but only in tests 1a and 2a does $A^{-1}$ exist as an everywhere defined bounded operator.

Correspondingly, four inverse linear problem of the form $Af = g$ are analysed with datum $g \in \text{ran}A$ given by the above explicit choice of the solution $f$. For each such problem conjugate gradient approximate solutions $f^{[N]}$ to $f$ are con-structed, namely $\xi$-iterates with $\xi = 1$, with initial guess $f^{[0]} = 0$ (the zero function on $\mathbb{R}$). Thus, each $f^{[N]}$ is searched for over the Krylov subspace $\mathcal{K}_N(A,g) = \mathrm{span}\{g, Ag, \ldots, A^{N-1}g\}$. $f^{[0]}$ is trivially smooth, and so are $f$ and $g$, therefore the 1-iterates are all well-defined. Owing to the injectivity of $A$ in all considered cases, necessarily $P_{\mathcal{S}} f^{[0]} = f$. The algorithm is well defined in all tests, as $f$ (and hence $g$) is a smooth function and is square-integrable together with all its derivatives (i.e., $g \in C^\infty(A)$).

Of course in practice the minimisation (3.3) is replaced with the standard, equivalent algebraic construction for the $f^{[N]}$'s [96, 73], so as to implement it as a routine in a symbolic computation software.

Iteratively the indicators

$$\begin{aligned}
\rho_0(f^{[N]}) &= \|f^{[N]} - f\|^2, \\
\rho_1(f^{[N]}) &= \langle f^{[N]} - f, A(f^{[N]} - f) \rangle, \\
\rho_2(f^{[N]}) &= \|Af^{[N]} - g\|^2
\end{aligned} \qquad (3.54)$$

(see (3.17) above) are evaluated, and their behaviour is monitored as $N$ increases.

The choice of the data $g$ in these tests is made so as tests 1a and 1b fall within the scope of Theorem 3.1, whereas tests 2a and 2b go *beyond* it. Indeed:

**Lemma 3.8** *With respect to the Hilbert space $L^2(\mathbb{R}, \mathrm{d}x)$ and the self-adjoint operators $A = -\frac{\mathrm{d}^2}{\mathrm{d}x^2} + 1$ or $A = -\frac{\mathrm{d}^2}{\mathrm{d}x^2}$ introduced above,*

*(i) the function $f = e^{-x^2}$, and hence $g = Af$ is analytic;*
*(ii) the function $f = (1+x^2)^{-1}$, and hence $g = Af$ is not quasi-analytic.*

*Proof* (i) It suffices to show that $f$ is $A$-analytic, the same conclusion for $g$ then follows from Lemma 3.6. By a standard criterion, $A$-analyticity of $f$ is tantamount as $f \in \mathcal{D}(e^{zT})$ for all $z \in \mathbb{C}$ with $|\Re z|$ small enough (see, e.g., [97, Corollary 7.9]). In terms of the Hilbert space isomorphism $L^2(\mathbb{R}, \mathrm{d}x) \xrightarrow{\cong} L^2(\mathbb{R}, \mathrm{d}p)$, $h \mapsto \widehat{h}$ induced by the Fourier transform, when $f = e^{-x^2}$ one has $\widehat{f} = \frac{1}{\sqrt{2}} e^{-\frac{1}{4}p^2}$ and the considered $A$'s are unitarily equivalent to the multiplication by $p^2 + 1$, or by $p^2$. Therefore,

$$\widehat{e^{zA}f} = \frac{1}{\sqrt{2}} e^{z(p^2+1)} e^{-\frac{1}{4}p^2} \in L^2(\mathbb{R}, \mathrm{d}p) \qquad \text{for } |\Re z| < \frac{1}{4},$$

and the same conclusion holds with $e^{zp^2}$ in place of $e^{z(p^2+1)}$, so both $A$'s are covered. Thus, $f = e^{-x^2}$ is indeed $A$-analytic. Incidentally the same criterion also shows that $f = (1+x^2)^{-1}$ is not analytic, because $\widehat{f} = \sqrt{\frac{\pi}{2}} e^{-|p|}$, and for no non-zero values of $\Re z$ is the function

$$\widehat{e^{zA}f} = \sqrt{\frac{\pi}{2}} e^{z(p^2+1)} e^{-|p|}$$

square-integrable on $\mathbb{R}$.

(ii) It suffices to show that $f$ is not quasi-analytic for $A$, the same conclusion for $g$ then follows from Lemma 3.4(i). Besides, it suffices to show the lack of quasi-analyticity of $f$ with respect to $A = -\frac{d^2}{dx^2}$: then, since for $n \in \mathbb{N}$

$$\|(-\tfrac{d^2}{dx^2} + 1)^n f\|_{L^2}^2 = \langle f, (-\tfrac{d^2}{dx^2}+1)^{2n} f\rangle_{L^2} \geqslant \langle f, (-\tfrac{d^2}{dx^2})^{2n} f\rangle_{L^2}$$
$$= \|(-\tfrac{d^2}{dx^2})^n f\|_{L^2}^2,$$

one concludes

$$\sum_{n\in\mathbb{N}} \|(-\tfrac{d^2}{dx^2}+1)^n f\|_{L^2}^{-\frac{1}{n}} \leqslant \sum_{n\in\mathbb{N}} \|(-\tfrac{d^2}{dx^2})^n f\|_{L^2}^{-\frac{1}{n}} < +\infty$$

(last inequality following from (3.19) and the fact that $f$ is not quasi-analytic for $-\frac{d^2}{dx^2}$), that is, the lack of quasi-analyticity of $f$ also for $-\frac{d^2}{dx^2} + 1$. So now let $A = -\frac{d^2}{dx^2}$, $f = (1+x^2)^{-1}$, and using $\widehat{f} = \sqrt{\frac{\pi}{2}} e^{-|p|}$ one computes

$$\|A^n f\|_{L^2}^2 = \int_{\mathbb{R}} \left(\sqrt{\frac{\pi}{2}} (p^2)^n e^{-|p|}\right)^2 dp = \pi \int_0^{+\infty} p^{4n} e^{-2p} dp$$
$$= \pi \frac{\Gamma(1+4n)}{2^{1+4n}}.$$

Using the known asymptotics [1, Eq. (6.1.37)] of the gamma function

$$\Gamma(t) \stackrel{t\to+\infty}{=} \sqrt{2\pi} e^{-t} t^{t-\frac{1}{2}} (1+O(t^{-1})),$$

one obtains, asymptotically as $n \to \infty$,

$$\|A^n f\|_{L^2}^2 = \frac{\pi\sqrt{2\pi}}{e} e^{-4n} (1+4n)^{\frac{1}{2}+4n} 2^{-(1+4n)} (1+O(n^{-1})),$$

whence

$$\|A^n f\|_{L^2}^{-\frac{1}{n}} = 4e^2 n^{-2} (1+O(n^{-1})).$$

This shows that the series $\sum_{n\in\mathbb{N}} \|A^n f\|_{L^2}^{-\frac{1}{n}}$ is asymptotic to $\sum_{n\in\mathbb{N}} n^{-2}$ and therefore converges. $f$ is not quasi-analytic for $A$. $\square$

On account of Lemma 3.8, Theorem 3.1 and Corollary 3.1 are applicable to tests 1a and 1b: since obviously $f^{[0]} \in \mathscr{C}_{A,g}(\sigma)$ $\forall \sigma \geqslant 0$, then Theorem 3.1 ensures that $\rho_\sigma(f^{[N]}) \to 0$ for any $\sigma \in [0,1]$. In particular, both the error ($\rho_0$) and the energy norm ($\rho_1$) are predicted to vanish as $N \to \infty$.

In the bounded case also the residual ($\rho_2$) would automatically vanish (Remark 3.2), but in tests 1a and 1b this indicator is not controlled by Theorem 3.1 and it is worth monitoring it.

A fourth meaningful quantity to monitor is $N^2 \rho_1(f^{[N]})$. Recall indeed that *if A were bounded* the energy norm would be predicted to vanish not slower than a rate of order $N^{-2}$ (as given by (3.25) with $\xi = 1$, $\sigma = 0$, $\sigma' = 1$). Thus, detecting now the possible failure of $N^2 \rho_1(f^{[N]})$ to stay bounded uniformly in $N$ is an immediate signature of the fact that one cannot apply to the unbounded-$A$ scenario the 'classical' quantitative convergence rate predicted by Nemirovskiy and Polyak for the bounded-$A$ scenario [81], which in fact was also proved to be optimal in that case [82].

The results of tests 1a and 1b are shown respectively in Figure 3.1 and 3.2.

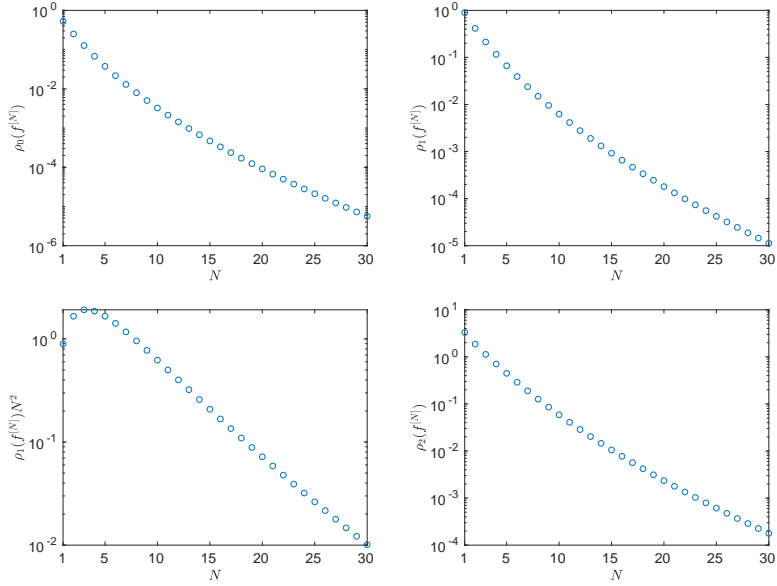

**Fig. 3.1** Numerical experiments for test 1a. From top left: $\rho_0(f^{[N]})$, $\rho_1(f^{[N]})$, $\rho_1(f^{[N]})N^2$, and $\rho_2(f^{[N]})$ indicators of convergence vs $N$.

Both tests 1a and 1b reveal that the iterates not only converge in the sense of the error and of the energy norm as predicted by Theorem 3.1, but also in the residual sense (not covered by Theorem 3.1). Of course in retrospect the error's vanishing in test 1a is consistent with the residual's vanishing, owing to the boundedness of $A^{-1}$ in test 1a: indeed, obviously,

$$\rho_0(f^{[N]}) = \|f^{[N]} - f\|^2 \leqslant \|A^{-1}\|_{\mathrm{op}}^2 \|Af^{[N]} - g\|^2 \leqslant \|A^{-1}\|_{\mathrm{op}}^2 \rho_2(f^{[N]}).$$

In addition, the classical Nemirovskiy-Polyak convergence rate for the energy norm is not violated in test 1a (Figure 3.1), whereas it appears to be violated in test 1b (Figure 3.2), where $A$ does not have a bounded inverse.

Heuristically, the slower vanishing rate of $\rho_0$, $\rho_1$, and $\rho_2$ in test 1b is due to the presence of zero in the spectrum of $A = -\frac{d^2}{dx^2}$: as we are approximating $A^{-1}g$ with polynomials $p_N(A)g$, the approximation to the inverse with polynomials is hampered about the "bad" spectral point $\lambda = 0$.

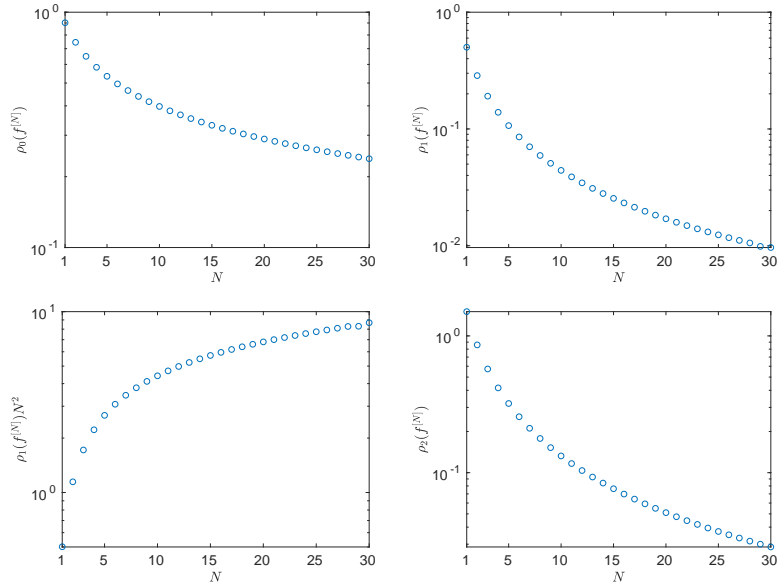

**Fig. 3.2** Numerical experiments for test 1b. From top left: $\rho_0(f^{[N]})$, $\rho_1(f^{[N]})$, $\rho_1(f^{[N]})N^2$, and $\rho_2(f^{[N]})$ indicators of convergence vs $N$.

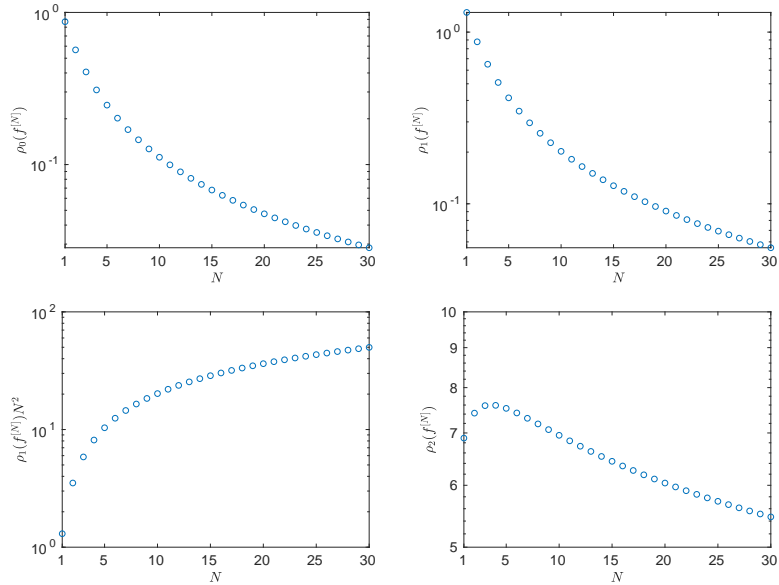

**Fig. 3.3** Numerical experiments for test 2a. From top left: $\rho_0(f^{[N]})$, $\rho_1(f^{[N]})$, $\rho_1(f^{[N]})N^2$, and $\rho_2(f^{[N]})$ indicators of convergence vs $N$.

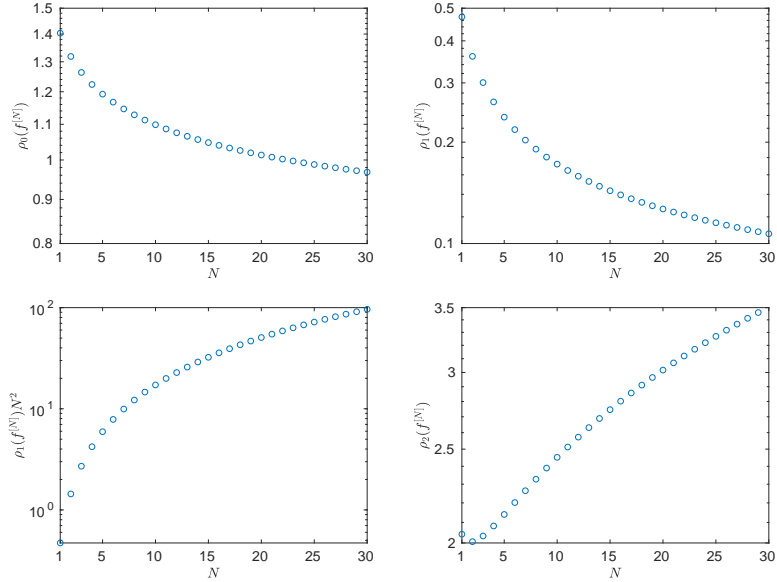

**Fig. 3.4** Numerical experiments for test 2b. From top left: $\rho_0(f^{[N]})$, $\rho_1(f^{[N]})$, $\rho_1(f^{[N]})N^2$, and $\rho_2(f^{[N]})$ indicators of convergence vs $N$.

As opposite to tests 1a and 1b, we know from Lemma 3.8(ii) that tests 2a and 2b, represented respectively in Figure 3.3 and 3.4, are *not covered* by Theorem 3.1 or Corollary 3.1, but for the fact that the quantities $\rho_0(f^{[N]})$ and $\rho_1(f^{[N]})$ are surely predicted to stay uniformly bounded in $N$ (Remark 3.15). Such uniform boundedness is confirmed numerically.

In test 2a, where $A$ has bounded inverse on the whole $\mathcal{H}$, numerics indicate that $\rho_0(f^{[N]}) \to 0$ and $\rho_1(f^{[N]}) \to 0$ as $N \to \infty$. That provides some practical evidence that there exist non-quasi-analytic data $g$ that still display "good behaviour", i.e., convergence at suitable $\rho_\sigma$-level. This is completely compatible with our Theorem 3.1 and Corollary 3.1: the use of (quasi-)analyticity that we made therein is solely localised in Proposition 3.4(v) in order to apply Carleman's criterion for the determinacy of the Hamburger moment problem, and that criterion is only a sufficient condition for the uniqueness of the $\nu_\xi$-measure. In fact, investigating the nebulous regime beyond quasi-analyticity would be of substantial relevance to understand what minimal assumptions on $g$ and $f^{[0]}$ guarantee the uniqueness of the $\nu_\xi$-measure (that surely fails in certain cases, as we saw in Proposition 3.5), or at least the convergence of the unbounded CG algorithm.

In comparison, in test 2b (unbounded $A^{-1}$) the decay rates of $\rho_0$ and $\rho_1$ appear to be slower than the counterpart 2a and it is unclear whether there is an actual vanishing, beside the evident decreasing behaviour.

The residual $\rho_2$ looks clearly decreasing in test 2a, with insufficient numerical evidence for vanishing, though, and instead manifestly divergent in test 2b. Here the solution $f = (1+x^2)^{-1}$ is not localised as the Gaussian of tests 1a and 1b, and has instead a long tail at large distances: the intuition suggests that this feature affects the convergence at higher regularity levels.

In either test 2a and 2b numerics give definite evidence of *violation* of the Nemirovskiy-Polyak convergence rate.

# Chapter 4
# Krylov solvability of unbounded inverse problems

## 4.1 Unbounded setting

We discuss in this Chapter an attempt of systematic extension of the analysis of Chapter 2 on the Krylov solvability of the inverse linear problem (1.1), switching now to the general case where the operator *A* is possibly *unbounded*. This means that all results presented here are applicable to the special case $A \in \mathcal{B}(\mathcal{H})$ (where in general stronger conclusions are possible), but their scope is discussed with explicit reference to unbounded operators on Hilbert space. In fact, in the latter class one can set various level of generality and abstraction, some of which depart excessively from the actual counterparts on applications: we shall therefore focus primarily on (unbounded) operators that are densely defined and closed.

As argued already in Chapter 3, an additional, obvious, and inescapable operational assumption is needed in the unbounded setting in order for the very notion of Krylov subspace to make sense. Given an operator *A* on a Hilbert space $\mathcal{H}$, we define the **N-th order Krylov subspace** and **Krylov subspace** (or **cyclic subspace**) relative to (or associated with) *A* and *g* to be the spaces

$$\mathcal{K}_N(A,g) := \mathrm{span}\{g, Ag, \ldots, A^{N-1}g\}, \qquad N \in \mathbb{N}, \tag{4.1}$$

and

$$\mathcal{K}(A,g) := \mathrm{span}\{A^k g \,|\, k \in \mathbb{N}_0\}, \tag{4.2}$$

*provided that g is a* **smooth vector for A** (or a **A-smooth vector**), meaning that

$$g \in C^\infty(A), \qquad C^\infty(A) := \bigcap_{k \in \mathbb{N}} \mathcal{D}(A^k). \tag{4.3}$$

For *g*'s outside the subspace of *A*-smooth vectors, definitions (4.1)-(4.2) are ill-posed, and on the other hand *A*-smooth vectors may constitute a tiny class, possibly the sole trivial set $\{0\}$! This circumstance indirectly restricts the interest to *A*'s and *g*'s admitting a satisfactorily large (for theoretical purposes and in applications) Krylov subspace.

*Example 4.1* On the Hilbert space $\mathcal{H} = L^2(\mathbb{R})$ consider the operator

$$\mathcal{D}(A) := \{\text{simple functions on } \mathbb{R}\}, \qquad A\varphi := x\varphi,$$

(**simple functions** in $L^2(\mathbb{R})$ being finite linear combinations of characteristic functions of Lebesgue-measurable subset of $\mathbb{R}$, [72, Section 1.17]). It is standard to see that $A$ is densely defined and essentially self-adjoint, yet obviously $\mathcal{D}(A) \cap \mathcal{D}(A^2) = \{0\}$, because the multiplication by $x$ destroys the step function structure. In this case, $C^\infty(A) = \{0\}$.

The object of this Chapter is therefore the analysis of Krylov solvability of the inverse problem

$$Af = g \tag{4.4}$$

on a given Hilbert space $\mathcal{H}$, in the unknown $f \in \mathcal{D}(A)$, under the working condition

$$g \in \operatorname{ran} A \cap C^\infty(A) \tag{4.5}$$

when the (possibly unbounded) operator $A$ is *densely defined* and *closed* in $\mathcal{H}$.

In fact, the analysis of Chapter 3 already produced a fundamental result for the Krylov solvability of inverse problems induced by (unbounded) positive self-adjoint operators on $\mathcal{H}$. Take indeed Theorem 3.1 with the special choice, in the notation therein, $f^{[0]} = 0$, $\sigma = 0$, $\xi = 1$: this yields the following property that for convenience we cast into a separate theorem. The subclass $\mathcal{D}^{qa}(A) \subset C^\infty(A)$ was defined in Section 3.2.

**Theorem 4.1** *With respect to a given Hilbert space $\mathcal{H}$ let $A = A^* \geqslant \mathbb{O}$ and let $g \in \operatorname{ran} A \cap \mathcal{D}^{qa}(A)$. Then the inverse problem $Af = g$ is Krylov solvable, and in particular one has*

$$\lim_{N \to \infty} \|f_N - f^\circ\|_{\mathcal{H}} = 0 \tag{4.6}$$

*along the sequence of the conjugate gradient iterates $f_N \in \mathcal{K}_N(A,g)$ defined in (3.3) with initial guess $f^{[0]} = 0$ and order $\xi = 1$, where $f^\circ$ is the minimal norm solution to the considered inverse problem.*

Here, and throughout this chapter, we shall keep the explicit notation $\|\ \|_{\mathcal{H}}$ for the norm of $\mathcal{H}$, not to generate confusion with the other frequently used norm in this discussion, the graph norm.

Theorem 4.1 is crucial for establishing Krylov solvability when $A$ is generically (unbounded and) self-adjoint or skew-adjoint, as the discussion in Section 4.2 shows.

Then in Section 4.3 we proceed on to the general case when $A$ is densely defined and closed on $\mathcal{H}$. We identify new obstructions in the problem of Krylov solvability, which are not present in the bounded case. A most serious one is the somewhat counter-intuitive phenomenon of 'Krylov escape', namely the possibility that vectors of $\overline{\mathcal{K}(A,g)}$ that also belong to the domain of $A$ are mapped by $A$ outside $\overline{\mathcal{K}(A,g)}$, whereas obviously $A\mathcal{K}(A,g) \subset \mathcal{K}(A,g)$. From a perspective that in fact we are not carrying over here, the possibility of Krylov escape adds further complication to the unbounded operator counterpart of the celebrated invariant subspace problem [119], at least when $g \neq 0$ and $g$ is not a cyclic vector for $A$, hence $\overline{\mathcal{K}(A,g)}$ is a proper closed subspace of $\mathcal{H}$.

In Section 4.3 we also determine that if the closures of $\mathcal{K}(A,g)$ in the Hilbert space norm and in the stronger $A$-graph norm are the same (up to intersection with $\mathcal{D}(A)$), an occurrence that we named 'Krylov-core condition', then Krylov escape is actually prevented.

This leads us to Section 4.4, where we demonstrate that, under assumptions like the Krylov-core condition (and, more generally, lack of Krylov escape), the intrinsic mechanisms of Krylov reducibility and triviality of the Krylov intersection play a completely analogous role as compared to the bounded case.

Last, in Section 4.5 we re-consider the (unbounded) self-adjoint scenario, already solved conceptually in Section 4.2, investigating Krylov solvability from the perspective of the abstract operator-theoretic mechanisms mentioned above. Noticeably, this is also a perspective rising up interesting open questions that certainly indicate future directions of investigation.

## 4.2 The general self-adjoint and skew-adjoint case

We have the following main result (let us recall once again that the class $\mathcal{D}^{qa}(A)$ of quasi-analytic vectors for $A$ was defined in Section 3.2).

**Theorem 4.2** *With respect to a given Hilbert space $\mathcal{H}$, let $A$ be a self-adjoint ($A^* = A$) or skew-adjoint ($A^* = -A$) operator on $\mathcal{H}$ and let $g \in \operatorname{ran} A \cap \mathcal{D}^{qa}(A)$. Then there exists a unique solution $f$ to $Af = g$ such that $f \in \overline{\mathcal{K}(A,g)}$. Thus, the inverse problem $Af = g$ is Krylov-solvable.*

***Proof*** <u>Existence</u>. In the self-adjoint case, since $A^2$ is self-adjoint and $A^2 \geqslant \mathbb{O}$, Theorem 4.1 implies that there exists $f \in \overline{\mathcal{K}(A^2, Ag)} \subset \overline{\mathcal{K}(A,g)}$ such that $A^2 f = Ag$. Analogously, in the skew-adjoint case, since $-A^2$ is self-adjoint and $-A^2 \geqslant \mathbb{O}$, then there exists $f \in \overline{\mathcal{K}(-A^2, -Ag)} = \overline{\mathcal{K}(A^2, Ag)} \subset \overline{\mathcal{K}(A,g)}$ such that $-A^2 f = -Ag$.

In either case, $f \in \overline{\mathcal{K}(A,g)}$, $f \in \mathcal{D}(A^2) \subset \mathcal{D}(A)$, and $A^2 f = Ag$, equivalently, $A(Af - g) = 0$. This shows that $Af - g \in \ker A$.

On the other hand, both $Af$ and $g$ belong to $\operatorname{ran} A$, whence $Af - g \in \operatorname{ran} A \subset (\ker A^*)^\perp = (\ker A)^\perp$, where the last identity is clearly valid both for in the self-adjoint and in the skew-adjoint case.

Then necessarily $Af - g = 0$, which proves that $f \in \overline{\mathcal{K}(A,g)}$ is a solution to the considered inverse problem.

<u>Uniqueness</u>. If $f_1, f_2 \in \overline{\mathcal{K}(A,g)}$ and $Af_1 = g = Af_2$, then $f_1 - f_2 \in \ker A \cap \overline{\mathcal{K}(A,g)}$. Moreover, $\ker A = \ker A^*$ and $\overline{\mathcal{K}(A,g)} \subset \overline{\operatorname{ran} A}$. Therefore, $f_1 - f_2 \in \ker A^* \cap \overline{\operatorname{ran} A}$. Thus, $f_1 = f_2$. □

*Remark 4.1* In view of the discussion made in Chapter 3 for Theorem 3.1, the proof of Theorem 4.2 shows that the actual Krylov-solution $f$ to $Af = g$ is the minimal norm solution to $A^2 f = Ag$ and admits approximants $f_N$, with $\|f_N - f\|_{\mathcal{H}} \to 0$ as $N \to \infty$, defined by

$$f_N := \arg\min_{h \in \mathcal{K}_N(A^2, Ag)} \|A(h-f)\|^2_{\mathcal{H}}.$$

Thus, the iterates of the CG algorithm applied to the auxiliary problem $A^2 f = Ag$ (interpreted as $-A^2 f = -Ag$ in the skew-adjoint case) converge precisely to the Krylov solution to $Af = g$.

*Remark 4.2* Unbounded skew-adjoint inverse problems are intimately related to inverse problems induced by so-called Friedrichs operators. These constitute a class of elliptic, parabolic, and hyperbolic differential operators, that can be also characterised as abstract operators on Hilbert space $\mathcal{H}$ [38, 5, 6, 7], having the typical (but not the only one) form $T = A + C$ where $A^* = -A$ and $C \in \mathcal{B}(\mathcal{H})$. Theorem 4.2 is applicable when $C$ is skew-adjoint itself.

## 4.3 New phenomena in the general unbounded case: 'Krylov escape', generalised Krylov reducibility, generalised Krylov intersection

Let us start in this section the analysis of Krylov solvability of the inverse problem (4.4) on a given Hilbert space $\mathcal{H}$, in the unknown $f \in \mathcal{D}(A)$, under the working condition (4.5) when the (possibly unbounded) operator $A$ is *densely defined* and *closed* in $\mathcal{H}$.

A number of substantial novelties, due to domain issues, emerge in this case as compared to the bounded case discussed in Chapter 2.

The first unavoidable difference concerns the invariance of $\overline{\mathcal{K}(A,g)}$ (respectively, $\mathcal{K}(A,g)^\perp$) under the action of $A$ (resp., of $A^*$). Indeed, the inclusions (2.21) certainly cannot be valid in general, because the above subspaces may well not be included, respectively, in $\mathcal{D}(A)$ and $\mathcal{D}(A^*)$.

*Example 4.2* The 'quantum mechanical creation operator' on $\mathcal{H} = L^2(\mathbb{R})$

$$A = \frac{1}{\sqrt{2}}\left(-\frac{\mathrm{d}}{\mathrm{d}x} + x\right)$$
$$\mathcal{D}(A) = \{h \in L^2(\mathbb{R}) \,|\, -h' + xh \in L^2(\mathbb{R})\}$$

is densely defined, unbounded, and closed, and has the well-known property that

$$\psi_{n+1} = \frac{1}{\sqrt{n+1}} A\psi_n \qquad n \in \mathbb{N}_0,$$

where $(\psi_n)_{n \in \mathbb{N}_0}$ is the orthonormal basis of $L^2(\mathbb{R})$ of the Hermite functions $\psi_n(x) = c_n H_n(x) e^{-x^2/2}$ (here $c_n$ is a normalisation factor and $H_n$ is the $n$-th Hermite polynomial). In particular, each $\psi_n$ is a $C^\infty(A)$-function. Choosing $g = \psi_1$ evidently yields $\overline{\mathcal{K}(A,g)} = \text{span}\{\psi_0\}^\perp$. But there are $L^2$-functions orthogonal to $\psi_0$ that do not belong to $\mathcal{D}(A)$.

It is then clear that only the possible invariance of $\overline{\mathcal{K}(A,g)} \cap \mathcal{D}(A)$ under $A$ and of $\mathcal{K}(A,g)^\perp \cap \mathcal{D}(A^*)$ under $A^*$ makes sense in general.

This naturally leads one to consider the operators $A|_{\overline{\mathcal{K}(A,g)} \cap \mathcal{D}(A)}$ (the **part of $A$ on $\overline{\mathcal{K}(A,g)}$**) and $A^*|_{\mathcal{K}(A,g)^\perp \cap \mathcal{D}(A^*)}$ (the **part of $A^*$ on $\mathcal{K}(A,g)^\perp$**). Noticeably, when $A$ is unbounded and $\mathcal{D}(A)$ is a proper dense subspace of $\mathcal{H}$, none of the two is densely defined in $\mathcal{H}$, unless $\overline{\mathcal{K}(A,g)} = \mathcal{H}$, as their domain is by construction the intersection of a proper dense and a proper closed subspace. Obviously, instead, $A|_{\overline{\mathcal{K}(A,g)} \cap \mathcal{D}(A)}$ is densely defined in the Hilbert subspace $\overline{\mathcal{K}(A,g)}$.

**Lemma 4.1** *For a given Hilbert space $\mathcal{H}$ let $A$ be a densely defined operator on $\mathcal{H}$ and let $g \in C^\infty(A)$. Then*

$$A^*\left(\mathcal{K}(A,g)^\perp \cap \mathcal{D}(A^*)\right) \subset \mathcal{K}(A,g)^\perp. \tag{4.7}$$

***Proof*** Let $z \in \mathcal{K}(A,g)^\perp \cap \mathcal{D}(A^*)$. For arbitrary $h \in \overline{\mathcal{K}(A,g)}$ let $(h_n)_{n \in \mathbb{N}}$ be a sequence in $\mathcal{K}(A,g)$ of norm-approximants of $h$. Then each $Ah_n \in \mathcal{K}(A,g)$, and therefore

$$\langle h, A^*z \rangle = \lim_{n \to \infty} \langle h_n, A^*z \rangle = \lim_{n \to \infty} \langle Ah_n, z \rangle = 0,$$

thus proving (4.7). $\square$

The counterpart inclusion to (4.7), namely $A(\overline{\mathcal{K}(A,g)} \cap \mathcal{D}(A)) \subset \overline{\mathcal{K}(A,g)}$ when $\overline{\mathcal{K}(A,g)}$ is only a proper closed subspace of $\mathcal{H}$, turns out to be considerably less trivial. In fact, as somewhat counter-intuitive as it appears, $A$ may indeed map vectors from $\overline{\mathcal{K}(A,g)} \cap \mathcal{D}(A)$ *outside* of $\overline{\mathcal{K}(A,g)}$. In the present context, we shall refer to this phenomenon, that has no analogue in the bounded case, as **Krylov escape**.

*Example 4.3 (Krylov escape)* Let $\mathcal{H}'$ be a Hilbert space and $T'$ be a self-adjoint operator in $\mathcal{H}'$ having a cyclic vector $g'$, meaning that there exists $g' \in \mathcal{D}(T')$ such that $\overline{K(T',g')} = \mathcal{H}'$. (It is straightforward to construct explicit examples for such a choice.) For any one-dimensional vector space $\mathcal{H}_0$, say, $\mathcal{H}_0 = \text{span}\{e_0\}$, set

$$\begin{aligned} \mathcal{H} &:= \mathcal{H}_0 \oplus \mathcal{H}', \\ T &:= \mathbb{O} \oplus T', \\ g &:= 0 \oplus g'. \end{aligned}$$

(The last condition is just an identification of $g$ as an element of $\mathcal{H}$.) Thus, $T$ is a self-adjoint operator in $\mathcal{H}$ such that $Te_0 = 0$ and $Tx' = T'x' \; \forall x' \in \mathcal{D}(T')$, and moreover $\overline{\mathcal{K}(T,g)} = \mathcal{H}'$. (Now the closure is taken with respect to $\mathcal{H}$.) Furthermore, let $x_0 \in \mathcal{H}$ such that $x_0 \in \mathcal{H}' \setminus \mathcal{D}(T')$. Then set

$$\begin{aligned} \mathcal{D}(A) &:= \mathcal{D}(T) \dotplus \text{span}\{x_0\}, \\ Ax_0 &:= e_0, \\ Ax &:= Tx \qquad \forall x \in \mathcal{D}(T). \end{aligned}$$

$A$ is meant to be defined by the above identities and extended by linearity on the whole $\mathcal{D}(T) \dotplus \mathrm{span}\{x_0\}$. The operator $A$ is densely defined in $\mathcal{H}$ by construction.

• <u>Closedness</u>. Let us check that if $x_n + \mu_n x_0 \to v$ and $A(x_n + \mu_n x_0) \to w$ in $\mathcal{H}$ as $n \to \infty$ for some vectors $v, w \in \mathcal{H}$, where $(x_n)_{n \in \mathbb{N}}$ is a generic sequence in $\mathcal{D}(T)$ and $(\mu_n)_{n \in \mathbb{N}}$ is a generic sequence in $\mathbb{C}$, then $v \in \mathcal{D}(A)$ and $Av = w$. First, we observe that it must be $\mu_n \to \mu$ for some $\mu \in \mathbb{C}$, for otherwise there would be no chance for the vectors $A(x_n + \mu_n x_0) = Tx_n + \mu_n e_0$ to converge as assumed, because $Tx_n \perp \mu_n e_0$. Next, since $\mu_n \to \mu$ and $x_n + \mu_n x_0 \to v$, then necessarily $x_n \to x$ for some $x \in \mathcal{H}$ satisfying $x + \mu x_0 = v$; analogously, since $Tx_n + \mu_n e_0 = A(x_n + \mu_n x_0) \to w$, then $Tx_n \to y$ for some $y \in \mathcal{H}$ satisfying $y + \mu e_0 = w$. As by construction $T$ is self-adjoint and hence closed, then necessarily $x \in \mathcal{D}(T)$ and $Tx = y$. In turn, this implies that $v \in \mathcal{D}(A)$ and $Av = w$. The conclusion is that $A$ is closed.

• <u>Occurrence of Krylov escape</u>. By the construction above, $\overline{\mathcal{K}(A, g)} = \mathcal{H}' = \mathrm{span}\{e_0\}^\perp$. Let us focus on the vector $e_0$. On the one hand, $e_0 = Ax_0$ and $x_0$ by definition belongs both to $\overline{\mathcal{K}(A, g)}$ and to $\mathcal{D}(A)$. Therefore $e_0 \in A(\overline{\mathcal{K}(A, g)} \cap \mathcal{D}(A))$. On the other hand, however, $e_0 \in (\mathcal{H}')^\perp = \overline{\mathcal{K}(A, g)}^\perp$, whence $e_0 \notin \overline{\mathcal{K}(A, g)}$. This provides a counterexample of a densely defined closed operator $A$ in $\mathcal{H}$ such that the inclusion

$$A\big(\overline{\mathcal{K}(A, g)} \cap \mathcal{D}(A)\big) \subset \overline{\mathcal{K}(A, g)}.$$

is *violated*.

Owing to the possible occurrence of the Krylov escape phenomenon, the invariance of $\overline{\mathcal{K}(A, g)} \cap \mathcal{D}(A)$ requires additional assumptions. A reasonable one is to assume further that the operator $A|_{\overline{\mathcal{K}(A,g)} \cap \mathcal{D}(A)}$ and its restriction $A|_{\mathcal{K}(A,g)}$ are in a sense as close as possible. To this aim, let us observe first the following.

**Lemma 4.2** *For a given Hilbert space $\mathcal{H}$ let $A$ be a densely defined and closed operator on $\mathcal{H}$ and let $g \in C^\infty(A)$. Then*

*(i) the operator $A|_{\overline{\mathcal{K}(A,g)} \cap \mathcal{D}(A)}$ is closed,*
*(ii) and the operator $A|_{\mathcal{K}(A,g)}$ is closable.*

**Proof** Obviously $A|_{\mathcal{K}(A,g)} \subset A|_{\overline{\mathcal{K}(A,g)} \cap \mathcal{D}(A)}$, in the sense of operator inclusion, so part (ii) follows at once from part (i). In turn, part (i) is true as is always the case when one restricts a *closed* operator to the intersection of its domain with a closed subspace: the restriction operator too is closed. Explicitly, let $((x_n, Ax_n))_{n \in \mathbb{N}}$ be an arbitrary $\mathcal{H} \oplus \mathcal{H}$-convergent sequence in the graph of $A|_{\overline{\mathcal{K}(A,g)} \cap \mathcal{D}(A)}$, that is, for some $x, y \in \mathcal{H}$ one has $\overline{\mathcal{K}(A, g)} \cap \mathcal{D}(A) \ni x_n \to x$ and $Ax_n \to y$ in $\mathcal{H}$. Then, by closedness of $A$, $x \in \mathcal{D}(A)$ and $Ax = y$. Moreover, since $x_n \in \overline{\mathcal{K}(A, g)}$ $\forall n$, also for the limit point one has $x \in \overline{\mathcal{K}(A, g)}$. Thus, $x \in \overline{\mathcal{K}(A, g)} \cap \mathcal{D}(A)$. This shows that the pair $(x, y)$ belongs to the graph of $A|_{\overline{\mathcal{K}(A,g)} \cap \mathcal{D}(A)}$, which is therefore closed in $\mathcal{H} \oplus \mathcal{H}$. □

*Remark 4.3* With a completely analogous argument one shows that the operator $A^*|_{\overline{\mathcal{K}(A,g)}^\perp \cap \mathcal{D}(A^*)}$ is closed too.

It is then natural to consider the case when the operator closure of $A|_{\mathcal{K}(A,g)}$ is precisely $A|_{\overline{\mathcal{K}(A,g)} \cap \mathcal{D}(A)}$.

We set this definition. For a given Hilbert space $\mathcal{H}$ let $A$ be a densely defined and closed operator on $\mathcal{H}$ and let $g \in C^\infty(A)$. Then the pair $(A,g)$ is said to satisfy the **Krylov-core condition** when the subspace $\mathcal{K}(A,g)$ is a core for $A|_{\overline{\mathcal{K}(A,g)} \cap \mathcal{D}(A)}$. Explicitly (see, e.g., [97, Section 1.1.2]), this is the requirement that

$$\overline{A|_{\mathcal{K}(A,g)}} = A|_{\overline{\mathcal{K}(A,g)} \cap \mathcal{D}(A)} \tag{4.8}$$

in the sense of operator closure, equivalently, it is the requirement that $\mathcal{K}(A,g)$ is dense in $\overline{\mathcal{K}(A,g)} \cap \mathcal{D}(A)$ in the **graph norm**

$$\|h\|_A := (\|h\|_{\mathcal{H}}^2 + \|Ah\|_{\mathcal{H}}^2)^{\frac{1}{2}}, \qquad h \in \mathcal{D}(A), \tag{4.9}$$

that is,

$$\overline{\mathcal{K}(A,g)}^{\|\ \|_A} = \overline{\mathcal{K}(A,g)} \cap \mathcal{D}(A). \tag{4.10}$$

*Remark 4.4* Since $A$ is a closed operator, the inclusion

$$\overline{\mathcal{K}(A,g)}^{\|\ \|_A} \subset \overline{\mathcal{K}(A,g)} \cap \mathcal{D}(A) \tag{4.11}$$

is always true, as one sees reasoning as for Lemma 4.2.

*Example 4.4* Let $\mathcal{H} = \ell^2(\mathbb{N}_0)$ and, in terms of the canonical orthonormal basis $(e_n)_{n \in \mathbb{N}_0}$, let $A$ be the densely defined and closed operator

$$\mathcal{D}(A) := \left\{ x \equiv (x_n)_{n \in \mathbb{N}_0} \ \Big| \ \sum_{n=0}^{\infty} (n+1)^2 |x_n|^2 < +\infty \right\}$$

$$Ax := (0, x_0, 2x_1, 3x_2, \ldots) \qquad \text{for } x \in \mathcal{D}(A)$$

(thus, in particular, $Ae_n = (n+1)e_{n+1}$ for any $n \in \mathbb{N}_0$). Obviously, we have $\mathcal{K}(A,e_0) = \mathrm{span}\{e_1,e_2,e_3,\ldots\}$ and $\overline{\mathcal{K}(A,e_0)} = \{e_0\}^\perp$. Let

$$x := (0, x_1, x_2, x_3, x_4, \ldots) \quad \text{with} \quad \sum_{n=1}^{\infty} (n+1)^2 |x_n|^2 < +\infty$$

be a generic vector in $\overline{\mathcal{K}(A,e_0)} \cap \mathcal{D}(A)$ and for each $N \in \mathbb{N}$ let

$$x^{(N)} := (0, x_1^{(N)}, \ldots, x_N^{(N)}, 0, 0, \ldots) \quad \text{with} \quad x_n^{(N)} = x_n + \frac{1}{n^2 N}.$$

Then $x^{(N)} \in \mathcal{K}(A, e_0)$ and

$$\|x - x^{(N)}\|_A^2 = \sum_{n=1}^{\infty} (1+(n+1)^2)|x_n - x_n^{(N)}|^2$$

$$= \frac{1}{N^2} \sum_{n=1}^{N} \frac{n^2 + 2n + 2}{n^4} + \sum_{n=N+1}^{\infty} (1+(n+1)^2)|x_n|^2 \xrightarrow{N \to \infty} 0,$$

whence $x \in \overline{\mathcal{K}(A,e_0)}^{\|\,\|_A}$. Thus, $\overline{\mathcal{K}(A,e_0)}^{\|\,\|_A} \supset \overline{\mathcal{K}(A,e_0)} \cap \mathcal{D}(A)$, which together with (4.11) shows that the pair $(A, e_0)$ does satisfy the Krylov-core condition (4.10).

The Krylov-core condition is indeed sufficient to finally ensure that $A$ maps $\overline{\mathcal{K}(A,g)} \cap \mathcal{D}(A)$ into $\overline{\mathcal{K}(A,g)}$.

**Lemma 4.3** *For a given Hilbert space $\mathcal{H}$ let $A$ be a densely defined and closed operator on $\mathcal{H}$ and let $g \in C^{\infty}(A)$.*

*(i) One has the inclusion*

$$A\big(\overline{\mathcal{K}(A,g)} \cap \mathcal{D}(A)\big) \subset \overline{A\mathcal{K}(A,g)} \tag{4.12}$$

*if and only if $A$ and $g$ satisfy the Krylov-core condition* (4.8).
*(ii) In particular, under the Krylov-core condition one has*

$$A\big(\overline{\mathcal{K}(A,g)} \cap \mathcal{D}(A)\big) \subset \overline{\mathcal{K}(A,g)}. \tag{4.13}$$

*Proof* Clearly (4.13) follows at once from (4.12) as $\overline{A\mathcal{K}(A,g)} \subset \overline{\mathcal{K}(A,g)}$. Let us focus then on the proof of part (i). Assume that (4.8) is satisfied and let $z \in \overline{\mathcal{K}(A,g)} \cap \mathcal{D}(A)$, the domain of $\overline{A|_{\overline{\mathcal{K}(A,g)} \cap \mathcal{D}(A)}}$. Since the latter operator is the closure of $A|_{\mathcal{K}(A,g)}$, then there exists $(z_n)_{n \in \mathbb{N}}$ in $\mathcal{K}(A,g)$ such that $z_n \to z$ and $A|_{\mathcal{K}(A,g)} z_n \to \overline{A|_{\overline{\mathcal{K}(A,g)} \cap \mathcal{D}(A)}} z$, i.e., $Az_n \to Az$. This shows that $Az \in \overline{A\mathcal{K}(A,g)}$. For the converse implication, assume now that $\overline{A|_{\overline{\mathcal{K}(A,g)} \cap \mathcal{D}(A)}}$ is a *proper* closed extension of $\overline{A|_{\mathcal{K}(A,g)}}$. Therefore there is some $z \in \overline{\mathcal{K}(A,g)} \cap \mathcal{D}(A)$ such that whatever sequence $(z_n)_{n \in \mathbb{N}}$ in $\mathcal{K}(A,g)$ of norm-approximants of $z$ is considered, one cannot have $Az_n \to Az$, because this would mean $\overline{A|_{\mathcal{K}(A,g)}} z_n \to \overline{A|_{\overline{\mathcal{K}(A,g)} \cap \mathcal{D}(A)}} z$, that is, $\overline{A|_{\mathcal{K}(A,g)}} z = \overline{A|_{\overline{\mathcal{K}(A,g)} \cap \mathcal{D}(A)}} z$, contrarily to the assumption $\overline{A|_{\mathcal{K}(A,g)}} \subsetneq \overline{A|_{\overline{\mathcal{K}(A,g)} \cap \mathcal{D}(A)}}$. Thus, $Az$ cannot have norm-approximants in $A\mathcal{K}(A,g)$ and hence (4.12) cannot be valid. □

*Remark 4.5* In general $\overline{A\mathcal{K}(A,g)} \subsetneq \overline{\mathcal{K}(A,g)}$. Example 4.2 shows that in that case $\overline{A\mathcal{K}(A,g)} = \text{span}\{\psi_0, \psi_1\}^{\perp}$, whereas $\overline{\mathcal{K}(A,g)} = \text{span}\{\psi_0\}^{\perp}$.

The two key notions of Krylov reducibility and Krylov intersection introduced in the bounded case in Section 2.2 can be now generalised to the unbounded case. Let $A$ be a densely defined and closed operator on a Hilbert space $\mathcal{H}$ and let $g \in C^{\infty}(A)$. $A$ is said to be **$\mathcal{K}(A,g)$-reduced** (in the generalised sense) when

$$\begin{aligned} A\big(\overline{\mathcal{K}(A,g)} \cap \mathcal{D}(A)\big) &\subset \overline{\mathcal{K}(A,g)}, \\ A\big(\mathcal{K}(A,g)^{\perp} \cap \mathcal{D}(A)\big) &\subset \mathcal{K}(A,g)^{\perp}. \end{aligned} \tag{4.14}$$

The subspace
$$\mathcal{I}(A,g) := \overline{\mathcal{K}(A,g)} \cap A\big(\mathcal{K}(A,g)^\perp \cap \mathcal{D}(A)\big) \tag{4.15}$$
is called the **(generalised) Krylov intersection** for the given $A$ and $g$.

As in the bounded case (Section 2.2), it is clear from (4.14) and (4.15) that also in the present, generalised sense Krylov reducibility implies triviality of the Krylov intersection.

*Remark 4.6* The first condition of the two requirements (4.14) for Krylov-reducibility is precisely the lack of Krylov escape (4.13). Thus this condition is matched if, for instance, $A$ and $g$ satisfy the Krylov-core condition (Lemma 4.3).

*Remark 4.7* Krylov reducibility in the generalised sense (4.14) for an *unbounded* operator differs from Krylov reducibility in the *bounded* case (Section 2.2), which is formulated as
$$A\overline{\mathcal{K}(A,g)} \subset \overline{\mathcal{K}(A,g)},$$
$$A\mathcal{K}(A,g)^\perp \subset \mathcal{K}(A,g)^\perp,$$
in that when $A$ is bounded the subspaces $\overline{\mathcal{K}(A,g)}$ and $\mathcal{K}(A,g)^\perp$ are *reducing* for $A$ and hence with respect to the Krylov decomposition $\mathcal{H} = \overline{\mathcal{K}(A,g)} \oplus \mathcal{K}(A,g)^\perp$ the operator $A$ decomposes as $A = A|_{\overline{\mathcal{K}(A,g)}} \oplus A|_{\mathcal{K}(A,g)^\perp}$ whereas if $A$ is unbounded and Krylov-reduced it is false in general that $A = A|_{\overline{\mathcal{K}(A,g)} \cap \mathcal{D}(A)} \oplus A|_{\mathcal{K}(A,g)^\perp \cap \mathcal{D}(A)}$.

## 4.4 Krylov solvability in the general unbounded case

Let us now discuss the Krylov solvability of the inverse problem (4.4) when the operator $A$ is densely defined and closed.

A crucial role in this matter turns out to be played by the lack of Krylov escape, namely the property (4.13), a feature that was automatically present in the bounded case. A first example is the following technical Lemma that will be useful in a moment.

**Lemma 4.4** *For a given Hilbert space $\mathcal{H}$ let $A$ be a densely defined and closed operator on $\mathcal{H}$, let $g \in \mathrm{ran}(A) \cap C^\infty(A)$, and let $f \in \mathcal{D}(A)$ satisfy $Af = g$. If in addition the pair $(A,g)$ satisfies the Krylov-core condition and $f \in \overline{\mathcal{K}(A,g)}$, then*
$$\overline{A\mathcal{K}(A,g)} = \overline{\mathcal{K}(A,g)}. \tag{4.16}$$

(We already observed in Remark 4.5 that in general $\overline{A\mathcal{K}(A,g)} \subsetneq \overline{\mathcal{K}(A,g)}$.)

***Proof (Proof of Lemma 4.4)*** By assumption, $f \in \overline{\mathcal{K}(A,g)} \cap \mathcal{D}(A)$, which is the same as $f \in \overline{\mathcal{K}(A,g)}^{\|\|\|_A}$, owing to the Krylov-core condition. Thus there exists a sequence $(f_n)_{n\in\mathbb{N}}$ in $\mathcal{K}(A,g)$ such that $f_n \xrightarrow{\|\|\|_A} f$ as $n \to \infty$, in particular $Af_n \xrightarrow{\|\|\|} Af =$ $g$, which implies that $g$ belongs to $\overline{A\mathcal{K}(A,g)}$. Clearly all vectors $Ag, A^2g, A^3g, \ldots$ belong to the same space too. Therefore,

$$\mathrm{span}\{A^k g \,|\, k \in \mathbb{N}_0\} \subset \overline{A\mathcal{K}(A,g)},$$

so that $\overline{\mathcal{K}(A,g)} \subset \overline{A\mathcal{K}(A,g)}$. The opposite inclusion is trivial, hence (4.16) follows. $\square$

For an inverse problem with unbounded injective $A$, the triviality of the Krylov intersection still implies Krylov solvability under additional assumptions that are automatically satisfied if $A$ is bounded. The first requires that the orthogonal projection onto $\overline{\mathcal{K}(A,g)}$ of the solution manifold $\mathcal{S}(A,g)$ is entirely contained in the domain of $A$. The other requirement is the lack of Krylov escape.

In fact, Proposition 2.3 admits the following analogues. For convenience we shall denote by $P_{\mathcal{K}}$ the orthogonal projection $P_{\mathcal{K}} : \mathcal{H} \to \mathcal{H}$ onto the Krylov subspace $\overline{\mathcal{K}(A,g)}$.

**Proposition 4.1** *For a given Hilbert space $\mathcal{H}$ let $A$ be a densely defined and closed operator on $\mathcal{H}$, let $g \in \mathrm{ran}(A) \cap C^\infty(A)$, and let $f \in \mathcal{D}(A)$ satisfy $Af = g$. Assume furthermore that*

*(a) $A$ is injective;*
*(b) $A(\overline{\mathcal{K}(A,g)} \cap \mathcal{D}(A)) \subset \overline{\mathcal{K}(A,g)}$ – this assumption holds true, for example, if the pair $(A,g)$ satisfies the Krylov-core condition (Lemma 4.3);*
*(c) $P_{\mathcal{K}} f \in \mathcal{D}(A)$;*
*(d) $\mathcal{I}(A,g) = \{0\}$,*

*or also, assume the more stringent assumptions*

*(a) $A$ is injective;*

*(b') $A$ is $\mathcal{K}(A,g)$-reduced;*
*(c) $P_{\mathcal{K}} f \in \mathcal{D}(A)$.*

*Under such assumptions, $f \in \overline{\mathcal{K}(A,g)}$.*

***Proof*** By assumption (c), $P_{\mathcal{K}} f \in \overline{\mathcal{K}(A,g)} \cap \mathcal{D}(A)$; by assumption (b), $AP_{\mathcal{K}} f \in \overline{\mathcal{K}(A,g)}$. Then $A(\mathbb{1} - P_{\mathcal{K}})f = g - AP_{\mathcal{K}} f \in \overline{\mathcal{K}(A,g)}$. On the other hand, again by assumption (c), $(\mathbb{1} - P_{\mathcal{K}})f \in \mathcal{D}(A)$, whence $A(\mathbb{1} - P_{\mathcal{K}})f \in A(\mathcal{K}(A,g)^\perp \cap \mathcal{D}(A))$. Thus, $A(\mathbb{1} - P_{\mathcal{K}})f \in \mathcal{I}(A,g)$. By assumptions (a) and (d), then $f = P_{\mathcal{K}} f$. $\square$

**Proposition 4.2** *For a given Hilbert space $\mathcal{H}$ let $A$ be a densely defined and closed operator on $\mathcal{H}$, let $g \in \mathrm{ran}(A) \cap C^\infty(A)$, and let $f \in \mathcal{D}(A)$ satisfy $Af = g$. Assume furthermore that*

*(a) $A$ is invertible with everywhere defined, bounded inverse on $\mathcal{H}$;*
*(b) the pair $(A,g)$ satisfies the Krylov-core condition.*

*Under such assumptions, if $f \in \overline{\mathcal{K}(A,g)}$, then $\mathcal{I}(A,g) = \{0\}$.*

***Proof*** Let $z \in \mathcal{I}(A,g)$. Then $z = Aw$ for some $w \in \mathcal{K}(A,g)^\perp \cap \mathcal{D}(A)$, and $z \in \overline{\mathcal{K}(A,g)}$. Owing to Lemma 4.4, $\overline{\mathcal{K}(A,g)} = \overline{A\mathcal{K}(A,g)}$, hence there is a sequence $(v_n)_{n \in \mathbb{N}}$ in $\mathcal{K}(A,g)$ such that $Av_n \to z = Aw$ as $n \to \infty$. Then also $v_n \to w$, because $\|v_n - w\|_\mathcal{H} \leqslant \|A^{-1}\|_{\mathrm{op}} \|Av_n - z\|_\mathcal{H}$. Since $v_n \perp w$ for each $n$, then

$$0 = \lim_{n \to \infty} \|v_n - w\|_\mathcal{H}^2 = \lim_{n \to \infty} \left( \|v_n\|_\mathcal{H}^2 + \|w\|_\mathcal{H}^2 \right) = 2\|w\|_\mathcal{H}^2,$$

whence $w = 0$ and also $z = Aw = 0$. □

When $A$ is not injective, from the perspective of the Krylov-solvability of the inverse problem (4.4) one immediately makes two observations. First, if $g = 0$, then trivially $\mathcal{K}(A,g) = \{0\}$ and therefore the Krylov space does not capture any of the non-zero solutions to (4.4), which all belong to $\ker A$. Second, if $g \neq 0$ and therefore the problem of Krylov-solvability is non-trivial, it is natural to ask whether a Krylov solution exists and whether it is unique.

The following two Propositions provide answers to these questions.

**Proposition 4.3** *For a given Hilbert space $\mathcal{H}$ let $A$ be a densely defined operator on $\mathcal{H}$ and let $g \in \mathrm{ran}(A) \cap C^\infty(A)$. If $\ker A \subset \ker A^*$ (in particular, if $A$ is normal), then there exists at most one $f \in \overline{\mathcal{K}(A,g)} \cap \mathcal{D}(A)$ such that $Af = g$.*

***Proof*** The argument is analogous to that used for the corresponding statement of Theorem 4.2. If $f_1, f_2 \in \overline{\mathcal{K}(A,g)}$ and $Af_1 = g = Af_2$, then $f_1 - f_2 \in \ker A \cap \overline{\mathcal{K}(A,g)}$. By assumption, $\ker A \subset \ker A^*$, and moreover obviously $\mathcal{K}(A,g) \subset \overline{\mathrm{ran}\,A}$. Therefore, $f_1 - f_2 \in \ker A^* \cap \overline{\mathrm{ran}\,A} = \{0\}$, whence $f_1 = f_2$. □

**Proposition 4.4** *For a given Hilbert space $\mathcal{H}$ let $A$ be a densely defined and closed operator on $\mathcal{H}$, let $g \in \mathrm{ran}(A) \cap C^\infty(A)$, and let $f \in \mathcal{D}(A)$ satisfy $Af = g$. Assume furthermore that*

*(a) $A(\overline{\mathcal{K}(A,g)} \cap \mathcal{D}(A)) \subset \overline{\mathcal{K}(A,g)}$;*
*(b) $P_\mathcal{K} f \in \mathcal{D}(A)$;*
*(c) $\mathcal{I}(A,g) = \{0\}$,*

*or also, assume the more stringent assumptions*

*(a') $A$ is $\mathcal{K}(A,g)$-reduced*
*(b') $P_\mathcal{K} f \in \mathcal{D}(A)$.*

*Then there exists $f_\circ \in \overline{\mathcal{K}(A,g)} \cap \mathcal{D}(A)$ such that $Af_\circ = g$.*

***Proof*** Let $f$ be a solution to $Af = g$ (it certainly exists, since $\mathcal{S}(A,g)$ is non-empty). Reasoning as in the proof of Proposition 4.1: $P_\mathcal{K} f \in \overline{\mathcal{K}(A,g)} \cap \mathcal{D}(A)$ (by (b)), $AP_\mathcal{K} f \in \overline{\mathcal{K}(A,g)}$ (by (a)), $(\mathbb{1} - P_\mathcal{K})f \in \mathcal{D}(A)$ (by (b)), whence $A(\mathbb{1} - P_\mathcal{K})f \in A(\mathcal{K}(A,g)^\perp \cap \mathcal{D}(A))$ and also $A(\mathbb{1} - P_\mathcal{K})f = g - AP_\mathcal{K} f \in \overline{\mathcal{K}(A,g)}$. Thus, $A(\mathbb{1} - P_\mathcal{K})f \in \mathcal{I}(A,g)$, whence $g = Af = AP_\mathcal{K} f$ (by (c)). Set

$$f_\circ := P_\mathcal{K} f \in \overline{\mathcal{K}(A,g)} \cap \mathcal{D}(A).$$

Now, if $g = 0$, then $\overline{\mathcal{K}(A,g)} = \{0\}$, whence $f_\circ = 0$: this (trivial) solution to the corresponding inverse problem is indeed a Krylov-solution. If instead $g \neq 0$, then necessarily $f_\circ \neq 0$. In either case $Af_\circ = g$. □

## 4.5 The self-adjoint case revisited: structural properties.

We can now re-consider the question of Krylov-solvability of unbounded, self-adjoint inverse problems. The analysis of Section 4.2 (Theorem 4.2) already provides an affirmative answer, based upon conjugate gradient arguments, yet it does not explain how the operator $A$ and the datum $g$ are connected, in the self-adjoint case, to the abstract operator-theoretic mechanisms for Krylov-solvability identified in Section 4.4, namely Krylov-reducibility, triviality of the Krylov intersection, stability of the solution manifold inside $\mathcal{D}(A)$ under the projection $P_\mathcal{K}$.

Of course, if $A = A^*$ and $g \in C^\infty(A)$, then the second of the two conditions (4.14) for Krylov reducibility is automatically true (owing to (4.7)) and therefore $\mathcal{I}(A,g)$ is always trivial.

Unlike the bounded case, however, in order for $A$ to be $\mathcal{K}(A,g)$-reduced no Krylov escape must occur, namely $A$ has to match also the first of the two conditions (4.14), and we have already observed (Remark 4.6) that an assumption on $A$ and $g$ such as the Krylov-core condition would indeed prevent the Krylov escape phenomenon.

As relevant as such issues are from a more abstract perspective in order to understand why 'structurally' a self-adjoint inverse problem is Krylov-solvable, and as deceptively simple the underlying mathematical questions appear, yet to our knowledge no complete answer is available in the affirmative (i.e., a proof) or in the negative (a counterexample) to the following questions:

**Q1**: When $A = A^*$ and $g \in C^\infty(A)$, is it true that $\overline{\mathcal{K}(A,g)}^{\|\,\|_A} = \overline{\mathcal{K}(A,g)} \cap \mathcal{D}(A)$, i.e., is the Krylov core condition satisfied by the pair $(A,g)$ ?

**Q2**: When $A = A^*$ and $g \in C^\infty(A)$, is it true that $A\big(\overline{\mathcal{K}(A,g)} \cap \mathcal{D}(A)\big) \subset \overline{\mathcal{K}(A,g)}$, i.e., is $A$ $\mathcal{K}(A,g)$-reduced in the generalised sense?

(Clearly in (**Q2**) it is tacitly understood that $\mathcal{K}(A,g)$ is not dense in $\mathcal{H}$.)

We can provide a partial answer in a vast class of cases, namely whenever the vector $g$ is 'bounded' for $A$.

Let us recall (see, e.g., [97, Section 7.4]) that a $C^\infty$-vector $g \in \mathcal{H}$ for a linear operator $A$ on a Hilbert space $\mathcal{H}$ is called a **bounded vector** for $A$ when there is a constant $B_g > 0$ such that

$$\|A^n g\|_\mathcal{H} \leqslant B_g^n \qquad \forall n \in \mathbb{N}. \tag{4.17}$$

As well-known (see, e.g., [97, Lemma 7.13]), the vector space of bounded vectors for $A$, when $A$ is self-adjoint, is dense in $\mathcal{H}$.

The special relevance of the assumption that $g$ be a bounded vector for $A$ emerges from Theorems 4.3 and 4.4 below. The former establishes a natural isomorphism between the corresponding Krylov space and $L^2(\mathbb{R}, d\mu_g^{(A)})$. Again, we refer to [97, Section 5.3] for the standard theory of the functional calculus of $A$.

**Theorem 4.3** *Let $A$ be a self-adjoint operator on the Hilbert space $\mathcal{H}$ and let $g \in \mathcal{H}$ be a bounded vector for $A$. Then there is the Hilbert space isomorphism*

$$L^2(\mathbb{R}, d\mu_g^{(A)}) \xrightarrow{\cong} \overline{\mathcal{K}(A, g)}, \qquad (4.18)$$
$$\varphi \longmapsto \varphi(A)g,$$

*where $\varphi(A)$ is the operator constructed with the measurable functional calculus of $A$.*

Observe that Theorem 4.3 is the counterpart in the unbounded self-adjoint setting of Proposition 2.8 in the bounded self-adjoint case. For its proof it is convenient to single out the following facts.

**Lemma 4.5** *Let $\mu$ be a positive finite measure on $\mathbb{R}$. Then the space of $\mathbb{R} \to \mathbb{C}$ functions that are restrictions of $\mathbb{C} \to \mathbb{C}$ entire functions and are square-integrable with respect to $\mu$ is dense in $L^2(\mathbb{R}, d\mu)$.*

*Proof* Let $f \in L^2(\mathbb{R}, d\mu)$. For every $\varepsilon > 0$ then there exists a continuous and $L^2(\mathbb{R}, d\mu)$-function $f_c$ such that $\|f - f_c\|_{L^2(\mathbb{R}, d\mu)} \leqslant \varepsilon$ (see, e.g., [91, Chapt. 3, Theorem 3.14]). In turn, by Carleman's theorem (see, e.g., [42, Chapt. 4, §3, Theorem 1]), there exists an entire function $f_e$ such that $|f_c(\lambda) - f_e(\lambda)| \leqslant \mathcal{E}(\lambda) \, \forall \lambda \in \mathbb{R}$, where $\mathcal{E}$ is an arbitrary error function. $\mathcal{E}$ can be therefore chosen so as $\|f_c - f_e\|_{L^2(\mathbb{R}, d\mu)} \leqslant \varepsilon$. This shows that $f_e \in L^2(\mathbb{R}, d\mu)$ and $\|f - f_e\|_{L^2(\mathbb{R}, d\mu)} \leqslant 2\varepsilon$. As $\varepsilon$ is arbitrary, this completes the proof. $\square$

**Lemma 4.6** *Let $A$ be a self-adjoint operator on the Hilbert space $\mathcal{H}$ and let $g \in \mathcal{H}$ be a bounded vector for $A$. If $f : \mathbb{C} \to \mathbb{C}$ is an entire function, then its restriction to the real line belongs to $L^2(\mathbb{R}, d\mu_g^{(A)})$ and $g \in \mathcal{D}(f(A))$.*

*Proof* As $f$ is entire, in particular $f(\lambda) = \sum_{k=0}^{\infty} \frac{f^{(k)}(0)}{k!} \lambda^k$ for every $\lambda \in \mathbb{R}$, where the series converges point-wise for every $\lambda$ and uniformly on compact subsets of $\mathbb{R}$. Thus, by functional calculus (see, e.g., [97, Proposition 4.12(v)]), for every $N \in \mathbb{N}$

$$\mathbf{1}_{[-N,N]}(A)f(A)g = \sum_{k=0}^{\infty} \frac{f^{(k)}(0)}{k!} (A\mathbf{1}_{[-N,N]}(A))^k g,$$

where $\mathbf{1}_\Omega$ denotes the characteristic function of $\Omega$ and the series in the r.h.s. above converges in the norm of $\mathcal{H}$. Again by standard properties of the functional calculus one then has

$$\left(\int_{-N}^{N}|f(\lambda)|^2\,d\mu_g^{(A)}(\lambda)\right)^{1/2} = \|\mathbf{1}_{[-N,N]}(A)f(A)g\|_{\mathcal{H}}$$

$$= \Big\|\sum_{k=0}^{\infty}\frac{f^{(k)}(0)}{k!}\mathbf{1}_{[-N,N]}(A)A^k g\Big\|_{\mathcal{H}}$$

$$\leqslant \sum_{k=0}^{\infty}\Big|\frac{f^{(k)}(0)}{k!}\Big|\,\|\mathbf{1}_{[-N,N]}(A)A^k g\|_{\mathcal{H}}$$

$$\leqslant \sum_{k=0}^{\infty}\Big|\frac{f^{(k)}(0)}{k!}\Big|\,\|A^k g\|_{\mathcal{H}}$$

$$\leqslant M_R\sum_{k=0}^{\infty}\Big(\frac{B_g}{R}\Big)^k = M_R(1-B_g/R)^{-1}\,,$$

where in the last inequality we used the bound $\|A^k g\|_{\mathcal{H}} \leqslant B_g^k$ for some $B_g > 0$ and we applied Cauchy's estimate $|\frac{f^{(k)}(0)}{k!}| \leqslant M_R/R^k$ for $R > B_g$ and $M_R := \max_{|z|=R}|f(x)|$. Taking the limit $N \to \infty$ yields the thesis. $\square$

**Lemma 4.7** *Let A be a self-adjoint operator on the Hilbert space $\mathcal{H}$ and let $g \in \mathcal{H}$ be a bounded vector for A. Let $f : \mathbb{C} \to \mathbb{C}$ be an entire function and let $f(z) = \sum_{k=0}^{\infty}\frac{f^{(k)}(0)}{k!}z^k$, $z \in \mathbb{C}$, be its Taylor expansion. Then,*

$$\lim_{n\to\infty}\Big\|f(A)g - \sum_{k=0}^{n}\frac{f^{(k)}(0)}{k!}A^k g\Big\|_{\mathcal{H}} = 0$$

*and therefore $f(A)g \in \overline{\mathcal{K}(A,g)}$.*

*Proof* Both $f$ and $\mathbb{C} \ni z \mapsto \sum_{k=0}^{n}\frac{f^{(k)}(0)}{k!}z^k$ are entire functions for each $n \in \mathbb{N}$, and so is their difference. Then, reasoning as in the proof of Lemma 4.6,

$$\Big\|f(A)g - \sum_{k=0}^{n}\frac{f^{(k)}(0)}{k!}A^k g\Big\|_{\mathcal{H}} = \Big\|\sum_{k=n+1}^{\infty}\frac{f^{(k)}(0)}{k!}A^k g\Big\|_{\mathcal{H}}$$

$$\leqslant \sum_{k=n+1}^{\infty}\Big|\frac{f^{(k)}(0)}{k!}\Big|\,\|A^k g\|_{\mathcal{H}}$$

$$\leqslant M_R\sum_{k=n+1}^{\infty}\Big(\frac{B_g}{R}\Big)^k$$

$$= M_R(B_g/R)^{n+1}(1-B_g/R)^{-1} \xrightarrow{n\to\infty} 0,$$

which completes the proof. $\square$

Let us finally prove Theorem 4.3.

*Proof (Proof of Theorem 4.3)* Let us denote by $\mathscr{E}(\mathbb{R})$ the space of $\mathbb{R} \to \mathbb{C}$ functions that are restrictions of $\mathbb{C} \to \mathbb{C}$ entire functions. Owing to Lemma 4.6, $\mathscr{E}(\mathbb{R}) \subset L^2(\mathbb{R},d\mu_g^{(A)})$, and owing to Lemma 4.5, $\mathscr{E}(\mathbb{R})$ is actually dense in $L^2(\mathbb{R},d\mu_g^{(A)})$.

Moreover, for each $\varphi \in \mathscr{E}(\mathbb{R})$, $\varphi(A)g \in \overline{\mathcal{K}(A,g)}$, as found in Lemma 4.7. As $\|\varphi(A)g\|_{\mathcal{H}} = \|\varphi\|_{L^2(\mathbb{R},\mathrm{d}\mu_g^{(A)})}$, the map $\varphi \mapsto \varphi(A)g$ is an isometry from $\mathscr{E}(\mathbb{R})$ to $\overline{\mathcal{K}(A,g)}$, which extends canonically to a $L^2(\mathbb{R},\mathrm{d}\mu_g^{(A)}) \to \overline{\mathcal{K}(A,g)}$ isometry, by standard linear bounded extension [88, Theorem I.7]. In fact, this map is also surjective and therefore is an isomorphism. Indeed, for a generic $w \in \overline{\mathcal{K}(A,g)}$ there exists a sequence $(p_n)_{n \in \mathbb{N}}$ of polynomials such that $p_n(A)g \to w$ in $\mathcal{H}$ as $n \to \infty$, and moreover $(p_n)_{n \in \mathbb{N}}$ is a Cauchy sequence in $L^2(\mathbb{R},\mathrm{d}\mu_g^{(A)})$ because $\|p_n - p_m\|_{L^2(\mathbb{R},\mathrm{d}\mu_g^{(A)})} = \|p_n(A)g - p_m(A)g\|_{\mathcal{H}}$ for any $n,m \in \mathbb{N}$. Therefore $p_n \to \varphi$ in $L^2(\mathbb{R},\mathrm{d}\mu_g^{(A)})$ for some $\varphi$ and consequently $w = \lim_{n \to \infty} p_n(A)g = \varphi(A)g$. $\square$

*Remark 4.8* The reasoning of the proof of Theorem 4.3 reveals that for *generic* $g \in C^\infty(A)$, not necessarily bounded for $A$, the map $h \mapsto h(A)g$ is an isomorphism

$$\overline{\{\text{polynomials } \mathbb{R} \to \mathbb{C}\}}^{\|\cdot\|_{L^2(\mathbb{R},\mathrm{d}\mu_g^{(A)})}} \xrightarrow{\cong} \overline{\mathcal{K}(A,g)}. \tag{4.19}$$

The assumption of *A*-boundedness for $g$ was used to show, concerning the l.h.s. above, that polynomials are dense in $L^2(\mathbb{R},\mathrm{d}\mu_g^{(A)})$.

In turn, following the same line of reasoning, Theorem 4.3 yields an alternative proof of Corollary 2.1: in this case we re-obtain Krylov solvability for the inverse problem (4.4) when $A$ is self-adjoint and injective.

***Proof (Alternative proof of Theorem 4.2, injective self-adjoint case)*** Let $\varphi$ be the $\mu_g^{(A)}$-measurable function defined by $\varphi(\lambda) := \lambda^{-1}$ for $\lambda \in \mathbb{R}$, and let $f = A^{-1}g$ the unique solution to $Af = g$ (by injectivity $A$ is invertible on its range). Since $\|f\|_{\mathcal{H}}^2 = \int_{\mathbb{R}} \lambda^{-2} \mathrm{d}\mu_g^{(A)}(\lambda)$, then $\varphi \in L^2(\mathbb{R},\mathrm{d}\mu_g^{(A)})$. Therefore, owing to Theorem 4.3, $f = \varphi(A)g \in \overline{\mathcal{K}(A,g)}$. $\square$

The other noticeable result based upon the boundedness of the datum $g$ is the following partial answer to question (**Q1**) above.

**Theorem 4.4** *Let A be a self-adjoint operator on the Hilbert space $\mathcal{H}$ and let $g \in \mathcal{H}$ be a bounded vector for A. Then the pair $(A,g)$ satisfies the Krylov-core condition and consequently A is $\mathcal{K}(A,g)$-reduced.*

As observed already (Remark 4.4), the thesis follows from the inclusion

$$\overline{\mathcal{K}(A,g)}^{\|\cdot\|_A} \supset \overline{\mathcal{K}(A,g)} \cap \mathcal{D}(A) \tag{4.20}$$

that we shall prove now.

Let us denote as usual by $E^{(A)}$ the spectral measure for $A$, and by $\mu_x^{(A)}$ the scalar spectral measure associated with $A$ and $x \in \mathcal{H}$, i.e., $\mu_x^{(A)}(\Omega) = \langle x, E^{(A)}(\Omega)x \rangle$ for every Borel subset $\Omega \subset \mathbb{R}$ [97, Section 5.2-5.3].

*Proof (Proof of Theorem 4.4)* Let $x \in \overline{\mathcal{K}(A,g)} \cap \mathcal{D}(A)$. We need to exhibit a sequence $(x_n)_{n \in \mathbb{N}}$ in $\mathcal{K}(A,g)$ such that $x_n \xrightarrow{\|\;\|_A} x$ as $n \to \infty$. Since $x \in \overline{\mathcal{K}(A,g)}$, there surely exists a sequence $(p_n(A)g)_{n \in \mathbb{N}}$, for some polynomials $p_n$, that converges to $x$ in the $\mathcal{H}$-norm: we shall show that convergence occurs in the stronger $\|\;\|_A$-norm too.

Theorem 4.3 implies that there exists some $\varphi \in L^2(\mathbb{R}, \mathrm{d}\mu_g^{(A)})$ such that $x = \varphi(A)g$. As

$$\|p_n(A)g - x\|_{\mathcal{H}}^2 = \int_{\mathbb{R}} |p_n(\lambda) - \varphi(\lambda)|^2 \, \mathrm{d}\mu_g^{(A)}, \tag{*}$$

the function $\varphi$ is the $L^2(\mathbb{R}, \mathrm{d}\mu_g^{(A)})$-limit of the polynomials $p_n$'s. Moreover, $x \in \mathcal{D}(A)$, equivalently, $\lambda \varphi \in L^2(\mathbb{R}, \mathrm{d}\mu_g^{(A)})$, whence

$$\|Ap_n(A)g - Ax\|_{\mathcal{H}}^2 = \int_{\mathbb{R}} |\lambda(p_n(\lambda) - \varphi(\lambda))|^2 \, \mathrm{d}\mu_g^{(A)}. \tag{**}$$

In fact, for arbitrary $n \in \mathbb{N}$ and $r \in (0,1]$ one has

$$\int_{\mathbb{R}} |\lambda(p_n(\lambda) - \varphi(\lambda))|^{2r} \, \mathrm{d}\mu_g^{(A)} < +\infty.$$

Indeed, setting $\Omega_n := \{\lambda \in \mathbb{R}; |\lambda(p_n(\lambda) - \varphi(\lambda))|^2 > 1\}$ and taking $r \in (0,1]$,

$$\begin{aligned}
\int_{\mathbb{R}} |\lambda(p_n(\lambda) - \varphi(\lambda))|^{2r} \, \mathrm{d}\mu_g^{(A)} &= \int_{\mathbb{R} \setminus \Omega_n} |\lambda(p_n(\lambda) - \varphi(\lambda))|^{2r} \, \mathrm{d}\mu_g^{(A)} \\
&\quad + \int_{\Omega_n} |\lambda(p_n(\lambda) - \varphi(\lambda))|^{2r} \, \mathrm{d}\mu_g^{(A)} \\
&\leqslant \int_{\mathbb{R} \setminus \Omega_n} \mathrm{d}\mu_g^{(A)} + \int_{\Omega_n} |\lambda(p_n(\lambda) - \varphi(\lambda))|^2 \, \mathrm{d}\mu_g^{(A)} \\
&\leqslant \int_{\mathbb{R}} \mathrm{d}\mu_g^{(A)} + \int_{\mathbb{R}} |\lambda(p_n(\lambda) - \varphi(\lambda))|^2 \, \mathrm{d}\mu_g^{(A)} \\
&= \|g\|_{\mathcal{H}}^2 + \|Ap_n(A)g - Ax\|_{\mathcal{H}}^2 < +\infty.
\end{aligned}$$

It is standard to see (see, e.g., [91, Chapter 3, Exercise 4(b)]) that for each $n \in \mathbb{N}$

$$f_n(r) := \int_{\mathbb{R}} |\lambda(p_n(\lambda) - \varphi(\lambda))|^{2r} \, \mathrm{d}\mu_g^{(A)}$$

defines a continuous function $f_n$ on the set $R_n := \{r \in (0, +\infty]; f_n(r) < +\infty\}$. For what shown above, $(0,1] \subset R_n$ for all $n \in \mathbb{N}$, whence

$$f_n(r) \xrightarrow{r \uparrow 1} f_n(1) = \int_{\mathbb{R}} |\lambda(p_n(\lambda) - \varphi(\lambda))|^2 \, \mathrm{d}\mu_g^{(A)} \qquad \forall n \in \mathbb{N}.$$

Each $f_n$ is estimated by means of the generalised Hölder inequality as

$$f_n(r) \leqslant \left( \int_{\mathbb{R}} |p_n(\lambda) - \varphi(\lambda)|^{2p} \, d\mu_g^{(A)} \right)^{\frac{r}{p}} \left( \int_{\mathbb{R}} |\lambda|^{2q} \, d\mu_g^{(A)} \right)^{\frac{r}{q}}$$

for $r \in (0,1]$ as fixed above and $p,q \in (0,+\infty)$ such that $p^{-1} + q^{-1} = r^{-1}$. The special choice $p = 1$ and $q \in \mathbb{N}$ guarantees that $r = \frac{q}{1+q} \in (0,1)$ and yields

$$f_n(r) \leqslant \|p_n(A)g - x\|_{\mathcal{H}}^{2r} B_g^{2r},$$

having used (*) and $\left( \int_{\mathbb{R}} |\lambda|^{2q} \, d\mu_g^{(A)} \right)^{\frac{r}{q}} = \|A^q g\|_{\mathcal{H}}^{\frac{2r}{q}} \leqslant B_g^{2r}$, where $B_g$ is the constant from (4.17) owing to the assumption of boundedness for $g$. Letting $q \to \infty$, and hence $r \uparrow 1$ in the above bound for $f_n(r)$, and using the continuity of $f_n$ on $R_n$ and (**), gives

$$\|A p_n(A)g - Ax\|_{\mathcal{H}}^2 \leqslant \|p_n(A)g - x\|_{\mathcal{H}}^2 B_g^2.$$

Thus,

$$\|p_n(A)g - x\|_A^2 \leqslant (1 + B_g^2) \|p_n(A)g - x\|_{\mathcal{H}}^2.$$

Therefore, the $\mathcal{H}$-norm convergence $p_n(A)g \to x$ is equivalent to the graph norm convergence. □

*Remark 4.9* In [19] we gave an alternative proof of Theorem 4.4 that, at the cost of being somewhat lengthier, identifies more constructively the actual sequence of graph norm approximants $p_n(A)g \in \mathcal{K}(A,g)$.

## 4.6 Remarks on rational Krylov subspaces and solvability of self-adjoint inverse problems

It is natural to supplement the analysis of the previous Section with a few additional remarks on the Krylov solvability with respect to *rational* Krylov subspaces.

Those defined as in (4.2)-(4.3) above are *standard* (also referred to as *polynomial*) Krylov subspaces, and are built using positive powers of the operator $A$ applied to a vector $g$. The concept of *rational* Krylov subspace, instead, uses the class of rational functions of $A$ applied to the vector $g$ to build up a vector space.

As a general definition, to a given closed operator $A$ acting on a Hilbert space $\mathcal{H}$, a vector $g \in \mathcal{H}$, and a sequence $\Xi \equiv (\xi_n)_{n \in \mathbb{N}}$ in the resolvent set $\rho(A)$ of $A$, one associates the **rational Krylov subspace**

$$\mathcal{K}^{(\Xi)}(A,g) := \mathrm{span}\left\{ g, \prod_{n=1}^{m}(A - \xi_n \mathbb{1})^{-1} g \,\bigg|\, m \in \mathbb{N} \right\} \tag{4.21}$$

as well as, for each $N \in \mathbb{N}$, the **$N$-th order rational Krylov subspace**

$$\mathcal{K}_N^{(\Xi)}(A,g) := \mathrm{span}\left\{ g, \prod_{n=1}^{m}(A - \xi_n \mathbb{1})^{-1} g \,\bigg|\, m \in \{1, \ldots N-1\} \right\}. \tag{4.22}$$

Rational Krylov subspaces were first explicitly introduced in the finite-dimensional setting [92] for solving the eigenvalue problem, and they now constitute an accepted and well-studied area of the numerical literature [36, 92, 46, 93, 94, 94, 79, 114, 85, 63, 32, 10, 54]. These spaces are particularly attractive in the study of numerical techniques of time-dependent partial differential equations [46, 79, 114, 32, 59] as they have applications in approximating time dependent functions used in the time-stepping of the solution. However, typically these studies have focussed on the approximation of certain operator functions and their eigenvalue eigenvector pairs, rather than the Krylov solvability of the linear inverse problem. Furthermore, these studies are mostly restricted to the finite-dimensional setting, and are also restricted to the class of bounded linear operators.

The results discussed in this Section show that also in the case of unbounded operators, from an approximation standpoint it may be advantageous to consider general rational approximations rather than standard polynomials. Additionally, the restriction that $g \in C^\infty(A)$ may be dropped. The practical drawback of such an approach is that there is an extra computational cost in calculating the operator-valued resolvent function $\rho(A) \to \mathcal{B}(\mathcal{H})$, $\xi \mapsto (A - \xi \mathbb{1})^{-1}$, although there are several examples of operators of great relevance in physics whose resolvent is a finite-rank perturbation of some reference 'free Hamiltonian' and hence has an explicit and particularly manageable expression [4, 76, 50, 45, 44, 43, 77, 78].

In this context, the solution $f$ to an inverse problem $Af = g$ as (4.4) for a given $g \in \mathcal{H}$ and a given closed operator $A$ acting on $\mathcal{H}$ is said to be a **rational Krylov solution** with respect to a given sequence $\Xi \equiv (\xi_n)_{n \in \mathbb{N}}$ in $\rho(A)$ when $f \in \overline{\mathcal{K}^{(\Xi)}(A,g)}$. Any such inverse problem admitting a rational Krylov solution is said **rational-Krylov solvable**.

While we do not have room in this monograph to make an exhaustive discussion of rational Krylov subspaces and rational solvability in general, in light of the previous Section it is worth highlighting a few important properties of the self-adjoint inverse problem.

**Lemma 4.8** *Let $A$ be a self-adjoint operator on the Hilbert space $\mathcal{H}$, let $g \in \mathcal{H}$, and let $(\xi_n)_{n \in \mathbb{N}}$ be a sequence in the resolvent set $\rho(A)$ which contains, for every element, also its complex conjugate. Denote by $\mathfrak{I}$ the collection of all finite index sets generated from $\mathbb{N}$, and for every $I \in \mathfrak{I}$ let $R_I : \sigma(A) \to \mathbb{C}$ be the rational function*

$$R_I(\lambda) := \prod_{n \in I} \frac{1}{\lambda - \xi_n}, \qquad \lambda \in \sigma(A). \tag{4.23}$$

*Then the linear space*

$$\mathfrak{R} := \mathrm{span}\{R_I \,|\, I \in \mathfrak{I}\} \tag{4.24}$$

*is dense in $L^2(\sigma(A), d\mu_g^{(A)})$, where $\mu_g^{(A)}$ is the scalar spectral measure associated with $A$ and $g$.*

***Proof*** The Stone-Weierstrass theorem for locally compact space is applicable to $\mathfrak{R}$. Indeed, $\mathfrak{R}$ is an involutive sub-algebra of $C_0(\sigma(A), \mathbb{C})$ (the complex-valued continuous functions vanishing at infinity, equipped with the supremum norm), as it contains, for each $R_I = \prod_{n \in I}(\lambda - \xi_n)^{-1}$, also $\overline{R_I} = \prod_{n \in I}(\lambda - \overline{\xi_n})^{-1}$, and is clearly closed under usual point-wise sums and products. Furthermore, $\sigma(A)$ is closed in $\mathbb{C}$ and therefore locally compact, and $\mathfrak{R}$ separates points in $\sigma(A)$ and vanishes nowhere on $\sigma(A)$. Thus, the Stone-Weierstrass theorem [88, Theorem IV.10] implies that $\mathfrak{R}$ is dense in $C_0(\sigma(A), \mathbb{C})$ in the supremum norm. On the other hand, $C_0(\sigma(A), \mathbb{C})$ is dense in $L^2(\sigma(A), d\mu_g^{(A)})$, because it contains the $L^2$-dense subspace $C_c(\sigma(A), \mathbb{C})$ [91, Theorem 3.14]. (Recall that $\mu_g^{(A)}$ is a regular Borel measure.) Thus, for any $\varphi \in L^2(\sigma(A), d\mu_g^{(A)})$ and arbitrary $\varepsilon > 0$ one finds some $h^{(\varepsilon)} \in C_0(\sigma(A), \mathbb{C})$ such that $\|\varphi - h^{(\varepsilon)}\|_{L^2} \leqslant \varepsilon$, and also $R_I^{(\varepsilon)} \in \mathfrak{R}$ such that $\|h^{(\varepsilon)} - R_I^{(\varepsilon)}\|_{L^\infty} \leqslant \varepsilon$. Therefore,

$$\begin{aligned}\|\varphi - R_I^{(\varepsilon)}\|_{L^2}^2 &\leqslant 2\|\varphi - h^{(\varepsilon)}\|_{L^2}^2 + 2\|h^{(\varepsilon)} - R_I^{(\varepsilon)}\|_{L^2}^2 \\ &\leqslant 2\varepsilon^2 + 2\|h^{(\varepsilon)} - R_I^{(\varepsilon)}\|_{L^\infty}^2 \int_{\sigma(A)} d\mu_g^{(A)}(\lambda) \\ &= 2\varepsilon^2(1 + \|g\|_{\mathcal{H}}^2), \end{aligned}$$

which proves that $\mathfrak{R}$ is dense in $L^2(\sigma(A), d\mu_g^{(A)})$. $\square$

*Remark 4.10* As an immediate corollary of (the proof of) Lemma 4.8, one has that the map $R_I \mapsto R_I(A)g$ induces the Hilbert space isomorphism

$$L^2(\mathbb{R}, d\mu_g^{(A)}) \cong L^2(\sigma(A), d\mu_g^{(A)}) \xrightarrow{\cong} \overline{\operatorname{span}\{R_I(A)g \,|\, I \in \mathfrak{I}\}} \\ \varphi \longmapsto \varphi(A)g \tag{4.25}$$

for the given sequence $\Xi$. This is completely analogous to the Hilbert space isomorphism $L^2(\mathbb{R}, d\mu_g^{(A)}) \cong \overline{\mathcal{K}(A, g)}$ (Theorem 4.3), as well as to the isomorphism (4.19) (Remark 4.8).

Among various relevant approximation results in the self-adjoint case, let us single out the following one.

**Proposition 4.5** *Let $A$ be a self-adjoint injective operator on the Hilbert space $\mathcal{H}$, let $g \in \operatorname{ran} A$, let $f \in \mathcal{D}(A)$ be the unique solution to $Af = g$, and let $\Xi \equiv (\xi_n)_{n \in \mathbb{N}}$ be a sequence in the resolvent set $\rho(A)$ which contains, for every element, also its complex conjugate. Then, in the notation of Lemma 4.8,*

$$f \in \overline{\operatorname{span}\{R_I(A)g \,|\, I \in \mathfrak{I}\}}. \tag{4.26}$$

A slightly modified version of the above approximation is also worth being highlighted.

**Proposition 4.6** *Let $A$ be a self-adjoint injective operator on the Hilbert space $\mathcal{H}$ admitting a real point in its resolvent, let $g \in \operatorname{ran} A$, let $f \in \mathcal{D}(A)$ be the unique solution to $Af = g$, and let $\xi \in \mathbb{R} \cap \rho(A)$. Then*

$$f \in \overline{\operatorname{span}\{(A - \xi \mathbb{1})^{-m} \,|\, m \in \mathbb{N}\}}. \tag{4.27}$$

*Proof (Proof of Proposition 4.5)* As the sequence $\Xi$ satisfies the assumptions of Lemma 4.8, the corresponding functional space $\mathfrak{R}$ defined in (4.24) is dense in $L^2\big(\sigma(A), \mathrm{d}\mu_g^{(A)}\big)$. On the other hand, since $A$ is injective, one has $f = A^{-1}g$: therefore, as argued in the alternative proof of Theorem 4.2 (Section 4.5), and as a consequence of the isomorphism established in Theorem 4.3, the $\mu_g^{(A)}$-measurable function defined by $\varphi(\lambda) := \lambda^{-1}$ for $\lambda \in \sigma(A)$ belongs to $L^2\big(\sigma(A), \mathrm{d}\mu_g^{(A)}\big)$. Thus, $\varphi$ is the $L^2$-limit of functions from $\mathfrak{R}$, implying that $f = \varphi(A)g$ is the $\mathcal{H}$-limit of finite linear combinations of vectors of the form $R_I(A)g$ with $R_I \in \mathfrak{R}$. □

*Proof (Proof of Proposition 4.6)* Applying Lemma 4.8, as in the proof of Proposition 4.5, with respect to the sequence $\Xi \equiv (\xi_n)_{n \in \mathbb{N}}$ defined by $\xi_n := \xi \ \forall n \in \mathbb{N}$, yields the conclusion that $f$ is the $\mathcal{H}$-limit of vectors of the form $\prod_{n=1}^m (A - \xi \mathbb{1})^{-1}g$, whence precisely (4.27). □

# Chapter 5
# Krylov solvability in a perturbative framework

## 5.1 Krylov solvability from a perturbative perspective

So far in this book the discussion has concerned inverse linear problems where both the datum $g$ and the linear operator $A$ are known *exactly*. A more general perspective is to allow for an amount of uncertainty affecting the knowledge of $A$, or $g$, or both. Or, from another point of view, instead of only focussing on the inverse problem of interest, one may consider also an auxiliary, possibly more tractable problem, close in some sense to the original one, which allows for useful approximate information.

Thus, in this Chapter we consider *perturbations* of the original problem $Af = g$ of the form $A'f' = g'$, where $A$ and $A'$, as well as $g$ and $g'$ are close in a controlled sense, and the goal is to study the *effect of the perturbation on the Krylov solvability*.

We keep regarding $A$ and $g$ as exactly known or, in principle, exactly accessible, but with the idea that close to the problem $Af = g$ there is a perturbed problem $A'f' = g'$ that serves as an auxiliary one, possibly more easily tractable, say, with Krylov subspace methods, in order to obtain conclusions on the Krylov solvability of the original problem. Or, conversely, we inquire under which conditions the nice property of Krylov solvability for $Af = g$ is stable enough to survive a small perturbation (that in applications could arise, again, from experimental or numerical uncertainties), or when instead Krylov solvability is washed out by even small inaccuracies in the precise knowledge of $A$ or $g$ – an occurrence in which Krylov subspace methods would prove to be unstable. And, more abstractly, we pose the question of a convenient notion of vicinity between the subspaces $\overline{\mathcal{K}(A,g)}$ and $\overline{\mathcal{K}(A',g')}$ when $A$ and $A'$ (respectively, $g$ and $g'$) are suitably close.

There exists a large amount of literature that accounts for perturbations in Krylov subspace methods, most of which approaches the problem from the point of view of '*inexact Krylov methods*' (see, e.g., [120, 115, 102, 101, 100, 118, 33]). The underlying idea is that the exact (typically non-singular) inverse problem in $\mathbb{C}^N$, $Af = g$, is perturbed in $A$ by a series of linear operators $(E_k)_{k=1}^N$ on $\mathbb{C}^N$ that may change at each step $k$ of the algorithm. Typical scenarios that could induce such perturbations at each step of the algorithm are, but not limited to, truncation and rounding errors in finite precision machines, or approximation errors from calculating complicated matrix-vector products. The main results include the convergence behaviour of the error and residual terms, in particular their rates, and typically bounds on how far these indicators of convergence are from the unperturbed setting at a given iteration number $k$.

Yet, these investigations are of a fundamentally different nature than what we discuss here. To begin with, the typical analysis of inexact Krylov methods is in the finite-dimensional setting, where there is already a fairly complete control on Krylov solvability and rates of convergence. Furthermore, that segment of the literature contains no discussion of how the underlying Krylov subspaces themselves change, as well as the richer phenomena pertaining to the Krylov solvability (or lack-of) of the inverse problem and its perturbations.

It is rather the latter kind of questions that the present Chapter focuses on. This is in fact a perspective that appeared to be essentially uncharted until our recent and partially systematic analysis [21], which the present Chapter is based upon. It is therefore beneficial in the first place to pose a set of general fundamental questions for the perturbative framework: this is the object of Section 5.2, which in a sense sets a meaningful, non-trivial agenda of future investigation.

Then, along this line, in Section 5.3 an overview is presented of typical phenomena that may occur to the Krylov solvability of an inverse problem $Af = g$ in terms of the Krylov solvability, or lack thereof, of auxiliary inverse problems where $A$ or $g$ or both are perturbed in a controlled sense. Noticeably, such a survey indicates that the sole control of the operator or of the data perturbation, in the respective operator and Hilbert norm, still leaves the possibility open to all phenomena such as the persistence, gain, or loss of Krylov solvability in the limit $A_n \to A$ or $g_n \to g$, where $A_n$ (and so $g_n$) is the generic element of a sequence of perturbed objects. The implicit explanation is that an information like $A_n \to A$ or $g_n \to g$ is not enough to account for a suitable vicinity of the corresponding Krylov subspaces. In Chapters 2 and 4 it was discussed how the Krylov solvability of the inverse problem $Af = g$ corresponds to certain structural properties of the subspace $\mathcal{K}(A,g)$, therefore one implicitly needs to monitor how the latter properties are preserved or altered under the perturbation.

This also suggests that performing the perturbation within certain sub-classes of operators may supplement further information on Krylov solvability. In this respect, Section 5.4 focuses on the operators of $\mathscr{K}$-class, previously considered in Section 2.5: the robustness and fragility of this class is discussed from the perturbative perspective of the induced inverse problems.

Then the second part of the Chapter, Sections 5.5-5.7, addresses more systematically the issue of vicinity of Krylov subspaces in a sense that be informative for the Krylov solvability of the corresponding inverse problems. This framework leads to appealing approximation results, as the inner approximability of Krylov subspaces established in Subsection 5.7.2. In particular, Proposition 5.3 is a prototype of the kind of perturbative results we have in mind, namely a control of the perturbation, formulated in terms of the perturbed and unperturbed Krylov subspaces, that *predicts* the persistence of Krylov solvability when the perturbation is removed.

In a sense, this Chapter is aimed at opening a perspective on a general and essentially new level of questions, that are operator-theoretic in nature, yet with direct motivations from Krylov approximation algorithms in numerical computation. The corpus of partial results presented here only scratch the surface of such a scenario. We collect more explicit conclusions and future perspectives in this respect in the final Section 5.8.

## 5.2 Fundamental perturbative questions

One easily realises that the question of the effects of perturbations on the Krylov solvability of an infinite-dimensional inverse problem takes a multitude of related, yet somewhat different formulations depending on the precise perspective one looks at it.

It is therefore instructive to organise the most relevant of such queries into a coherent scheme – which is the goal of this Section. This serves both as a reference for the results and explicit partial answers of this Chapter, as well as an ideal road map for future studies.

In practice, let us discuss the following main categories of connected problems. The first three of them are set for $A \in \mathcal{B}(\mathcal{H})$, $\mathcal{H}$ being as usual the underlying complex infinite-dimensional Hilbert space.

**I. Comparison between "close" Krylov subspaces.** This is the abstract problem of providing a meaningful comparison between $\overline{\mathcal{K}(A,g)}$ and $\overline{\mathcal{K}(A',g')}$, as two closed subspaces of $\mathcal{H}$, for given $A, A' \in \mathcal{B}(\mathcal{H})$ and $g, g' \in \mathcal{H}$ such that in some convenient sense $A$ and $A'$, as well as $g$ and $g'$ are close. As a priori such subspaces might only have a trivial intersection, the framework is rather that of comparison of subspaces of a normed space, in practice introducing convenient topologies or metric distances.

A more application-oriented version of the same problem is the following. Given $A \in \mathcal{B}(\mathcal{H})$ and $g \in \mathcal{H}$, one considers approximants of one or the other (or both), say, sequences $(A_n)_{n \in \mathbb{N}}$ and $(g_n)_{n \in \mathbb{N}}$ respectively in $\mathcal{B}(\mathcal{H})$ and $\mathcal{H}$, such that $\|A_n - A\|_{\mathrm{op}} \to 0$ and $\|g_n - g\| \to 0$ as $n \to \infty$. Then the question is whether a meaningful notion of limit $\overline{\mathcal{K}(A_n, g_n)} \to \overline{\mathcal{K}(A,g)}$ can be defined.

**II. Perturbations preserving/creating Krylov solvability.** This question is inspired by the possibility that, given $A \in \mathcal{B}(\mathcal{H})$ and $g \in \mathrm{ran}A$, instead of solving the "difficult" inverse problem $Af = g$ one solves a convenient perturbed problem $A'f' = g'$, with $A' \in \mathcal{B}(\mathcal{H})$ and $g' \in \mathrm{ran}A'$ close respectively to $A$ and $g$, which is "easily" Krylov solvable, and the Krylov solution of which provides approximate information to the original solution $f$.

Here is an explicit formulation. Assume that one finds $(A_n)_{n \in \mathbb{N}}$ in $\mathcal{B}(\mathcal{H})$ and $(g_n)_{n \in \mathbb{N}}$ in $\mathcal{H}$ such that the inverse problems $A_n f_n = g_n$ are all Krylov solvable and $\|A_n - A\|_{\mathrm{op}} \to 0$ and $\|g_n - g\| \to 0$ as $n \to \infty$. Is $Af = g$ Krylov solvable too? And if at each perturbed level $n$ there is a unique Krylov solution $f_n$, does one have $\|f_n - f\| \to 0$ where $f$ is a (Krylov) solution to $Af = g$?

One scenario of applications is that for $A_n f_n = g_n$ Krylov solvability comes with a much more easily (say, faster) solvable solution algorithm, so that $f$ is rather determined as $f = \lim_{n\to\infty} f_n$ instead of directly approaching the problem $Af = g$.

Another equally relevant scenario is that the possible Krylov solvability of the problem of interest $Af = g$ is initially *unknown*, and prior to launching resource-consuming Krylov algorithms for solving the problem, one wants to be guaranteed that a Krylov solution indeed exists. To this aim one checks the Krylov solvability for $A_n f_n = g_n$ for each $n$, and the convergence $A_n \to A$, $g_n \to g$, thus coming to an affirmative answer.

**III. Perturbations destroying Krylov solvability.** The opposite occurrence has to be monitored as well, namely the possibility that a small perturbation of the Krylov solvable problem $Af = g$ produces a non-Krylov solvable problem $A'f' = g'$. Say, if in the above setting *none* of the problems $A_n f_n = g_n$ are Krylov solvable and yet $A_n \to A$ and $g_n \to g$, under what conditions does one gain Krylov solvability in the limit for the problem $Af = g$? A comprehension of this phenomenon would be of great relevance to identify those circumstances when Krylov methods are intrinsically unstable, in the sense that even a tiny uncertainty in the knowledge of $A$ and/or $g$ yields a perturbed problem $A'f' = g'$ for which, unlike the exact problem of interest $Af = g$, Krylov methods are not applicable.

In the *unbounded* case, more precisely when the operator $A$ is closed and unbounded on $\mathcal{H}$, in principle all the above questions have their own counterpart, except that the fundamental condition $g \in C^\infty(A)$ required to have a meaningful notion of $\mathcal{K}(A,g)$ is highly unstable under perturbations, and one has to ensure case by case that certain problems are well posed.

Yet, for its evident relevance let us highlight the following additional class of questions.

**IV. Perturbations-regularisations exploiting Krylov solvability.** For the problem of interest $Af = g$ one might well have $g \in \mathrm{ran}A$ *but* $g \notin C^\infty(A)$. In this case Krylov methods are not applicable: there is no actual notion of Krylov subspace associated to $A$ and $g$, hence no actual Krylov approximants to utilise iteratively (Section 4.1). Assume though that one finds a sequence $(g_n)_{n\in\mathbb{N}}$ entirely in $\mathrm{ran}A \cap C^\infty(A)$ with $\|g_n - g\| \to 0$. This occurrence is quite typical: if $A$ is a differential operator on $L^2(\mathbb{R}^d)$, everyone is familiar with sequences $(g_n)_{n\in\mathbb{N}}$ of functions that all have high regularity, uniformly in $n$, and for which the $L^2$-limit $g_n \to g$ produces a rough function $g$. Assume further that each problem $Af_n = g_n$ is Krylov solvable. For example, as we showed in Theorem 4.2, for a vast class of self-adjoint or skew-adjoint $A$'s, possibly unbounded, there exists a unique solution $f_n \in \overline{\mathcal{K}(A,g_n)}$. This leads to the following questions. First, do the $f_n$'s have a limit $f$ and does $f$ solve $Af = g$? And, more abstractly speaking, is there a meaningful notion of the limit $\lim_{n\to\infty} \mathcal{K}(A,g_n)$, *irrespectively* of the approximant sequence $(g_n)_{n\in\mathbb{N}}$, that could then be interpreted as a replacement for the non-existing Krylov subspace associated to $A$ and $g$? The elements of such limit subspace would provide exploitable approximants for the solution to the original problem $Af = g$.

One last remark concerns the topologies underlying all the questions above. We explicitly formulated them in terms of the operator norm and Hilbert norm, but al-ternatively there is a variety of weaker notions of convergence that are still highly informative for the solution to the considered inverse problem – a perspective discussed in Appendix A. Thus, the "weaker" counterpart of the above questions represents equally challenging and potentially useful problems to address.

## 5.3 Gain or loss of Krylov solvability under perturbations

Let us survey examples of different behaviours that may occur in those perturbative scenarios belonging to the categories II and III contemplated in the previous Section.

In practice we are comparing here the "unperturbed" inverse linear problem $Af = g$ with "perturbed" problems of the form $Af_n = g_n$, or $A_n f_n = g$, along a sequence $(g_n)_{n\in\mathbb{N}}$ such that $\|g_n - g\| \to 0$, or along a sequence $(A_n)_{n\in\mathbb{N}}$ such that $\|A_n - A\|_{\mathrm{op}} \to 0$. Our particular focus is the Krylov solvability, namely its preservation, or gain, or loss in the limit $n \to \infty$. Next to operator perturbations ($A_n \to A$) and data perturbations ($g_n \to g$), it is also natural to consider simultaneous perturbations of the both of them.

The purpose here is two-fold: we want to convey a concrete flavour of how inverse problems behave under controlled perturbations of the operator or of the datum, as far as having Krylov solutions is concerned, and we also want to highlight the emerging, fundamental, and a priori unexpected lesson. Which is going to be, in short: the sole control that $A_n \to A$ or $g_n \to g$ is not enough to predict whether Krylov solvability is preserved, or gained, or lost in the limit, i.e., each such behaviour can actually occur. The immediate corollary of this conclusion is: one must describe the perturbation of the problem $Af = g$ by means of additional information, say, by restricting to particular sub-classes of inverse problems, or by introducing suitable notions of vicinity of Krylov subspaces, in order to control the effect of the perturbation on Krylov solvability. It is this latter consideration that motivates the more specific discussion of Section 5.4 and of Sections 5.5-5.7.

### 5.3.1 Operator perturbations

*Example 5.1* Let $R$ be the weighted (compact) right shift operator

$$R := \sum_{k=1}^{\infty} \frac{1}{k^2} |e_{k+1}\rangle\langle e_k|$$

on the Hilbert space $\ell^2(\mathbb{N})$ (analogously to Example 2.3), and define

$$R_n := \sum_{k=1}^{n-1} \frac{1}{k^2} |e_{k+1}\rangle\langle e_k| + \frac{1}{n^2} |e_1\rangle\langle e_n|, \qquad n \in \mathbb{N},\ n \geqslant 2.$$

As

$$R - R_n = \sum_{k=n}^{\infty} \frac{1}{k^2} |e_{k+1}\rangle\langle e_k| - \frac{1}{n^2} |e_1\rangle\langle e_n|,$$

$$\|R - R_n\|_{\mathrm{op}} \leqslant \sum_{k=n}^{\infty} \frac{1}{k^2} + \frac{1}{n^2},$$

then $R_n \to R$ in operator norm as $n \to \infty$. For $g := e_2$ and any $n \geqslant 2$, the inverse problem induced by $R_n$ and with datum $g$ has unique solution $f_n := f := e_1$, and so does the inverse problem induced by $R$ and with the same datum, i.e., $R_n f_n = g$ and $Rf = g$. On the other hand,

$$\overline{\mathcal{K}(R_n, g)} = \mathcal{K}(R_n, g) = \mathrm{span}\{e_1, \dots, e_n\},$$
$$\overline{\mathcal{K}(R, g)} = \{e_1\}^\perp.$$

Thus, the inverse problem $R_n f_n = g$ is Krylov solvable, and obviously $f_n \to f$ in norm, yet the inverse problem $Rf = g$ is not.

*Example 5.2* Let $R = \sum_{k=1}^{\infty} |e_{k+1}\rangle\langle e_k|$ be the right shift operator on the Hilbert space $\ell^2(\mathbb{N})$ (Example 2.1), and let

$$A_n := |e_2\rangle\langle e_2| + \frac{1}{n} R, \qquad n \in \mathbb{N},$$
$$A := |e_2\rangle\langle e_2|,$$
$$g := e_2.$$

Clearly, $A_n \to A$ in operator norm as $n \to \infty$. The inverse problem $A_n f_n = g$ has unique solution $f_n = n e_1$; as $(A_n)^k g = \sum_{j=0}^{k} n^{-j} e_{j+2}$ for any $k \in \mathbb{N}_0$, and therefore $\overline{\mathcal{K}(A_n, g)} = \{e_1\}^\perp$, such solution is not a Krylov solution. Instead, passing to the limit, the inverse problem $Af = g$ has unique solution $f = e_2$, which is a Krylov solution since $\overline{\mathcal{K}(A, g)} = \mathrm{span}\{e_2\}$. Observe also that $f_n$ does not converge to $f$.

### 5.3.2 Data perturbations

*Example 5.3* Let $R : \ell^2(\mathbb{Z}) \to \ell^2(\mathbb{Z})$ be the usual right shift operator (Example 2.2). As seen, $R$ is unitary with $R^* = R^{-1} = L$, the left shift operator. Moreover $R$ admits a dense subset $\mathcal{C} \subset \ell^2(\mathbb{Z})$ of cyclic vectors and a dense subset $\mathcal{N} \subset \ell^2(\mathbb{Z})$, consisting of all finite linear combinations of canonical basis vectors, such that the solution $f$ to the inverse problem $Rf = g$ does not belong to $\overline{\mathcal{K}(R, g)}$. All vectors in $\mathcal{N}$ are non-cyclic for $R$.

(i) (Loss of Krylov solvability.) For a datum $g \in \mathcal{N}$, the inverse problem $Rf = g$ admits a unique solution $f$, and $f$ is not a Krylov solution. Yet, by density, there exists a sequence $(g_n)_{n \in \mathbb{N}}$ in $\mathcal{C}$ with $g_n \to g$ (in $\ell^2$-norm) as $n \to \infty$, and each perturbed inverse problem $Rf_n = g_n$ is Krylov solvable with unique solution $f_n = R^{-1}g_n = Lg_n$, and with $f_n \to f$ as $n \to \infty$. Krylov solvability is lost in the limit, still with the approximant Krylov solutions converging to the solution $f$ to the original problem.

(ii) (Gain of Krylov solvability.) For a datum $g \in \mathcal{C}$, the inverse problem $Rf = g$ is obviously Krylov solvable, as $\overline{\mathcal{K}(R,g)} = \ell^2(\mathbb{Z})$, owing to the cyclicity of $g$. Yet, by density, there exists a sequence $(g_n)_{n \in \mathbb{N}}$ in $\mathcal{N}$ with $g_n \to g$, and each perturbed inverse problem $Rf_n = g_n$ is not Krylov solvable. Krylov solvability is absent along the perturbations and only emerges in the limit, still with the solution approximation $f_n \to f$.

*Example 5.4* With respect to the Hilbert space orthogonal sum $\mathcal{H} = \mathcal{H}_1 \oplus \mathcal{H}_2$, let $A = A_1 \oplus A_2$ with $A^{(j)} \in \mathcal{B}(\mathcal{H}_j)$, and $g^{(j)} \in \operatorname{ran} A^{(j)}$, $j \in \{1,2\}$, such that the problem $A^{(1)} f^{(1)} = g^{(1)}$ is Krylov solvable in $\mathcal{H}_1$, with Krylov solution $f^{(1)}$ (for instance, $\mathcal{H}_1 = L^2[0,1]$, $A^{(1)} = V$, the Volterra operator from Example 2.5, $g^{(1)} = x$, $f^{(1)} = \mathbf{1}$), and the problem $A^{(2)} f^{(2)} = g^{(2)}$ is not Krylov solvable in $\mathcal{H}_2$ (for instance, $\mathcal{H}_2 = \ell^2(\mathbb{N})$, $A^{(2)} = R$, the right shift from Example 2.1, $g^{(2)} = e_2$, $f^{(2)} = e_1$).

(i) (Lack of Krylov solvability persists in the limit.) The inverse problems $Af_n = g_n$, $n \in \mathbb{N}$, with $g_n := \left(\frac{1}{n}g^{(1)}\right) \oplus g^{(2)}$, are all non-Krylov solvable, with solution(s) $f_n = \left(\frac{1}{n}f^{(1)}\right) \oplus f^{(2)}$. In the limit, $g_n \to g := 0 \oplus g^{(2)}$ in $\mathcal{H}$, whence also $f_n \to 0 \oplus f^{(2)} =: f$. The inverse problem $Af = g$ has solution $f$ (modulo $\ker A^{(1)} \oplus \{0\}$), but $f$ is not a Krylov solution.

(ii) (Krylov solvability emerges in the limit.) The inverse problems $Af_n = g_n$, $n \in \mathbb{N}$, with $g_n := g^{(1)} \oplus \left(\frac{1}{n}g^{(2)}\right)$, are all non-Krylov solvable, with solution(s) $f_n = f^{(1)} \oplus \left(\frac{1}{n}f^{(2)}\right)$. In the limit, $g_n \to g := g^{(1)} \oplus 0$ in $\mathcal{H}$, whence also $f_n \to f^{(1)} \oplus 0 =: f$. The inverse problem $Af = g$ has solution $f$ (modulo $\{0\} \oplus \ker A^{(2)}$), and $f$ is a Krylov solution.

### 5.3.3 Simultaneous perturbations of operator and data

*Example 5.5* Same setting as in Example 5.4.

(i) (Lack of Krylov solvability persists in the limit.) The inverse problems $A_n f_n = g_n$, $n \in \mathbb{N}$, with $A_n := \left(\frac{1}{n}A^{(1)}\right) \oplus A^{(2)}$ and $g_n := \left(\frac{1}{n}g^{(1)}\right) \oplus g^{(2)}$, are all non-Krylov solvable, with solutions $f_n = f^{(1)} \oplus f^{(2)}$. In the limit, $A_n \to A := \mathbb{O} \oplus A^{(2)}$ in operator norm and $g_n \to g := 0 \oplus g^{(2)}$ in $\mathcal{H}$. The inverse problem $Af = g$ has solution $f := 0 \oplus f^{(2)}$ (modulo $\ker A^{(1)} \oplus \{0\}$), which is not a Krylov solution. Moreover in general $f_n$ does not converge to $f$.

(ii) (Krylov solvability emerges in the limit.) The inverse problems $A_n f_n = g_n$, $n \in \mathbb{N}$, with $A_n := A^{(1)} \oplus \left(\frac{1}{n}A^{(2)}\right)$ and $g_n := g^{(1)} \oplus \left(\frac{1}{n}g^{(2)}\right)$, are all non-Krylov solvable, with solutions $f_n = f^{(1)} \oplus f^{(2)}$. In the limit, $A_n \to A := A^{(1)} \oplus \mathbb{O}$ in operator norm and $g_n \to g := g^{(1)} \oplus 0$ in $\mathcal{H}$. The inverse problem $Af = g$ has solution $f := f^{(1)} \oplus 0$ (modulo $\{0\} \oplus \ker A^{(2)}$), which is a Krylov solution. Moreover in general $f_n$ does not converge to $f$.

## 5.4 Krylov solvability along perturbations of $\mathscr{K}$-class

In Section 2.5 operators of $\mathscr{K}$-class were introduced that in retrospect display relevant behaviour as far as Krylov solvability along perturbations is concerned. In this Section we elaborate further on that class, in view of the scheme of general questions presented in Section 5.2.

Let us recall that a linear operator $A$ acting on a complex Hilbert space $\mathcal{H}$ is in the $\mathscr{K}$-class when $A$ is everywhere defined and bounded, and there exists a bounded open $\mathcal{W} \subset \mathbb{C}$ containing the spectrum $\sigma(A)$ and such that $0 \notin \overline{\mathcal{W}}$ and $\mathbb{C} \setminus \mathcal{W}$ is connected. In particular, a $\mathscr{K}$-class operator has everywhere defined bounded inverse.

We have this result.

**Theorem 5.1** *Let $A$ be a $\mathscr{K}$-class operator on a complex Hilbert space $\mathcal{H}$.*

(i) *For every $g \in \mathcal{H}$ the inverse problem $Af = g$ is Krylov solvable, with unique solution $f = A^{-1}g$.*
(ii) *The $\mathscr{K}$-class is open in $\mathcal{B}(\mathcal{H})$. In particular, there is $\varepsilon_A > 0$ such that for any other operator $A' \in \mathcal{B}(\mathcal{H})$ with $\|A' - A\|_{\mathrm{op}} < \varepsilon_A$ the inverse problem $A'f' = g$ has a unique solution, $f' = {A'}^{-1}g$, which is also a Krylov solution.*
(iii) *When in addition (and without loss of generality) $\|A - A'\|_{\mathrm{op}} \leqslant (2\|A^{-1}\|_{\mathrm{op}})^{-1}$, then $f$ and $f'$ from (i) and (ii) satisfy*

$$\|f - f'\| \leqslant 2\|g\| \|A^{-1}\|_{\mathrm{op}}^2 \|A' - A\|_{\mathrm{op}}.$$

Theorem 5.1 addresses questions of type II from the general scheme of Section 5.2: it provides a framework where Krylov solvability is preserved under perturbations of the linear operator inducing the inverse problem. Indeed, an obvious consequence of Theorem 5.1 is: if a sequence $(A_n)_{n \in \mathbb{N}}$ in $\mathcal{B}(\mathcal{H})$ satisfies $A_n \to A$ in operator norm for some $\mathscr{K}$-class operator $A$, then eventually in $n$ the $A_n$'s are all of $\mathscr{K}$-class, the associated inverse problems $A_n f_n = g$ are Krylov solvable with unique solution $f_n = A_n^{-1}g$, and moreover $f_n \to f$ in $\mathcal{H}$, where $f = A^{-1}g$ is the unique and Krylov solution to $Af = g$.

*Proof (Proof of Theorem 5.1)* Part (i) was proved in Corollary 2.2.

Concerning (ii), we use the fact that $\sigma(A)$ is an upper semi-continuous function of $A \in \mathcal{B}(\mathcal{H})$ (see, e.g., [55, Problem 103] and [67, Theorem IV.3.1 and Remark IV.3.3]), meaning that for every bounded open set $\Omega \subset \mathbb{C}$ with $\sigma(A) \subset \Omega$ there exists $\varepsilon_A > 0$ such that if $A' \in \mathcal{B}(\mathcal{H})$ with $\|A' - A\|_{\mathrm{op}} < \varepsilon_A$, then $\sigma(A') \subset \Omega$. Applying this to $\Omega = \mathcal{W}$ we deduce that any such $A'$ is again of $\mathscr{K}$-class. The remaining part of the thesis then follows from (i).

As for (iii), clearly

$$\|f - f'\| \leq \|g\| \|A^{-1} - A'^{-1}\|_{\mathrm{op}} = \|g\| \|A'^{-1}(A - A')A^{-1}\|_{\mathrm{op}}$$
$$\leq \|g\| \|A'^{-1}\|_{\mathrm{op}} \|A^{-1}\|_{\mathrm{op}} \|A - A'\|_{\mathrm{op}},$$

and

$$A'^{-1} = A^{-1} \sum_{n=0}^{\infty} \left((A - A')A^{-1}\right)^n \quad \text{when } \|A - A'\|_{\mathrm{op}} < \|A^{-1}\|_{\mathrm{op}}^{-1},$$

whence, when additionally $\|A - A'\|_{\mathrm{op}} \leq \left(2\|A^{-1}\|_{\mathrm{op}}\right)^{-1}$, $\|A'^{-1}\|_{\mathrm{op}} \leq 2\|A^{-1}\|_{\mathrm{op}}$. Plugging the latter inequality into the above estimate for $\|f - f'\|$ yields the conclusion. $\square$

Theorem 5.1 only scratches the surface of expectedly relevant features of $\mathscr{K}$-class operators, in view of the study of perturbations preserving Krylov solvability (see questions of type II in Section 5.2).

That the issue is non-trivial, however, is demonstrated by important difficulties that one soon encounters when trying to extend the scope of Theorem 5.1. Let us discuss here one point in particular: in the same spirit of questions of type II, it is natural to inquire whether $\mathscr{K}$-class operators allow to establish Krylov solvability in the limit when the perturbation is removed.

To begin with, if a sequence $(A_n)_{n \in \mathbb{N}}$ of $\mathscr{K}$-class operators on $\mathcal{H}$ converges in operator norm, the limit $A$ fails in general to be of $\mathscr{K}$-class. Indeed:

**Lemma 5.1** *The $\mathscr{K}$-class is not closed in $\mathcal{B}(\mathcal{H})$.*

***Proof*** It suffices to consider a positive, compact operator $A : \mathcal{H} \to \mathcal{H}$, with zero in its spectrum, and its perturbations $A_n := A + n^{-1}\mathbb{1}$, $n \in \mathbb{N}$. Then each $A_n$ is of $\mathscr{K}$-class and $\|A_n - A\|_{\mathrm{op}} \to 0$ as $n \to \infty$, but by construction $A$ is not of $\mathscr{K}$-class. $\square$

One might be misled to believe that the general mechanism for such failure is the appearance of zero in the spectrum of the limit operator $A$, and that therefore a uniform separation of $\sigma(A_n)$ from zero as $n \to \infty$ would produce a limit $A$ still in the $\mathscr{K}$-class. To show that this is not the case either, let us work out the following example.

*Example 5.6* Let $A$ and $A_n$, $n \in \mathbb{N}$, be the operators on the Hilbert space $L^2[0,1]$ defined by

$$(Af)(x) := e^{2\pi i x} f(x),$$
$$(A_n f)(x) := \begin{cases} e^{2\pi i x} f(x), & \text{if } x \in (\frac{1}{2\pi n}, 1], \\ (1 + \frac{1}{n}) e^{2\pi i x} f(x), & \text{if } x \in [0, \frac{1}{2\pi n}], \end{cases}$$

for $f \in L^2[0,1]$ and a.e. $x \in [0,1]$. Clearly,

$$\|A\|_{\mathrm{op}} = 1, \qquad \|A_n\|_{\mathrm{op}} = 1 + \frac{1}{n}, \qquad \|A - A_n\|_{\mathrm{op}} = \frac{1}{n} \|A\|_{\mathrm{op}} \xrightarrow{n \to \infty} 0.$$

Moreover, each $A_n$ is a $\mathscr{K}$-class operator: its spectrum $\sigma(A_n)$ covers the unit circle except for a 'lid' arc, corresponding to the angle $x \in [0, \frac{1}{n}]$, which lies on the circle of larger radius $1+\frac{1}{n}$, therefore it is possible to include $\sigma(A_n)$ into a suitable bounded open $\mathcal{W}$ separated from zero and with connected complement in $\mathbb{C}$. In the limit $n \to \infty$ the 'lid' closes the unit circle: $\sigma(A)$ is indeed the whole unit circle. Thus, even if the limit operator $A$ satisfies $0 \notin \sigma(A)$, $A$ fails to belong to the $\mathscr{K}$-class. In addition, not only the $\mathscr{K}$-class condition is lost in the limit, but so too is the Krylov solvability. Indeed, by means of the Hilbert space isomorphism

$$L^2[0,1] \xrightarrow{\cong} \ell^2(\mathbb{Z}), \qquad e^{2\pi i k x} \longmapsto e_k$$

(namely with respect to the orthonormal bases $(e^{2\pi i k x})_{k \in \mathbb{Z}}$ and $(e_k)_{k \in \mathbb{Z}}$), $A$ is unitarily equivalent to the right shift operator on $\ell^2(\mathbb{Z})$: thus, any choice $g \cong e_k$ for some $k \in \mathbb{Z}$ produces a non-Krylov solvable inverse problem $Af = g$.

In conclusion, the $\mathscr{K}$-class proves to be an informative sub-class of operators that is very robust and preserves Krylov solvability under perturbations of an unperturbed $\mathscr{K}$-class inverse problem (Theorem 5.1), but on the contrary is very fragile when from a sequence of approximating inverse problems of $\mathscr{K}$-class one wants to extract information on the Krylov solvability of the limit problem (Lemma 5.1, Example 5.6).

## 5.5 Weak gap-metric for weakly closed parts of the unit ball

Let us introduce now a convenient indicator of vicinity of closed subspaces of a given Hilbert space, which turns out to possess convenient properties when comparing (closures of) Krylov subspaces, and to provide a rigorous language to express and control limits of the form $\overline{\mathcal{K}(A_n, g_n)} \to \overline{\mathcal{K}(A, g)}$. Even though such an indicator is not optimal, in that it lacks other desired properties that would make it fully informative, we discuss it in depth here as a first attempt towards an efficient measurement of vicinity and convergence of Krylov subspaces under perturbations.

One natural motivation is provided by the *failure* of describing the intuitive convergence $\mathcal{K}_N(A, g) \xrightarrow{N \to \infty} \overline{\mathcal{K}(A, g)}$, where

$$\mathcal{K}_N(A, g) := \mathrm{span}\{g, Ag, \ldots, A^{N-1}g\}, \qquad N \in \mathbb{N} \tag{5.1}$$

is the $N$-th order Krylov subspace, by means of the ordinary **gap metric** between closed subspaces of the underlying Hilbert space.

Let us recall (see, e.g., [67, Chapt. 4, §2]) that given a Hilbert space $\mathcal{H}$ and two closed subspaces $U, V \subset \mathcal{H}$, the **gap** and the **gap distance** between them are, respectively, the quantities

$$\widehat{\delta}(U, V) := \max\{\delta(U, V), \delta(V, U)\} \tag{5.2}$$

and

$$\widehat{d}(U, V) := \max\{d(U, V), d(V, U)\}, \tag{5.3}$$

where

$$\delta(U,V) := \sup_{\substack{u \in U \\ \|u\|=1}} \inf_{v \in V} \|u-v\|,$$
$$d(U,V) := \sup_{\substack{u \in U \\ \|u\|=1}} \inf_{\substack{v \in V \\ \|v\|=1}} \|u-v\|, \qquad (5.4)$$

and with the tacit definitions $\delta(\{0\},V) := 0$, $d(\{0\},V) := 0$, $d(U,\{0\}) := 2$ for $U \neq \{0\}$, when one of the two entries is the trivial subspace.

The short-hands $B_\mathcal{H}$ for the closed unit ball of $\mathcal{H}$, $S_\mathcal{H}$ for the closed unit sphere, $B_U := U \cap B_\mathcal{H}$, $S_U := U \cap S_\mathcal{H}$, and $\mathrm{dist}(x,C)$ for the norm distance of a point $x \in \mathcal{H}$ from the closed subset $C \subset \mathcal{H}$ will be used throughout.

Thus,

$$\delta(U,V) = \sup_{u \in S_U} \mathrm{dist}(u,V), \qquad d(U,V) = \sup_{u \in S_U} \mathrm{dist}(u,S_V). \qquad (5.5)$$

As a matter of fact, on the set of all closed subspaces of $\mathcal{H}$ both $\widehat{\delta}$ and $\widehat{d}$ are two equivalent metrics, with

$$\widehat{\delta}(U,V) \leqslant \widehat{d}(U,V) \leqslant 2\widehat{\delta}(U,V), \qquad (5.6)$$

and the resulting metric space is complete [52].

The construction recalled here is for the Hilbert space setting and was introduced first in [70] as '*opening*' between (closed) subspaces (i.e., the operator norm distance between their orthogonal projections). It also applies to the more general case when $\mathcal{H}$ is a Banach space, a generalisation originally discussed in in [71], and [3, §34] (except that in the non-Hilbert case the gap $\widehat{\delta}$ is not a metric, even though it still satisfies (5.6) and hence induces the same topology as the metric $\widehat{d}$). Let us also recall that by linearity the closedness of the above subspaces $U$ and $V$ can be equivalently formulated in the norm or in the weak topology of $\mathcal{H}$.

Now, given $A \in \mathcal{B}(\mathcal{H})$ and $g \in \mathcal{H}$, for the closed subspaces $\mathcal{K} := \overline{\mathcal{K}(A,g)}$ and $\mathcal{K}_N := \mathcal{K}_N(A,g)$, $N \in \mathbb{N}$, of $\mathcal{H}$ one obviously has $\mathcal{K}_N \subset \mathcal{K}$ and hence $\delta(\mathcal{K}_N,\mathcal{K}) = 0$; on the other hand, if $\dim \mathcal{K} = \infty$, one can find for every $N$ a vector $u \in S_\mathcal{K}$ such that $u \perp \mathcal{K}_N$, thus with $\mathrm{dist}(u,\mathcal{K}_N) = 1$, and hence $\delta(\mathcal{K},\mathcal{K}_N) \geqslant 1$. This shows that $\widehat{d}(\mathcal{K}_N,\mathcal{K}) \geqslant \widehat{\delta}(\mathcal{K}_N,\mathcal{K}) \geqslant 1$, therefore the sequence $(\mathcal{K}_N)_{N \in \mathbb{N}}$ fails to converge to $\mathcal{K}$ in the $\widehat{d}$-metric. In this respect, the $\widehat{d}$-metric is certainly not a convenient tool to monitor the vicinity of Krylov subspaces, for it cannot accommodate the most intuitive convergence $\mathcal{K}_N \to \mathcal{K}$.

With this observation in mind, it is natural to weaken the ordinary gap distance $\widehat{d}$-metric so as to encompass a larger class of limits. To do so, we exploit the fact (see, e.g., [13, Theorem 3.29]) that in any separable Hilbert space $\mathcal{H}$ the norm-closed unit ball $B_\mathcal{H}$ is metrisable in the Hilbert space weak topology. More precisely, there exists a norm $\|\cdot\|_w$ on $\mathcal{H}$ (and hence a metric $\rho_w(x,y) := \|x-y\|_w$) such that $\|x\|_w \leqslant \|x\|$ and whose metric topology restricted to $B_\mathcal{H}$ is precisely the Hilbert space weak topology. For concreteness one may define

$$\|x\|_w := \sum_{n=1}^{\infty} \frac{1}{2^n} |\langle \xi_n, x \rangle|$$

for a dense countable collection $(\xi_n)_{n\in\mathbb{N}}$ in $B_{\mathcal{H}}$ which identifies the norm $\|\cdot\|_w$. On the other hand, since a Hilbert space is reflexive, $B_{\mathcal{H}}$ is compact in the weak topology (see, e.g., [13, Theorem 3.16]), and hence in the $\rho_w$-metric. Being $(B_{\mathcal{H}}, \rho_w)$ a metric space, its compactness is equivalent to the property of being simultaneously complete and totally bounded (see, e.g., [80, Theorem 45.1]). In conclusion, the metric space $(B_{\mathcal{H}}, \rho_w)$ is compact and complete, and its metric topology is the Hilbert space weak topology (restricted to $B_{\mathcal{H}}$). In fact, the construction that follows, including Theorem 5.2 below, is applicable to the more general case where $\mathcal{H}$ is a reflexive Banach space with separable dual: indeed, the same properties above for $(B_{\mathcal{H}}, \rho_w)$ hold.

In $(B_{\mathcal{H}}, \rho_w)$ we denote the **relative weakly open balls** (namely the $\rho_w$-*open* balls of $\mathcal{H}$ intersected with $B_{\mathcal{H}}$) as

$$\mathfrak{B}_w(x_0, \varepsilon) := \{x \in B_{\mathcal{H}} \mid \|x - x_0\|_w < \varepsilon\} \tag{5.7}$$

for given $x_0 \in B_{\mathcal{H}}$ and $\varepsilon > 0$. Observe that any such open ball $\mathfrak{B}_w(x_0, \varepsilon)$ always contains points of the unit sphere (not all, if $\varepsilon$ is small enough); thus, at fixed $x_0 \in B_{\mathcal{H}}$, and along a sequence of radii $\varepsilon_n \downarrow 0$, one can select a sequence $(y_n)_{n\in\mathbb{N}}$ with $\|y_n\| = 1$ and $\|y_n - x_0\|_w < \varepsilon_n$, whence the conclusion $y_n \xrightarrow{\rho_w} x_0$, which reproduces, in the metric space language, the topological statement that $y_n \rightharpoonup x_0$, i.e., that the unit ball is the weak closure of the unit sphere.

Based on the weak (and metric) topology $B_{\mathcal{H}}$ it is natural to weaken the gap distance $\widehat{d}$ considered before, as we shall do in a moment, except that dealing now with weak limits instead of norm limits one has to expect possible "discontinuous jumps", say, in the form of sudden expansions or contractions of the limit object as compared to its approximants (in the same spirit of taking the closure of the unit sphere $S_{\mathcal{H}}$: the norm-closure gives again $S_{\mathcal{H}}$, the weak closure gives the whole $B_{\mathcal{H}}$). For this reason we set up the new notion of **weak gap metric** in the more general class

$$\mathcal{C}_w(\mathcal{H}) := \{\text{non-empty and weakly closed subsets of } B_{\mathcal{H}}\}, \tag{5.8}$$

instead of the subclass of unit balls of closed subspaces of $\mathcal{H}$.

For $U, V \in \mathcal{C}_w(\mathcal{H})$ let us then set

$$\begin{aligned} d_w(U, V) &:= \sup_{u \in U} \inf_{v \in V} \|u - v\|_w, \\ \widehat{d_w}(U, V) &:= \max\{d_w(U, V), d_w(V, U)\}. \end{aligned} \tag{5.9}$$

We shall now establish the fundamental properties of the map $\widehat{d_w}$. They are summarised as follows.

**Theorem 5.2** *Let $\mathcal{H}$ be a separable Hilbert space.*

(i) $\widehat{d}_w$ is a metric on $\mathcal{C}_w(\mathcal{H})$.
(ii) The metric space $(\mathcal{C}_w(\mathcal{H}), \widehat{d}_w)$ is complete.
(iii) If $\widehat{d}_w(U_n, U) \xrightarrow{n \to \infty} 0$ for an element $U$ and a sequence $(U_n)_{n \in \mathbb{N}}$ in $\mathcal{C}_w(\mathcal{H})$, then

$$U = \{u \in B_\mathcal{H} \,|\, u_n \rightharpoonup u \text{ for a sequence } (u_n)_{n \in \mathbb{N}} \text{ with } u_n \in U_n\}. \quad (5.10)$$

(iv) The metric space $(\mathcal{C}_w(\mathcal{H}), \widehat{d}_w)$ is compact.
(v) If $\widehat{d}_w(U_n, U) \xrightarrow{n \to \infty} 0$ for an element $U$ and a sequence $(U_n)_{n \in \mathbb{N}}$ in $\mathcal{C}_w(\mathcal{H})$, then $\widehat{d}_w(f(U_n), f(U)) \xrightarrow{n \to \infty} 0$ for any weakly closed and weakly continuous map $f : \mathcal{H} \to \mathcal{H}$ such that $f(B_\mathcal{H}) \subset B_\mathcal{H}$.

We shall also write $U_n \xrightarrow{\widehat{d}_w} U$ as an alternative to $\widehat{d}_w(U_n, U) \to 0$.

*Remark 5.1* The completeness and the compactness result of Theorem 5.2 are in a sense folk knowledge in the context of the **Hausdorff distance**. In fact, the gap distance $\widehat{d}(U, V)$ introduced in (5.3)-(5.4) is, apart from zero-sets, the Hausdorff distance between $U$ and $V$ as subsets of the metric (normed) space $(\mathcal{H}, \|\cdot\|)$, and our modified weak gap distance $\widehat{d}_w(U, V)$ defined in (5.9) between elements of $\mathcal{C}_w(\mathcal{H})$ is the Hausdorff distance between sets in the metric space $(B_H, \rho_w)$. The completeness and the compactness of $(B_H, \rho_w)$ then lift, separately, to the completeness and compactness of $(\mathcal{C}_w(\mathcal{H}), \widehat{d}_w)$ – they are actually equivalent (see, e.g., [109, Theorem 5.38]). We present here both the completeness and the compactness result, and their proofs, in detail in order to expose reasonings, tailored on the weak topology setting, which will be repeatedly used in the proof of the various statements of the following Sections. Also, having the proof of completeness of $(\mathcal{C}_w(\mathcal{H}), \widehat{d}_w)$ fully laid down is of further help in understanding the *failure* of completeness of the metric space $(\mathscr{S}(\mathcal{H}), \widehat{d}_w)$ that we will consider in the next Section, for applications to Krylov subspaces.

There are additional technical properties of $d_w$ and $\widehat{d}_w$ that are worth being singled out. Let us collect them in the following lemmas.

**Lemma 5.2** *Let $U, V, Z \in \mathcal{C}_w(\mathcal{H})$ for some separable Hilbert space $\mathcal{H}$. Then*

$$d_w(U, V) = 0 \Leftrightarrow U \subset V, \quad (5.11)$$
$$\widehat{d}_w(U, V) = 0 \Leftrightarrow U = V, \quad (5.12)$$
$$d_w(U, Z) \leqslant d_w(U, V) + d_w(V, Z), \quad (5.13)$$
$$\widehat{d}_w(U, Z) \leqslant \widehat{d}_w(U, V) + \widehat{d}_w(V, Z). \quad (5.14)$$

**Lemma 5.3** *Given a separable Hilbert space $\mathcal{H}$ and a collection $(U_n)_{n \in \mathbb{N}}$ in $\mathcal{C}_w(\mathcal{H})$, the set*

$$\mathcal{U} := \{x \in B_\mathcal{H} \,|\, u_n \rightharpoonup x \text{ for a sequence } (u_n)_{n \in \mathbb{N}} \text{ with } u_n \in U_n\} \quad (5.15)$$

*is closed in the weak topology of $\mathcal{H}$.*

Given $U \in \mathcal{C}_w(\mathcal{H})$ we define its **weakly open $\varepsilon$-expansion** in $B_\mathcal{H}$ as

$$U(\varepsilon) := \bigcup_{u \in U} \mathfrak{B}_w(u, \varepsilon). \tag{5.16}$$

Observe that $U(\varepsilon)$ is a weakly open subset of $B_\mathcal{H}$.

**Lemma 5.4** *Let $U, V \in \mathcal{C}_w(\mathcal{H})$ for some separable Hilbert space $\mathcal{H}$ and let $\varepsilon > 0$. Then:*

*(i) $d_w(U,V) < \varepsilon \Rightarrow V \cap \mathfrak{B}_w(u, \varepsilon) \neq \emptyset \; \forall u \in U$;*
*(ii) $d_w(U,V) < \varepsilon \Leftrightarrow U \subset V(\varepsilon)$;*
*(iii) $\left.\begin{array}{l} U \subset V(\varepsilon) \\ U \cap \mathfrak{B}_w(v, \varepsilon) \neq \emptyset \; \forall v \in V \end{array}\right\} \Rightarrow \widehat{d}_w(V, U) < \varepsilon.$*

The remaining part of this Section is devoted to proving the above statements.

***Proof (Proof of Lemma 5.2)*** For (5.11), the inclusion $U \subset V$ implies $\inf_{v \in V} \|u - v\|_w = 0$ for every $u \in U$, whence $d_w(U,V) = 0$; conversely, if $0 = d_w(U,V) = \sup_{u \in U} \inf_{v \in V} \|u-v\|_w$, then $\inf_{v \in V} \|u-v\|_w = 0$ for every $u \in U$, whence the fact, by weak closedness of $V$, that any such $u$ belongs also to $V$. As for the property (5.12), it follows from (5.11) exploiting separately both inclusions $U \subset V$ and $U \supset V$. Last, let us prove the triangular inequalities (5.13)-(5.14). Let $u_0 \in U$: then, owing to the weak compactness of $V$ (as a closed subset of the compact metric space $(B_\mathcal{H}, \rho_w)$), $\inf_{v \in V} \|u_0 - v\|_w = \|u_0 - v_0\|_w$ for some $v_0 \in V$, whence

$$\|u_0 - v_0\|_w = \inf_{v \in V} \|u_0 - v\|_w \leqslant \sup_{u \in U} \inf_{v \in V} \|u - v\|_w = d_w(U,V) \leqslant \widehat{d}_w(U,V).$$

As a consequence,

$$\begin{aligned}\inf_{z \in Z} \|u_0 - z\|_w &\leqslant \|u_0 - v_0\|_w + \inf_{z \in Z} \|v_0 - z\|_w \\ &\leqslant d(U,V) + d(U,Z) \\ &\leqslant \widehat{d}_w(U,V) + \widehat{d}_w(V,Z),\end{aligned}$$

having used the triangular inequality of the $\|\cdot\|_w$-norm in the first inequality. By the arbitrariness of $u_0 \in U$, thus taking the supremum over all such $u_0$'s,

$$\begin{aligned} d_w(U,Z) &\leqslant d_w(U,V) + d_w(V,Z), \\ d_w(U,Z) &\leqslant \widehat{d}_w(U,V) + \widehat{d}_w(V,Z).\end{aligned}$$

With the first inequality above we proved (5.13). Next, let us combine the second inequality above with the corresponding bound for $d_w(V,U)$, which is established in a similar manner: let now $z_0 \in Z$, and again by weak compactness there exists $v_0 \in V$ with $\|v_0 - z_0\|_w = \inf_{v \in V} \|v - z_0\|_w \leqslant d_w(Z,V) \leqslant \widehat{d}_w(V,Z)$, and also $\inf_{u \in U} \|u - v_0\|_w \leqslant d(V,U) = \widehat{d}(U,V)$, whence

$$\inf_{u \in U} \|u - z_0\|_w \leqslant \inf_{u \in U} \|u - v_0\|_w + \|v_0 - z_0\|_w \leqslant \widehat{d}_w(U,V) + \widehat{d}_w(V,Z).$$

Taking the supremum over all $z_0 \in Z$ yields

$$d_w(Z,U) \leq \widehat{d}_w(U,V) + \widehat{d}_w(V,Z).$$

Combining the above estimates for $d_w(U,Z)$ and $d_w(Z,U)$ yields the conclusion. $\square$

*Proof (Proof of Theorem 5.2(i))* It follows directly from (5.12)-(5.14) of Lemma 5.2. $\square$

*Proof (Proof of Lemma 5.3)* Let $x \in \overline{\mathcal{U}}^w$, the closure of $\mathcal{U}$ in the weak topology, and let us construct a sequence $(u_n)_{n\in\mathbb{N}}$ with $u_n \in U_n$ and $u_n \rightharpoonup x$, thereby showing that $x \in \mathcal{U}$.

By assumption $\exists x_2 \in \mathcal{U}$ with $\|x - x_2\|_w \leq \frac{1}{4}$ and $\|x_2 - u_n^{(2)}\|_w \xrightarrow{n\to\infty} 0$ for a sequence $(u_n^{(2)})_{n\in\mathbb{N}}$ with $u_n^{(2)} \in U_n$. In particular, there is $N_2 \in \mathbb{N}$ with $\|x_2 - u_n^{(2)}\|_w \leq \frac{1}{4}$ $\forall n \geq N_2$. For the integer $N_3$ with $N_3 > N_2 + 1$ to be fixed in a moment, set

$$u_n := u_n^{(2)}, \qquad n \in \{N_2,\ldots,N_3-1\}.$$

By construction, $\|x - u_n\|_w \leq \|x - x_n\|_w + \|x_n - u_n\|_w \leq \frac{1}{2}$ $\forall n \in \{N_2,\ldots,N_3-1\}$.

Next, for each integer $k \geq 3$ one identifies recursively a sequence $(N_k)_{k=3}^\infty$ in $\mathbb{N}$ with $N_k > N_{k-1} + 1$, and vectors $u_n \in U_n$ for $n \in \{N_k,\ldots,N_{k+1}-1\}$ as follows. By assumption $\exists x_k \in \mathcal{U}$ with $\|x - x_k\|_w \leq \frac{1}{2k}$ and $\|x_k - u_n^{(k)}\|_w \xrightarrow{n\to\infty} 0$ for a sequence $(u_n^{(k)})_{n\in\mathbb{N}}$ with $u_n^{(k)} \in U_n$. In particular, it is always possible to find $N_k \in \mathbb{N}$ with $N_k > N_{k-1} + 1$ such that $\|x_k - u_n^{(k)}\|_w \leq \frac{1}{2k}$ $\forall n \geq N_k$. For the integer $N_{k+1} > N_k + 1$ set

$$u_n := u_n^{(k)}, \qquad n \in \{N_k,\ldots,N_{k+1}-1\}.$$

By construction, $\|x - u_n\|_w \leq \|x - x_n\|_w + \|x_n - u_n\|_w \leq \frac{1}{k}$ $\forall n \in \{N_k,\ldots,N_{k+1}-1\}$.

This yields a sequence $(u_n)_{n\in\mathbb{N}}$ (having added, if needed, finitely many irrelevant vectors $u_1,\ldots,u_{N_2-1}$) with $u_n \in U_n$ and such that, for any integer $k \geq 2$,

$$\|x - u_n\|_w \leq \frac{1}{k} \quad \forall n \geq N_k.$$

Hence $\|x - u_n\|_w \xrightarrow{n\to\infty} 0$, thus $x \in \mathcal{U}$. $\square$

*Proof (Proof of Theorem 5.2(ii) and (iii))* One needs to show that given $(U_n)_{n\in\mathbb{N}}$, Cauchy sequence in $\mathcal{C}_w(\mathcal{H})$, there exists $U \in \mathcal{C}_w(\mathcal{H})$ with $\widehat{d}_w(U_n,U) \xrightarrow{n\to\infty} 0$, and that $U$ has precisely the form (5.10). Moreover, as $(\mathcal{C}_w(\mathcal{H}), \widehat{d}_w)$ is a metric space, it suffices to establish the above statement for one subsequence of $(U_n)_{n\in\mathbb{N}}$.

By the Cauchy property, $\widehat{d}_w(U_n,U_m) \xrightarrow{n,m\to\infty} 0$. Up to extracting a subsequence, henceforth denoted again with $(U_n)_{n\in\mathbb{N}}$, one can further assume that

$$\widehat{d}_w(U_n,U_m) \leq \frac{1}{2^n} \quad \forall m \geq n.$$

We shall establish the $\widehat{d}_w$-convergence of such (sub-)sequence.

First of all, fixing any $n \in \mathbb{N}$ and any $u_n \in U_n$, we construct a sequence with

- an (irrelevant) choice of vectors $u_1, \ldots, u_{n-1}$ in the first $n-1$ positions, such that $u_1 \in U_1, \ldots, u_{n-1} \in U_{n-1}$,
- precisely the considered vector $u_n$ in position $n$,
- and an infinite collection $u_{n+1}, u_{n+2}, u_{n+3}, \ldots$ determined recursively so that, given $u_k \in U_k$ ($k \geqslant n$), the next $u_{k+1}$ is that element of $U_{k+1}$ satisfying $\inf_{v \in U_{k+1}} \|u_k - v\|_w = \|u_k - u_{k+1}\|_w$ – a choice that is always possible, owing to the weak compactness of $U_{k+1}$ as a closed subset of the compact metric space $(B_{\mathcal{H}}, \rho_w)$.

Let us refer to such $(u_k)_{k \in \mathbb{N}}$ as the sequence 'originating from the given $u_n$' (tacitly understanding that it is one representative of infinitely many sequences with the same property, owing to the irrelevant choice of the first $n-1$ vectors). When the originating vector need be indicated, we shall write $\big(u_k^{(u_n)}\big)_{k \in \mathbb{N}}$: thus, $u_n^{(u_n)} \equiv u_n$.

By construction, for any $k \geqslant n$,

$$\|u_k - u_{k+1}\|_w = \inf_{v \in U_{k+1}} \|u_k - v\|_w \leqslant \sup_{z \in U_k} \inf_{v \in U_{k+1}} \|z - v\|_w$$

$$= d_w(U_k, U_{k+1}) \leqslant \widehat{d}_w(U_k, U_{k+1}) \leqslant \frac{1}{2^k},$$

whence, for any $m > n$,

$$\|u_n - u_m\|_w \leqslant \sum_{k=n}^{m-1} \|u_k - u_{k+1}\|_w \leqslant \sum_{k=n}^{m-1} \frac{1}{2^k} \leqslant \frac{1}{2^{n-1}}.$$

This implies that the sequence $(u_k)_{k \in \mathbb{N}}$ originating from the considered $u_n \in U_n$ is a Cauchy sequence in $(B_{\mathcal{H}}, \rho_w)$ and we denote its weak limit as $u_\infty^{(u_n)} \in B_{\mathcal{H}}$. The same construction can be repeated for any $n \in \mathbb{N}$ and starting the sequence from any $u_n \in U_n$: the collection of all possible limit points is

$$U_\infty := \big\{ u \in B_{\mathcal{H}} \,|\, u = u_\infty^{(u_n)} \text{ for some } n \in \mathbb{N} \text{ and some 'starting' } u_n \in U_n \big\}.$$

Compare now the set $U_\infty$ with the set

$$\widetilde{U} := \{ u \in B_{\mathcal{H}} \,|\, u_n \rightharpoonup u \text{ for a sequence } (u_n)_{n \in \mathbb{N}} \text{ with } u_n \in U_n \}$$

$\widetilde{U}$ is weakly closed (Lemma 5.3), and obviously $U_\infty \subset \widetilde{U}$. We claim that

$$\widetilde{U} = \overline{U_\infty}^{\|\,\|_w} \qquad \text{(the weak closure of } U_\infty\text{)}.$$

For arbitrary $u \in \widetilde{U}$ and $\varepsilon > 0$ there is $n_\varepsilon \in \mathbb{N}$ and $u_{n_\varepsilon} \in U_{n_\varepsilon}$ with $\|u - u_{n_\varepsilon}\|_w \leqslant \varepsilon$. Non-restrictively, $n_\varepsilon \to \infty$ as $\varepsilon \downarrow 0$. For a sequence $\big(u_k^{(u_{n_\varepsilon})}\big)_{k \in \mathbb{N}}$ originating from $u_{n_\varepsilon}$ and for its weak limit $u_\infty^{(u_{n_\varepsilon})} \in U_\infty$, there is $k_\varepsilon \in \mathbb{N}$ with $k_\varepsilon > n_\varepsilon$ satisfying both $\|u_{n_\varepsilon} - u_{k_\varepsilon}^{(u_{n_\varepsilon})}\|_w \leqslant 2^{-(n_\varepsilon-1)}$ (because of the above property of the sequences originating from one element) and $\|u_{k_\varepsilon}^{(u_{n_\varepsilon})} - u_\infty^{(u_{n_\varepsilon})}\|_w \leqslant \varepsilon$ (because of the convergence $u_k^{(u_{n_\varepsilon})} \rightharpoonup$

$u_\infty^{(u_{n_\varepsilon})}$). Thus,

$$\|u - u_\infty^{(u_{n_\varepsilon})}\|_w \leqslant \|u - u_{n_\varepsilon}\|_w + \|u_{n_\varepsilon} - u_{k_\varepsilon}^{(u_{n_\varepsilon})}\|_w + \|u_{k_\varepsilon}^{(u_{n_\varepsilon})} - u_\infty^{(u_{n_\varepsilon})}\|_w$$
$$\leqslant 2^{-(n_\varepsilon - 1)} + 2\varepsilon.$$

Taking $\varepsilon \downarrow 0$ shows that $u$ indeed belongs to the weak closure of $U_\infty$.

It remains to prove that $\widehat{d}_w(U_n, \widetilde{U}) \xrightarrow{n \to \infty} 0$. Let us control $d_w(U_n, \widetilde{U})$ first. Pick $n \in \mathbb{N}$ and $u_n \in U_n$. For the sequence $(u_k)_{k \in \mathbb{N}}$ originating from $u_n$ (thus, $u_k \rightharpoonup u_\infty^{(u_n)}$) and for arbitrary $\varepsilon > 0$ there is $k_\varepsilon \in \mathbb{N}$ with $k_\varepsilon > n$ such that $\|u_\infty^{(u_n)} - u_{k_\varepsilon}\|_w \leqslant \varepsilon$ and $\|u_{k_\varepsilon} - u_n\|_w \leqslant 2^{-(n-1)}$. Thus,

$$\inf_{u \in \widetilde{U}} \|u - u_n\|_w \leqslant \|u_\infty^{(u_n)} - u_n\|_w \leqslant \|u_\infty^{(u_n)} - u_{k_\varepsilon}\|_w + \|u_{k_\varepsilon} - u_n\|_w \leqslant \varepsilon + \frac{1}{2^{n-1}},$$

whence also

$$d_w(U_n, \widetilde{U}) = \sup_{u_n \in U_n} \inf_{u \in \widetilde{U}} \|u - u_n\|_w \leqslant \varepsilon + \frac{1}{2^{n-1}}.$$

This implies that $\limsup_n d_w(U_n, \widetilde{U}) \leqslant \varepsilon$, and owing to the arbitrariness of $\varepsilon$, finally $d_w(U_n, \widetilde{U}) \xrightarrow{n \to \infty} 0$. The other limit $d_w(\widetilde{U}, U_n) \xrightarrow{n \to \infty} 0$ is established in much the same way, exploiting additionally the density of $U_\infty$ in $\widetilde{U}$. Pick arbitrary $n \in \mathbb{N}$, $\varepsilon > 0$, and $u \in \widetilde{U}$. We already argued, for the proof of the identity $\widetilde{U} = \overline{U_\infty}^{\|\cdot\|_w}$, that there is $n_\varepsilon \in \mathbb{N}$ with $n_\varepsilon \to \infty$ as $\varepsilon \downarrow 0$, and there is $u_{n_\varepsilon} \in U_{n_\varepsilon}$, such that

$$\|u - u_\infty^{(u_{n_\varepsilon})}\|_w \leqslant 2^{-(n_\varepsilon - 1)} + 2\varepsilon.$$

Non-restrictively, $n_\varepsilon > n$. In turn, as $u_k^{(u_{n_\varepsilon})} \rightharpoonup u_\infty^{(u_{n_\varepsilon})}$, there is $k_\varepsilon \in \mathbb{N}$ with $k_\varepsilon \geqslant n_\varepsilon > n$ such that $\|u_\infty^{(u_{n_\varepsilon})} - u_{k_\varepsilon}^{(n_\varepsilon)}\|_w \leqslant \varepsilon$ and $\|u_{k_\varepsilon}^{(n_\varepsilon)} - u_n^{(n_\varepsilon)}\| \leqslant 2^{-(n-1)}$. Therefore,

$$\inf_{v \in U_n} \|u - w\|_w \leqslant \|u - u_\infty^{(u_{n_\varepsilon})}\|_w + \inf_{v \in U_n} \|u_\infty^{(u_{n_\varepsilon})} - v\|_w$$
$$\leqslant \|u - u_\infty^{(u_{n_\varepsilon})}\|_w + \|u_\infty^{(u_{n_\varepsilon})} - u_{k_\varepsilon}^{(n_\varepsilon)}\|_w + \|u_{k_\varepsilon}^{(n_\varepsilon)} - u_n^{(n_\varepsilon)}\|$$
$$\leqslant \frac{1}{2^{n_\varepsilon - 1}} + 3\varepsilon + \frac{1}{2^{n-1}},$$

whence also

$$d_w(U, U_n) = \sup_{u \in U} \inf_{u_n \in U_n} \|u - u_n\|_w \leqslant \frac{1}{2^{n_\varepsilon - 1}} + 3\varepsilon + \frac{1}{2^{n-1}}.$$

As above, the limit $\varepsilon \downarrow 0$ and the arbitrariness of $n$ imply $d_w(U, U_n) \xrightarrow{n \to \infty} 0$ and finally $\widehat{d}_w(U_n, U) \xrightarrow{n \to \infty} 0$. □

***Proof (Proof of Lemma 5.4)*** (i) If, for contradiction, $V \cap \mathfrak{B}_w(u_0, \varepsilon) = \emptyset$ for some $u_0 \in U$, then the weak metric distance $(\rho_w)$ of $u_0$ from $V$ is at least $\varepsilon$, meaning that

$$d_w(U,V) \geqslant \inf_{v \in V} \|u_0 - v\|_w \geqslant \varepsilon.$$

(ii) Assume that $d_w(U,V) < \varepsilon$. On account of the weak compactness of $V$ (as a closed subset of the compact metric space $(B_{\mathcal{H}}, \rho_w)$), for any $u \in U$ there is $v_u \in V$ with

$$\|u - v_u\|_w = \inf_{v \in V} \|u - v\|_w \leqslant d_w(U,V) < \varepsilon,$$

meaning that $u \in \mathfrak{B}_w(v_u, \varepsilon)$. Thus, $U \subset V(\varepsilon)$. Conversely, if $U \subset V(\varepsilon)$, then any $u \in U$ belongs to a ball $\mathfrak{B}_w(v_u, \varepsilon)$ for some $v_u \in V$, whence

$$f(u) := \inf_{v \in V} \|u - v\|_w \leqslant \|u - v_u\|_w < \varepsilon.$$

The function $f : U \to \mathbb{R}$ is continuous on the weak compact set $U$, hence it attains its maximum at a point $u = u_0$ and

$$d_w(U,V) = \sup_{u \in U} f(u) = \inf_{v \in V} \|u_0 - v\|_w < \varepsilon.$$

(iii) As by assumption $U \subset V(\varepsilon)$, we know from (ii) that $d_w(U,V) < \varepsilon$. In addition, for any $v \in V$ it is assumed that $\mathfrak{B}_w(v, \varepsilon)$ is not disjoint from $U$, meaning that there is $u_v \in U$ with $\|u_v - v\|_w < \varepsilon$. Therefore,

$$g(v) := \inf_{u \in U} \|u - v\|_w < \varepsilon,$$

and from the continuity of $g : V \to \mathbb{R}$ on the weak compact $V$,

$$d_w(V,U) = \sup_{v \in V} g(v) = g(v_0) < \varepsilon,$$

where $v_0$ is some point of maximum for $g$. In conclusion, $\widehat{d}_w(V,U) < \varepsilon$. □

*Proof (Proof of Theorem 5.2(iv))* As the metric space $(\mathcal{C}_w(\mathcal{H}), \widehat{d}_w)$ is complete, compactness follows if one proves that for any $\varepsilon > 0$ the set $\mathcal{C}_w(\mathcal{H})$ can be covered by finitely many $\widehat{d}_w$-open balls of radius $\varepsilon$ (total boundedness and completeness indeed imply compactness for a metric space).

To this aim, let us observe first that, owing to the compactness of $(B_{\mathcal{H}}, \rho_w)$, for any $\varepsilon > 0$ we may cover it with finitely many open balls $\mathfrak{B}_w(x_1, \varepsilon), \ldots, \mathfrak{B}_w(x_M, \varepsilon)$ for some $x_1, \ldots x_M \in B_{\mathcal{H}}$ and $M \in \mathbb{N}$ all depending on $\varepsilon$. Each $\mathfrak{B}_w(x_n, \varepsilon)$ is the $\varepsilon$-expansion of the weakly closed set $\{x_n\}$, hence

$$Z(\varepsilon) = \bigcup_{x \in Z} \mathfrak{B}_w(x, \varepsilon) \qquad \forall Z \subset \mathcal{Z}_M := \{x_1, \ldots, x_M\}.$$

Let us now show that the finitely many $\widehat{d}_w$-open balls of the form

$$\{U \in \mathcal{C}_w(\mathcal{H}) \,|\, \widehat{d}_w(U,Z) < \varepsilon\},$$

centred at some $Z \subset \mathcal{Z}_M$, actually cover $\mathcal{C}_w(\mathcal{H})$. Pick $U \in \mathcal{C}_w(\mathcal{H})$: as $U \subset B_\mathcal{H}$, $U$ intersects some of the balls $\mathfrak{B}_w(x_n, \varepsilon)$, so let $Z_U \subset \mathcal{Z}_M$ be the collection of the corresponding centres of such balls. Thus, $U \subset Z_U(\varepsilon)$ and $U \cap \mathfrak{B}_w(x, \varepsilon) \neq \emptyset$ for any $x \in Z_U$. The last two properties are precisely the assumption of Lemma 5.4(iii), that then implies $\widehat{d}_w(U, Z_U) < \varepsilon$. In conclusion, each $U \in \mathcal{C}_w(\mathcal{H})$ belongs to the $\widehat{d}_w$-open ball centred at $Z_U$ and with radius $\varepsilon$, and irrespectively of $U$ the number of such balls is finite, thus realising a finite cover of $\mathcal{C}_w(\mathcal{H})$. □

*Proof (Proof of Theorem 5.2(v))* Both $f(U_n)$ and $f(U)$ are weakly closed, hence also weakly compact subsets of $B_\mathcal{H}$. In particular it makes sense to evaluate $\widehat{d}_w(f(U_n), f(U))$.

We start with proving that $d_w(f(U_n), f(U)) \xrightarrow{n \to \infty} 0$. Let $\varepsilon > 0$. The weakly open $\varepsilon$-expansion $f(U)(\varepsilon)$ of $f(U)$ (see (5.16) above) is weakly open in $B_\mathcal{H}$, namely open in the relative topology of $B_\mathcal{H}$ induced by the weak topology of $\mathcal{H}$. By weak continuity, $f^{-1}(f(U)(\varepsilon))$ too is weakly open in $B_\mathcal{H}$, and in fact it is a relatively open neighbourhood of $U$, for $f(U) \subset f(U)(\varepsilon) \Rightarrow U \subset f^{-1}(f(U)(\varepsilon))$. The set $B_\mathcal{H} \setminus f^{-1}(f(U)(\varepsilon))$ is therefore weakly closed and hence weakly compact in $B_\mathcal{H}$, implying that from any point $u \in U$ one has a notion of weak metric distance between $u$ and $B_\mathcal{H} \setminus f^{-1}(f(U)(\varepsilon))$. So set

$$\widetilde{\varepsilon} := \inf_{u \in U} \inf \left\{ \|u - z\|_w \,\Big|\, z \in B_\mathcal{H} \setminus f^{-1}(f(U)(\varepsilon)) \right\}.$$

It must be $\widetilde{\varepsilon} > 0$, otherwise there would be a common point in $U$ and $B_\mathcal{H} \setminus f^{-1}(f(U)(\varepsilon))$ (owing to the weak closedness of the latter). Thus, any weakly open expansion of $U$ up to $U(\widetilde{\varepsilon})$ is surely contained in $f^{-1}(f(U)(\varepsilon))$, whence also $f(U(\widetilde{\varepsilon})) \subset f(U)(\varepsilon)$. Now, as $U_n \xrightarrow{\widehat{d}_w} U$, there is $n_\varepsilon \in \mathbb{N}$ (in fact depending on $\widetilde{\varepsilon}$, and therefore on $\varepsilon$) such that $d_w(U_n, U) < \widetilde{\varepsilon}$ for all $n \geqslant n_\varepsilon$: then (Lemma 5.4(ii)) $U_n \subset U(\widetilde{\varepsilon})$ for all $n \geqslant n_\varepsilon$. As a consequence, for all $n \geqslant n_\varepsilon$, $f(U_n) \subset f(U(\widetilde{\varepsilon})) \subset f(U)(\varepsilon)$. Using again Lemma 5.4(ii), $d_w(f(U_n), f(U)) < \varepsilon$ for all $n \geqslant n_\varepsilon$, meaning that $d_w(f(U_n), f(U)) \xrightarrow{n \to \infty} 0$.

Let us now turn to proving that $d_w(f(U), f(U_n)) \xrightarrow{n \to \infty} 0$. Assume for contradiction that, up to passing to a subsequence, still denoted with $(U_n)_{n \in \mathbb{N}}$, there is $\varepsilon_0 > 0$ such that $d_w(f(U), f(U_n)) \geqslant \varepsilon_0 \,\forall n \in \mathbb{N}$. With respect to such $\varepsilon_0$, as proved in the first part, there is $n_{\varepsilon_0} \in \mathbb{N}$ such that $f(U_n) \subset f(U)(\varepsilon_0) \,\forall n \geqslant n_{\varepsilon_0}$. For any such $n \geqslant n_{\varepsilon_0}$, on account of Lemma 5.4(iii) one deduces from the latter two properties, namely $d_w(f(U), f(U_n)) \geqslant \varepsilon_0$ and $f(U_n) \subset f(U)(\varepsilon_0)$, that there is $y_n \in f(U)$ such that $f(U_n) \cap \mathfrak{B}_w(y_n, \varepsilon_0) = \emptyset$, whence also

$$U_n \cap f^{-1}(\mathfrak{B}_w(y_n, \varepsilon_0)) = \emptyset.$$

From this condition we want now to construct a sufficiently small weak open ball of a point $u \in U$ that is disjoint from all the $U_n$'s as well. The sequence $(u^{(n)})_{n=n_{\varepsilon_0}}^\infty$ with each $u^{(n)} \in U$ such that $f(u^{(n)}) = y_n$, owing to the weak compactness of $U$, has a weakly convergent subsequence to some $u \in U$. (The superscript in $u^{(n)}$ is to warn that each $u^{(n)}$ belongs to $U$, not to $U_n$.) So, up to further refinement, $u^{(n)} \rightharpoonup u$ in $U$, and by weak continuity $y_n = f(u^{(n)}) \rightharpoonup f(u) =: y$. The latter convergence implies that, eventually in $n$, say, $\forall n \geqslant m_{\varepsilon_0}$ for some $m_{\varepsilon_0} \in \mathbb{N}$, $\mathfrak{B}_w(y, \frac{1}{2}\varepsilon_0) \subset \mathfrak{B}_w(y_n, \varepsilon_0)$. In view of the disjointness condition above, one then deduces

$$U_n \cap f^{-1}(\mathfrak{B}_w(y, \tfrac{1}{2}\varepsilon_0)) = \emptyset \qquad \forall n \geqslant m_{\varepsilon_0}\,.$$

As $f^{-1}(\mathfrak{B}_w(y, \frac{1}{2}\varepsilon_0))$ above is an open neighbourhood of $u \in U$ in the relative weak topology of $B_\mathcal{H}$ (weak continuity of $f$), it contains a ball $\mathfrak{B}_w(u, \varepsilon_1)$ around $u$ for some radius $\varepsilon_1 > 0$, whence

$$U_n \cap \mathfrak{B}_w(u, \varepsilon_1) = \emptyset \qquad \forall n \geqslant m_{\varepsilon_0}\,.$$

On account of Lemma 5.4(i), this implies $d_w(U, U_n) \geqslant \varepsilon_1$ $\forall n \geqslant m_{\varepsilon_0}$. However, this contradicts the assumption $\widehat{d}_w(U, U_n) \to 0$. $\square$

## 5.6 Weak gap metric for linear subspaces

Our primary interest is to exploit the $\widehat{d}_w$-convergence for closed subspaces of $\mathcal{H}$, and ultimately for Krylov subspaces, in the sense of the convergence naturally induced by the convergence of the corresponding unit balls as elements of $(\mathcal{C}_w(\mathcal{H}), \widehat{d}_w)$.

In other words, given two closed subspaces $U, V \subset \mathcal{H}$, by definition we identify

$$\widehat{d}_w(U, V) \equiv \widehat{d}_w(B_U, B_V) \tag{5.17}$$

with the r.h.s. defined in (5.9), since $B_U, B_V \in \mathcal{C}_w(\mathcal{H})$. Analogously, given $U$ and a sequence $(U_n)_{n\in\mathbb{N}}$, all closed subspaces of $\mathcal{H}$, we write $U_n \xrightarrow{\widehat{d}_w} U$ to mean that $B_{U_n} \xrightarrow{\widehat{d}_w} B_U$ in the sense of the definition given in the previous Section. This provides a metric topology and a notion of convergence on the set

$$\mathscr{S}(\mathcal{H}) := \{\text{closed linear subspaces of } \mathcal{H}\}\,. \tag{5.18}$$

By linearity, the closedness of each subspace of $\mathcal{H}$ is equivalently meant in the $\mathcal{H}$-norm or in the weak topology. (Recall, however, that the weak topology on $\mathcal{H}$ is *not* induced by the norm $\|\,\|_w$, as this is only the case in $B_\mathcal{H}$.)

**Lemma 5.5** *The set $(\mathscr{S}(\mathcal{H}), \widehat{d}_w)$ is a metric space.*

***Proof*** Positivity and triangular inequality are obvious from (5.17) and Lemma 5.2. Last, to deduce from $\widehat{d}_w(U, V) = 0$ that $U = V$, one observes that $\widehat{d}_w(B_U, B_V) = 0$ and hence $B_U = B_V$. If $u \in U$, then $u/\|u\| \in B_U = B_V$, whence by linearity $u \in V$, thus, $U \subset V$. Exchanging the role of the two subspaces, also $V \subset U$. $\square$

The metric space $(\mathscr{S}(\mathcal{H}), \widehat{d}_w)$ contains in particular the closures of Krylov subspaces, and monitoring the distance between two such subspaces in the $\widehat{d}_w$-metric turns out to be informative in many respects. Unfortunately there is a major drawback, for:

**Lemma 5.6** *The metric space* $(\mathscr{S}(\mathcal{H}), \widehat{d}_w)$ *is not complete.*

*Proof* It is enough to provide an example of $\widehat{d}_w$-Cauchy sequence in $\mathscr{S}(\mathcal{H})$ that does not converge in $\mathscr{S}(\mathcal{H})$. So take $\mathcal{H} = \ell^2(\mathbb{N})$, with the usual canonical orthonormal basis $(e_n)_{n \in \mathbb{N}}$. For $n \in \mathbb{N}$, set $U_n := \mathrm{span}\{e_1 + e_n\} \subset \mathscr{S}(\mathcal{H})$.

Let us show first of all that the sequence $(U_n)_{n \in \mathbb{N}}$ is $\widehat{d}_w$-Cauchy, i.e., that the corresponding unit balls form a Cauchy sequence $(B_{U_n})_{n \in \mathbb{N}}$ in the metric space $(\mathcal{C}_w, \widehat{d}_w)$. A generic $u \in B_{U_n}$ has the form $u = \alpha_u(e_1 + e_n)$ for some $\alpha_u \in \mathbb{C}$ with $|\alpha_u| \leqslant \frac{1}{\sqrt{2}}$. Therefore,

$$\inf_{v \in B_{U_m}} \|u - v\|_w \;\leqslant\; \|\alpha_u(e_1+e_n) - \alpha_u(e_1+e_m)\|_w \;\leqslant\; \frac{1}{\sqrt{2}} \|e_n - e_m\|_w,$$

the first inequality following from the concrete choice $v = \alpha_u(e_1 + e_m) \in B_{U_m}$. Using the above estimate and the fact that $e_n \rightharpoonup 0$, and hence $(e_n)_{n \in \mathcal{H}}$ is Cauchy in $\mathcal{H}$, one deduces

$$d_w(B_{U_n}, B_{U_m}) \;=\; \sup_{u \in U_n} \inf_{v \in B_{U_m}} \|u - v\|_w \;\leqslant\; \frac{1}{\sqrt{2}} \|e_n - e_m\|_w \;\xrightarrow{n,m \to \infty}\; 0.$$

Inverting $n$ and $m$ one also finds $d_w(B_{U_m}, B_{U_n}) \xrightarrow{n,m \to \infty} 0$. The Cauchy property is thus proved.

On account of the completeness of $(\mathcal{C}_w, \widehat{d}_w)$ (Theorem 5.2(ii)), $B_{U_n} \xrightarrow{\widehat{d}_w} B$ for some $B \in \mathcal{C}_w$. Next, let us show that there is no closed subspace $U \subset \mathcal{H}$ with $B_U = B$, which prevents the sequence $(U_n)_{n \in \mathbb{N}}$ to converge in $(\mathscr{S}(\mathcal{H}), \widehat{d}_w)$. To this aim, we shall show that although the line segment

$$\{\beta e_1 \,|\, \beta \in \mathbb{C}, \, |\beta| \leqslant \tfrac{1}{\sqrt{2}}\}$$

is entirely contained in $B$, however $e_1 \notin B$: this clearly prevents $B$ to be the unit ball of a linear subspace. Assume for contradiction that $e_1 \in B$; then, owing to Theorem 5.2(iii) (see formula (5.10) therein), $e_1 \leftharpoonup u_n$ for a sequence $(u_n)_{n \in \mathbb{N}}$ with $u_n \in B_{U_n}$. In fact, weak approximants from $B_{\mathcal{H}}$ of points of the unit sphere $S_{\mathcal{H}}$ are necessarily also norm approximants: explicitly, owing to weak convergence, the sequence $(u_n)_{n \in \mathbb{N}}$ is norm lower semi-continuous, thus,

$$1 \;=\; \|e_1\| \;\leqslant\; \liminf_{n \to \infty} \|u_n\| \;\leqslant\; 1, \qquad \text{whence} \qquad \lim_{n \to \infty} \|u_n\| \;=\; 1;$$

then, since $u_n \rightharpoonup e_1$ and $\|u_n\| \to \|e_1\|$, one has $u_n \to x$ in the $\mathcal{H}$-norm. As a consequence, writing $u_n = \alpha_n(e_1 + e_n)$ for a suitable $\alpha_n \in \mathbb{C}$ with $|\alpha_n| \leqslant \frac{1}{\sqrt{2}}$, one has $|\alpha_n|\sqrt{2} = \|u_n\| \to 1$, whence $|\alpha_n| \to \frac{1}{\sqrt{2}}$. This implies though that

$$\|u_n - e_1\|^2 \;=\; \|\alpha_n(e_1 + e_n) - e_1\|^2 \;=\; |1 - \alpha_n|^2 + |\alpha_n|^2$$

cannot vanish as $n \to \infty$, a contradiction. Therefore, $e_1 \notin B$.

On the other hand, for any $\beta \in \mathbb{C}$ with $|\beta| \leqslant \frac{1}{\sqrt{2}}$, $B_{U_n} \ni \beta(e_1 + e_n) \rightharpoonup \beta e_1$, which by Theorem 5.2(iii) means that $\beta e_1 \in B$. □

Despite the lack of completeness, the metric $\widehat{d_w}$ in $\mathscr{S}(\mathcal{H})$ displays useful properties for our purposes. The first is the counterpart of Theorem 5.2(iii).

**Proposition 5.1** *Let $\mathcal{H}$ be a separable Hilbert space and assume that $U_n \xrightarrow{\widehat{d_w}} U$ as $n \to \infty$ for some $(U_n)_{n \in \mathbb{N}}$ and $U$ in $\mathscr{S}(\mathcal{H})$. Then*

$$U = \{u \in \mathcal{H} \,|\, u_n \rightharpoonup u \text{ for a sequence } (u_n)_{n \in \mathbb{N}} \text{ with } u_n \in U_n\}. \qquad (5.19)$$

*Proof* Call temporarily

$$\widehat{U} := \{u \in \mathcal{H} \,|\, u_n \rightharpoonup u \text{ for a sequence } (u_n)_{n \in \mathbb{N}} \text{ with } u_n \in U_n\},$$

so that the proof consists of showing that $\widehat{U} = U$. From Theorem 5.2(iii) we know that

$$B_{U_n} \xrightarrow{\widehat{d_w}} B_U = \{\widetilde{u} \in B_{\mathcal{H}} \,|\, \widetilde{u}_n \rightharpoonup \widetilde{u} \text{ for a sequence } (\widetilde{u}_n)_{n \in \mathbb{N}} \text{ with } \widetilde{u}_n \in B_{U_n}\}.$$

So now if $u \in U$, then $B_U \ni u/\|u\| \leftharpoonup \widetilde{u}_n$ for some $(\widetilde{u}_n)_{n \in \mathbb{N}}$ with $\widetilde{u}_n \in B_{U_n}$, whence $u \leftharpoonup \|u\|\widetilde{u}_n \in U_n$, meaning that $u \in \widehat{U}$. Conversely, if $u \in \widehat{U}$, and hence $u_n \rightharpoonup u$ for some $(u_n)_{n \in \mathbb{N}}$ with $u_n \in U_n$, then by uniform boundedness $\|u_n\| \leqslant \kappa \,\forall n \in \mathbb{N}$ for some $\kappa > 0$, and by lower semi-continuity of the norm along the limit $\|u\| \leqslant \kappa$ as well. As a consequence, $B_{U_n} \ni \kappa^{-1} u_n \rightharpoonup \kappa^{-1} u$, meaning that $\kappa^{-1} u \in B_U$ and therefore $u \in U$. □

Other relevant features of the $\widehat{d_w}$-metric will be worked out in the next Section in application to Krylov subspaces.

## 5.7 Krylov perturbations in the weak gap-metric

We are mainly concerned with controlling how close two (closures of) Krylov subspaces $\mathcal{K} \equiv \overline{\mathcal{K}(A, g)}$ and $\mathcal{K}' \equiv \overline{\mathcal{K}(A', g')}$ are within the metric space $(\mathscr{S}(\mathcal{H}), \widehat{d_w})$ of closed subspaces of the separable Hilbert space $\mathcal{H}$ with the weak gap metric $\widehat{d_w}$, for given $A, A' \in \mathcal{B}(\mathcal{H})$ and $g, g' \in \mathcal{H}$.

### 5.7.1 Some technical features of the vicinity of Krylov subspaces

A first noticeable feature, that closes the problem left open as one initial motivation in Section 5.5, is the $\widehat{d_w}$-convergence of the finite-dimensional Krylov subspace to the corresponding closed Krylov subspace.

**Lemma 5.7** *Let $A \in \mathcal{B}(\mathcal{H})$ and $g \in \mathcal{H}$ for an infinite-dimensional, separable Hilbert space $\mathcal{H}$. Then*

$$\widehat{d}_w\big(\mathcal{K}_N(A,g), \overline{\mathcal{K}(A,g)}\big) \xrightarrow{N \to \infty} 0, \qquad (5.20)$$

*with the two spaces defined, respectively, in* (5.1) *and* (2.2).

This means that the $\widehat{d}_w$-metric provides the appropriate language to measure the distance between $\mathcal{K}_N(A,g)$ and $\overline{\mathcal{K}(A,g)}$ in an informative way: $\mathcal{K}_N(A,g) \xrightarrow{\widehat{d}_w} \overline{\mathcal{K}(A,g)}$, whereas we saw that it is false in general that $\mathcal{K}_N(A,g) \xrightarrow{\widehat{d}} \overline{\mathcal{K}(A,g)}$.

For the simple proof of this fact, and for later purposes, it is convenient to work out the following useful construction.

**Lemma 5.8** *Let $\mathcal{H}$ be a separable Hilbert space, and let $A \in \mathcal{B}(\mathcal{H})$, $g \in \mathcal{H}$. Set $\mathcal{K} := \overline{\mathcal{K}(A,g)}$.*

*(i) For every $x \in B_{\mathcal{K}}$ there exists a sequence $(u_n)_{n \in \mathbb{N}}$ in $\mathcal{K}(A,g)$ such that $\|u_n\| < 1$ $\forall n \in \mathbb{N}$ and $u_n \xrightarrow[n \to \infty]{\|\cdot\|} x$.*

*(ii) For every $\varepsilon > 0$ there is a cover of $B_{\mathcal{K}}$ consisting of finitely many weakly open balls $\mathfrak{B}_w(x_1, \varepsilon), \ldots, \mathfrak{B}_w(x_M, \varepsilon)$ for some $x_1, \ldots, x_M \in B_{\mathcal{K}}$ and $M \in \mathbb{N}$ all depending on $\varepsilon$. Moreover, each centre $x_j \in B_{\mathcal{K}}$ has an approximant $p_j(A)g$ for some polynomial $p_j$ on $\mathbb{R}$ with $\|p_j(A)g\| < 1$ and $\|p_j(A)g - x_j\| \leqslant \varepsilon$ $\forall j \in \{1, \ldots, M\}$.*

***Proof*** As $x \in \overline{\mathcal{K}(A,g)}$, then $\widetilde{u}_n \xrightarrow{\|\cdot\|} x$ for some sequence $(\widetilde{u}_n)_{n \in \mathbb{N}}$ in $\mathcal{K}(A,g)$. In particular, $\|\widetilde{u}_n\| \to \|x\|$, and it is not restrictive to assume $\|\widetilde{u}_n\| > 0$ for all $n$. Therefore,

$$u_n := \frac{n-1}{n} \frac{\|x\|}{\|\widetilde{u}_n\|} \widetilde{u}_n \in \mathcal{K}(A,g), \qquad \text{and} \qquad \|u_n\| < \|\widetilde{u}_n\| \leqslant 1,$$

and obviously $u_n \xrightarrow[n \to \infty]{\|\cdot\|} x$. This proves part (i). Concerning part (ii), the existence of such cover follows from the compactness of $(B_{\mathcal{H}}, \rho_w)$ and the closure of $B_{\mathcal{K}}$ in $B_{\mathcal{H}}$. The approximants $p_n(A)g$ are then found based on part (i). $\square$

***Proof (Proof of Lemma 5.7)*** Let us use the shorthand $\mathcal{K}_N \equiv \mathcal{K}_N(A,g)$ and $\mathcal{K} \equiv \overline{\mathcal{K}(A,g)}$. As $\mathcal{K}_N \subset \mathcal{K}$, then $d_w(\mathcal{K}_N, \mathcal{K}) = 0$ (see (5.11), Lemma 5.2 above), so only the limit $d_w(\mathcal{K}, \mathcal{K}_N) \equiv d_w(B_{\mathcal{K}}, B_{\mathcal{K}_N}) \to 0$ is to be checked.

For $\varepsilon > 0$ take the finite open $\varepsilon$-cover of $B_{\mathcal{K}}$ constructed in Lemma 5.8 with centres $x_1, \ldots, x_M$ and Krylov approximants $p_1(A)g, \ldots, p_M(A)g$. Let $N_0$ be the largest degree of the $p_j$'s, thus ensuring that $p_j(A)g \in B_{\mathcal{K}_N}$ $\forall j \in \{1, \ldots, M\}$ and $\forall N \geqslant N_0 + 1$.

Now consider an arbitrary integer $N \geqslant N_0 + 1$ and an arbitrary $u \in B_{\mathcal{K}}$. The vector $u$ clearly belongs to at least one of the balls of the finite open cover above: up to renaming the centres, it is non-restrictive to claim that $u \in \mathfrak{B}_w(x_1, \varepsilon)$, and consider the above approximant $p_1(A)g \in B_{\mathcal{K}_N}$ of the ball's centre $x_1$. Thus, $\|u - x_1\|_w < \varepsilon$ and $\|x_1 - p_1(A)g\|_w \leqslant \|x_1 - p_1(A)g\| \leqslant \varepsilon$. Then

$$\inf_{v \in B_{\mathcal{K}_N}} \|u - v\|_w \;\leqslant\; \|u - x_1\|_w + \|x_1 - p_1(A)g\|_w + \inf_{v \in B_{\mathcal{K}_N}} \|p_1(A)g - v\|_w$$
$$< 2\varepsilon$$

whence also $d_w(B_{\mathcal{K}}, B_{\mathcal{K}_N}) = \sup_{u \in B_{\mathcal{K}}} \inf_{v \in B_{\mathcal{K}_N}} \|u - v\|_w \leqslant 2\varepsilon.$   □

Despite the encouraging property stated in Lemma 5.7, one soon learns that the sequences of (closures of) Krylov subspaces with good convergence properties of the Krylov data $A$ and/or $g$ display in general quite a diverse (including non-convergent) behaviour in the $\widehat{d_w}$-metric. This suggests that an efficient control of $\widehat{d_w}$-convergence of Krylov subspaces is only possible under suitable restrictive assumptions.

Lemma 5.9, Example 5.7 and Example 5.8 below are meant to shed some light on this scenario. In particular, Lemma 5.9 establishes that the convergence $g_n \to g$ in $\mathcal{H}$ is sufficient to have $d_w(\mathcal{K}(A,g), \mathcal{K}(A,g_n)) \to 0$.

**Lemma 5.9** *Given a separable Hilbert space $\mathcal{H}$ and $A \in \mathcal{B}(\mathcal{H})$, assume that $g_n \xrightarrow[n \to \infty]{\|\,\|} g$ for vectors $g, g_n \in \mathcal{H}$. Set $\mathcal{K}_n \equiv \overline{\mathcal{K}(A, g_n)}$ and $\mathcal{K} \equiv \overline{\mathcal{K}(A, g)}$.*

*(i) One has $d_w(\mathcal{K}, \mathcal{K}_n) \xrightarrow{n \to \infty} 0$.*
*(ii) From a sequence $(B_{\mathcal{K}_n})_{n \in \mathbb{N}}$ extract, by compactness of $(\mathcal{C}_w(\mathcal{H}), \widehat{d_w})$, a convergent subsequence to some $B \in \mathcal{C}_w(\mathcal{H})$. Then $B_{\mathcal{K}} \subset B$.*

*Proof*
(i) For $\varepsilon > 0$ take the finite open $\varepsilon$-cover of $B_{\mathcal{K}}$ constructed in Lemma 5.8 with centres $x_1, \ldots, x_M$ and Krylov approximants $p_1(A)g, \ldots, p_M(A)g$. In view of the finitely many conditions $\|p_j(A)g\| < 1$ and $p_j(A)g_n \xrightarrow[n \to \infty]{\|\,\|} p_j(A)g$, $j \in \{1, \ldots, M\}$, there is $n_\varepsilon \in \mathbb{N}$ such that $\|p_j(A)g_n - p_j(A)g\| \leqslant \varepsilon$ and $\|p_j(A)g_n\| < 1$ for all $n \geqslant n_\varepsilon$ and $j \in \{1, \ldots, M\}$.

Take $u \in B_{\mathcal{K}}$. Up to re-naming the centres of the cover's balls, $\|u - x_1\|_w < \varepsilon$, $\|x_1 - p_1(A)g\|_w \leqslant \varepsilon$, $\|p_1(A)g_n - p_1(A)g\| \leqslant \varepsilon$, and $\|p_1(A)g_n\| < 1$ $\forall n \geqslant n_\varepsilon$. Then, for any $n \geqslant n_\varepsilon$,

$$\inf_{v \in B_{\mathcal{K}_n}} \|u - v\|_w \;\leqslant\; \|u - x_1\|_w + \|x_1 - p_1(A)g\|_w + \inf_{v \in B_{\mathcal{K}_n}} \|p_1(A)g - v\|_w$$
$$\leqslant 2\varepsilon + \|p_1(A)g - p_1(A)g_n\|_w \leqslant 3\varepsilon,$$

whence also, for $n \geqslant n_\varepsilon$,

$$d_w(\mathcal{K}, \mathcal{K}_n) \equiv d_w(B_{\mathcal{K}}, B_{\mathcal{K}_n}) = \sup_{u \in B_{\mathcal{K}}} \inf_{v \in B_{\mathcal{K}_n}} \|u - v\|_w \leqslant 3\varepsilon.$$

This means precisely that $d_w(\mathcal{K}, \mathcal{K}_n) \to 0$.

(ii) Rename the extracted subsequence again as as $(B_{\mathcal{K}_n})_{n \in \mathbb{N}}$, so that $B_{\mathcal{K}_n} \xrightarrow{\widehat{d_w}} B$. On account of (5.13) (Lemma 5.2),

$$d_w(B_{\mathcal{K}}, B) \leqslant d_w(B_{\mathcal{K}}, B_{\mathcal{K}_n}) + d_w(B_{\mathcal{K}_n}, B).$$

Since $d_w(B_{\mathcal{K}_n}, B) \leqslant \widehat{d}_w(B_{\mathcal{K}_n}, B) \to 0$ by assumption, and $d_w(B_{\mathcal{K}}, B_{\mathcal{K}_n}) \to 0$ as established in part (i), then $d_w(B_{\mathcal{K}}, B) = 0$. Owing to (5.11) (Lemma 5.2), this implies $B_{\mathcal{K}} \subset B$. □

*Example 5.7* In general, the assumptions of Lemma 5.9 are *not enough* to guarantee that also $d_w(\mathcal{K}_n, \mathcal{K}) \to 0$ and hence $\mathcal{K}_n \xrightarrow{\widehat{d}_w} \mathcal{K}$. Consider for instance $\mathcal{H} = \ell^2(\mathbb{N})$ and the right shift operator $A \equiv R$, acting as $Re_k = e_{k+1}$ on the canonical basis $(e_k)_{k \in \mathbb{N}}$. As in Example 5.3, $R$ admits a dense of cyclic vectors, as well as a dense of non-cyclic vectors: so, with respect to the general setting of Lemma 5.9, take now $g$ to be *non-cyclic*, say, $g = e_2$, and $(g_n)_{n \in \mathbb{N}}$ to be a sequence of $\mathcal{H}$-norm approximants of $g$ that are all *cyclic*. Concerning the subspaces $\mathcal{K}_n := \overline{\mathcal{K}(R, g_n)}$ and $\mathcal{K} := \overline{\mathcal{K}(R, g)}$, $\mathcal{K}_n = \mathcal{H}$ $\forall n \in \mathbb{N}$ by cyclicity, and $\mathcal{K} = \{e_1\}^\perp \subsetneq \mathcal{H}$. As $\mathcal{K} \subset \mathcal{K}_n$, then $d_w(\mathcal{K}, \mathcal{K}_n) = 0$, a conclusion consistent with Lemma 5.9, for $(\mathcal{K}_n)_{n \in \mathbb{N}}$ is obviously $\widehat{d}_w$-Cauchy and Lemma 5.9 implies $d_w(\mathcal{K}, \mathcal{K}_n) \to 0$. On the other hand,

$$d_w(B_{\mathcal{K}_n}, B_{\mathcal{K}}) = \sup_{u \in B_{\mathcal{K}_n}} \inf_{v \in B_{\mathcal{K}}} \|u - v\|_w \geqslant \inf_{v \in B_{\mathcal{K}}} \|e_1 - v\|_w > 0,$$

which prevents $d_w(\mathcal{K}_n, \mathcal{K})$ to vanish with $n$.

*Example 5.8* In general, with respect to the setting of Lemma 5.9 and Example 5.7, the sole convergence $g_n \xrightarrow{\|\cdot\|} g$ is *not enough* to guarantee that $(\mathcal{K}_n)_{n \in \mathbb{N}}$ be $\widehat{d}_w$-Cauchy. For, again with the right shift $R$ on $\mathcal{H} = \ell^2(\mathbb{N})$, take now a sequence $(\widetilde{g}_n)_{n \in \mathbb{N}}$ of cyclic vectors for $R$ such that $\widetilde{g}_n \xrightarrow{\|\cdot\|} e_2$, and set

$$g_n := \begin{cases} \widetilde{g}_n & \text{for even } n, \\ e_2 & \text{for odd } n. \end{cases}$$

Thus, $g_n \xrightarrow{\|\cdot\|} g := e_2$. For even $n$, $\mathcal{K}_n := \overline{\mathcal{K}(R, g_n)} = \mathcal{H}$ and $\mathcal{K}_{n+1} = \{e_1\}^\perp$, whence

$$d_w(B_{\mathcal{K}_n}, B_{\mathcal{K}_{n+1}}) = \sup_{u \in B_{\mathcal{K}_n}} \inf_{v \in B_{\mathcal{K}_{n+1}}} \|u - v\|_w \geqslant \inf_{v \in B_{\mathcal{K}_{n+1}}} \|e_1 - v\|_w > 0,$$

which prevents $(\mathcal{K}_m)_{m \in \mathbb{N}}$ to be $\widehat{d}_w$-Cauchy.

### 5.7.2 Existence of $\widehat{d}_w$-limits. Krylov inner approximability.

Based on the examples discussed above, one is to expect a variety sufficient conditions ensuring the convergence of a sequence of (closures of) Krylov subspaces to a (closure of) Krylov subspace. In this Subsection we discuss one mechanism of convergence that is meaningful in our context of Krylov perturbations.

**Proposition 5.2** *Let $\mathcal{H}$ be a separable Hilbert space, $A \in \mathcal{B}(\mathcal{H})$, and $g \in \mathcal{H}$. Assume further that there is a sequence $(g_n)_{n \in \mathbb{N}}$ such that*

$$g_n \in \overline{\mathcal{K}(A,g)} \;\; \forall n \in \mathbb{N} \qquad \text{and} \qquad g_n \xrightarrow[n\to\infty]{\|\;\|} g. \qquad (5.21)$$

*Then $\overline{\mathcal{K}(A,g_n)} \xrightarrow{\widehat{d_w}} \overline{\mathcal{K}(A,g)}$.*

***Proof*** Let us use the shorthand $\mathcal{K}_n \equiv \overline{\mathcal{K}(A,g_n)}$, $\mathcal{K} \equiv \overline{\mathcal{K}(A,g)}$. As $\mathcal{K}_n \subset \mathcal{K}$, then $d_w(\mathcal{K}_n, \mathcal{K}) = 0$. As $g_n \to g$ in $\mathcal{H}$, then $d_w(\mathcal{K}, \mathcal{K}_n) \to 0$ (Lemma 5.9). Thus, $\mathcal{K}_n \xrightarrow{\widehat{d_w}} \mathcal{K}$. □

The above convergence $\overline{\mathcal{K}(A,g_n)} \xrightarrow{\widehat{d_w}} \overline{\mathcal{K}(A,g)}$, in view of condition (5.21), expresses the **inner approximability** of $\overline{\mathcal{K}(A,g)}$. In fact, $\overline{\mathcal{K}(A,g_n)} \subset \overline{\mathcal{K}(A,g)}$. Condition (5.21) includes also the case of approximants $g_n$ from $\mathcal{K}(A,g)$ or also from $\mathcal{K}_n(A,g)$ (the $n$-th order Krylov subspace (5.1)). For instance, set

$$g_n := \sum_{k=0}^{n-1} \frac{1}{n^{2k} \|A\|_{\mathrm{op}}^k} A^k g \in \mathcal{K}_n(A,g), \qquad n \in \mathbb{N},$$

and as

$$\|g - g_n\| \leqslant \sum_{k=1}^{n-1} \frac{\|A^k g\|}{n^{2k} \|A\|_{\mathrm{op}}^k} \leqslant \|g\| \sum_{k=1}^{n-1} \frac{1}{n^{2k}} \leqslant \frac{\|g\|}{n},$$

then $\overline{\mathcal{K}(A,g)} \supset \mathcal{K}(A,g) \supset \mathcal{K}_n(A,g) \ni g_n \to g$ in $\mathcal{H}$.

### 5.7.3 Krylov solvability along $\widehat{d_w}$-limits

Let us finally scratch the surface of a very central question for the present investigation, namely how a perturbation of a given inverse linear problem, that is small in $\widehat{d_w}$-sense for the corresponding Krylov subspaces, affects the Krylov solvability.

Far from a general answer, we have at least the tools to control the following class of cases. The proof is fast, but it relies on two non-trivial toolboxes.

**Proposition 5.3** *Let $\mathcal{H}$ be a separable Hilbert space. The following be given:*

- *an operator $A \in \mathcal{B}(\mathcal{H})$ with inverse $A^{-1} \in \mathcal{B}(\mathcal{H})$;*
- *a sequence $(g_n)_{n \in \mathbb{N}}$ in $\mathcal{H}$ such that for each $n$ the (unique) solution $f_n := A^{-1} g_n$ to the inverse problem $A f_n = g_n$ is a Krylov solution;*
- *a vector $g \in \mathcal{H}$ such that $\overline{\mathcal{K}(A,g_n)} \xrightarrow{\widehat{d_w}} \overline{\mathcal{K}(A,g)}$ as $n \to \infty$.*

*Then the (unique) solution $f := A^{-1} g$ to the inverse problem $Af = g$ is a Krylov solution. If in addition $g_n \to g$, respectively $g_n \rightharpoonup g$, then $f_n \to f$, respectively $f_n \rightharpoonup f$.*

*Proof* As $A$ is a bounded bijection of $\mathcal{H}$ with bounded inverse, $A$ is a strongly continuous and closed $\mathcal{H} \to \mathcal{H}$ (linear) map, and therefore it also weakly continuous and weakly closed. Up to a non-restrictive scaling one may assume that $\|A\|_{\mathrm{op}} \leqslant 1$, implying that $A$ maps $B_\mathcal{H}$ into itself. The conditions of Theorem 5.2(v) are therefore matched. Thus, from $\overline{\mathcal{K}(A,g_n)} \xrightarrow{\widehat{d_w}} \overline{\mathcal{K}(A,g)}$ one deduces $A\overline{\mathcal{K}(A,g_n)} \xrightarrow{\widehat{d_w}} A\overline{\mathcal{K}(A,g)}$. On the other hand, based on a result that we proved in Proposition 2.2(ii), the assumption that $f_n \in \overline{\mathcal{K}(A,g_n)}$ is equivalent to $A\mathcal{K}(A,g_n) = \mathcal{K}(A,g_n)$. Thus, $\overline{\mathcal{K}(A,g_n)} = A\overline{\mathcal{K}(A,g_n)} \xrightarrow{\widehat{d_w}} A\overline{\mathcal{K}(A,g)}$. The $\widehat{d_w}$-limit being unique, $A\overline{\mathcal{K}(A,g)} = \overline{\mathcal{K}(A,g)}$. Then, again on account of Proposition 2.2(ii), $f \in \overline{\mathcal{K}(A,g_n)}$. This proves the main statement; the additional convergences of $(f_n)_{n\in\mathbb{N}}$ to $f$ are obvious. □

*Remark 5.2* It is worth stressing that the the control of the perturbation in Proposition 5.3, namely the assumption $\overline{\mathcal{K}(A,g_n)} \xrightarrow{\widehat{d_w}} \overline{\mathcal{K}(A,g)}$, does not necessarily correspond to some $\mathcal{H}$-norm vicinity between $g_n$ and $g$ (in Proposition 5.2, instead, we had discussed a case where $\overline{\mathcal{K}(A,g_n)} \xrightarrow{\widehat{d_w}} \overline{\mathcal{K}(A,g)}$ is a consequence of $g_n \to g$ in $\mathcal{H}$). The following example elucidates the situation. With respect to the general setting of Proposition 5.3, consider $\mathcal{H} = \ell^2(\mathbb{N})$, $A = \mathbb{1}$, $g = 0$, $g_n = e_n$ (the $n$-th canonical basis vector), and hence

$$\mathcal{K}_n := \overline{\mathcal{K}(A,g_n)} = \mathrm{span}\{e_n\},$$
$$\mathcal{K} := \overline{\mathcal{K}(A,g)} = \{0\}.$$

Obviously $B_K \subset B_{K_n}$, whence $d_w(K,K_n) = 0$, on account of (5.11) and (5.17). On the other hand, a generic $u \in B_{K_n}$ has the form $u = \alpha e_n$ for some $|\alpha| \leqslant 1$. Therefore,

$$d_w(\mathcal{K}_n, \mathcal{K}) = \sup_{u \in B_{\mathcal{K}_n}} \inf_{v \in B_{\mathcal{K}}} \|u-v\|_w = \sup_{u \in B_{\mathcal{K}_n}} \|u\|_w \leqslant \|e_n\|_w \xrightarrow{n \to \infty} 0.$$

s This shows that $\overline{\mathcal{K}(A,g_n)} \xrightarrow{\widehat{d_w}} \overline{\mathcal{K}(A,g)}$. Thus, all assumptions of Proposition 5.3 are matched. However, it is false that $g_n$ converges to $g$ in norm: in this case it is only true that $g_n \rightharpoonup g$ (weakly in $\mathcal{H}$), indeed $e_n \rightharpoonup 0$.

## 5.8 Perspectives on the perturbation framework

In retrospect, a few concluding observations are in order.

It was already argued in Section 5.1 that the main perspective here is to regard a perturbed inverse problem as a potentially "easier" source of information, including Krylov solvability, for the original, unperturbed problem, and conversely to understand when a given inverse problem looses Krylov solvability under small perturbations, that in practice would correspond to uncertainties of various sort, thus making Krylov subspace methods potentially unstable.

The evidence from Section 5.3 is that a controlled vicinity of the perturbed operator or the perturbed datum is not sufficient, alone, to decide on the above questions, for Krylov solvability may well persist, disappear, or appear in the limit when the perturbation is removed. And the idea inspiring Section 5.4 is that constraining the perturbation within certain classes of operators may provide the additional information needed.

In Sections 5.5-5.7 the inverse problem perturbation is formulated in terms of a convenient topology that allows one to predict whether Krylov solvability persists or is washed out. On a conceptual footing this *is* the appropriate approach, because we know from Chapters 2 and 4 that Krylov solvability is essentially a structural property of the Krylov subspace $\mathcal{K}(A,g)$, therefore it is natural to compare Krylov subspaces in a meaningful sense.

The weak gap metric for linear subspaces of $\mathcal{H}$, while being encouraging in many respects ($\mathcal{K}_N \to \mathcal{K}$, inner approximability, stability under perturbations in the sense of Proposition 5.3), suffers various limitations (indirectly due to the lack of completeness of the $\widehat{d}_w$-metric out of the Hilbert closed unit ball, in turn due to the lack of metrisability of the weak topology out of the unit ball).

It is plausible to expect that the informative control of the inverse problem perturbation, as far as Krylov solvability is concerned, is a combination of an efficient distance between Krylov subspaces, vicinity of operators and of data, and restriction to classes of distinguished operators.

# Appendix A
# Outlook on general projection methods and weaker convergence

## A.1 Standard projection methods and beyond

As anticipated in Section 1.2, we detour from the main line of reasoning therein and devote this Appendix to broadening the main point of view of this book – the solvability of inverse problems within a Krylov-based scheme – so as to include general projection methods.

It was already argued in Section 1.2 that Krylov subspace methods are indeed a popular subclass of a much wider spectrum of theoretical and numerical schemes for finite- and infinite-dimensional inverse problems, where the truncation is performed with respect to a variety of different solution and trial spaces.

Our key perspective here is the *convergence* of finite-dimensional approximants, in particular the strength (in fact, the topology) of convergence. In this book the investigation of *Krylov solvability* of the abstract inverse problem $Af = g$ concerned the convergence $f_N \xrightarrow{N\to\infty} f$ of the approximants $f_N \in \mathcal{K}(A,g)$ (for concreteness, $f_N \in \mathcal{K}_N(A,g)$) to a solution $f$ *in the Hilbert space norm*. More generally, approximating solutions in norm is the purpose of other standard projection methods, each of which is precisely designed for and applicable to those classes of inverse problems for which a satisfactory norm convergence theory can be demonstrated.

This leaves the question open of analysing the scenario where the finite-dimensional truncation is performed under weakened, yet practically reasonable (or unavoidable) conditions that do *not* ensure norm convergence – a problem to which we cannot but attribute far-reaching conceptual relevance at this abstract level, let alone practical importance in applications, and yet lacks at present a comprehensive investigation.

In this respect, the possible failure of Krylov solvability is an example of the impossibility of having certain finite-dimensional approximants (the Krylov iterates) to converge in norm to an actual solution.

In the same spirit, the focus of this Appendix is on *generic convergence phenomena* for approximate solutions to linear inverse problems, in particular (weak) convergence mechanisms not considered in typical (strong) approximation schemes.

Our long term ambition would be to systematically describe a much wider variety of weaker, yet fully informative ways to control the finite-dimensional truncations of an infinite-dimensional inverse problem and the possible asymptotic convergence of the approximated solutions. At this initial stage we cannot be so comprehensive and we content ourselves to outline certain general features. As a compensation for such an initial lack of systematicness, we supplement our discussion with an amount of instructing examples that challenge the intuition.

We already sketched a description of standard '*projection methods*' in Section 1.2 (other truncation schemes include [104, 105, 69, 113, 110, 111, 65]: their typical approach consists of studying the discrete convergence using a system of connection operators, this way encompassing a vast class of approximation schemes). Surely the most popular are the celebrated and already-mentioned '*Petrov-Galerkin projection methods*' [69, Chapter 4], [96, Chapter 5], [9, Chapter 9]. Such truncation schemes are customarily implemented under special assumptions that at this abstract level can be presented as follows:

- the spanning vectors $(u_n)_{n\in\mathbb{N}}$ and $(v_n)_{n\in\mathbb{N}}$ constitute orthonormal *bases* of $\mathcal{H}$;
- the linear operator and the chosen bases are such that both the infinite-dimensional problem (1.1) and each truncated problem (1.3) admit a unique solution, respectively $f$ and $\widehat{f^{(N)}}$, satisfying $\|\widehat{f^{(N)}} - f\|_{\mathcal{H}} \to 0$ as $N \to \infty$.

All this is very familiar for certain classes of boundary value problems on $L^2(\Omega)$, where $\Omega$ is some domain in $\mathbb{R}^d$, the typical playground for Galerkin and Petrov-Galerkin finite element methods [37, 87]. In these cases $A$ is an *unbounded* operator, say, of elliptic type [37, Chapter 3], [87, Chapter 4], of Friedrichs type [37, Section 5.2], [38, 6, 8], of parabolic type [37, Chapter 6], [87, Chapter 5], of 'mixed' (i.e., inducing saddle-point problems) type [37, Section 2.4 and Chapter 4], etc. Such $A$'s are assumed to satisfy (and so do they in applications) some kind of coercivity (as is the case for elliptic problems [87, Section 4.2.3]), or more generally one among the various classical conditions that ensure the corresponding problem (1.1) to be well-posed, such as the Banach-Nečas-Babuška Theorem or the Lax-Milgram Lemma [37, Chapter 2]. When the differential operator is non-coercive, additional sufficient conditions have been studied for the stability of the truncated problem and for the quasi optimality of the discretization scheme [17, 16, 15]. All this makes the finite-dimensional truncation and the infinite-dimensional error analysis widely studied and well understood. In that context, in order for the finite-dimensional solutions to converge in the Hilbert norm, one requires stringent yet often plausible conditions [37, Section 2.2-2.4], [87, Section 4.2]:

(a) the truncation spaces need to approximate suitably well the ambient space $\mathcal{H}$ ('*approximability*', e.g., the interpolation capability of finite elements),
(b) the reduced problems admit solutions that are uniformly controlled by the data ('*uniform stability*'),
(c) finite-dimensional solutions are suitably good approximants for the original (infinite-dimensional) problem ('*asymptotic consistency*'),
(d) one has some suitable boundedness of the problem in appropriate topologies ('*uniform continuity*').

As plausible as the above conditions are, they are *not* matched by several other types of inverse problems of applied interest.

Mathematically this is the case, for instance, whenever *A* does not have a 'good' inverse, say, when *A* is a compact operator on $\mathcal{H}$ with arbitrarily small singular values. Another meaningful example is precisely the lack of Krylov solvability, as commented above.

For these reasons, we refer in this appendix to '*general projection methods*', thus moving *outside* the standard working assumptions of the Petrov-Galerkin framework, in the sense that, while still projecting the infinite-dimensional problem to finite-dimensional truncations, yet

(i) the problem (1.1) is only assumed to be solvable;
(ii) no a priori conditions are assumed that guarantee the truncated problems (1.3) to be well-defined (they may thus have a multiplicity of solutions), let alone solvable (they may have no solution at all);
(iii) not all the standard conditions (a)-(d) above are assumed (unlike for Petrov-Galerkin schemes), which ensured a convergence theory in the Hilbert norm of the error and/or the residual along the sequence of approximate solutions.

In particular, with reference to (iii), we carry on the point of view that error and residual may be controlled in a still informative way in some *weaker* sense than the expected norm topology of the Hilbert space.

To this aim, we identify practically plausible sufficient conditions for the error or the residual to be small in such generalised senses and we discuss the *mechanisms* why the same indicators may actually fail to vanish in norm. We also argue the genericity of such situations. We specialise our observations first on the case when *A* is compact (Section A.3) and then when *A* is generically bounded (Section A.4).

Let us stress once again that this is only scratching the surface of a line of reasoning that surely deserves further investigation.

Let us denote, throughout this Appendix, the (infinite-dimensional) Hilbert norm and the finite-dimensional norm with explicit symbols, $\|\ \|_{\mathcal{H}}$ and $\|\ \|_{\mathbb{C}^N}$, so as to avoid confusion.

## A.2 Finite-dimensional truncation

### A.2.1 Set up and notation

Let us start with revisiting, in the present abstract terms, the setting and the formalism for a generic projection scheme – in the framework of **Galerkin** and **Petrov-Galerkin methods** this is customarily referred to as the **approximation setting** [37, Section 2.2.1].

Let $(u_n)_{n\in\mathbb{N}}$ and $(v_n)_{n\in\mathbb{N}}$ be two orthonormal *systems* of the considered Hilbert space $\mathcal{H}$. They need not be orthonormal *bases*, although their completeness is often crucial for the goodness of the approximation. Their choice depends on the specific approach. In the framework of finite element methods they can be taken to be the global shape functions of the interpolation scheme [37, Chapter 1]. For Krylov subspace methods they are just the orthonormalisaton of the spanning vectors of the associated Krylov subspace (Section 2.1 and 4.1).

Correspondingly, for each $N \in \mathbb{N}$, the orthonormal projections in $\mathcal{H}$ respectively onto $\mathrm{span}\{u_1,\ldots,u_N\}$ and $\mathrm{span}\{v_1,\ldots,v_N\}$ shall be

$$P_N := \sum_{n=1}^{N} |u_n\rangle\langle u_n|, \qquad Q_N := \sum_{n=1}^{N} |v_n\rangle\langle v_n|. \tag{A.1}$$

Associated to a given inverse problem

$$Af = g \tag{A.2}$$

in $\mathcal{H}$, one considers the finite-dimensional truncations induced by $P_N$ and $Q_N$, that is, for each $N$, the problem to find solutions $\widehat{f^{(N)}} \in P_N\mathcal{H}$ to the equation

$$(Q_N A P_N)\widehat{f^{(N)}} = Q_N g. \tag{A.3}$$

In (A.3) $Q_N g = \sum_{n=1}^{N}\langle v_n, g\rangle v_n$ is the datum and $\widehat{f^{(N)}} = \sum_{n=1}^{N}\langle u_n, \widehat{f^{(N)}}\rangle u_n$ is the unknown, and the **compression** $Q_N A P_N$ is only non-trivial as a map from $P_N\mathcal{H}$ to $Q_N\mathcal{H}$, its kernel containing at least the subspace $(\mathbb{1} - P_N)\mathcal{H}$.

There is an obvious and irrelevant degeneracy (which is infinite when $\dim \mathcal{H} = \infty$) in (A.3) when it is regarded as a problem on the whole $\mathcal{H}$. The actual interest towards (A.3) is the problem resulting from the identification $P_N\mathcal{H} \cong \mathbb{C}^N \cong Q_N\mathcal{H}$, in terms of which $P_N f \in \mathcal{H}$ and $Q_N g \in \mathcal{H}$ are canonically identified with the vectors

$$f_N = \begin{pmatrix} \langle u_1, f\rangle \\ \vdots \\ \langle u_N, f\rangle \end{pmatrix} \in \mathbb{C}^N, \qquad g_N = \begin{pmatrix} \langle v_1, g\rangle \\ \vdots \\ \langle v_N, g\rangle \end{pmatrix} \in \mathbb{C}^N, \tag{A.4}$$

and $Q_N A P_N$ with a $\mathbb{C}^N \to \mathbb{C}^N$ linear map represented by the $N \times N$ matrix $A_N = (A_{N;ij})_{i,j\in\{1,\ldots,N\}}$

$$A_{N;ij} = \langle v_i, Q_N A P_N u_j\rangle. \tag{A.5}$$

The matrix $A_N$ is what in the framework of finite element methods for partial differential equations is customarily referred to as the **stiffness matrix**.

We shall call the linear inverse problem

$$A_N f^{(N)} = g_N \tag{A.6}$$

with datum $g_N \in \mathbb{C}^N$ and unknown $f^{(N)} \in \mathbb{C}^N$, and matrix $A_N$ defined by (A.5), the **$N$-dimensional truncation** of the original problem (A.2).

Let us stress the meaning of the present notation.

- $Q_N A P_N$, $P_N f$, and $Q_N g$ are objects (one operator and two vectors) referred to the whole Hilbert space $\mathcal{H}$, whereas $A_N$, $f^{(N)}$, $f_N$, and $g_N$ are the analogues referred now to the space $\mathbb{C}^N$.
- The subscript in $A_N$, $f_N$, and $g_N$ indicates that the components of such objects are precisely the corresponding components, up to order $N$, respectively of $A$, $f$, and $g$, with respect to the tacitly declared bases $(u_n)_{n\in\mathbb{N}}$ and $(v_n)_{n\in\mathbb{N}}$, through formulas (A.4)-(A.5).
- As opposite, the superscript in $f^{(N)}$ indicates that the components of the $\mathbb{C}^N$-vector $f^{(N)}$ are not necessarily to be understood as the first $N$ components of the $\mathcal{H}$-vector $f$ with respect to the basis $(u_n)_{n\in\mathbb{N}}$, and in particular for $N_1 < N_2$ the components of $f^{(N_1)}$ are not a priori equal to the first $N_1$ components of $f^{(N_2)}$. In fact, if $f \in \mathcal{H}$ is a solution to $Af = g$, it is evident from obvious counterexamples that in general the truncations $A_N$, $f_N$, $g_N$ do *not* satisfy the identity $A_N f_N = g_N$, whence the notation $f^{(N)}$ for the unknown in (A.6).
- For a $\mathbb{C}^N$-vector $f^{(N)}$ the notation $\widehat{f^{(N)}}$ indicates a vector in $\mathcal{H}$ whose first $N$ components, with respect to the basis $(u_n)_{n\in\mathbb{N}}$, are precisely those of $f^{(N)}$, all others being zero. Thus, as pedantic as it looks, $f^{(N)} = (\widehat{f^{(N)}})_N$ and $f_N = (\widehat{f_N})_N$, and of course in general $f \neq \widehat{f_N}$.

With $A$, $g$, $(u_n)_{n\in\mathbb{N}}$, and $(v_n)_{n\in\mathbb{N}}$ explicitly known, the truncated problem (A.6) is explicitly formulated and, being finite-dimensional, it is suited for numerical algorithms.

This poses the general question on *whether the truncated problem itself is solvable, and whether its exact or approximate solution $f^{(N)}$ is close to the exact solution $f$ and in which (possibly quantitative) sense*.

### A.2.2 Singularity of the truncated problem

It is clear, first of all, that the question of the singularity of the truncated problem (A.6) makes sense here *eventually in N*, meaning for all $N$'s that are large enough. For a *fixed* value of $N$ the truncation might drastically alter the problem so as to make it manifestly non-informative as compared to $Af = g$, such an alteration then disappearing for larger values.

Yet, even when the solvability of $A_N f^{(N)} = g_N$ is inquired eventually in $N$, it is no surprise that the answer is generically negative.

*Example A.1* That the matrix $A_N$ may remain singular for every $N$, even when the operator $A$ is injective, can be seen for example from the truncation of the weighted right-shift operator (Example 2.3) $R_\sigma = \sum_{n\in\mathbb{Z}} \sigma_n |e_{n+1}\rangle\langle e_n|$ on $\ell^2(\mathbb{Z})$ onto the span$\{e_{-N}, e_{-N+1}, \ldots, e_{N-1}, e_N\}$, where $(e_n)_{n\in\mathbb{Z}}$ is the canonical orthonormal basis of $\ell^2(\mathbb{Z})$ and $\sigma \equiv (\sigma_n)_{n\in\mathbb{Z}}$ is sequence derived from a given bounded sequence $(\widetilde{\sigma}_n)_{n\in\mathbb{N}_0}$ with $0 < \widetilde{\sigma}_{n+1} < \widetilde{\sigma}_n$ $\forall n \in \mathbb{N}_0$, $\lim_{n\to\infty} \widetilde{\sigma}_n = 0$, and $\sigma_n = \widetilde{\sigma}_{|n|}$ $\forall n \in \mathbb{Z}$. $R_\sigma$ is indeed compact and injective. However,

$$R_{\sigma,2N+1} := \begin{pmatrix} 0 & \cdots & \cdots & \cdots & 0 \\ \sigma_{-N} & 0 & \cdots & \cdots & 0 \\ 0 & \sigma_{-N+1} & 0 & \cdots & 0 \\ \vdots & \vdots & \ddots & \ddots & 0 \\ 0 & 0 & \cdots & \sigma_{N-1} & 0 \end{pmatrix} \tag{A.7}$$

is singular irrespectively of $N \in \mathbb{N}$, with $\ker R_{\sigma,2N+1} = \operatorname{span}\{e_N\}$. (See Lemma A.1 below for a more general perspective on such an example.) Noticeably, $\sigma(R_{\sigma,2N+1}) = \{0\}$ for each $N$, (explicitly, the vector with all null components but the last one is precisely the eigenvector spanning the kernel of $R_{\sigma,2N+1}$): the erroneous eigenvalue zero that arises in the limit $N \to \infty$ is a typical manifestation of the phenomenon known as **spectral pollution**.

*Remark A.1* In the same spirit, it is not difficult to produce variations of the above example where the matrix $A_N$ is, say, alternately singular and non-singular as $N \to \infty$. Of course, on the other hand, it may also well happen that the truncated matrix is always non-singular: the truncation of the multiplication operator $M = \sum_{n=1}^{\infty} \frac{1}{n} |e_n\rangle\langle e_n|$ on $\ell^2(\mathbb{N})$ with respect to $(e_n)_{n \in \mathbb{N}}$ (Example 2.4) yields the matrix $M_N = \operatorname{diag}(1, \frac{1}{2}, \ldots, \frac{1}{N})$, which is a $\mathbb{C}^N \to \mathbb{C}^N$ bijection for every $N$.

In fact, 'bad' truncations are always possible, as the following mechanism shows.

**Lemma A.1** *Let $\mathcal{H}$ be a Hilbert space with $\dim \mathcal{H} = \infty$, and let $A \in \mathcal{B}(\mathcal{H})$. There always exist two orthonormal systems $(u_n)_{n \in \mathbb{N}}$ and $(v_n)_{n \in \mathbb{N}}$ of $\mathcal{H}$ such that the corresponding truncated matrix $A_N$ defined as in (A.5) is singular for every $N \in \mathbb{N}$.*

*Proof* Let us pick an arbitrary orthonormal system $(u_n)_{n \in \mathbb{N}}$ and construct the other system $(v_n)_{n \in \mathbb{N}}$ inductively. When $N = 1$, it suffices to choose $v_1$ such that $v_1 \perp Au_1$ and $\|v_1\|_{\mathcal{H}} = 1$. Let now $(v_n)_{n \in \{1,\ldots,N-1\}}$ be an orthonormal system in $\mathcal{H}$ satisfying the thesis up to the order $N-1$ and let us construct $v_N$ so that $(v_n)_{n \in \{1,\ldots,N\}}$ satisfies the thesis up to order $N$. To this aim, let us show that a choice of $v_N$ is always possible so that the final row in the matrix $A_N$ has all zero entries. In fact, $(A_N)_{ij} = (Q_N A P_N)_{ij} = \langle v_i, A u_j \rangle$ for $i \in \{1,\ldots,N-1\}$ and $j \in \{1,\ldots,N\}$ and in order for $\langle v_N, Au_j \rangle = 0$ for $j \in \{1, \cdots, N\}$ it suffices to take

$$v_N \perp \operatorname{ran}(AP_N), \qquad v_N \perp \operatorname{ran} Q_{N-1}, \qquad \|v_N\|_{\mathcal{H}} = 1,$$

where $P_N$ and $Q_{N-1}$ are the orthogonal projections defined in A.1. Since $\operatorname{ran}(AP_N)$ and $\operatorname{ran} Q_{N-1}$ are finite-dimensional subspaces of $\mathcal{H}$, there is surely a vector $v_N \in \mathcal{H}$ with the above properties. □

The occurrence described by Lemma A.1 may happen both with an orthogonal and with an oblique projection scheme, namely both when $P_N = Q_N$ and when $P_N \neq Q_N$ eventually in $N$.

As commented already, in the standard framework of (Petrov-)Galerkin methods such an occurrence is *explicitly prevented* by suitable assumptions on $A$, typically coercivity [37, Section 2.2], [87, Section 4.1], and by similar assumptions at the truncated level (e.g., [37, Section 3.2, Section 4.2]), so as to ensure the truncated problem (A.6) to be solvable for all $N$.

As in our discussion we do not exclude a priori such an occurrence, we are compelled to regard $f^{(N)}$ as an approximate solution to the truncated problem, in the sense that

$$A_N f^{(N)} = g_N + \varepsilon^{(N)} \qquad \text{for some } \varepsilon^{(N)} \in \mathbb{C}^N. \tag{A.8}$$

(We write $\varepsilon^{(N)}$ and not $\varepsilon_N$ because there is no reason to claim that the residual $\varepsilon^{(N)}$ in the $N$-dimensional problem is the actual truncation for every $N$ of the same infinite-dimensional vector $\varepsilon \in \mathcal{H}$.)

It would be desirable that $\varepsilon^{(N)}$ is indeed small and asymptotically vanishing with $N$, or even that $\varepsilon^{(N)} = 0$ for $N$ large enough, as is the case in some applications. Morally (up to passing to the weak formulation of the inverse problem), this is the assumption of *asymptotic consistency* naturally made for approximations by Petrov-Galerkin methods [37, Definition 2.15 and Theorem 2.24].

In the present abstract context it is worth remarking that an assumption of that sort it is motivated by the following property, whose proof is postponed to Section A.4.

**Lemma A.2** *Let $A \in \mathcal{B}(\mathcal{H})$ and $g \in \operatorname{ran} A$. Let $A_N$ and $g_N$ be defined as in (A.4)-(A.5) above, for some orthonormal bases $(u_n)_{n \in \mathbb{N}}$, $(v_n)_{n \in \mathbb{N}}$ of $\mathcal{H}$. Then there always exists a sequence $(f^{(N)})_{N \in \mathbb{N}}$ such that*

$$f^{(N)} \in \mathbb{C}^N \qquad \text{and} \qquad \lim_{N \to \infty} \|A_N f^{(N)} - g_N\|_{\mathbb{C}^N} = 0.$$

In other words, there do exist approximate solutions $f^{(N)}$ to (A.6) actually satisfying (A.8) with $\|\varepsilon^{(N)}\|_{\mathbb{C}^N} \to 0$ as $N \to \infty$, so that (A.8) is asymptotically consistent.

### A.2.3 Convergence of the truncated problems

Depending on the context, the vanishing of the displacement $f - \widehat{f^{(N)}}$ as $N \to \infty$ is monitored in various forms with respect to the Hilbert norm of $\mathcal{H}$, the two most typical ones are the norms of the **infinite-dimensional error** $\mathscr{E}_N$ and the **infinite-dimensional residual** $\mathfrak{R}_N$, defined respectively as

$$\begin{aligned} \mathscr{E}_N &:= f - \widehat{f^{(N)}}, \\ \mathfrak{R}_N &:= A(f - \widehat{f^{(N)}}) = g - A\widehat{f^{(N)}}. \end{aligned} \tag{A.9}$$

The nomenclature 'infinite-dimensional', that we shall drop when no confusion arises, is to distinguish them from the error and residual at fixed $N$, which may be indexed by the number of steps in an iterative algorithm.

A first evident obstruction to the actual vanishing of $\mathscr{E}_N$ when $\dim \mathcal{H} = \infty$ is the use of a *non-complete* orthonormal system $(u_n)_{n \in \mathbb{N}}$, that is, such that $\operatorname{span}\{u_n \,|\, n \in$ $\mathbb{N}\}$ is not dense in $\mathcal{H}$. Truncations with respect to a potentially non-complete orthonormal system might appear unwise, but in certain contexts are natural: in the framework of Krylov-based truncations this was the case in Examples 2.1, 2.2, 2.3, 2.6(ii), 2.11, 2.12, 2.13, 2.16, among others. In standard (Petrov-)Galerkin methods such an obstruction is ruled out by an ad hoc '*approximability*' assumption [37, Definition 2.14 and Theorem 2.24] that can be rephrased as the request that $(u_n)_{n\in\mathbb{N}}$ is indeed an orthonormal *basis* of $\mathcal{H}$. Besides, the approximability property is known to fail in situations of engineering interest, as is the case for the failure of the Lagrange finite elements in differential problems for electromagnetism [37, Section 2.3.3].

Even when (complete) orthonormal bases of $\mathcal{H}$ are employed for the truncation, another consequence of the infinite dimensionality is the possibility that error and residual are asymptotically small only in some weaker sense than the customary norm topology of $\mathcal{H}$. On a related scenario, special classes of linear ill-conditioned problems (rank-deficient and discrete ill-posed problems) can be treated with regularisation methods in which the solution is stabilised [108, 57]. The most notable regularisation methods, namely the Tikhonov-Phillips method, the Landweber-Fridman iteration method, and the truncated singular value decomposition, produce indeed a strongly vanishing error [53, 75]. Yet, when the linear inverse problem $Af = g$ is governed by an infinite-rank compact operator $A$, it can be seen that the conjugate gradient method, as well as $\alpha$-processes (in particular, the method of steepest descent), although having regularising properties, may have strongly divergent error in the presence of noise [34] and one is forced to consider weaker forms of convergence. In fact, in [34] the presence of component-wise convergence is also alluded to.

It is therefore natural to contemplate, next to the strong ($\mathcal{H}$-norm) convergence vanishing $\|\mathfrak{R}_N\|_{\mathcal{H}} \to 0$, resp., $\|\mathscr{E}_N\|_{\mathcal{H}} \to 0$, where obviously

$$\|\mathfrak{R}_N\|_{\mathcal{H}} \leqslant \|A\|_{\mathrm{op}} \|\mathscr{E}_N\|_{\mathcal{H}}, \tag{A.10}$$

also the **weak convergence** [88, Section IV.5] $\mathfrak{R}_N \rightharpoonup 0$ or $\mathscr{E}_N \rightharpoonup 0$, or even the **component-wise convergence**, namely the vanishing of each component of the vector $\mathfrak{R}_N$ or $\mathscr{E}_N$ with respect to the considered basis.

As we shall elaborate further on: *a strong control such as $\|\mathfrak{R}_N\|_{\mathcal{H}} \to 0$ or $\|\mathscr{E}_N\|_{\mathcal{H}} \to 0$ is not generic in truncation schemes that do not obey the more stringent assumptions of Petrov-Galerkin projection methods, and only holds under specific a priori conditions on the linear inverse problem.*

On the other hand, even a mere component-wise vanishing of $\mathscr{E}_N$ is in many respects already satisfactorily informative, for in this case each component of $\widehat{f^{(N)}}$ (with respect to the basis $(u_n)_{n\in\mathbb{N}}$) approximates the corresponding component of the exact solution $f$.

## A.3 The compact linear inverse problem

When the operator $A$ is compact on $\mathcal{H}$, it is standard to refer to its **singular value decomposition** [88, Theorem VI.17] as the expansion

$$A = \sum_{n \in J} \sigma_n |\psi_n\rangle\langle\varphi_n|, \tag{A.11}$$

where $n$ runs over $J := \{1, \ldots, M\}$, with $M < \infty$ or $M = \infty$, $(\varphi_n)_n$ and $(\psi_n)_n$ are two orthonormal systems of $\mathcal{H}$, $\sigma_n \geqslant \sigma_{n+1} > 0$ for all $n$, and $\sigma_n \to 0$ if $\mathrm{rank}\, A = \infty$, so that the above series, if infinite, converges in operator norm.

Let us recall that the injectivity of $A$ is tantamount as $(\varphi_n)_{n\in\mathbb{N}}$ being an orthonormal basis and that $\overline{\mathrm{ran}\, A} = \mathcal{H}$ if an only if $(\psi_n)_{n\in\mathbb{N}}$ is an orthonormal basis.

Let us consider the inverse problem (A.2) for compact *and injective $A$*, $g \in \mathrm{ran}\, A$, and $\dim \mathcal{H} = \infty$. The problem is well-defined: there exists a unique $f \in \mathcal{H}$ such that $Af = g$.

The compactness of $A$ has two noticeable consequences here. First, since $\dim \mathcal{H} = \infty$, $A$ is invertible on its range only, and cannot have an everywhere defined bounded inverse: $\mathrm{ran}\, A$ can be dense in $\mathcal{H}$, as in the case of the Volterra operator on $L^2[0,1]$ (Example 2.5), or dense in a closed proper subspace of $\mathcal{H}$, as for the weighted right-shift on $\ell^2(\mathbb{N})$ (Example 2.3).

Furthermore, $A$ and its compression (in the usual meaning of Section A.2.1) are close in a robust sense, as the following standard Lemma shows.

**Lemma A.3** *With respect to an infinite-dimensional separable Hilbert space $\mathcal{H}$, let $A : \mathcal{H} \to \mathcal{H}$ be a compact operator and let $(u_n)_{n\in\mathbb{N}}$ and $(v_n)_{n\in\mathbb{N}}$ be two orthonormal bases of $\mathcal{H}$. Then*

$$\|A - Q_N A P_N\|_{\mathrm{op}} \xrightarrow{N \to \infty} 0, \tag{A.12}$$

*$P_N$ and $Q_N$ being as usual the orthogonal projections* (A.1).

*Proof* Upon splitting

$$A - Q_N A P_N = (A - Q_N A) + Q_N (A - A P_N),$$

it suffices to prove that $\|A - AP_N\|_{\mathrm{op}} \xrightarrow{N \to \infty} 0$ and $\|A - Q_N A\|_{\mathrm{op}} \xrightarrow{N \to \infty} 0$. Let us prove the first limit (the second being completely analogous).

Clearly, it is enough to prove that $\|A - AP_N\|_{\mathrm{op}}$ vanishes assuming further that $A$ has finite rank. Indeed, the difference $(A - AP_N) - (\widetilde{A} - \widetilde{A}P_N)$, where $\widetilde{A}$ is a finite-rank approximant of the compact operator $A$, is controlled in operator norm by $2\|A - \widetilde{A}\|_{\mathrm{op}}$ and hence can be made arbitrarily small.

Thus, we consider non-restrictively $A = \sum_{k=1}^{M} \sigma_k |\psi_k\rangle\langle\varphi_k|$ for some integer $M$, where $(\varphi_k)_{k=1}^{M}$ and $(\psi_k)_{k=1}^{M}$ are two orthonormal systems, and $0 < \sigma_M \leqslant \cdots \leqslant \sigma_1$. Now, for a generic $\xi = \sum_{n=1}^{\infty} \xi_n v_n \in \mathcal{H}$ one has

$$\|(A-AP_N)\xi\|_{\mathcal{H}}^2 = \Bigg\|\sum_{k=1}^{M}\sigma_k\Bigg(\sum_{n=N+1}^{\infty}\xi_n\langle\varphi_k,v_n\rangle\Bigg)\psi_k\Bigg\|_{\mathcal{H}}^2$$

$$= \sum_{k=1}^{M}\sigma_k^2\Bigg|\sum_{n=N+1}^{\infty}\xi_n\langle\varphi_k,v_n\rangle\Bigg|^2 \leqslant \|\xi\|_{\mathcal{H}}^2\sum_{k=1}^{M}\sigma_k^2\,\|(\mathbb{1}-P_N)\varphi_k\|_{\mathcal{H}}^2,$$

therefore

$$\|A-AP_N\|_{\mathrm{op}}^2 \leqslant M\sigma_1^2 \cdot \max_{k\in\{1,\dots,M\}}\|(\mathbb{1}-P_N)\varphi_k\|_{\mathcal{H}}^2 \xrightarrow{N\to\infty} 0,$$

since the above maximum is taken over $M$ (hence, finitely many) quantities, each of which vanishes as $N\to\infty$. $\square$

In the following Theorem we describe the generic behaviour of the well-defined compact inverse problem.

**Theorem A.1** *Consider*

- *the linear inverse problem $Af = g$ in a separable Hilbert space $\mathcal{H}$ for some compact and injective $A: \mathcal{H} \to \mathcal{H}$ and some $g \in \mathrm{ran}A$;*
- *the finite-dimensional truncation $A_N$ obtained by compression with respect to the orthonormal bases $(u_n)_{n\in\mathbb{N}}$ and $(v_n)_{n\in\mathbb{N}}$ of $\mathcal{H}$.*

*Let $(f^{(N)})_{N\in\mathbb{N}}$ be a sequence of approximate solutions to the truncated problems in the quantitative sense*

$$A_N f^{(N)} = g_N + \varepsilon^{(N)}, \qquad f^{(N)},\varepsilon^{(N)} \in \mathbb{C}^N, \qquad \|\varepsilon^{(N)}\|_{\mathbb{C}^N} \xrightarrow{N\to\infty} 0$$

*for every (sufficiently large) $N$. If $\widehat{f^{(N)}}$ is $\mathcal{H}$-norm bounded uniformly in $N$, then*

$$\|\mathfrak{R}_N\|_{\mathcal{H}} \to 0 \quad\text{and}\quad \mathscr{E}_N \rightharpoonup 0 \quad\text{as } N\to\infty.$$

*Proof* We split

$$\begin{aligned}
A\widehat{f^{(N)}} - g &= (A - Q_N A P_N)\widehat{f^{(N)}} \\
&\quad + Q_N A P_N \widehat{f^{(N)}} - Q_N g \\
&\quad + Q_N g - g.
\end{aligned} \qquad (*)$$

By assumption, $\|Q_N g - g\|_{\mathcal{H}} \xrightarrow{N\to\infty} 0$ and

$$\|Q_N A P_N \widehat{f^{(N)}} - Q_N g\|_{\mathcal{H}} = \|A_N f^{(N)} - g_N\|_{\mathbb{C}^N} = \|\varepsilon^{(N)}\|_{\mathbb{C}^N} \xrightarrow{N\to\infty} 0.$$

Moreover, Lemma A.3 and the uniform boundedness of $\widehat{f^{(N)}}$ imply

$$\|(A - Q_N A P_N)\widehat{f^{(N)}}\|_{\mathcal{H}} \leqslant \|A - Q_N A P_N\|_{\mathrm{op}} \|\widehat{f^{(N)}}\|_{\mathcal{H}} \xrightarrow{N\to\infty} 0$$

Plugging the three limits above into (*) proves $\|\mathfrak{R}_N\|_{\mathcal{H}} \to 0$.

Next, in terms of the singular value decomposition (A.11) of $A$, where now $(\varphi_n)_{n\in\mathbb{N}}$ is an orthonormal basis of $\mathcal{H}$, $(\psi_n)_{n\in\mathbb{N}}$ is an orthonormal system, and $0 < \sigma_{n+1} \leqslant \sigma_n$ for all $n \in J = \mathbb{N}$, we write

$$\widehat{f^{(N)}} = \sum_{n\in\mathbb{N}} f_n^{(N)} \varphi_n, \qquad f = \sum_{n\in\mathbb{N}} f_n \varphi_n,$$

whence

$$0 = \lim_{N\to\infty} \|A\widehat{f^{(N)}} - g\|_{\mathcal{H}}^2 = \lim_{N\to\infty} \sum_{n\in J} \sigma_n^2 |f_n^{(N)} - f_n|^2.$$

Then necessarily $\widehat{f^{(N)}}$ converges to $f$ component-wise:

$$\langle u_n, \widehat{f^{(N)}} \rangle \xrightarrow{N\to\infty} \langle u_n, f \rangle \quad \forall n \in \mathbb{N}.$$

On the other hand, $\widehat{f^{(N)}}$ is uniformly bounded in $\mathcal{H}$, $\sup_{N\in\mathbb{N}} \|\widehat{f^{(N)}}\|_{\mathcal{H}} < +\infty$. As well known, latter two conditions are equivalent to the fact that $\widehat{f^{(N)}}$ converges to $f$ weakly ($\mathscr{E}_N \rightharpoonup 0$). □

Theorem A.1 provides sufficient conditions for some form of vanishing of the error and the residual. The key assumptions are:

- *injectivity* of $A$,
- *asymptotic solvability of the truncated problems*, i.e., asymptotic smallness of the finite-dimensional residual $A_N f^{(N)} - g_N$,
- *uniform boundedness of the approximate solutions* $f^{(N)}$.

In fact, injectivity was only used in the analysis of the error in order to conclude $\mathscr{E}_N \rightharpoonup 0$; instead, the conclusion $\|\mathfrak{R}_N\|_{\mathcal{H}} \to 0$ follows irrespectively of injectivity.

It is also worth observing that both the assumption of asymptotic solvability of the truncated problems (which is inspired to Lemma A.2, as commented already) and the assumption of uniform boundedness of the $f^{(N)}$'s are finite-dimensional in nature: their check, from the practical point of view of scientific computing, requires a control of explicit $N$-dimensional vectors for a 'numerically satisfactory' amount of integers $N$.

To further understand the impact of such assumptions, a few remarks are in order.

*Remark A.2 (Genericity)*

Under the conditions of Theorem A.1, the occurrence of the *strong* vanishing of the residual ($\|\mathfrak{R}_N\|_{\mathcal{H}} \to 0$) and the *weak* vanishing of the error ($\mathscr{E}_N \rightharpoonup 0$) as $N \to \infty$ is indeed a *generic behaviour*. For example, the compact inverse problem $Rf = 0$ in $\ell^2(\mathbb{N})$ associated with the weighted right-shift $R$ (Example 2.1) has exact solution $f = 0$. The truncated problem $R_N f^{(N)} = 0$ with respect to $\mathrm{span}\{e_1,\ldots,e_N\}$, where

$$R_N = \begin{pmatrix} 0 & \cdots & \cdots & \cdots & 0 \\ \sigma_1 & 0 & \cdots & \cdots & 0 \\ 0 & \sigma_2 & 0 & \cdots & 0 \\ \vdots & \vdots & \ddots & \ddots & 0 \\ 0 & 0 & \cdots & \sigma_{N-1} & 0 \end{pmatrix}, \qquad (A.13)$$

is solved by the $\mathbb{C}^N$-vectors whose first $N-1$ components are zero, i.e., $\widehat{f^{(N)}} = e_N$. The sequence $(\widehat{f^{(N)}})_{N\in\mathbb{N}} \equiv (e_N)_{N\in\mathbb{N}}$ converges weakly to zero in $\ell^2(\mathbb{N})$, whence indeed $\mathscr{E}_N \rightharpoonup 0$, and also, by compactness, $\|\mathfrak{R}_N\|_{\mathcal{H}} \to 0$. However, $\|\mathscr{E}_N\|_{\mathcal{H}} = 1$ for every $N$, thus the error cannot vanish in the $\mathcal{H}$-norm.

*Remark A.3 ('Bad' approximate solutions)*

The example considered in Remark A.2 is also instructive to understand that generically one may happen to select 'bad' approximate solutions $\widehat{f^{(N)}}$ such that, despite the 'good' strong convergence $\|A_N f^{(N)} - g_N\|_{\mathbb{C}^N} \to 0$, have the unsatisfactory feature $\|f^{(N)}\|_{\mathbb{C}^N} = \|\widehat{f^{(N)}}\|_{\mathcal{H}} \to +\infty$: this is the case if one chooses, for instance, $\widehat{f^{(N)}} = N e_N$. Thus, the uniform boundedness of $\widehat{f^{(N)}}$ in $\mathcal{H}$ required in Theorem A.1 is *not* redundant. (This is also consistent with the fact that whereas by compactness $\widehat{f^{(N)}} \rightharpoonup f$ implies $\|A\widehat{f^{(N)}} - Af\| \to 0$, yet the opposite implication is not true in general.)

*Remark A.4 (The density of* ran$A$ *does not help)* The genericity of convergence discussed in Remarks A.2 and A.3 is not improved by additionally requiring that $\overline{\mathrm{ran} A} = \mathcal{H}$. For instance, the inverse problem considered in Example A.1 involves an operator that is compact, injective, and with dense range, but the finite-dimensional projections produce *for every N* the matrix (A.7) that is singular and for which, therefore, all the considerations of Remarks A.2 and A.3 can be repeated verbatim.

*Remark A.5 ('Bad' truncations and 'good' truncations)*

We saw in Lemma A.1 that 'bad' truncations (i.e., leading to matrices $A_N$ that are, eventually in $N$, all singular) are always possible. On the other hand, there always exists a "good" choice for the truncation – although such a choice might not be identifiable explicitly – which makes the infinite-dimensional residual and error vanish in norm and with no extra assumption of uniform boundedness on the approximate solutions. For instance, in terms of the singular value decomposition (A.11) of $A$, it is enough to choose

$$(u_n)_{n\in\mathbb{N}} = (\varphi_n)_{n\in\mathbb{N}}, \qquad (v_n)_{n\in\mathbb{N}} = (\psi_n)_{n\in\mathbb{N}},$$

in which case $Q_N A P_N = \sum_{n=1}^N \sigma_n |\psi_n\rangle\langle\varphi_n|$ and $A_N = \mathrm{diag}(\sigma_1, \ldots, \sigma_N)$, and for given $g = \sum_{n\in\mathbb{N}} g_n \psi_n$ one has $\widehat{f^{(N)}} = \sum_{n=1}^N \frac{g_n}{\sigma_n} \varphi_n$, where the sequence $(\frac{g_n}{\sigma_n})_{n\in\mathbb{N}}$ belongs to $\ell^2(\mathbb{N})$ owing to the assumption $g \in \mathrm{ran} A$, whence

$$\|f - \widehat{f^{(N)}}\|_{\mathcal{H}}^2 = \sum_{n=N+1}^\infty \left|\frac{g_n}{\sigma_n}\right|^2 \xrightarrow{N\to\infty} 0.$$

## A.4 The bounded linear inverse problem

When $\dim\mathcal{H}=\infty$ and a generic bounded linear operator $A:\mathcal{H}\to\mathcal{H}$ is compressed (in the usual sense of Section A.2) between the spans of the first $N$ vectors of the orthonormal bases $(u_n)_{n\in\mathbb{N}}$ and $(v_n)_{n\in\mathbb{N}}$, then surely $Q_N A P_N \to A$ as $N\to\infty$ in the strong operator topology, that is, $\|Q_N A P_N \psi - A\psi\|_{\mathcal{H}} \xrightarrow{N\to\infty} 0 \ \forall \psi\in\mathcal{H}$, yet the convergence *may fail to occur in the operator norm*.

The first statement is an obvious consequence of the inequality

$$\|(A - Q_N A P_N)\psi\|_{\mathcal{H}} \leqslant \|(\mathbb{1} - Q_N)A\psi\|_{\mathcal{H}} + \|A\|_{\mathrm{op}}\|\psi - P_N\psi\|_{\mathcal{H}}$$

valid for any $\psi\in\mathcal{H}$. The lack of operator norm convergence is also clear for any non-compact $A$: indeed, the operator norm limit of finite-rank operators can only be compact.

For this reason, the control of the infinite-dimensional inverse problem in terms of its finite-dimensional truncated versions is in general less strong.

As a counterpart of Theorem A.1 above, let us discuss the following generic behaviour of *well-posed* bounded linear inverse problems.

**Theorem A.2** *Consider*

- *the linear inverse problem $Af = g$ in a Hilbert space $\mathcal{H}$ for some bounded and injective $A:\mathcal{H}\to\mathcal{H}$ and some $g\in\mathcal{H}$;*
- *the finite-dimensional truncation $A_N$ obtained by compression with respect to the orthonormal bases $(u_n)_{n\in\mathbb{N}}$ and $(v_n)_{n\in\mathbb{N}}$ of $\mathcal{H}$.*

*Let $(f^{(N)})_{N\in\mathbb{N}}$ be a sequence of approximate solutions to the truncated problems in the quantitative sense*

$$A_N f^{(N)} = g_N + \varepsilon^{(N)}, \qquad f^{(N)}, \varepsilon^{(N)} \in \mathbb{C}^N, \qquad \|\varepsilon^{(N)}\|_{\mathbb{C}^N} \xrightarrow{N\to\infty} 0$$

*for every (sufficiently large) $N$.*

*(i) The residuals $\mathfrak{R}_N$ vanish strongly, or weakly or component-wise as $N\to\infty$ if and only if so do the vectors $(A - Q_N A P_N)\widehat{f^{(N)}}$.*

*(ii) In particular, if $\widehat{f^{(N)}}$ converges strongly in $\mathcal{H}$, equivalently, if*

$$\|f^{(N)} - f^{(M)}\|_{\mathbb{C}^{\max\{N,M\}}} \xrightarrow{N,M\to\infty} 0,$$

*then*

$$\|\mathscr{E}_N\|_{\mathcal{H}} \to 0 \qquad \text{and} \qquad \|\mathfrak{R}_N\|_{\mathcal{H}} \to 0 \qquad \text{as } N\to\infty.$$

*Proof* Since

$$\begin{aligned} A\widehat{f^{(N)}} - g &= (A - Q_N A P_N)\widehat{f^{(N)}} \\ &\quad + Q_N A P_N \widehat{f^{(N)}} - Q_N g \\ &\quad + Q_N g - g, \end{aligned} \qquad (*)$$

and since by assumption $\|Q_N g - g\|_{\mathcal{H}} \xrightarrow{N\to\infty} 0$ and

$$\|Q_N A P_N \widehat{f^{(N)}} - Q_N g\|_{\mathcal{H}} = \|A_N f^{(N)} - g_N\|_{\mathbb{C}^N} = \|\varepsilon^{(N)}\|_{\mathbb{C}^N} \xrightarrow{N\to\infty} 0,$$

then the strong, weak, or component-wise vanishing of $A\widehat{f^{(N)}} - g$ is tantamount as, respectively, the strong, weak, or component-wise vanishing of $(A - Q_N A P_N)\widehat{f^{(N)}}$.

Since in addition $\|\widehat{f^{(N)}} - \widetilde{f}\|_{\mathcal{H}} \xrightarrow{N\to\infty} 0$ for some $\widetilde{f} \in \mathcal{H}$, then

$$\|(A - Q_N A P_N)\widehat{f^{(N)}}\|_{\mathcal{H}} \leqslant \|(A - Q_N A P_N)\widetilde{f}\|_{\mathcal{H}} + 2\|A\|_{\mathrm{op}} \|\widetilde{f} - \widehat{f^{(N)}}\|_{\mathcal{H}}$$
$$\xrightarrow{N\to\infty} 0$$

(the first summand in the r.h.s. above vanishing due to the operator strong convergence $Q_N A P_N \to A$), and (*) thus implies $\|\mathfrak{R}_N\|_{\mathcal{H}} = \|A\widehat{f^{(N)}} - g\|_{\mathcal{H}} \xrightarrow{N\to\infty} 0$.

Thus, part (i) is proved. Under the assumption of part (ii), one has $A\widehat{f^{(N)}} \to g$ (owing to (i)) and $A\widehat{f^{(N)}} \to A\widetilde{f}$ (by continuity), whence $A\widetilde{f} = g = Af$ and also (by injectivity) $f = \widetilde{f}$. This shows that $\|\mathscr{E}_N\|_{\mathcal{H}} = \|f - \widehat{f^{(N)}}\|_{\mathcal{H}} = \|\widetilde{f} - \widehat{f^{(N)}}\|_{\mathcal{H}} \to 0$. □

We observe that also here injectivity was only used in the analysis of the error, whereas it is not needed to conclude that $\|\mathfrak{R}_N\|_{\mathcal{H}} \to 0$.

As compared to Theorem A.1, Theorem A.2 now relies on the following hypotheses:

- *injectivity* of $A$,
- *asymptotic solvability of the truncated problems*,
- *convergence of the approximate solutions* $f^{(N)}$.

Once again, the last two assumptions are finite-dimensional in nature and from a numerical perspective their check involves a control of explicit $N$-dimensional vectors for satisfactorily many $N$'s.

*Remark A.6* When passing from a (well-defined) *compact* to a generic (well-defined) *bounded* inverse problem, one has to strengthen the hypothesis of uniform boundedness of the $\widehat{f^{(N)}}$'s to their actual strong convergence, in order for the residual $\mathfrak{R}_N$ to vanish strongly (in which case, as a by-product, also the error $\mathscr{E}_N$ vanishes strongly). In the compact case, $A - Q_N A P_N \to \mathbb{O}$ in operator norm (Lemma A.3), and it suffices that the $\widehat{f^{(N)}}$'s be uniformly bounded (or, in principle, have increasing norm $\|\widehat{f^{(N)}}\|_{\mathcal{H}}$ compensated by the vanishing of $\|A - Q_N A P_N\|_{\mathrm{op}}$), in order for $\|\mathfrak{R}_N\|_{\mathcal{H}} \to 0$.

When instead the sequence of the $\widehat{f^{(N)}}$'s *does not* converge strongly, in general one has to expect only weak vanishing of the residual, $\mathfrak{R}_N \rightharpoonup 0$, which in turn prevents the error to vanish strongly – for otherwise $\|\mathscr{E}_N\|_{\mathcal{H}} \to 0$ would imply $\|\mathfrak{R}_N\|_{\mathcal{H}} \to 0$, owing to (A.10). The following example shows such a possibility.

*Example A.2* For the right-shift $R = \sum_{n=1}^{\infty} |e_{n+1}\rangle\langle e_n|$ on $\ell^2(\mathbb{N})$ (Example 2.1), the inverse problem $Rf = g$ with $g = 0$ admits the unique solution $f = 0$. The truncated finite-dimensional problems induced by the bases $(u_n)_{n\in\mathbb{N}} = (v_n)_{n\in\mathbb{N}} = (e_n)_{n\in\mathbb{N}}$ is governed by the sub-diagonal matrix $R_N$ given by (A.13) with $\sigma_j = 1 \ \forall j \in \mathbb{N}$. Let us consider the sequence $(\widehat{f^{(N)}})_{N\in\mathbb{N}}$ with $\widehat{f^{(N)}} := e_N$ for each $N$. Then:

- $R_N f^{(N)} = 0 = g_N$ (the truncated problems are solved exactly),
- $\widehat{f^{(N)}} \rightharpoonup 0$ (only weakly, not strongly),
- $\mathfrak{R}_N = g - R\widehat{f^{(N)}} = -e_{N+1} \rightharpoonup 0$ (only weakly, not strongly).

Of course, what discussed so far highlights features of generic bounded inverse problems (as compared to compact ones). Ad hoc analyses for special classes of bounded inverse problems are available and complement the picture of Theorem A.2. This is the case, to mention one, when $A$ is an *algebraic operator*, namely $p(A) = \mathbb{O}$ for some polynomial $p$ (which includes finite-rank $A$'s) and one treats the inverse problem with the generalised minimal residual method (GMRES) [47].

In retrospect, based on the reasoning of this Section we prove Lemma A.2.

*Proof (Proof of Lemma A.2)* Let $f$ satisfy $Af = g$. The sequence $(f^{(N)})_{N\in\mathbb{N}}$ defined by

$$f^{(N)} := (P_N f)_N = f_N \qquad (\text{that is, } \widehat{f^{(N)}} = P_N f)$$

does the job, as a straightforward consequence of the fact, argued already at the beginning of this Section, that $Q_N A P_N \to A$ strongly in the operator topology. Indeed,

$$\begin{aligned}\|A_N f^{(N)} - g_N\|_{\mathbb{C}^N} &= \|Q_N A P_N \widehat{f^{(N)}} - Q_N g\|_{\mathcal{H}} \\ &\leqslant \|(Q_N A P_N - A)f\|_{\mathcal{H}} + \|(1 - Q_N)Af\|_{\mathcal{H}}.\end{aligned}$$

The strong limit yields the conclusion. $\square$

## A.5 Effects of changing the truncation basis: numerical evidence

Let us examine some of the features discussed theoretically so far through a few numerical tests concerning different choices of the truncation bases. We employed a Legendre, complex Fourier, and a Krylov basis to truncate the problems.

The two model operators that we considered are the already mentioned Volterra operator $V$ in $L^2[0,1]$ and the self-adjoint multiplication operator $M : L^2[1,2] \to L^2[1,2]$, $\psi \mapsto x\psi$. We examined the following two inverse problems.

**First problem:** $Vf_1 = g_1$, with $g_1(x) = \frac{1}{2}x^2$. The problem has unique solution

$$f_1(x) = x, \qquad \|f_1\|_{L^2[0,1]} = \frac{1}{\sqrt{3}} \simeq 0.5774 \qquad (A.14)$$

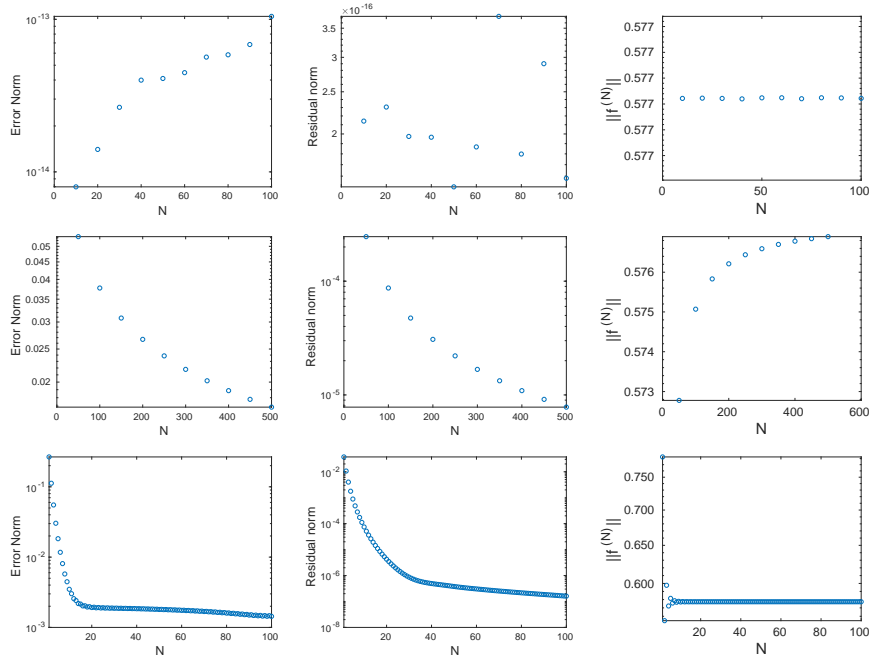

**Fig. A.1** Norm of the infinite-dimensional error and residual, and of the approximated solution for the Volterra inverse problem truncated with the Legendre (top row), complex Fourier (middle row), and Krylov bases (bottom row).

and $f_1$ is a Krylov solution, i.e., $f_1 \in \overline{\mathcal{K}(V,g)}$, although $f_1 \notin \mathcal{K}(V,g)$ (Examples 2.5 and 2.14).

**Second problem:** $Mf_2 = g_2$, with $g_2(x) = x^2$. The problem has unique solution

$$f_2(x) = x, \qquad \|f_2\|_{L^2[1,2]} = \sqrt{\frac{7}{3}} \simeq 1.5275 \qquad (A.15)$$

and $f_2$ is a Krylov solution. Indeed, $\mathcal{K}(M,g) = \{x^2 p \mid p \text{ is a polynomial on } [1,2]\}$ and $\overline{\mathcal{K}(M,g)} = \{x^2 h(x) \mid h \in L^2[1,2]\} = L^2[1,2]$, whence $f_2 \in \overline{\mathcal{K}(M,g)}$ and $f_2 \notin \mathcal{K}(M,g)$.

Both problems are treated with three different orthonormal bases: the Legendre polynomials and the complex Fourier modes (on the intervals $[0,1]$ or $[1,2]$, depending on the problem) solved using the QR factorisation algorithm, and the Krylov basis generated using the GMRES algorithm.

Incidentally, computationally speaking, generating accurate representations of the Legendre polynomials is quite demanding and accuracy can be lost rather soon due to their highly oscillatory nature, particularly at the end points. For this reason computations are limited up to $N = 100$ when considering the Legendre basis, and

A.5 Effects of changing the truncation basis: numerical evidence    131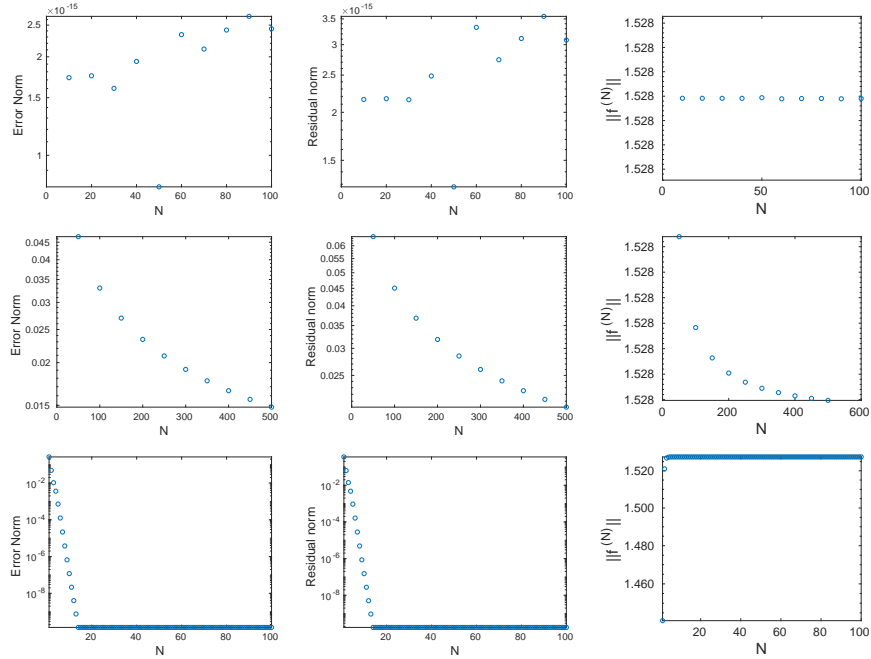

**Fig. A.2** Norm of the infinite-dimensional error, residual, and approximated solution for the $M$-multiplication inverse problem truncated with the Legendre (top row), complex Fourier (middle row), and Krylov bases (bottom row).

$N = 500$ when considering the complex Fourier basis. It is expected that there is no significant numerical error from the computation of the Legendre basis, as the $L^2[0,1]$ and $L^2[1,2]$ norms of the basis polynomials have less than 1% error compared to their exact unit value.

For each problem and each choice of the basis, three norms are monitored: infinite-dimensional error $\|\mathscr{E}_N\|_{L^2} = \|f - \widehat{f^{(N)}}\|_{L^2}$ ($f = f_1$ or $f_2$), infinite-dimensional residual $\|\mathfrak{R}_N\|_{L^2} = \|g - A\,\widehat{f^{(N)}}\|_{L^2}$ ($g = g_1$ or $g_2$; $A = V$ or $M$), and approximated solution $\|\widehat{f^{(N)}}\|_{L^2} = \|f^{(N)}\|_{\mathbb{C}^N}$.

Figure A.1 highlights the difference between the computation in the three bases for the Volterra operator.

- In the Legendre basis, $\|\mathscr{E}_N\|_{L^2}$ and $\|\mathfrak{R}_N\|_{L^2}$ are almost zero. $\|\widehat{f^{(N)}}\|_{L^2}$ stays bounded and constant with $N$ and matches the expected value (A.14). The approximated solutions reconstruct the exact solution $f_1$ at any truncation number.
- In the complex Fourier basis, both $\|\mathscr{E}_N\|_{L^2}$ and $\|\mathfrak{R}_N\|_{L^2}$ are some orders of magnitude *larger* than in the Legendre basis and decrease monotonically with $N$; in fact, $\|\mathscr{E}_N\|_{L^2}$ and $\|\mathfrak{R}_N\|_{L^2}$ display an evident convergence to zero, however attaining values that are more than ten orders of magnitude larger than the correspond-

ing error and residual norms for the same $N$ in the Legendre case. $\|\widehat{f^{(N)}}\|_{L^2}$, on the other hand, increases monotonically and appears to approach the theoretical value (A.14).

- In the Krylov basis $\|\mathscr{E}_N\|_{L^2}$ and $\|\mathfrak{R}_N\|_{L^2}$ decrease monotonically, relatively fast for small $N$'s, then rather slowly with $N$. Such quantities are smaller than in the Fourier basis. $\|\widehat{f^{(N)}}\|_{L^2}$ displays some initial highly oscillatory behaviour, but quickly approaches the theoretical value (A.14).

In contrast, Figure A.2 highlights the difference between the computation in the three bases for the $M$-multiplication operator.

- In the Legendre basis, $\|\mathscr{E}_N\|_{L^2}$ and $\|\mathfrak{R}_N\|_{L^2}$ are again almost zero. $\|\widehat{f^{(N)}}\|_{L^2}$ is constant with $N$ at the expected value (A.15). The approximated solutions reconstruct the exact solution $f_2$ at any truncation number.
- In the Fourier basis the behaviour of the above indicators is again qualitatively the same, and again with a much milder convergence rate in $N$ to the asymptotic values as compared with the Legendre case. $\|\mathscr{E}_N\|_{L^2}$ and $\|\mathfrak{R}_N\|_{L^2}$ still display an evident convergence to zero.
- The Krylov basis displays a fast initial decrease of both $\|\mathscr{E}_N\|_{L^2}$ and $\|\mathfrak{R}_N\|_{L^2}$ to the tolerance level of $10^{-10}$ that was set for the residual. $\|\widehat{f^{(N)}}\|_{L^2}$ also increases rapidly and remains constant at the expected value (A.15).

All this gives numerical evidence that the choice of the truncation basis *does* affect the sequence of solutions. The Legendre basis is best suited to these problems as $f_1$, $f_2$, $g_1$ and $g_2$ can be very well represented by the first few basis vectors.

# Index